\documentclass[reqno]{amsart}
\usepackage{amsfonts}
\usepackage{graphicx}
\usepackage{amsfonts,amsmath, amssymb}
\usepackage{bbm,mathrsfs,mathscinet}
\usepackage{marginnote}
\usepackage{color}
\usepackage{placeins}
\usepackage{verbatim}
\usepackage{tikz}
\usepackage{graphicx}
\usepackage{subfigure}
\usepackage{float}
\usetikzlibrary{decorations,arrows}
\allowdisplaybreaks
\counterwithin{figure}{section}
\numberwithin{figure}{section}

\numberwithin{equation}{section}

\usepackage[colorlinks,linkcolor=blue,anchorcolor=red,citecolor=blue]{hyperref}

\usepackage{cite}
\usepackage{cleveref}

\newcommand{\R}{\mathbb{R}}

\makeatletter
\newcommand{\rmnum}[1]{\romannumeral #1}
\newcommand{\Rmnum}[1]{\expandafter\@slowromancap\romannumeral #1@}
\makeatother

\newtheorem{theorem}{Theorem}[section]
\newtheorem{corollary}[theorem]{Corollary}
\newtheorem{lemma}[theorem]{Lemma}

\newtheorem{remark}{Remark}[section]

\newtheorem{definition}{Definition}[section]

\def\v{\varepsilon}

\def\la{\langle}
\def\ra{\rangle}

\begin{document}
	\title[The Boltzmann equation for 2D Taylor-Couette Flow]{The Boltzmann equation for 2D Taylor-Couette Flow}
	
	\author[R. J. Duan]{Renjun Duan}
	\address[R. J. Duan]{Department of Mathematics, The Chinese University of Hong Kong, Shatin, Hong Kong.}
	\email[R.J. Duan]{rjduan@math.cuhk.edu.hk}

	\author[W. Q. Wang]{Weiqiang Wang}
	\address[W. Q. Wang]{Department of Mathematics, University of Pittsburgh, Pittsburgh, PA 15260, USA}
	\email{wew179@pitt.edu}
	
	\author[Y. Wang]{Yong Wang}
	\address[Y. Wang]{Academy of Mathematics and Systems Science, Chinese Academy of Sciences, Beijing 100190, China; School of Mathematical Sciences, University of Chinese Academy of Sciences, Beijing 100049, China.}
	\email{yongwang@amss.ac.cn}

\begin{abstract}
In this paper, we investigate the existence of 2-D Taylor-Couette flow for a rarefied gas between two coaxial rotating cylinders, characterized by differing angular velocities at the outer boundary $\{r=1\}$ and the inner boundary $\{r=r_{1}>0\}$, with a small relative strength denoted by $\alpha$. We formulate the problem using the steady Boltzmann equation in polar coordinates and seek a solution invariant under rotation. We assume that the steady state has the specific form $F(r,v_{r},v_{\phi}-\alpha\frac{r-r_{1}}{1-r_{1}},v_{z})$, where the translation angular velocity $\alpha\frac{r-r_{1}}{1-r_{1}}$ is linearly sheared along the radial direction. With this ansatz, the problem is reduced to solve the nonlinear steady Boltzmann equation with geometric correction, subject to an external shear force of strength $\alpha$ and the homogeneous non-moving diffuse reflection boundary condition. We establish the existence of a non-equilibrium steady solution for any small enough shear rate $\alpha$ through Caflisch's decomposition, complemented by careful uniform estimates based on Guo's $L^{\infty} \cap L^2$ framework. The steady profile displays a polynomial tail behavior at large velocities. For the proof, we develop a delicate double-parameter $(\epsilon,\sigma)$-approximation argument for the construction of solutions. In particular, we obtain uniform macroscopic dissipation estimates in the absence of mass conservation for $\sigma \in [0,1)$ getting close to 1. Additionally, due to the non-trivial geometric effects, we develop subtle constructions of test functions by solving second-order ODEs with geometric corrections to establish macroscopic dissipation. Furthermore, we justify the non-negativity of the steady profile by demonstrating its large-time asymptotic stability with an exponential convergence rate under radial perturbations.
\end{abstract}


\date{\today}
\maketitle
	
\setcounter{tocdepth}{2}
\tableofcontents
	
\thispagestyle{empty}
	
\section{Introduction}
	
\subsection{Background} The Boltzmann equation is a fundamental model in kinetic theory which can be used to describe the motion of a rarefied gas when the Knudsen number is finite, cf. \cite{Ce1988, GS-2003,Ko1969,Sone07}. The stability near the global Maxwellian equilibrium is a natural result of the celebrated $H-$theorem and the intrinsic structure of the Boltzmann collision operator. The physical importance of studying the non-equilibrium dynamics of a rarefied gas which admits a non-Maxwellian steady state has been discussed in a recent review article \cite{EM-JSP} and references therein. In most situations, the appearance of non-equilibrium states is caused by inhomogeneous data at physical boundaries or at far fields with variable temperature and velocity, or possibly caused by a given external force. In fact, based on perturbation theory, some non-equilibrium steady states are non-Maxwellians but are still close to a global Maxwellian and decay exponentially for large velocities, see \cite{EGKM-13,EGKM-18,EGM-2018,UA-1983,UA-1986}.

However, there also exists a class of non-equilibrium states inherent to the Boltzmann dynamics that exhibit a polynomial tail for large velocities. Indeed, the numerical analysis \cite{BNV-2019,CC-1994,JNV-ARMA,Ke1,Ke2,LL-1961,OSA-1989,O2016,STO-1990,SY-1970} for kinetic plane Couette flow governed by the steady Boltzmann equation indicates that the solution should exhibit a polynomial tail, see also the monograph \cite{GS-2003} and the survey \cite{NV-2023}. Recently, the first author of this paper together with his collaborators \cite{DL-2020,DL-2022,DLY-2021} gave a rigorous justification of the existence and positivity of the 1D uniform shear flow and the Couette flow confined by two parallel infinite plates moving relative to each other with opposite velocities; see also \cite{MW-2024} for the asymptotic stability of such 1D Couette flow for the steady Boltzmann equation under 3D perturbation. The diffuse limit of such Couette flow for the Boltzmann equation was also rigorously justified in \cite{DLSY-arXiv}, see also \cite{DE-96, ELM-94, ELM-95, Guo-2006} for the Navier-Stokes limit of the Boltzmann equation. 

Closely related to the plane Couette flow, the Taylor-Couette flow between two coaxial rotating cylinders provides another example for studying the existence and stability issues on the Couette flow solution in the cylindrical framework. In fact, the Taylor-Couette problem is a classical problem in fluid dynamics, cf.~\cite{ChIo94}, and has been extensively studied in many aspects; see some mathematical progress \cite{AHL-NS, AHL24, CaLe24, KT20, Sch98} and references therein. The same problem for a rarefied gas flow governed by the nonlinear Boltzmann equation or related models is also fundamental and important in the kinetic theory of gases and has been investigated from various points of view. Most studies are based on direct numerical analysis that focuses on the effect of different types of boundary conditions, cf.~\cite{CerSer67, Ko, YA}. The main mathematical study of the stationary Boltzmann equation in a Taylor-Couette flow setting was made by Arkeryd-Nouri \cite{AN06,AN05} in the hydrodynamic regime; see also a recent work \cite{AEMN} for further investigation of the ghost effect induced by the boundary curvature. To the best of our knowledge, so far there have been few results on the existence and stability of the Taylor-Couette flow using kinetic theory. 
	
The current study of Taylor-Couette flow for the Boltzmann equation  between two coaxial rotating cylinders is motivated by \cite{DLY-2021}. The finite Knudsen number and diffuse reflection boundary condition are assumed. We aim at showing existence of the Taylor-Couette type non-equilibrium state for the steady Boltzmann equation driven by a weakly shearing effect due to the difference of rotating velocities at outer and inner boundaries, and also obtaining the exponential time asymptotic stability of the steady Taylor-Couette profile under small initial perturbation for the initial-boundary value problem.

\subsection{Reformulation of problem}
We consider the stationary motion of a rarefied gas flow contained in an annular domain 
$
\Omega:=\{{x}=(x_{1},x_{2})\,|\,0<r_{1}<|x|^2=|x_{1}|^2+|x_{2}|^2<1\}
$
on the horizontal plane with slab symmetry in the vertical $x_{3}$-axis direction. Let the non-negative unknown function $F(x,v)$ be the time-independent density distribution function of gas particles with velocity $v=(v_{1},v_{2},v_{3})\in \R^3$ located at position $x=(x_{1},x_{2})\in \Omega$. We assume that the motion of such rarefied gas flow can be governed by the steady Boltzmann equation
	\begin{equation}\label{1.1}
    v_{1}\partial_{x_{1}}F+v_{2}\partial_{x_{2}}F=\frac{1}{\mathrm{K\!n}}Q(F,F),
	\end{equation}
	subject to the completely diffuse reflection boundary condition
	\begin{equation}\label{1.2}
		\begin{aligned}
			&F(x,v)\vert_{\{|x|=1, v\cdot \vec{n}_{1}<0\}}=\sqrt{2\pi}\mathcal{M}_{1}\int_{u\cdot \vec{n}_{1}>0}F(x,u)|u\cdot \vec{n}_{1}|\,\mathrm{d}u,\\
			&F(x,v)\vert_{\{|x|=r_{1}, v\cdot \vec{n}_{r_1}<0\}}=\sqrt{2\pi}\mathcal{M}_{r_1}\int_{u\cdot \vec{n}_{r_1}>0} F(x,u) |u\cdot \vec{n}_{r_1}|\,\mathrm{d}u,
		\end{aligned}
	\end{equation}
    as well as a given total mass
		\begin{equation}\label{tm}
			\int_{\Omega}\int_{\R^3}F(x,v)\,\mathrm{d}v\mathrm{d}x=M_{0},
		\end{equation}
    for some constant $M_{0}>0$. 
	Here, the non-dimensional parameter $\mathrm{K\!n}>0$ is the Knudsen
number given by the ratio of the mean free path to the typical length, $\vec{n}_{1}=(x_{1},x_{2},0)$ and $\vec{n}_{r_1}=(-x_1/r_1,-x_2/r_1,0)$ are the outward unit normal vectors at the outer boundary $|x|=1$ and the inner boundary $|x|=r_1$, respectively, and the incoming boundary equilibrium states are given by
	$$
	\mathcal{M}_{1}=\frac{1}{(2\pi)^{\frac{3}{2}}}e^{-\frac{|v-\mathfrak{u}_{1}|^2}{2}},\quad \mathcal{M}_{r_1}=\frac{1}{(2\pi)^{\frac{3}{2}}}e^{-\frac{|v-\mathfrak{u}_{r_{1}}|^2}{2}},
	$$
    with $\mathfrak{u}_{1}$ and $\mathfrak{u}_{r_1}$ being the angular velocities of the outer and inner circles, respectively. The bilinear Boltzmann collision operator $Q(\cdot,\cdot)$ is defined as
	\begin{align}\label{Q-op}
		Q(F_1,F_2)&=\int_{\R^3}\int_{\mathbb{S}^2}B(\omega,v-v_\ast)[F_1(v_\ast')F_2(v')-F_1(v_\ast)F_2(v)]\,\mathrm{d}\omega\mathrm{d}v_{\ast}=:Q_{+}(F_{1},F_{2})-Q_{-}(F_{1},F_{2}),
	\end{align}
	where we have denoted $v'=v+[(v_\ast-v)\cdot\omega]\omega$ and $v_\ast'=v_\ast-[(v_\ast-v)\cdot\omega]\omega$ with $\omega\in \mathbb{S}^2$ as the velocities after collision, which can be solved from the conservation laws 
	$$
	v_\ast+v=v_\ast'+v',\quad  |v_\ast|^2+|v|^2=|v_\ast'|^2+|v'|^2.
	$$ 
	Throughout this paper, we consider the case of general hard potentials under the Grad's angular cutoff assumption, that is, the collision kernel $B(\omega,v-v_{*})$ is given by
	\begin{align*}
		 B(\omega,v-v_\ast)=|v-v_\ast|^\gamma B_0(\cos\theta),\ \cos\theta=\omega\cdot\frac{v-v_\ast}{|v-v_\ast|},\ \omega\in\mathbb{S}^2,
	\end{align*}
	with $0\leq\gamma\leq1$ and $\ 0\leq B_0(\cos\theta)\leq C|\cos\theta|$. 

In the following, without loss of generality, we assume
	$$
	\mathfrak{u}_{1}=\alpha\, \Big(-\frac{x_{2}}{\sqrt{x_{1}^2+x_{2}^2}},\frac{x_{1}}{\sqrt{x_{1}^2+x_{2}^2}},0\Big),\quad \mathfrak{u}_{r_1}=(0,0,0),
	$$
	where $\alpha=|\mathfrak{u}_{1}|>0$ represents the relative strength of the rotation rate. We remark that the general boundary angular velocities $\mathfrak{u}_{1}$ and $\mathfrak{u}_{r_{1}}$ with a small relative strength $|\mathfrak{u}_{1}-\mathfrak{u}_{r_{1}}|$ can also be treated in a similar way. For convenience, we introduce the polar coordinates
	$$
	\left\{\begin{aligned}
		&x_{1}=r\cos \phi,\\
		&x_{2}=r\sin \phi,\\
	\end{aligned}
	\right.\qquad r\in (r_{1},1),\quad \phi\in [-\pi,\pi).
	$$
Then, one has $\mathfrak{u}_{1}=\alpha(-\sin\phi, \cos\phi,0)$, $\vec{n}_{1}=(\cos\phi,\sin\phi,0)$, and $\vec{n}_{r_{1}}=(-\cos\phi,-\sin\phi,0)$. Due to the effect of boundary curvature, it is also convenient to study \eqref{1.1} in the polar coordinates. To this end, as in \cite{AN05,AN06,Wu-2016}, we shall make some changes of variables to reformulate the problem \eqref{1.1}--\eqref{1.2}. For simplicity of presentation, we set $\mathrm{K\!n}=1$ throughout the paper.
	
	\textbf{Substitution 1:} Let
	\begin{equation*}
		\left\{\begin{aligned}
			&x_{1}=r\cos\phi,\\
			&x_{2}=r\sin\phi,\\
			&v=v,
		\end{aligned}
		\right.
	\end{equation*}
	then \eqref{1.1}-\eqref{1.2} becomes
	\begin{equation}\label{1.4} 
		(v_{1}\cos\phi+v_{2}\sin\phi)\frac{\partial F}{\partial r}+\frac{1}{r}(-v_{1}\sin\phi+v_{2}\cos\phi)\frac{\partial F}{\partial \phi}=Q(F,F),
	\end{equation}
	and
	\begin{equation}\label{1.5}
		\begin{aligned} 
			&F(1,\phi,v)\vert_{v\cdot \vec{n}_{1}<0}=\sqrt{2\pi}\mathcal{M}_{1}\int_{u\cdot \vec{n}_{1}>0}F(1,\phi,u)|u\cdot \vec{n}_{1}|\,\mathrm{d}u,\\
			&F(r_{1},\phi,v)\vert_{v\cdot \vec{n}_{r_1}<0}=\sqrt{2\pi}\mathcal{M}_{r_1}\int_{u\cdot \vec{n}_{r_{1}}>0}F(r_{1},\phi,u)|u\cdot \vec{n}_{r_1}|\,\mathrm{d}u.
		\end{aligned}
	\end{equation}

	\textbf{Substitution 2:} Let
	\begin{equation*}
		\left\{\begin{aligned}
			&\eta=1-r,\\
			&\phi=\phi,\\
			&v=v,
		\end{aligned}
		\right.
	\end{equation*}
	then \eqref{1.4}--\eqref{1.5} becomes
	\begin{equation}\label{1.7} 
		-(v_{1}\cos\phi+v_{2}\sin\phi)\frac{\partial F}{\partial \eta}+\frac{1}{1-\eta}(-v_{1}\sin\phi+v_{2}\cos\phi)\frac{\partial F}{\partial \phi}=Q(F,F),
	\end{equation}
	and
	\begin{equation}\label{1.8}
		\begin{aligned} 
			&F(0,\phi,v)\vert_{v\cdot \vec{n}_{1}<0}=\sqrt{2\pi}\mathcal{M}_{1}\int_{u\cdot \vec{n}_{1}>0}F(0,\phi,u)|u\cdot \vec{n}_{1}|\,\mathrm{d}u,\\
			&F(\eta_{1},\phi,v)\vert_{v\cdot \vec{n}_{r_1}<0}=\sqrt{2\pi}\mathcal{M}_{r_1}\int_{u\cdot \vec{n}_{r_1}>0}F(\eta_{1},\phi,u)|u\cdot \vec{n}_{r_1}|\,\mathrm{d}u,
		\end{aligned}
	\end{equation}
	where we have denoted $\eta_{1}=1-r_1$.
	
	\textbf{Substitution 3:} Let 
	\begin{equation*}
		\left\{\begin{aligned}
			&\eta=\eta,\\
			&\phi=\phi,\\
			&\tilde{v}=-v,
		\end{aligned}
		\right.
	\end{equation*}
	then \eqref{1.7}--\eqref{1.8} becomes
	\begin{equation}\label{1.10} 
		(\tilde{v}_{1}\cos\phi+\tilde{v}_{2}\sin\phi)\frac{\partial F}{\partial \eta}-\frac{1}{1-\eta}(-\tilde{v}_{1}\sin\phi+\tilde{v}_{2}\cos\phi)\frac{\partial F}{\partial \phi}=Q(F,F),
	\end{equation}
	and
	\begin{equation}\label{1.11}
		\begin{aligned}
			&F(0,\phi,\tilde{v})\vert_{\tilde{v}\cdot \vec{n}_{1}>0}=\sqrt{2\pi}\tilde{\mathcal{M}}_{1}\int_{u\cdot \vec{n}_{1}<0}F(0,\phi,u)|u\cdot \vec{n}_{1}|\,\mathrm{d}u,\\
			&F(\eta_{1},\phi,\tilde{v})\vert_{\tilde{v}\cdot \vec{n}_{r_1}>0}=\sqrt{2\pi}\tilde{\mathcal{M}}_{r_1}\int_{u\cdot \vec{n}_{r_{1}}<0}F(\eta_{1},\phi,u)|u\cdot \vec{n}_{r_1}|\,\mathrm{d}u,
		\end{aligned}
	\end{equation}
	where
	$$
	\tilde{\mathcal{M}}_{1}=\frac{1}{(2\pi)^{\frac{3}{2}}}e^{-\frac{|\tilde{v}+\mathfrak{u}_{1}|^2}{2}},\quad \tilde{\mathcal{M}}_{r_1}=\frac{1}{(2\pi)^{\frac{3}{2}}}e^{-\frac{|\tilde{v}|^2}{2}}.
	$$
	
	\textbf{Substitution 4:} Let
	\begin{equation*}
		\left\{\begin{aligned}
			&\eta=\eta,\\
			&\phi=\phi,\\
			&v_{\eta}=\tilde{v}_{1}\cos\phi+\tilde{v}_{2}\sin \phi,\\
			&v_{\phi}=-\tilde{v}_{1}\sin\phi+\tilde{v}_{2}\cos \phi,\\
            &v_{z}=\tilde{v}_{3},
		\end{aligned}
		\right.
	\end{equation*}
	then \eqref{1.10}-\eqref{1.11} can be rewritten as
	\begin{equation}\label{1.13} 
		v_{\eta}\frac{\partial F}{\partial \eta}-\frac{1}{1-\eta}\Big(v_{\phi}\frac{\partial F}{\partial \phi}+v_{\phi}^2\frac{\partial F}{\partial v_{\eta}}-v_{\eta}v_{\phi}\frac{\partial F}{\partial v_{\phi}}\Big)=Q(F,F),
	\end{equation}
	and
	\begin{equation}\label{1.14}
		\begin{aligned}
			&F(0,\phi,v_{\eta},v_{\phi},v_{z})\vert_{v_{\eta}>0}=\sqrt{2\pi}\tilde{\mathcal{M}}_{1}\int_{u_{\eta}<0}F(0,\phi,u_{\eta},u_{\phi},u_{z})|u_{\eta}|\,\mathrm{d}u_{\eta}\mathrm{d}u_{\phi}{\rm d}u_{z},\\
			&F(\eta_{1},\phi,v_{\eta},v_{\phi},v_{z})\vert_{v_{\eta}<0}=\sqrt{2\pi}\tilde{\mathcal{M}}_{r_1}\int_{u_{\eta}>0}F(\eta_{1},\phi,u_{\eta},u_{\phi},u_{z})|u_{\eta}|\,\mathrm{d}u_{\eta}\mathrm{d}u_{\phi}{\rm d}u_{z},
		\end{aligned}
	\end{equation}
	where
	\begin{equation}\label{def.bdym12}
	\tilde{\mathcal{M}}_{1}=\frac{1}{(2\pi)^{\frac{3}{2}}}e^{-\frac{|v_{\eta}|^2+|v_{\phi}+\alpha|^2+|v_{z}|^2}{2}},\quad \tilde{\mathcal{M}}_{r_1}=\frac{1}{(2\pi)^{\frac{3}{2}}}e^{-\frac{|v_{\eta}|^2+|v_{\phi}|^2+|v_{z}|^2}{2}}.
	\end{equation}
	
	For \eqref{1.13}--\eqref{1.14}, we assume that the distribution function $F$ is invariant under rotations, that is, we can drop the dependence on $\phi$. With this ansatz, we search the solution of the form
	$$
	F(\eta,\phi, v_{\eta},v_{\phi},v_{z})=G(\eta,v_{\eta}+\mathfrak{u}_{\eta},v_{\phi}+\mathfrak{u}_{\phi},v_{z}),
	$$
	where the macroscopic translation velocity $\mathfrak{u}=(\mathfrak{u}_{\eta},\mathfrak{u}_{\phi},\mathfrak{u}_{z})$ is chosen as the classical Taylor-Couette type flow:
	$$
	\mathfrak{u}=(\mathfrak{u}_{\eta},\mathfrak{u}_{\phi},\mathfrak{u}_{z})=:(0,\frac{\alpha(\eta_{1}-\eta)}{\eta_{1}},0).
	$$
	Note that the angular velocity $\mathfrak{u}_{\phi}=\alpha (1-\eta/\eta_1)$ is linear in the radial variable $\eta \in (0,\eta_1)$, connecting the boundary angular velocities $\alpha$ and $0$ as given in \eqref{def.bdym12}. 
    Then \eqref{1.13}--\eqref{1.14} can be rewritten as
	\begin{equation}\label{1.16} 
		v_{\eta}\frac{\partial G}{\partial \eta}-\frac{\alpha}{\eta_{1}}v_{\eta}\frac{\partial G}{\partial v_{\phi}}-\frac{1}{1-\eta}\Big(v_{\phi}^2\frac{\partial G}{\partial v_{\eta}}-v_{\eta}v_{\phi}\frac{\partial G}{\partial v_{\phi}}\Big)=Q(G,G),
	\end{equation}
	with the rest diffuse reflection boundary condition:
	\begin{equation}\label{1.17}
		\begin{aligned} 
			&G(0,v_{\eta},v_{\phi},v_{z})\vert_{v_{\eta}>0}=\sqrt{2\pi}\mu\int_{u_{\eta}<0}G(0,u_{\eta},u_{\phi},u_{z})|u_{\eta}|\,\mathrm{d}u,\\
			&G(\eta_{1},v_{\eta},v_{\phi},v_{z})\vert_{v_{\eta}<0}=\sqrt{2\pi}\mu\int_{u_{\eta}>0}G(\eta_{1},u_{\eta},u_{\phi},u_{z})|u_{\eta}|\,\mathrm{d}u,
		\end{aligned}
	\end{equation}
	where
	$$
	\mu=\mu(v_{\eta},v_{\phi},v_{z})=\frac{1}{(2\pi)^{\frac{3}{2}}}e^{-\frac{|v_{\eta}|^2+|v_{\phi}|^2+|v_{z}|^2}{2}}.
	$$
	
\subsection{Asymptotic analysis}
	Observe that if the shear strength $\alpha$ tends to zero in \eqref{1.16}, then the trivial equilibrium state $G=C\mu$ for a suitable constant $C$ with \eqref{tm} satisfied is the stationary solution of the boundary-value problem \eqref{1.16}--\eqref{1.17}. Motivated by this observation, we let
	$$
	G=\tilde{M}_{0}\mu+\alpha \sqrt{\mu} g,
	$$
    with $\tilde{M}_{0}$ being the constant determined by the total mass:
    \begin{align}\label{2.0}
2\pi\tilde{M}_{0}\int_{0}^{\eta_{1}}\int_{\R^3}(1-\eta)\mu\,{\rm d}v{\rm d}\eta=M_{0},
    \end{align}
	then $g$ satisfies
	\begin{equation}\label{2.1}
		\begin{aligned}
			&v_{\eta}\frac{\partial g}{\partial \eta}-\frac{1}{1-\eta}\Big(v_{\phi}^2\frac{\partial g}{\partial v_{\eta}}-v_{\eta}v_{\phi}\frac{\partial g}{\partial v_{\phi}}\Big)-\frac{\alpha}{\eta_{1}}\Big( v_{\eta}\frac{\partial g}{\partial v_{\phi}}-\frac{v_{\eta}v_{\phi}}{2}g\Big)+\mathbf{L}g=-\frac{v_{\eta}v_{\phi}\sqrt{\mu}}{\eta_{1}}+\alpha \Gamma (g,g),
		\end{aligned}
	\end{equation}
	supplemented with the zero mass condition
    $$
    \int_{0}^{\eta_{1}}\int_{\R^3}(1-\eta)\sqrt{\mu}g\,{\rm d}v{\rm d}\eta=0,
    $$
 and the following boundary conditions:
	\begin{equation*}
		\begin{aligned}
			&g(0,v_{\eta},v_{\phi},v_{z})\vert_{v_{\eta}>0}=\sqrt{2\pi\mu}\int_{u_{\eta}<0}|u_{\eta}| (\sqrt{\mu}g)(0,u_{\eta},u_{\phi},u_{z})\,\mathrm{d}u,\\
			&g(\eta_{1},u_{\eta},u_{\phi},u_{z})\vert_{v_{\eta}<0}=\sqrt{2\pi\mu}\int_{u_{\eta}>0}|u_{\eta}| (\sqrt{\mu}g)(\eta_{1},u_{\eta},u_{\phi},u_{z})\,\mathrm{d}u.
		\end{aligned}
	\end{equation*}
Here, $\mathbf{L}=\nu-K$ is the linearized Boltzmann operator defined by
	$$
	\mathbf{L}g=\nu g-Kg=\nu g-\int_{\R^3}k(v,u)g(u)\,\mathrm{d}u,
	$$
	with
	$$
	\nu=\nu(v)=\int_{\R^3}\int_{\mathbb{S}^{2}}|v-v_{*}|^{\gamma}\mu(v_{*})B_{0}(\cos\theta)\,{\rm d}\omega{\rm d}v_{*},
	$$
and $k(v,u)$ being the symmetric integral kernel of $K$. The nonlinear term is given as
	$$
	\Gamma(f,g)=\frac{1}{\sqrt{\mu}}Q(\sqrt{\mu}f,\sqrt{\mu}g).
	$$
	It is easy to check that $\nu(v)$ depends only on $|v|$ and is an even function with respect to $v_{\eta}$, $v_{\phi}$ and $v_{z}$. Due to the trouble term $\frac{\alpha}{2\eta_{1}}v_{\eta}v_{\phi}g$ induced by the shear effect in \eqref{2.1}, we cannot directly get $\|g\|_{L^2}$ using the energy estimate. To overcome this difficulty, motivated by \cite{DL-2020, DL-2022, DLY-2021}, we use the following {\bf Caflisch's decomposition}
	$$
	\sqrt{\mu}g=g_{1}+\sqrt{\mu}g_{2},
	$$
	and require $g_{1}$ to satisfy
	\begin{equation}\label{2.3}
		\left\{\begin{aligned}
			&v_{\eta}\frac{\partial g_{1}}{\partial \eta}-\frac{1}{1-\eta}\Big(v_{\phi}^2\frac{\partial g_{1}}{\partial v_{\eta}}-v_{\eta}v_{\phi}\frac{\partial g_{1}}{\partial v_{\phi}}\Big)-\frac{\alpha}{\eta_{1}} v_{\eta}\frac{\partial g_{1}}{\partial v_{\phi}}+\frac{\alpha}{2\eta_{1}}v_{\eta}v_{\phi}\sqrt{\mu}g_{2}+\nu g_{1}\\
			&\quad =(1-\chi_{M})\mathcal{K}g_{1}-\frac{1}{\eta_{1}}v_{\eta}v_{\phi}\mu+\alpha Q(\sqrt{\mu}g,\sqrt{\mu}g),\\
			&g_{1}(0,v_{\eta},v_{\phi},v_{z})\vert_{v_{\eta}>0}=0,\quad g_{1}(\eta_{1},v_{\eta},v_{\phi},v_{z})\vert_{v_{\eta}<0}=0,
		\end{aligned}\right.
	\end{equation}
	and $g_{2}$ to satisfy 
	\begin{equation}\label{2.4}
		\left\{\begin{aligned}
			&v_{\eta}\frac{\partial g_{2}}{\partial \eta}-\frac{1}{1-\eta}\Big(v_{\phi}^2\frac{\partial g_{2}}{\partial v_{\eta}}-v_{\eta}v_{\phi}\frac{\partial g_{2}}{\partial v_{\phi}}\Big)-\frac{\alpha}{\eta_{1}} v_{\eta}\frac{\partial g_{2}}{\partial v_{\phi}}+\nu g_{2}=Kg_{2}+\chi_{M}\mu^{-\frac{1}{2}}\mathcal{K}g_{1},\\
			&g_{2}(0,v_{\eta},v_{\phi},v_{z})\vert_{v_{\eta}>0}=\sqrt{2\pi\mu}\int_{u_{\eta}<0}|u_{\eta}|(g_{1}+\sqrt{\mu}g_{2})(0,u_{\eta},u_{\phi},u_{z})\,\mathrm{d}u,\\
			&g_{2}(\eta_{1},v_{\eta},v_{\phi},v_{z})\vert_{v_{\eta}<0}=\sqrt{2\pi\mu}\int_{u_{\eta}>0}|u_{\eta}|(g_{1}+\sqrt{\mu}g_{2})(\eta_{1},u_{\eta},u_{\phi},u_{z})\,\mathrm{d}u,
		\end{aligned}\right.
	\end{equation}
	where $\chi_{M}(v)$ is a smooth cut-off function satisfying
	\begin{equation*}
		\chi_{M}(v)=\left\{
		\begin{aligned}
			&1,\quad |v|\leq M,\\
			&0,\quad |v|\geq M+1,
		\end{aligned}
		\right.
	\end{equation*}
	and 
	$$
	\mathcal{K}f:=\sqrt{\mu}K(\frac{f}{\sqrt{\mu}})=\int_{\R^3}\mu^{\frac{1}{2}}(v)k(v,u)\mu^{-\frac{1}{2}}(u)f(u)\,\mathrm{d}u.
	$$
	Note that the boundary conditions for $g_{1}$ are set to be zero for incoming velocities.

	\subsection{Main results} Let $w^{\ell}(v)=(1+|v|^2)^{\frac{\ell}{2}}$ for $\ell>0$. Now, we state the main result for the stationary problem \eqref{1.1}--\eqref{1.2}.
	
	\begin{theorem}\label{thm1}
		Let $0\leq\gamma\leq 1$ and $\ell_{\infty}\gg 4$. For any given total mass $M_{0}>0$, there exist positive constants $\alpha_{*}$ and $C_{0}$ such that for any $\alpha\in (0,\alpha_{*})$, the steady boundary value problem \eqref{1.1}--\eqref{1.2} ${\rm (}$equivalently \eqref{1.16}--\eqref{1.17}${\rm )}$ and \eqref{tm} admits a unique solution in the form of
		\begin{align*}
		F_{st}(x_{1},x_{2},v_{1},v_{2},v_{3})&=F_{st}(\eta,\phi,v_{\eta},v_{\phi},v_{z})=G_{st}(\eta,v_{\eta},v_{\phi}+\mathfrak{u}_{\phi},v_{z})\\
&=\tilde{M}_{0}\mu(v_{\eta},v_{\phi}+\mathfrak{u}_{\phi},v_{z})+\alpha (g_{1}+\sqrt{\mu}g_{2})(\eta,v_{\eta},v_{\phi}+\mathfrak{u}_{\phi},v_{z}),
		\end{align*}
		where $\mathfrak{u}_{\phi}=\frac{\alpha(\eta_{1}-\eta)}{\eta_{1}}$, $\tilde{M}_{0}$ is determined by \eqref{2.0},
		and $(g_{1},g_{2})(\eta,v_{\eta},v_{\phi},v_{z})$ is the unique solution to the coupled system \eqref{2.3}--\eqref{2.4} satisfying
		$$
		\|w^{\ell_{\infty}}(g_{1}, g_{2})\|_{L_{\eta,v}^{\infty}}\leq C_{0}.
		$$
	\end{theorem}
	
\begin{remark}
    Compared with \cite{DL-2020, DL-2022,DLY-2021}, it is very difficult to obtain the estimates of the derivatives of the remainder term with respect to the velocity variable due to the appearance of geometric correction terms in \eqref{2.1} arising from the boundary curvature, cf.~\cite{Wu-2016}. This also restricts us to expand $G$ only up to the first order of $\alpha$.
\end{remark}
	
To show the non-negativity of the steady profile $F_{st}(x,v)$ established in Theorem \ref{thm1}, we shall consider the following time-evolutionary problem
\begin{equation}\label{2.6}
\left\{\begin{aligned}
	&\partial_{t}G+v_{\eta}\frac{\partial G}{\partial \eta}-\frac{\alpha}{\eta_{1}}v_{\eta}\frac{\partial G}{\partial v_{\phi}}-\frac{1}{1-\eta}(v_{\phi}^2\frac{\partial G}{\partial v_{\eta}}-v_{\eta}v_{\phi}\frac{\partial G}{\partial v_{\phi}})=Q(G,G),\\
	&G(0,\eta,v_{\eta},v_{\phi},v_{z})=G_{0}(\eta,v_{\eta},v_{\phi},v_{z}),\\
	&G(t,0,v_{\eta},v_{\phi},v_{z})\vert_{v_{\eta}>0}=\sqrt{2\pi}\mu\int_{u_{\eta}<0}|u_{\eta}|G(t,0,u_{\eta},u_{\phi},u_{z})\,{\rm d}u,\\
	&G(t,\eta_{1},v_{\eta},v_{\phi},v_{z})\vert_{v_{\eta}<0}=\sqrt{2\pi}\mu\int_{u_{\eta}>0}|u_{\eta}|G(t,\eta_{1},u_{\eta},u_{\phi},u_{z})\,{\rm d}u.
\end{aligned}
\right.
\end{equation}
If $G_{0}(\eta,v_{\eta},v_{\phi},v_{z})$ is very close to $G_{st}(\eta,v_{\eta},v_{\phi},v_{z})$, one may expect that the solution of time-dependent problem \eqref{2.6} converges to the solution $G_{st}(\eta,v_{\eta},v_{\phi},v_{z})$ of the steady problem \eqref{1.16}--\eqref{1.17} as the time goes to infinity. For this, the second result is concerned with the large time asymptotic stability of the stationary solution $F_{st}$ established in Theorem \ref{thm1}, which implies the non-negativity of $F_{st}$.

\begin{theorem}\label{thm2}
Let $F_{st}(x,v)$ be the steady state obtained in Theorem \ref{thm1}, that is,
\begin{align*}
F_{st}(x,v)&=G_{st}(\eta,v_{\eta}+\mathfrak{u}_{\eta},v_{\phi}+\mathfrak{u}_{\phi},v_{z}),
\end{align*}
where $G_{st}(\eta,v_{\eta},v_{\phi},v_{z})$ is the steady solution of problem \eqref{1.16}--\eqref{1.17}. There are constants $\alpha_{*}>0$, $\varepsilon_{0}>0$, $\lambda_{0}>0$ and $C>0$ such that for any $\alpha\in (0,\alpha_{*})$, if the initial data $F_{0}(x,v)=G_{0}(\eta,v_{\eta}+\mathfrak{u}_{\eta},v_{\phi}+\mathfrak{u}_{\phi},v_{z})\geq 0$ satisfy
$$
\|w^{\ell_{\infty}}[G_{0}-G_{st}]\|_{L_{\eta,v}^{\infty}}\leq \varepsilon_{0},
$$
with
\begin{align}\label{M1}
\int_{0}^{\eta_{1}}\int_{\R^3}(1-\eta)[G_{0}-G_{st}](\eta,v_{\eta},v_{\phi},v_{z})\,{\rm d}v{\rm d}\eta=0,
\end{align}
then the initial-boundary value problem \eqref{2.6} admits a unique solution $G(t,\eta,v)\geq 0$ satisfying the following exponential decay estimate
\begin{align}\label{M2}
  \|w^{\ell_{\infty}}[G-G_{st}]\|_{L_{\eta,v}^{\infty}}\leq Ce^{-\lambda_{0}t}\|w^{\ell_{\infty}}[G_{0}-G_{st}]\|_{L_{\eta,v}^{\infty}},  
\end{align}
for any $t\geq 0$.
\end{theorem}

\medskip

Now, we briefly emphasize the key points in the proof of main results above. 

\medskip
{(\rmnum{1})} The requirement of zero boundary conditions for \eqref{2.3} is essential to obtain weighted $L_{\eta,v}^{\infty}$ estimates on $g_{1}$. In fact, if we impose the following diffusive boundary conditions on $g_{1}$:
    \begin{align*}
    &g_{1}(0,v_{\eta},v_{\phi},v_{z})\vert_{v_{\eta}>0}=\sqrt{2\pi}\mu\int_{u_{\eta}<0}g_{1}(0,u_{\eta},u_{\phi},u_{z})|u_{\eta}|\,{\rm d}u,\\ &
    g_{1}(\eta_{1},v_{\eta},v_{\phi})\vert_{v_{\eta}<0}=\sqrt{2\pi}\mu\int_{u_{\eta}>0}g_{1}(\eta_{1},u_{\eta},u_{\phi},u_{z})|u_{\eta}|\,{\rm d}u,
    \end{align*}
    we have to iterative many times along the backward trajectory to utilize the smallness of diffusion of boundary data to close the estimates on $\|w^{\ell}g_{1}\|_{L_{\eta,v}^{\infty}}$, cf. \cite{Guo-2010}. However, for each iteration process, one has to control the terms
    \begin{align}\label{I1}
    \begin{aligned}
   &\sqrt{2\pi} w^{\ell}\mu\int_{u_{\eta}<0}(1-\chi_{M})\mathcal{K}_{w}(w^{\ell}g_{1})(0,u)(w^{-\ell}(u)|u_{\eta}|)\,{\rm d}u,\\
   &\sqrt{2\pi} w^{\ell}\mu\int_{u_{\eta}>0}(1-\chi_{M})\mathcal{K}_{w}(w^{\ell}g_{1})(\eta_{1},u)(w^{-\ell}(u)|u_{\eta}|)\,{\rm d}u,
   \end{aligned}
\end{align}
arising from the non-homogeneous diffusive boundary conditions. Here, $\mathcal{K}_{w}f=:w^{\ell}\mathcal{K}(w^{-\ell}f)$. As mentioned in \cite{DL-2020,DL-2022,DLY-2021}, the smallness of the term $(1-\chi_{M})\mathcal{K}_{w}(w^{\ell}g_{1})$ comes mainly from the weight $\ell$ for large velocities (see Lemma \ref{K}). For example, for $\gamma=0$, we only have
$$
\|(1-\chi_{M})\mathcal{K}_{w}(w^{\ell}g_{1})\|_{L_{\eta,v}^{\infty}}\lesssim \frac{1}{\ell}\|w^{\ell}g_{1}\|_{L_{\eta,v}^{\infty}},
$$
for $\ell$ and $M$ sufficiently large. Then one may only bound \eqref{I1} as
\begin{align*}
\frac{\|w^{\ell}\mu\|_{L_{v}^{\infty}}}{\ell}\|w^{\ell}g_{1}\|_{L_{\eta,v}^{\infty}}.
\end{align*}
Noting that $\|w^{\ell}\mu\|_{L_{v}^{\infty}}$ is finite but depends on $\ell$, and $\lim\limits_{\ell\to \infty}\|w^{\ell}\mu\|_{L_{v}^{\infty}}=\lim\limits_{\ell\to \infty}\sqrt{e}\big(\sqrt{\frac{\ell}{e}}\big)^{\ell}=\infty$, one fails to get any smallness in front of $\|w^{\ell}g_{1}\|_{L_{\eta,v}^{\infty}}$ so that it is impossible to close the estimates of $\|w^{\ell}g_{1}\|_{L_{\eta,v}^{\infty}}$. A similar treatment is also applied to solve the unsteady problem \eqref{2.6}.

\medskip
(\rmnum{2}) Although the zero boundary conditions for $g_{1}$ is very helpful in establishing $\|w^{\ell}g_{1}\|_{L_{\eta,v}^{\infty}}$, it makes the boundary condition for $g_{2}$ coupled with the one for $g_{1}$, which results in many difficulties in establishing the existence of $g_{2}$. Indeed, we have to construct the solution for $(g_{1},g_{2})$ through studying the following coupled approximate systems with two parameters $\epsilon>0$ and $0\leq \sigma\leq 1$:
\begin{align}\label{I2-0}
	\left\{\begin{aligned}
		&\epsilon g_{1}+v_{\eta}\frac{\partial g_{1}}{\partial \eta}-\frac{1}{1-\eta}\Big(v_{\phi}^2\frac{\partial g_{1}}{\partial v_{\eta}}-v_{\eta}v_{\phi}\frac{\partial g_{1}}{\partial v_{\phi}}\Big)-\frac{\alpha}{\eta_{1}}v_{\eta}\frac{\partial g_{1}}{\partial v_{\phi}}+\frac{\alpha}{2\eta_{1}}v_{\eta}v_{\phi}\sqrt{\mu}g_{2}+\nu g_{1}\\
		&\quad =\sigma(1-\chi_{M})\mu^{-\frac{1}{2}}\mathcal{K}g_{1}+\mathcal{S}_{1},\\
		&g_{1}(0,v_{\eta},v_{\phi},v_{z})\vert_{v_{\eta}>0}=g_{1}(\eta_{1},v_{\eta},v_{\phi},v_{z})\vert_{v_{\eta}<0}=0,
	\end{aligned}
	\right.
\end{align}
\begin{align}\label{I2}
\left\{\begin{aligned}
    &\epsilon g_{2}+v_{\eta}\frac{\partial g_{2}}{\partial \eta}-\frac{1}{1-\eta}\Big(v_{\phi}^2\frac{\partial g_{2}}{\partial v_{\eta}}-v_{\eta}v_{\phi}\frac{\partial g_{2}}{\partial v_{\phi}}\Big)-\frac{\alpha}{\eta_{1}}v_{\eta}\frac{\partial g_{2}}{\partial v_{\phi}}+\nu g_{2}\\
    &\quad =\sigma Kg_{2}+\sigma\chi_{M}\mu^{-\frac{1}{2}}\mathcal{K}g_{1}+\mathcal{S}_{2},\\
    &g_{2}(0,v_{\eta},v_{\phi},v_{z})\vert_{v_{\eta}>0}=\sqrt{2\pi\mu}\int_{u_{\eta}<0}|u_{\eta}|(g_{2}+\sigma\sqrt{\mu}g_{1})(0,u)\,{\rm d}u+\mathcal{S}_{2,b}^{(1)},\\
    &g_{2}(\eta_{1},v_{\eta},v_{\phi},v_{z})\vert_{v_{\eta}<0}=\sqrt{2\pi\mu}\int_{u_{\eta}>0}|u_{\eta}|(g_{2}+\sigma\sqrt{\mu}g_{1})(0,u)\,{\rm d}u+\mathcal{S}_{2,b}^{(2)}.
    \end{aligned}
    \right.
\end{align}
For given $\epsilon>0$, one has to apply the Leray-Schauder fixed point theorem to lift $\sigma$ from $0$ to $1$. One of the most important ingredients in applying the Leray-Schauder fixed point theorem is to establish uniform estimates with respect to $\sigma$. However, due to the coupling of $g_{1}$ and $g_{2}$, the weighted $L_{\eta,v}^{\infty}$ estimates may only lead to
\begin{align}\label{I2-1}
\|w^{\ell}g_{1}\|_{L_{\eta,v}^{\infty}}&\leq C\alpha\|w^{\ell}g_{2}\|_{L_{\eta,v}^{\infty}}+C\|\nu^{-1}w^{\ell}\mathcal{S}_{1}\|_{L_{\eta,v}^{\infty}},\nonumber\\
\|w^{\ell}g_{2}\|_{L_{\eta,v}^{\infty}}&\leq C_{\ell}\|w^{\ell}g_{1}\|_{L_{\eta,v}^{\infty}}+C_{\ell}\|\nu^{-1}w^{\ell}\mathcal{S}_{1}\|_{L_{\eta,v}^{\infty}}+C_{\ell}|w^{\ell}[\mathcal{S}_{2,b}^{(0)},\mathcal{S}_{2,b}^{(1)}]|_{L^{\infty}(\gamma_{-})}+C_{\ell}\|g_{2}\|_{L_{\eta,v}^2},
\end{align}
see Lemma \ref{lem2.2} for details. Therefore, to close the estimate, one has to establish $L_{\eta,v}^2$ estimates on $g_{2}$. Applying the direct energy estimates and the trace Lemma \ref{tt} to \eqref{I2} leads to
\begin{align}\label{I3}
&\epsilon\|g_{2}\|_{L_{\eta,v}^2}^2+(1-\sigma)\|g_{2}\|_{L_{\nu}^2}^2+\sigma\|(\mathbf{I-P})g_{2}\|_{L_{\nu}^2}^2+|(I-P_{\gamma})g_{2}|_{L^2(\gamma_{+})}^2\nonumber\\
&\leq C\delta\big[\|g_{2}\|_{L_{\nu}^2}^2+|(I-P_{\gamma})g_{2}|_{L^2(\gamma_{+})}^2\big]+C_{\delta}\big[\|w^{\ell}g_{1}\|_{L_{\eta,v}^{\infty}}^2+\|\nu^{-1}w^{\ell}\mathcal{S}_{2}\|_{L_{\eta,v}^{\infty}}^2+|w^{\ell}(\mathcal{S}_{2,b}^{(1)},\mathcal{S}_{2,b}^{(2)})|_{L^{\infty}(\gamma_{-})}^2\big],
\end{align}
see \eqref{E1}--\eqref{E5-1} for details. Then, when $\sigma$ is arbitrarily close to $1$ but never equals $1$, we can only obtain from \eqref{I3} that
\begin{align}\label{I4}
\|g_{2}\|_{L_{\eta,v}^2}^2\leq C_{\epsilon}\big[\|w^{\ell}g_{1}\|_{L_{\eta,v}^{\infty}}^2+\|\nu^{-1}w^{\ell}\mathcal{S}_{2}\|_{L_{\eta,v}^{\infty}}^2+|w^{\ell}(\mathcal{S}_{2,b}^{(1)},\mathcal{S}_{2,b}^{(2)})|_{L^{\infty}(\gamma_{-})}^2\big],
\end{align}
by taking $\delta\ll \epsilon$. Substituting \eqref{I4} into \eqref{I2-1}, one has
\begin{align*}
    &\|w^{\ell}g_{1}\|_{L_{\eta,v}^{\infty}}\leq C\alpha\|w^{\ell}g_{1}\|_{L_{\eta,v}^{\infty}}+C\|\nu^{-1}w^{\ell}\mathcal{S}_{1}\|_{L_{\eta,v}^{\infty}},\\
    &\|w^{\ell}g_{2}\|_{L_{\eta,v}^{\infty}}\leq C_{\epsilon}\big[\|w^{\ell}g_{1}\|_{L_{\eta,v}^{\infty}}+\|\nu^{-1}w^{\ell}\mathcal{S}_{2}\|_{L_{\eta,v}^{\infty}}+|w^{\ell}(\mathcal{S}_{2,b}^{(1)},\mathcal{S}_{2,b}^{(2)})|_{L^{\infty}(\gamma_{-})}\big].
\end{align*}
To close the above estimate, one has to require $\alpha\ll \epsilon$. However, this results in $\alpha\to 0$ when we take the limit $\epsilon\to 0$. Therefore, we need to establish uniform estimate on $\|\mathbf{P}g_{2}\|_{L_{\eta,v}^2}$ with respect to $\epsilon$ and $\sigma$, especially when $\sigma$ is very close to 1.

The classical method to establish $\|\mathbf{P}g_{2}\|_{L_{\eta,v}^2}$ is to apply some appropriate test functions to \eqref{I2} so that macroscopic estimates can be recovered from the term $v_{\eta}\frac{\partial g_{2}}{\partial \eta}$. It should be pointed out that the conservation of mass is very essential in constructing appropriate test functions to make the subsequent boundary terms disappear or be controlled. As in \cite{DL-2020,DL-2022,DLY-2021}, the conservation of mass does not hold for \eqref{I2} since we have used the Caflisch decomposition, so one has to consider the original function $g=\mu^{-\frac{1}{2}}g_{1}+g_{2}$ and attempt to get $\|\mathbf{P}g\|_{L_{\eta,v}^2}$. Then use $\|\mathbf{P}g\|_{L_{\eta,v}^{2}}$ and $\|w^{\ell}g_{1}\|_{L_{\eta,v}^{\infty}}$ to control $\|\mathbf{P}g\|_{L_{\eta,v}^2}$.

However, mass conservation still fails for $g$ due to the appearance of $\sigma$, resulting in many difficulties in establishing the estimate on $\|\mathbf{P}g\|_{L_{\eta,v}^2}$. In fact, a direct calculation shows that $g$ satisfies 
\begin{align*}
    &\epsilon g+v_{\eta}\frac{\partial g}{\partial \eta}-\frac{1}{1-\eta}\Big(v_{\phi}^2\frac{\partial g}{\partial v_{\eta}}-v_{\eta}v_{\phi}\frac{\partial g}{\partial v_{\phi}}\Big)-\frac{\alpha}{\eta_{1}}v_{\eta}\frac{\partial g}{\partial v_{\phi}}+\frac{\alpha}{2\eta_{1}}v_{\eta}v_{\phi}\sqrt{\mu}g+\nu g= \sigma Kg+(\mu^{-\frac{1}{2}}\mathcal{S}_{1}+\mathcal{S}_{2}).
\end{align*}
Let
$$
\mathbf{P}g=\big[a(\eta)+b_{1}(\eta)v_{\eta}+b_{2}(\eta)v_{\phi}+b_{3}(\eta)v_{z}+c(\eta)(|v|^2-3)\big]\sqrt{\mu}.
$$
Multiplying \eqref{I2} by $(1-\eta)\sqrt{\mu}$ and integrating over $[0,\eta_{1}]\times \R^3$ leads to
\begin{align}\label{I6}
[\epsilon +(1-\sigma)\beta_{0}]\int_{0}^{\eta_{1}}(1-\eta)a(\eta)\,{\rm d}\eta
&=-(1-\sigma)\Big((1-\eta_{1})\int_{v_{\eta}>0}g_{1}(\eta_{1})|v_{\eta}|\,{\rm d}v-\int_{v_{\eta}<0}g_{1}(0)|v_{\eta}|\,{\rm d}v\Big)\nonumber\\
&\quad -(1-\eta_{1})\int_{v_{\eta}<0}v_{\eta}\sqrt{\mu}\mathcal{S}_{2,b}^{(2)}\,{\rm d}v+\int_{v_{\eta}>0}v_{\eta}\sqrt{\mu}\mathcal{S}_{2,b}^{(1)}\,{\rm d}v\nonumber\\
&\quad -(1-\sigma)\beta_{\gamma}\int_{0}^{\eta_{1}}(1-\eta)c(\eta)\,{\rm d}\eta
+\int_{0}^{\eta_{1}}(1-\eta)\int_{\R^3}(\mathcal{S}_{1}+\sqrt{\mu}\mathcal{S}_{2})\,{\rm d}v{\rm d}\eta\nonumber\\
&\quad -(1-\sigma)\int_{0}^{\eta_{1}}(1-\eta)\int_{\R^3}\nu(v)(\mathbf{I-P})g\sqrt{\mu}\,{\rm d}v{\rm d}\eta,
\end{align}
where $\beta_{0}=\int_{\R^3}\nu(v)\mu\,{\rm d}v$ and $\beta_{\gamma}=\int_{\R^3}\nu(v)(|v|^2-3)\mu\,{\rm d}v$ are two constants.
Therefore, $\int_{0}^{\eta_{1}}(1-\eta)a(\eta)\,{\rm d}\eta$ does not vanish and even depends on $c(\eta)$, which leads that the subsequent boundary terms does not vanish and is hard to be controlled after applying the test function. Fortunately, we can obtain from \eqref{I6} that
\begin{align}\label{I6-1}
\int_{0}^{\eta_{1}}(1-\eta)a(\eta)\,{\rm d}\eta&\lesssim \frac{1}{\beta_{0}}\|w^{\ell}g_{1}\|_{L_{\eta,v}^{\infty}}+\frac{1}{c_{0}}\|(\mathbf{I-P})g\|_{L_{\nu}^2}+\frac{\beta_{\gamma}}{\beta_{0}}\|c\|_{L^2}\nonumber\\
&\quad +\frac{1}{\epsilon}\|\nu^{-1}w^{\ell}(\mathcal{S}_{1},\mathcal{S}_{2})\|_{L_{\eta,v}^{\infty}}+\frac{1}{\epsilon}|w^{\ell}(\mathcal{S}_{2,b}^{(1)},\mathcal{S}_{2,b}^{(2)})|_{L^{\infty}(\gamma_{-})},
\end{align}
which leads to the following estimates on the boundary terms:
\begin{align*}
   &\int_{\R^3}v_{\eta}(g\Psi_{a})(0,v)\,{\rm d}v-\int_{\R^3}v_{\eta}(g\Psi_{a})(\eta_{1},v)\,{\rm d}v\nonumber\\
   &\leq \kappa_{1}|P_{\gamma}g|_{L^2(\gamma_{+})}^2+\kappa_{1}\|a\|_{L^2}^2+C_{\kappa_{1}}\|c\|_{L^2}^2+C_{\kappa_{1}}|(I-P_{\gamma})g|_{L^2(\gamma_{+})}^2 +C_{\kappa_{1}}\|w^{\ell}g_{1}\|_{L_{\eta,v}^{\infty}}^2\nonumber\\
&\quad +C_{\epsilon,\kappa_{1}}|w^{\ell}(\mathcal{S}_{2,b}^{(1)},\mathcal{S}_{2,b}^{(2)})|_{L^{\infty}(\gamma_{-})}^2+C_{\epsilon,\kappa_{1}}\|\nu^{-1}w^{\ell}(\mathcal{S}_{1},\mathcal{S}_{2})\|_{L_{\eta,v}^{\infty}}^2,
\end{align*}
after applying the test function $\Psi_{a}$. Then, applying the trace Lemma \ref{tt} to control the $|P_{\gamma}g|_{L^{\infty}(\gamma_{+})}$ and appropriate test functions to control the other terms, we can obtain 
\begin{align*}
	\|a\|_{L^2}^2&\lesssim \alpha\|a\|_{L^2}\|b_{2}\|_{L^2}+(1-\sigma_{0})\|a\|_{L^2}\|b_{1}\|_{L^2}+\kappa_{1}\|g_{2}\|_{L_{\nu}^2}^2+C_{\kappa_{1}}\|(b_{1},c)\|_{L^2}^2\nonumber\\
	 &\quad +C_{\kappa_{1}}|(I-P_{\gamma})g|_{L^2(\gamma_{+})}^2+C_{\kappa_{1}}\|(\mathbf{I-P})f\|_{L_{\nu}^2}^2+C_{\kappa_{1}}|w^{\ell}g_{1}|_{L^{\infty}(\gamma_{+})}^2\nonumber\\
	 &\quad +C_{\epsilon,\kappa_{1}}|w^{\ell}(\mathcal{S}_{2,b}^{(1)},\mathcal{S}_{2,b}^{(2)})|_{L^{\infty}(\gamma_{-})}^2+C_{\epsilon,\kappa_{1}}\|\nu^{-1}w^{\ell}(\mathcal{S}_{1},\mathcal{S}_{2})\|_{L_{\eta,v}^{\infty}}^2,
\end{align*}
see \eqref{a}--\eqref{E27} for details. Even though there is no small coefficient in front of $\|(b_{1},c)\|_{L^2}$ on the right hand side of the above estimates, we can close them by constructing some appropriate test functions in treating $\|(b_{1},c)\|_{L^2}$ such that the estimate of $\|(b_{1},c)\|_{L^2}$ is decoupled with $\|a\|_{L^2}$. Here, it is also worth pointing out $\epsilon-$dependence only appears in front of the source terms, and vanishes in front of $\|w^{\ell}g_{1}\|_{L_{\eta,v}^{\infty}}$. Then we are able to use it to close the uniform estimates \eqref{I2-1} without the requirement $\alpha\ll \epsilon$. Moreover, the $\epsilon-$dependence does not matter in front of the source terms since we can use a bootstrap argument to lift $\sigma$ from $0$ to $1$ for any fixed $\epsilon$. As for the subsequent limit $\epsilon\to 0$, the sources terms are observed to automatically satisfy $\int_{0}^{\eta_{1}}\int_{\R^3}(\mathcal{S}_{1}+\sqrt{\mu}\mathcal{S}_{2})\,{\rm d}v{\rm d}\eta=0$ so that the conversation of mass holds. On the other hand, we would emphasize that it is essential to introduce the penalty terms $\epsilon g_{1}$ and $\epsilon g_{2}$ in \eqref{I2-0}--\eqref{I2}; otherwise, there must be a $\sigma-$dependent coefficient in front of the source terms (see from \eqref{I6-1}), which will prevent us from establishing uniform estimates for any $\sigma\in [0,1]$ and taking the limit $\sigma\to 1$.

\medskip
(\rmnum{3}) Since the domain under consideration has nonzero curvature, a geometric correction term $\frac{1}{1-\eta}$ appears in equations \eqref{2.3}--\eqref{2.4}. On one hand, the presence of such geometric terms complicates the characteristic lines compared to those in \cite{DL-2022, DLY-2021}, necessitating more extensive analysis in establishing the {\it a priori} $L_{\eta,v}^{\infty}$ estimates and the trace lemma (Lemma \ref{tt}). On the other hand, the geometric effect must be carefully accounted for when deriving the macroscopic dissipation estimate for $\|\mathbf{P}g\|_{L^2}$, particularly in the non-steady case. In fact, to control the temporal derivative term $\langle g,\partial_{t}\Psi\rangle$ that arises after applying a test function, we have to construct suitable test functions by solving a second-order ODE with variable coefficients with zero Dirichlet boundary conditions, and then analyze carefully the weak formulation the equations of $\partial_{t}\mathbf{P}g$, see Subsection \ref{sec7.3} for details.

\medskip
The rest of this paper is organized as follows. In Section \ref{sec2}, we reformulate the linearized problem and construct the corresponding linearized approximate systems. The {\it a priori} weighted $L^{\infty}$ estimates and $L^2$ estimates for linearized approximate systems will be given in Sections \ref{sec3} and \ref{sec4} respectively. In Section \ref{sec5}, with the help of {\it a priori} weighted $L^{\infty}$ estimates, we will apply the Leray-Schauder fixed point theorem to establish the existence and uniqueness of the solutions to the linearized problem. Then, the existence and uniqueness of solutions to the nonlinear problem and the proof of Theorem \ref{thm1} will be given in Section \ref{sec6}. In Section \ref{sec7}, we will consider the time-evolutionary problem. The local-in-time existence of non-negative solution will be established through a well-designed iteration scheme in Section \ref{sec7.1}, and the {\it a priori} $L^{\infty}\cap L^2$ uniform estimate will be shown in Sections \ref{sec7.2}--\ref{sec7.3}. Then, we will establish the exponential time asymptotic stability of the stationary solution so that the non-negativity of steady solution follows. In Appendix \ref{AppendixA}, some useful estimates, such as the smallness of $\mathcal{K}$ for large velocity, the trace lemma and estimates of iterated integrals, will be proved, which is used extensively in establishing the {\it a priori} estimates.

\smallskip

{\it Notations}. We list some notations and norms used in the paper. Throughout this paper, $C$ denotes some generic positive (generally large) constant independent of any parameters, and $C_{a}$ denotes the constant depending on $a$. $D\lesssim E$ means that there is a generic constant $C>0$ such that $D\leq CE$. ${\bf 1}_{A}$ denotes the characteristic function on the set $A$. Sometimes without confusion, we will use $\|\cdot\|_{L^2}$ to denote either the $L^{2}((0,\eta_{1})\times \R^3)-$norm or the $L^2(0,\eta_{1})-$norm or the $L^2(\R^3)-$norm, and use $\|\cdot\|_{L^{\infty}}$ to denote either the $L^{\infty}([0,\eta_{1}]\times \R^3)-$norm or the $L^{\infty}(\R^3)-$norm. Moreover, $\langle\,\cdot\, \rangle$ denotes the $L^2$ inner product in $(0,\eta_{1})\times \R^3$ with the $L^2-$norm $\|\cdot\|_{L^2}$, and $\|\cdot\|_{L_{w}^2}$ denotes the weighted $L^2-$norm with weight $w$. We denote by
\begin{align*}
&\gamma_{+}=\{(0,v)\,|\,v_{\eta}<0\}\cup \{(\eta_{1},v)\,|\,v_{\eta}>0\},\quad \text{outgoing boundary},\\
&\gamma_{-}=\{(0,v)\,|\,v_{\eta}>0\}\cup \{(\eta_{1},v)\,|\,v_{\eta}<0\},\quad \text{incoming boundary},\\
&\gamma_{0}=\{(0,v)\,|\,v_{\eta}=0\}\cup \{(\eta_{1},v)\,|\,v_{\eta}=0\},\quad\,\,\text{grazing boundary}.
\end{align*}
Furthermore, $|f|_{L^2(\gamma_{\pm})}$ and $|f|_{L^{\infty}(\gamma_{\pm})}$ represent the $L^2-$norm and $L^{\infty}-$norm at the boundaries $(\eta,v)\in \gamma_{\pm}$, respectively.

\section{The linearized approximate systems}\label{sec2}
	In this section, we consider the following linearized system for $(f_{1},f_{2})$:
	\begin{equation}\label{3.1}
		\left\{\begin{aligned}
			&v_{\eta}\frac{\partial f_{1}}{\partial \eta}-\frac{1}{1-\eta}\Big(v_{\phi}^2\frac{\partial f_1}{\partial v_{\eta}}-v_{\eta}v_{\phi}\frac{\partial f_{1}}{\partial v_{\phi}}\Big)-\frac{\alpha}{\eta_{1}} v_{\eta}\frac{\partial f_{1}}{\partial v_{\phi}}+\frac{\alpha}{2\eta_{1}}v_{\eta}v_{\phi}\sqrt{\mu}f_{2}+\nu f_{1}=(1-\chi_{M})\mathcal{K}f_{1}+\mathcal{S},\\
			&f_{1}(0,v_{\eta},v_{\phi},v_{z})\vert_{v_{\eta}>0}=0,\quad f_{1}(\eta_{1},v_{\eta},v_{\phi},v_{z})\vert_{v_{\eta}<0}=0,
		\end{aligned}\right.
	\end{equation}
	and
	\begin{equation}\label{3.2}
		\left\{\begin{aligned}
			&v_{\eta}\frac{\partial f_{2}}{\partial \eta}-\frac{1}{1-\eta}\Big(v_{\phi}^2\frac{\partial f_2}{\partial v_{\eta}}-v_{\eta}v_{\phi}\frac{\partial f_2}{\partial v_{\phi}}\Big)-\frac{\alpha}{\eta_{1}} v_{\eta}\frac{\partial f_2}{\partial v_{\phi}}+\nu f_2=Kf_2+\chi_{M}\mu^{-\frac{1}{2}}\mathcal{K}f_1,\\
			&f_2(0,v_{\eta},v_{\phi},v_{z})\vert_{v_{\eta}>0}=\sqrt{2\pi\mu}\int_{u_{\eta}<0}|u_{\eta}|(f_{1}+\sqrt{\mu}f_2)(0,u)\,\mathrm{d}u,\\
			&f_2(\eta_{1},v_{\eta},v_{\phi},v_{z})\vert_{v_{\eta}<0}=\sqrt{2\pi\mu}\int_{u_{\eta}>0}|u_{\eta}|(f_{1}+\sqrt{\mu}f_2)(\eta_{1},u)\,\mathrm{d}u,
		\end{aligned}\right.
	\end{equation}
where $\int_{0}^{\eta_{1}}\int_{\R^3}(1-\eta)\mathcal{S}\,{\rm d}v\,{\rm d}\eta=0$ and $\|\nu^{-1}w^{\ell}\mathcal{S}\|_{L^{\infty}}<\infty$. We denote $\mathcal{I}:=(0,\eta_{1})$ and $\gamma=\partial \mathcal{I} \times \R^3=\{0,\eta_{1}\}\times \R^3$ as the phase boundary of $\mathcal{I}\times \R^3$. We divide $\gamma$ into three disjoint parts, the outgoing boundary $\gamma_{+}$, the incoming boundary $\gamma_{-}$, and the singular boundary $\gamma_{0}$ for grazing velocities:
	\begin{align*}
		&\gamma_{+}=\{(\eta,v)\in \partial\mathcal{I}\times \R^3: \vec{n}(\eta)\cdot v>0\}=\{(0,v)\,|\,v_{\eta}<0\}\cup \{(\eta_{1},v)\,|\,v_{\eta}>0\},\\
		&\gamma_{-}=\{(\eta,v)\in \partial\mathcal{I}\times \R^3: \vec{n}(\eta)\cdot v<0\}=\{(0,v)\,|\,v_{\eta}>0\}\cup \{(\eta_{1},v)\,|\,v_{\eta}<0\},\\
		&\gamma_{0}=\{(\eta,v)\in \partial \mathcal{I}\times \R^3: \vec{n}(\eta)\cdot v=0\}=\{(0,v)\,|\,v_{\eta}=0\}\cup \{(\eta_{1},v)\,|\,v_{\eta}=0\},
	\end{align*}
	where $\vec{n}(\eta)$ is the outward unit normal vector on $\partial\mathcal{I}=\{0,\eta_{1}\}$. It is clear to see that $\partial \mathcal{I}=\{0,\eta_{1}\}$, $\vec{n}(0)=(-1,0,0)$ and $\vec{n}(\eta_{1})=(1,0,0)$. For later use, we denote
	\begin{equation}\label{3.3}
		P_{\gamma}f(\eta,v):=\sqrt{2\pi\mu}\int_{u\cdot \vec{n}(\eta)>0}|u_{\eta}|(\sqrt{\mu}f)
		(\eta,u)\,\mathrm{d}u\qquad \text{for }(\eta,v)\in \gamma_-.
	\end{equation}
	Let
	$$
	h_{1}(\eta,v)=w^{\ell}(v)f_{1}(\eta,v),\qquad h_{2}(\eta,v)=w^{\ell}(v)f_{2}(\eta,v),
	$$
	then we obtain from \eqref{3.1}--\eqref{3.2} that
	\begin{equation}\label{3.6}
		\left\{\begin{aligned}
			&v_{\eta}\frac{\partial h_{1}}{\partial \eta}-\frac{1}{1-\eta}\Big(v_{\phi}^2\frac{\partial h_1}{\partial v_{\eta}}-v_{\eta}v_{\phi}\frac{\partial h_{1}}{\partial v_{\phi}}\Big)-\frac{\alpha}{\eta_{1}} v_{\eta}\frac{\partial h_{1}}{\partial v_{\phi}}+\frac{\alpha}{2\eta_{1}}v_{\eta}v_{\phi}\sqrt{\mu}h_{2}+\nu h_{1}+\frac{\alpha}{\eta_{1}}\frac{\ell v_{\eta}v_{\phi}}{1+|v|^2}h_{1}\\
			&\quad =(1-\chi_{M})\mathcal{K}_{w}h_{1}+w^{\ell}\mathcal{S},\\
			&h_{1}(0,v_{\eta},v_{\phi},v_{z})\vert_{v_{\eta}>0}=0,\quad h_{1}(\eta_{1},v_{\eta},v_{\phi},v_{z})\vert_{v_{\eta}<0}=0,
		\end{aligned}\right.
	\end{equation}
	and
	\begin{equation}\label{3.7}
		\left\{\begin{aligned}
			&v_{\eta}\frac{\partial h_{2}}{\partial \eta}-\frac{1}{1-\eta}\Big(v_{\phi}^2\frac{\partial h_2}{\partial v_{\eta}}-v_{\eta}v_{\phi}\frac{\partial h_2}{\partial v_{\phi}}\Big)-\frac{\alpha}{\eta_{1}} v_{\eta}\frac{\partial h_2}{\partial v_{\phi}}+\nu h_2+\frac{\alpha}{\eta_{1}}\frac{\ell v_{\eta}v_{\phi}}{1+|v|^2}h_{2}\\
			&\quad =K_{w}h_2+\chi_{M}\mu^{-\frac{1}{2}}\mathcal{K}_{w}h_1,\\
			&h_2(0,v_{\eta},v_{\phi},v_{z})\vert_{v_{\eta}>0}=\frac{1}{\tilde{w}(v)}\int_{u_{\eta}<0}\big(\frac{h_{1}}{\sqrt{\mu}}+h_2\big)(0,u)\tilde{w}(u)\,\mathrm{d}\sigma,\\
			&h_2(\eta_{1},v_{\eta},v_{\phi},v_{z})\vert_{v_{\eta}<0}=\frac{1}{\tilde{w}(v)}\int_{u_{\eta}>0}\big(\frac{h_{1}}{\sqrt{\mu}}+h_2\big)(\eta_{1},u)\tilde{w}(u)\,\mathrm{d}\sigma,
		\end{aligned}\right.
	\end{equation}
	where
	\begin{equation*}
		\tilde{w}(v)=(\sqrt{2\pi\mu}w^{\ell})^{-1}\,\,\text{ and }\,\,\mathrm{d}\sigma=\sqrt{2\pi}\mu(v)|v_{\eta}|\,\mathrm{d}v_{\eta}\mathrm{d}v_{\phi}{\rm d}v_{z}.
	\end{equation*}
	Hereafter, we denote $K_{w}h=w^{\ell}K(\frac{h}{w^{\ell}})$ and $\mathcal{K}_{w}h=w^{\ell}\mathcal{K}(\frac{h}{w^{\ell}})$, that is,
	\begin{align}
		&K_{w}h(v)=\int_{\R^3}k_{w}(v,u)h(u)\,\mathrm{d}u\quad \text{with }k_{w}(v,u)=w^{\ell}(v)k(v,u)w^{-\ell}(u),\label{3.9}\\
		&\mathcal{K}_{w}h(v)=\int_{\R^3}\kappa_{w}(v,u)h(u)\,\mathrm{d}u\quad \text{with }\kappa_{w}(v,u)=w^{\ell}(v)\mu^{\frac{1}{2}}(v)k(v,u)\mu^{-\frac{1}{2}}(u)w^{-\ell}(u).\label{3.10}
	\end{align}
 Some properties for $K_{w}$ and $\mathcal{K}_{w}$ are listed in Lemma \ref{K}.
 
Motivated by \cite{DL-2022,DLY-2021,Huang-Wang-2022}, to show the existence of \eqref{3.7}, we shall consider the following coupled approximate systems
	\begin{equation}\label{3.6-1}
		\left\{\begin{aligned}
			&\epsilon h_{1}+v_{\eta}\frac{\partial h_{1}}{\partial \eta}-\frac{1}{1-\eta}\Big(v_{\phi}^2\frac{\partial h_1}{\partial v_{\eta}}-v_{\eta}v_{\phi}\frac{\partial h_{1}}{\partial v_{\phi}}\Big)-\frac{\alpha}{\eta_{1}} v_{\eta}\frac{\partial h_{1}}{\partial v_{\phi}}+\frac{\alpha}{2\eta_{1}}v_{\eta}v_{\phi}\sqrt{\mu}h_{2}+\nu h_{1}\\
			&\quad =-\frac{\alpha}{\eta_{1}}\frac{\ell v_{\eta}v_{\phi}}{1+|v|^2}h_{1}\sigma(1-\chi_{M})\mathcal{K}_{w}h_{1}+w^{\ell}\mathcal{S}_{1},\\
			&h_{1}(0,v)\vert_{v_{\eta}>0}=0,\quad h_{1}(\eta_{1},v)\vert_{v_{\eta}<0}=0,\\
		\end{aligned}\right.
	\end{equation}
	and
	\begin{equation}\label{3.7-1}
		\left\{\begin{aligned}
			&\epsilon h_{2}+v_{\eta}\frac{\partial h_{2}}{\partial \eta}-\frac{1}{1-\eta}\Big(v_{\phi}^2\frac{\partial h_2}{\partial v_{\eta}}-v_{\eta}v_{\phi}\frac{\partial h_2}{\partial v_{\phi}}\Big)-\frac{\alpha}{\eta_{1}} v_{\eta}\frac{\partial h_2}{\partial v_{\phi}}+\nu h_2+\frac{\alpha}{\eta_{1}}\frac{\ell v_{\eta}v_{\phi}}{1+|v|^2}h_{2}\\
			&\quad =\sigma K_{w}h_2+\sigma \chi_{M}\mu^{-\frac{1}{2}}\mathcal{K}_{w}h_1+w^{\ell}\mathcal{S}_{2},\\
			&h_2(0,v)\vert_{v_{\eta}>0}=\frac{1}{\tilde{w}(v)}\int_{u_{\eta}<0}h_2(0,u)\tilde{w}(u)\,\mathrm{d}\sigma+\frac{\sigma}{\tilde{w}(v)}\int_{u_{\eta<0}}\frac{h_{1}}{\sqrt{\mu}}(0,u)\tilde{w}(u)\,{\rm d}\sigma+w^{\ell}\mathcal{S}_{2,b}^{(1)},\\
			&h_2(\eta_{1},v)\vert_{v_{\eta}<0}=\frac{1}{\tilde{w}(v)}\int_{u_{\eta}>0}h_2(\eta_{1},u)\tilde{w}(u)\,\mathrm{d}\sigma+\frac{\sigma}{\tilde{w}(v)}\int_{u_{\eta}>0}\frac{h_{1}}{\sqrt{\mu}}(\eta_{1},u)\tilde{w}(u)\,{\rm d}\sigma+w^{\ell}\mathcal{S}_{2,b}^{(2)},
		\end{aligned}\right.
	\end{equation}
with parameters $\epsilon,\sigma\in [0,1]$ and
$$
\|\nu^{-1}w^{\ell}(\mathcal{S}_{1},\mathcal{S}_{2})\|_{L^{\infty}}<\infty,\,\,|w^{\ell}(\mathcal{S}_{2,b}^{(1)},\mathcal{S}_{2,b}^{(2)})|_{L^{\infty}(\gamma_{-})}<\infty.
$$
We shall establish the {\it a priori} $L^{\infty}$ estimate for the solution of the coupled systems \eqref{3.6-1}--\eqref{3.7-1} for any given $\sigma$ based on $L^{\infty}\cap L^2$ framework developed in \cite{Guo-2010} in Sections \ref{sec3} and \ref{sec4}. Then we will  use the {\it a priori} estimates to take the limits $\sigma\to 1$ and $\epsilon\to 0$ for establishing the existence of the linearized approximate systems \eqref{3.1}--\eqref{3.2} in Section \ref{sec5}.

\section{{\it A priori} $L^{\infty}$ estimate for approximate systems \texorpdfstring{\eqref{3.6-1}--\eqref{3.7-1}}{（2.9)--(2.10)}}\label{sec3}

For the coupled systems \eqref{3.6-1}--\eqref{3.7-1}, we are going to establish the {\it a priori} $L^{\infty}$ estimate by employing the characteristic method. In the following, we introduce a uniform parameter $t\in \R$ and regard $[h_{1},h_{2}](\eta,v)=[h_{1},h_{2}](t,\eta,v)$.
	\begin{definition}
		Given $(t,\eta,v)=(t,\eta,v_{\eta},v_{\phi},v_{z})$, let $$
		[X(s),V_{\eta}(s),V_{\phi}(s),V_{z}(s)]=:[X(s;t,\eta,v),V_{\eta}(s;t,\eta,v),V_{\phi}(s;t,\eta,v),V_{z}(s;t,\eta,v)]
		$$
		be the backward characteristics passing through $(t,\eta,v)$ for \eqref{3.6}--\eqref{3.7}, which are defined by
		\begin{equation}\label{3.11}
			\left\{
			\begin{aligned}
				&\frac{\mathrm{d}X(s)}{\mathrm{d}s}=V_{\eta}(s),\quad \frac{\rm{d}V_{z}(s)}{\rm{d} s}=0,\\
				&\frac{\mathrm{d}V_{\eta}(s)}{\mathrm{d}s}=-\frac{1}{1-X(s)}V_{\phi}^2(s),\\
				&\frac{\mathrm{d}V_{\phi}(s)}{\mathrm{d}s}=\frac{1}{1-X(s)}V_{\eta}(s)V_{\phi}(s)-\frac{\alpha}{\eta_{1}}V_{\eta}(s),\\
				&[X(s),V_{\eta}(s),V_{\phi}(s),V_{z}(s)]\vert_{s=t}=[\eta,v_{\eta},v_{\phi},v_{z}].
			\end{aligned}
			\right.
		\end{equation}
	\end{definition}
	
\smallskip
It follows from $\eqref{3.11}_{1}$ and $\eqref{3.11}_{4}$ that
\begin{align*}
V_{z}(s)\equiv v_{z}.
\end{align*}
Moreover, we obtain from $\eqref{3.11}_{1}$--$\eqref{3.11}_{3}$ that
	$$
	\frac{\mathrm{d}((1-X(s))V_{\phi}(s))}{\mathrm{d}s}=-\frac{\alpha}{\eta_{1}}(1-X(s))\frac{\mathrm{d}X(s)}{\mathrm{d}s},
	$$
	which, together with $\eqref{3.11}_{4}$, yields that
	\begin{equation}\label{3.12}
		V_{\phi}(s)=V_{\phi}(s; t,\eta, v)=\frac{1-\eta}{1-X(s)}v_{\phi}+\frac{\alpha}{2\eta_{1}}(1-X(s))-\frac{\alpha}{2\eta_{1}}\frac{(1-\eta)^2}{1-X(s)}.
	\end{equation}
	Then it follows from $\eqref{3.11}_{1}$--$\eqref{3.11}_{2}$ that
	$$
	\frac{\mathrm{d}V_{\eta}^2(s)}{\mathrm{d}X(s)}=-\frac{2}{1-X(s)}\Big[\frac{1-\eta}{1-X(s)}v_{\phi}+\frac{\alpha}{2\eta_{1}}(1-X(s))-\frac{\alpha}{2\eta_{1}}\frac{(1-\eta)^2}{1-X(s)}\Big]^2,
	$$
	which implies that
	\begin{equation}\label{3.13}
		V_{\eta}(s)^2=V_{\eta}^2(s;t,\eta,v)=v_{\eta}^2+\int_{X(s)}^{\eta}\frac{2}{1-y}\Big[\frac{1-\eta}{1-y}v_{\phi}+\frac{\alpha}{2\eta_{1}}(1-y)-\frac{\alpha}{2\eta_{1}}\frac{(1-\eta)^2}{1-y}\Big]^2\,\mathrm{d}y.
	\end{equation}
Noting that for any $\eta\in [0,\eta_{1}]$ and $v=(v_{\eta}, v_{\phi},v_{z})\in \R^3$, 
	\begin{equation*}
		v_{\eta}^2<\int_{\eta}^{1}\frac{2}{1-y}\Big[\frac{1-\eta}{1-y}v_{\phi}+\frac{\alpha}{2\eta_{1}}(1-y)-\frac{\alpha}{2\eta_{1}}\frac{(1-\eta)^2}{1-y}\Big]^2\,\mathrm{d}y=+\infty,
	\end{equation*}
there must exist a position $\eta_{+}(\eta,v)\in [\eta,1)$ such that 
	$$
	v_{\eta}^2=\int_{\eta}^{\eta_{+}}\frac{2}{1-y}\Big[\frac{1-\eta}{1-y}v_{\phi}+\frac{\alpha}{2\eta_{1}}(1-y)-\frac{\alpha}{2\eta_{1}}\frac{(1-\eta)^2}{1-y}\Big]^2\,\mathrm{d}y.
	$$
 
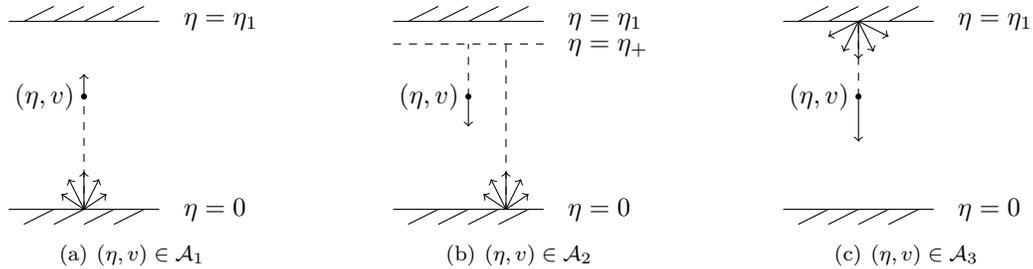
\begin{figure}[H]
				\centering
				\subfigure[ $(\eta,v)\in \mathcal{A}_{1}$]{
					\begin{tikzpicture}
						\draw (0,0.5)--(2,0.5);
						\draw (0,3)--(2,3);
						
						\draw (0.2,3)--(0.6,3.2);
						\draw (0.6,3)--(1.0,3.2);
						\draw (1.0,3)--(1.4,3.2);
						\draw (1.4,3)--(1.8,3.2);
						\draw (0.2,0.3)--(0.6,0.5);
						\draw (0.6,0.3)--(1.0,0.5);
						\draw (1.0,0.3)--(1.4,0.5);
						\draw (1.4,0.3)--(1.8,0.5);
						\draw[dashed] (1,0.5)--(1,2);
						\fill[black] (1,2)circle (.04);
						\node [left] at (1,2) {$(\eta,v)$};
						\node [right] at (2.2,0.5) {$\eta=0$};
						\node [right] at (2.2,3) {$\eta=\eta_{1}$};						
				\draw[black,->] (1,2)--(1,2.3);
				\draw[->] (1,0.5)--(1,1);
				\draw[->] (1,0.5)--(1.2,0.9);
				\draw[->] (1,0.5)--(0.8,0.9);
				\draw[->] (1,0.5)--(0.7,0.7);
				\draw[->] (1,0.5)--(1.3,0.7);
				
				\end{tikzpicture}}
				\qquad\qquad
			     \subfigure[$(\eta,v)\in \mathcal{A}_{2}$]{
					\begin{tikzpicture}
						\draw (0,0.5)--(2,0.5);
						\draw (0,3)--(2,3);
						\draw (0.2,3)--(0.6,3.2);
						\draw (0.6,3)--(1.0,3.2);
						\draw (1.0,3)--(1.4,3.2);
						\draw (1.4,3)--(1.8,3.2);
						\draw (0.2,0.3)--(0.6,0.5);
						\draw (0.6,0.3)--(1.0,0.5);
						\draw (1.0,0.3)--(1.4,0.5);
						\draw (1.4,0.3)--(1.8,0.5);
						
			    \draw[black,->] (1,2)--(1,1.6);
						\draw[dashed] (1,2)--(1,2.7);
						
				\draw[dashed] (0,2.7)--(2,2.7);
				
				\draw[dashed] (1.5,2.7)--(1.5,0.5);

				\draw[black,->] (1.5,0.5)--(1.5,1);
						
				\draw[black,->] (1.5,0.5)--(1.3,0.9);
				
				\draw[black,->] (1.5,0.5)--(1.7,0.9);
				
				\draw[black,->] (1.5,0.5)--(1.2,0.7);
				
				\draw[black,->] (1.5,0.5)--(1.8,0.7);
						\fill[black] (1,2)circle (.04);
						\node [left] at (1,2) {$(\eta,v)$};
			    \node [right] at (2.2,0.5) {$\eta=0$};
						\node [right] at (2.2,3) {$\eta=\eta_{1}$};
			    \node [right] at (2.2,2.65) {$\eta=\eta_{+}$};	\end{tikzpicture}}
			\qquad\qquad	\subfigure[$(\eta,v)\in \mathcal{A}_{3}$]{
					\begin{tikzpicture}
						\draw (0,0.5)--(2,0.5);
						\draw (0,3)--(2,3);
						\draw (0.2,3)--(0.6,3.2);
						\draw (0.6,3)--(1.0,3.2);
						\draw (1.0,3)--(1.4,3.2);
						\draw (1.4,3)--(1.8,3.2);
						\draw (0.2,0.3)--(0.6,0.5);
						\draw (0.6,0.3)--(1.0,0.5);
						\draw (1.0,0.3)--(1.4,0.5);
						\draw (1.4,0.3)--(1.8,0.5);
						
			    \draw[black,->] (1,2)--(1,1.4);
						\draw[dashed] (1,2)--(1,3);
				
				\draw[black,->] (1,3)--(1,2.5);
						
				\draw[black,->] (1,3)--(0.8,2.6);
				
				\draw[black,->] (1,3)--(1.2,2.6);
				
				\draw[black,->] (1,3)--(0.6,2.8);
				
				\draw[black,->] (1,3)--(1.4,2.8);
						\fill[black] (1,2)circle (.04);
						\node [left] at (1,2) {$(\eta,v)$};
					
				\node [right] at (2.2,3) {$\eta=\eta_{1}$};
				
				\node [right] at (2.2,0.5) {$\eta=0$};
				\end{tikzpicture}}
				\caption{The schematic of backward characteristic lines}
                \label{fig1}
			\end{figure}

	Now for each $(\eta, v)\in \bar{\mathcal{I}}\times \R^3$, we define its backward exit time $t_{ *}(t,\eta,v)\geq 0$ to be 
$$
	t_{*}(t,\eta, v)=\sup\{s\geq 0:\big(X(t-s;t,\eta, v),V_{\eta}(t-s;t,\eta,v),V_{\phi}(t-s;t,\eta,v),V_{z}(t-s;t,\eta,v)\big)\in \gamma_{-}\}.
	$$
	We also define
	$$
	x_{*}(t,\eta, v)=X(t-t_{*};t,\eta,v)\in\{0,\eta_{1}\}.
	$$
Denote
	$$
	\begin{aligned}
		&\mathcal{A}_{1}:=\{(y,u)\in \bar{\mathcal{I}}\times \R^3 :u_{\eta}\geq 0\},\\
		&\mathcal{A}_{2}:=\{(y,u)\in \bar{\mathcal{I}}\times \R^3:u_{\eta}<0,\,\eta_{+}(y, u)\leq \eta_{1}\},\\
		&\mathcal{A}_{3}:=\{(y,u)\in \bar{\mathcal{I}}\times \R^3 :u_{\eta}<0,\,\eta_{+}(y, u)> \eta_{1}\}.
	\end{aligned}
	$$
Then,
\begin{align*}
\begin{aligned}
	&V_{\eta}(s)=\sqrt{v_{\eta}^2+\int_{X(s)}^{\eta}\frac{2}{1-y}\Big[\frac{1-\eta}{1-y}v_{\phi}+\frac{\alpha}{2\eta_{1}}(1-y)-\frac{\alpha}{2\eta_{1}}\frac{(1-\eta)^2}{1-y}\Big]^2\,\mathrm{d}y},\\
		&X(s)=\eta-\int_{s}^{t}V_{\eta}(\tau)\,{\rm d}\tau,\quad x_{*}(t,\eta,v)=0,
	\end{aligned}\qquad\,\,\,\,\text{if }(\eta,v)\in \mathcal{A}_{1},
\end{align*}
and
\begin{align*}
\begin{aligned}
&V_{\eta}(s)=-\sqrt{v_{\eta}^2+\int_{X(s)}^{\eta}\frac{2}{1-y}\Big[\frac{1-\eta}{1-y}v_{\phi}+\frac{\alpha}{2\eta_{1}}(1-y)-\frac{\alpha}{2\eta_{1}}\frac{(1-\eta)^2}{1-y}\Big]^2\,\mathrm{d}y},\\
&X(s)=\eta-\int_{s}^{t}V_{\eta}(\tau)\,{\rm d}\tau,\quad x_{*}(t,\eta,v)=\eta_{1},
	\end{aligned}\qquad\text{if }(\eta,v)\in \mathcal{A}_{3}.
\end{align*}
Finally, if $(\eta,v)\in \mathcal{A}_{2}$, then there exists a moment $t_{*,0}\in [t-t_{*},t]$ such that $X(t_{*,0})=\eta_{+}(\eta,v)$ and $V_{\eta}(t_{*,0})=0$. Thus,
\begin{align*}
	&
	\begin{aligned}
		&V_{\eta}(s)=-\mathbf{1}_{\{[t_{*,0},t]\}}(s)\sqrt{v_{\eta}^2+\int_{X(s)}^{\eta}\frac{2}{1-y}\Big[\frac{1-\eta}{1-y}v_{\phi}+\frac{\alpha}{2\eta_{1}}(1-y)-\frac{\alpha}{2\eta_{1}}\frac{(1-\eta)^2}{1-y}\Big]^2\,\mathrm{d}y}\\
		&\qquad\quad\,\, +\mathbf{1}_{\{[t_{*},t_{*,0})\}}(s)\sqrt{\int_{X(s)}^{\eta_{+}}\frac{2}{1-y}\Big[\frac{1-\eta}{1-y}v_{\phi}+\frac{\alpha}{2\eta_{1}}(1-y)-\frac{\alpha}{2\eta_{1}}\frac{(1-\eta)^2}{1-y}\Big]^2\,\mathrm{d}y},\\
		&X(s)=\eta-\int_{s}^{t}V_{\eta}(\tau)\,{\rm d}\tau,\quad x_{*}(t,\eta,v)=0.
	\end{aligned}\quad\text{if }(\eta,v)\in \mathcal{A}_{2}.
\end{align*}
See Figure \ref{fig1}.

We denote
	$$
	(t_{0},x_{0},v_{0}):=(t,\eta,v),
	$$
	and define
	\begin{equation*}
		(t_{k+1},x_{k+1},v_{k+1})=(t_k-t_{*}(t_{k},x_{k},v_{k}),x_{*}(t_{k},x_{k},v_{k}),v_{k+1})
	\end{equation*}
   for $k\geq 0$, where
	\begin{equation*}
		v_{k+1}\in \mathcal{V}_{k+1}:=\{v_{k+1}\in \R^3:v_{k+1}\cdot \vec{n}(x_{k+1})>0\}.
	\end{equation*}
	Then the back-time cycle can be defined as
	
\begin{equation}\label{3.16}
		\left\{
		\begin{aligned}
			&\mathscr{X}(s)=\sum_{k}{\bf 1}_{[t_{k+1},t_{k})}(s)X(s;t_{k},x_{k},v_{k}),\\
			& \mathscr{V}_{\eta}(s)=\sum_{k}{\bf 1}_{[t_{k+1},t_{k})}(s)\Big[({\bf 1}_{\mathcal{A}_{1}}-{\bf 1}_{\mathcal{A}_{3}})(x_{k},v_{k})|V_{\eta}(s;t_{k},x_{k},v_{k})|\\
			&\qquad \qquad \qquad \qquad \qquad \quad
			-{\bf 1}_{\mathcal{A}_{2}}(x_{k},v_{k}){\bf 1}_{[t_{k,0},t_{k})}(s)|V_{\eta}(s;t_{k},x_{k},v_{k})|\\
			&\qquad \qquad \qquad \qquad \qquad \quad+{\bf 1}_{\mathcal{A}_{2}}(x_{k},v_{k}){\bf 1}_{[t_{k+1},t_{k,0})}(s)|V_{\eta}(s;t_{k},x_{k},v_{k})|\Big],\\
			&\mathscr{V}_{\phi}(s)=\sum_{k}{\bf 1}_{[t_{k+1},t_{k})}(s)V_{\phi}(s;t_{k},x_{k},v_{k}),\quad \mathscr{V}_{z}(s)= \sum_{k}{\bf 1}_{[t_{k+1},t_{k})}(s)V_{z}(s;t_{k},x_{k},v_{k}),
		\end{aligned}
		\right.
	\end{equation}
	where $t_{k,0}\in [t_{k+1},t_{k})$ is the moment when $V_{\eta}(t_{k,0};t_{k},x_{k},v_{k})=0$ for $(x_{k},v_{k})\in \mathcal{A}_{2}$.	
	
	\begin{lemma}\label{lem2.3}
		{\rm {1)}} If $x_{k}=0$, $v_{k}\in \mathcal{V}_{k}$ and $\eta_{+}(x_{k},v_{k})\leq \eta_{1}$, then it holds that
		\begin{align}\label{V1}
			\mathscr{V}_{\phi}(t_{k+1})=v_{k,\phi},\quad \mathscr{V}_{\eta}(t_{k+1})=-v_{k,\eta}.
		\end{align}
		{\rm {2)}} If $x_{k}=0$, $v_{k}\in \mathcal{V}_{k}$ and $\eta_{+}(x_{k},v_{k})>\eta_{1}$, then it holds that
		\begin{align}\label{V2}
			\frac{1+|\mathscr{V}(t_{k+1})|^2}{1+|v_{k}|^2}\leq 1+C\alpha.
		\end{align}
		{\rm {3)}} If $x_{k}=\eta_{1}$ and $v_{k}\in \mathcal{V}_{k}$, then it holds that
		\begin{align}\label{V3}
			1\leq \frac{1+|\mathscr{V}(t_{k+1})|^2}{1+|v_{k}|^2}\leq 1+C\alpha.
		\end{align}
	\end{lemma}
	
\noindent\textbf{Proof.} 
1. If $x_{k}=0$, $v_{k}\in \mathcal{V}_{k}$ and $\eta_{+}(x_{k},v_{k})\leq \eta_{1}$, then
	$$
	v_{k,\eta}<0,\quad x_{k+1}=X(t_{k+1})=0.
	$$
It follows from \eqref{3.12} and \eqref{3.16} that
\begin{align*}
	&\mathscr{V}_{\phi}(t_{k+1})=V_{\phi}(t_{k+1};t_{k},x_{k},v_{k})=\frac{1}{1-X(t_{k+1})}v_{k,\phi}+\frac{\alpha}{2\eta_{1}}(1-X(t_{k+1}))-\frac{\alpha}{2\eta_{1}}\frac{1}{1-X(t_{k+1})}=v_{k,\phi},\\
	&\mathscr{V}_{\eta}(t_{k+1})=|V_{\eta}(t_{k+1};t_{k},x_{k},v_{k})|
	\\
	&\qquad \qquad =\sqrt{\int_{X(t_{k+1})}^{\eta_{+}}\frac{2}{1-y}\Big[\frac{1-\eta}{1-y}v_{\phi}+\frac{\alpha}{2\eta_{1}}(1-y)-\frac{\alpha}{2\eta_{1}}\frac{(1-\eta)^2}{1-y}\Big]^2\,\mathrm{d}y}\\
	&\qquad\qquad=\sqrt{\int_{0}^{\eta_{+}}\frac{2}{1-y}\Big[\frac{1-\eta}{1-y}v_{\phi}+\frac{\alpha}{2\eta_{1}}(1-y)-\frac{\alpha}{2\eta_{1}}\frac{(1-\eta)^2}{1-y}\Big]^2\,\mathrm{d}y}=-v_{k,\eta},
\end{align*}
which concludes \eqref{V1}.

2. If $x_{k}=0$, $v_{k}\in \mathcal{V}_{k}$ and $\eta_{+}(0,v_{k})>\eta_{1}$, then $x_{k+1}=x_{*}(t_{k},x_{k},v_{k})=\eta_{1}$, and
	\begin{align*}
		&\mathscr{V}_{\phi}(t_{k+1})=V_{\phi}(t_{k+1};t_{k},0,v_{k})=\frac{1}{1-\eta_{1}}v_{k,\phi}+\frac{\alpha}{2\eta_{1}}(1-\eta_{1})-\frac{\alpha}{2\eta_{1}}\frac{1}{1-\eta_{1}},\\ &|\mathscr{V}_{\eta}(t_{k+1})|^2=|V_{\eta}(t_{k+1};t_{k},0,v_{k})|^2=v_{k,\eta}^2-\int_{0}^{\eta_{1}}\frac{2}{1-y}\Big[\frac{1}{1-y}v_{k,\phi}+\frac{\alpha}{2\eta_{1}}(1-y)-\frac{\alpha}{2\eta_{1}}\frac{1}{1-y}\Big]^2\,{\rm d}y\nonumber\\
		&\qquad \qquad\,\,\,\,\,  =v_{k,\eta}^2+v_{k,\phi}^2-\frac{1}{(1-\eta_{1})^2}v_{k,\phi}^2+\frac{\alpha}{\eta_{1}}v_{k,\phi}\big[2\ln(1-\eta_{1})+\frac{1}{(1-\eta_{1})^2}-1\big]\nonumber\\
		\
		&\qquad \qquad \qquad -\frac{\alpha^2}{4\eta_{1}^2}\big[2-(1-\eta_{1})^2-\frac{1}{(1-\eta_{1})^2}+4\ln(1-\eta_{1})\big].
	\end{align*}
	Then a direct calculation shows that
	$$
	(1-C\alpha)(1+|v_{k}|^2)\leq 1+|\mathscr{V}(t_{k+1})|^2\leq (1+C\alpha)(1+|v_{k}|^2),
	$$
which yields \eqref{V2}. The proof of \eqref{V3} is similar. Hence, the proof of Lemma \ref{lem2.3} is completed. $\hfill\square$

\medskip
We define the iterated integral
	\begin{equation*}
		\int_{\prod_{j=1}^{k-1}\mathcal{V}_{j}}\prod_{j=1}^{k-1}\,{\rm d}\sigma_{j}=\int_{\mathcal{V}_{1}}\,{\rm d}\sigma_{1}\int_{\mathcal{V}_{2}}\,{\rm d}\sigma_{2}\cdots\int_{\mathcal{V}_{k-1}}\,{\rm}d\sigma_{k-1},
	\end{equation*}
	where ${\rm d}\sigma_{j}=\sqrt{2\pi}\mu(v_{j})|v_{j,\eta}|\,{\rm d}v_{j}$ for $j=1,2,\cdots, k-1$.
	It is direct to know that ${\rm d}\sigma_{j}$ is a probability measure on $\mathcal{V}_{j}$, i.e.,
	\begin{align*}
	    \int_{\mathcal{V}_{j}}\,{\rm d}\sigma_{j}=1.
	\end{align*}

	\begin{lemma}\label{lem2.2}
   For any given $\epsilon,\sigma\in [0,1]$. Let $(h_{1},h_{2})$ be the $L^{\infty}-$solution of the coupled systems \eqref{3.6-1}--\eqref{3.7-1}. Then
    \begin{align}
   &\|h_{1}\|_{L^{\infty}}\leq C\alpha\|h_{2}\|_{L^{\infty}}+C\|\nu^{-1}w^{\ell}\mathcal{S}_{1}\|_{L^{\infty}},\label{3.17-1}\\
   &\|h_{2}\|_{L^{\infty}}\leq C_{\ell}\|h_{1}\|_{L^{\infty}}+C_{\ell}\|\nu^{-1}w^{\ell}\mathcal{S}_{2}\|_{L^{\infty}}+C_{\ell}\|f_{2}\|_{L^2}+C_{\ell}|w^{\ell}(\mathcal{S}_{2,b}^{(1)},\mathcal{S}_{2,b}^{(2)})|_{L^{\infty}(\gamma_{-})},\label{3.17-2}
   \end{align}
   where $C_{\ell}$ is a positive constant depending on the weight $\ell$, and $C$ is a generic constant independent of $\ell$, $\sigma$ and $\epsilon$.
    \end{lemma}
    
    \noindent\textbf{Proof}. We divide the proof into four steps.
    
    Step 1. {\it Representation of $h_{1}$ along the characteristic line}. For $(\eta,v)\in (\bar{\mathcal{I}}\times \R^3\backslash (\gamma_{0}\cup \gamma_{-}))$. Recall the characteristic line $[\mathscr{X}(s),\mathscr{V}(s)]$ defined in \eqref{3.16}. We can represent $h_{1}$ as 
    \begin{align}\label{3.17-4}
h_{1}(\eta,v)&=\int_{t_1}^{t}e^{-\int_{s}^{t}\mathcal{A}(\mathscr{V}(\tau))\,{\rm d}\tau}[-\frac{\alpha}{2\eta_{1}}\mathscr{V}_{\eta}\mathscr{V}_{\phi}\sqrt{\mu}h_{2}+\sigma (1-\chi_{M})\mathcal{K}_{w}h_{1}+w^{\ell}\mathcal{S}_{1}](\mathscr{X}(s),\mathscr{V}(s))\,{\rm d}s\nonumber\\
&=:I_{1}+I_{2}+I_{3},
\end{align}
with
\begin{align*}
	\mathcal{A}(v):=\epsilon+\nu(v)+\frac{\alpha}{\eta_{1}}\frac{\ell v_{\eta}v_{\phi}}{1+|v|^2}.
\end{align*}

Step 2. {\it Estimates on $\|h_{1}\|_{L^{\infty}}$}. Noting that the collision frequency $\nu(v)\cong (1+|v|)^{\gamma}$, we have 
$$
\mathcal{A}(v)=\epsilon+\nu(v)+\frac{\alpha}{\eta_{1}}\frac{\ell v_{\eta}v_{\phi}}{1+|v|^2}\geq \frac{1}{2}\nu(v)\geq \nu_{0}>0,
$$
provided $\alpha$ small enough. Here $\nu_{0}$ is a positive constant. Then, by direct calculations, one has
\begin{align}\label{3.18-4}
|I_{1}|+|I_{3}|\leq C\alpha \|h_{2}\|_{L^{\infty}}+C\|\nu^{-1}w^{\ell}\mathcal{S}_{1}\|_{L^{\infty}}.
\end{align}
For $I_{2}$, it follows from Lemma \ref{K} that
\begin{align}\label{3.18-6-1}
  |I_{2}|\leq C\Big[\frac{1}{\ell}{\bf 1}_{\{\gamma=0\}}+\big[(1+M)^{-\gamma}+\zeta\big]{\bf 1}_{\{0<\gamma\leq 1\}}\Big]\|h_{1}\|_{L^{\infty}},   
\end{align} 
where $M=M(\ell)>0$ is sufficiently large and $\zeta=\zeta(\ell)$ is sufficiently small for $\ell$ large enough.

Combining \eqref{3.18-4}--\eqref{3.18-6-1} together, we conclude
\begin{align}\label{3.18-7-1}
\|h_{1}\|_{L^{\infty}}\leq C\big[\frac{1}{\ell}+(1+M)^{-\gamma}+\zeta\big]\|h_{1}\|_{L^{\infty}}+C\alpha\|h_{2}\|_{L^{\infty}}+C\|\nu^{-1}w^{\ell}\mathcal{S}_{1}\|_{L^{\infty}}.
\end{align}
Taking $\ell$, $M=M(\ell)$ large enough and $\zeta=\zeta(\ell)$ small enough, we obtain from \eqref{3.18-7-1} that
\begin{align*}
\|h_{1}\|_{L^{\infty}}\leq C\alpha\|h_{2}\|_{L^{\infty}}+C\|\nu^{-1}w^{\ell}\mathcal{S}_{1}\|_{L^{\infty}},
\end{align*}
where $C$ is independent of $\ell$, $\epsilon$ and $\sigma$. This concludes the proof of \eqref{3.17-1}.

Step 3. {\it Representation of $h_{2}$ along the characteristic line}. We divide it into three cases.

Case 1. For the set $\mathcal{C}_{1}: (\eta,v)\in (\bar{\mathcal{I}}\times \R^3\backslash (\gamma_{0}\cup \gamma_{-}))\cap \mathcal{A}_{1}$, it holds that
		\begin{equation}\label{3.18-3}
			h_{2}(\eta,v)=\sum\limits_{j=1}^{4}J_{j}+\sum\limits_{j=5}^{15}{\bf 1}_{\{t_{1}>0\}}J_{j},
		\end{equation}
		with
		\begin{align*}
			J_{1}&={\bf 1}_{\{t_{1}\leq 0\}}e^{-\int_{0}^{t}\mathcal{A}(\mathscr{V}(\tau))\,{\rm d}\tau}h_{2}(\mathscr{X}(0),\mathscr{V}(0)),\\
			J_{2}+J_{3}+J_{4}&=\int_{\max\{t_{1},0\}}^{t}e^{-\int_{s}^{t}\mathcal{A}(\mathscr{V}(\tau))\,{\rm d}\tau}(\sigma\chi_{M}\mu^{-\frac{1}{2}} \mathcal{K}_{w}h_{1}+\sigma K_{w}h_{2}+w^{\ell}\mathcal{S}_{2})(\mathscr{X}(s),\mathscr{V}(s))\,{\rm d}s,\\
		J_{5}&=\sigma\frac{e^{-\int_{t_{1}}^{t}\mathcal{A}(\mathscr{V}(\tau))\,{\rm d}\tau}}{\tilde{w}(\mathscr{V}(t_{1}))}\int_{\prod_{j=1}^{k-1}\mathcal{V}_{j}}\sum\limits_{l=1}^{k-1}{\bf 1}_{\{t_{l}>0\}}(\frac{h_{1}}{\sqrt{\mu}})(\mathscr{X}(t_{l}),\mathscr{V}(t_{l}))\,{\rm d}\Sigma_{l}(t_{l}),\\
		J_{6}&=e^{-\int_{t_{1}}^{t}\mathcal{A}(\mathscr{V}(\tau))\,{\rm d}\tau}(w^{\ell}\mathcal{S}_{2,b}^{(1)})(\mathscr{V}(t_{1})),\\
	J_{7}&=\frac{e^{-\int_{t_{1}}^{t}\mathcal{A}(\mathscr{V}(\tau))\,{\rm d}\tau}}{\tilde{w}(\mathscr{V}(t_{1}))}\Big\{\int_{\prod_{j=1}^{k-1}\mathcal{V}_{j}}\sum\limits_{l=1}^{k-1}{\bf 1}_{\mathcal{A}_{2}}(\mathscr{X}(t_{2\ell-1}),\mathscr{V}(t_{2l-1})){\bf 1}_{\{t_{2l}>0\}}(w^{\ell}\mathcal{S}_{2,b}^{(1)})\,{\rm d}\Sigma_{l}(t_{2l})\\
		&\qquad\qquad\qquad\qquad\,\,
		+\int_{\prod_{j=1}^{k-1}\mathcal{V}_{j}}\sum\limits_{l=1}^{k-1}{\bf 1}_{\mathcal{A}_{3}}(\mathscr{X}(t_{2\ell-1}),\mathscr{V}(t_{2l-1})){\bf 1}_{\{t_{2l+1}>0\}}(w^{\ell}\mathcal{S}_{2,b}^{(1)})\,{\rm d}\Sigma_{l}(t_{2l+1})\\
		&\qquad\qquad\qquad\qquad\,\, +\int_{\prod_{j=1}^{k-1}\mathcal{V}_{j}}\sum\limits_{l=1}^{k-1}{\bf 1}_{\mathcal{A}_{3}}(\mathscr{X}(t_{2\ell-1}),\mathscr{V}(t_{2l-1})){\bf 1}_{\{t_{2l>0}\}}(w^{\ell}\mathcal{S}_{2,b}^{(2)})\,{\rm d}\Sigma_{l}(t_{2l})\Big\},\\
		J_{8}+J_{9}&=\sigma\frac{e^{-\int _{t_{1}}^{t}\mathcal{A}(\mathscr{V}(\tau))\,{\rm d}\tau}}{\tilde{w}(\mathscr{V}(t_{1}))}\int_{\prod_{j=1}^{k-1}\mathcal{V}_{j}}\sum\limits_{l=1}^{k-1}{\bf 1}_{\{t_{l+1}>0\}}\int_{t_{l+1}}^{t_{l}}(\chi_{M}\mu^{-\frac{1}{2}}\mathcal{K}_{w}h_{1}+ K_{w}h_{2})(s)\,{\rm d}\Sigma_{l}(s){\rm d}s,\\
		 J_{10}+J_{11}&=\sigma\frac{e^{-\int _{t_{1}}^{t}\mathcal{A}(\mathscr{V}(\tau))\,{\rm d}\tau}}{\tilde{w}(\mathscr{V}(t_{1}))}\int_{\prod_{j=1}^{k-1}\mathcal{V}_{j}}\sum\limits_{l=1}^{k-1}{\bf 1}_{\{t_{l+1}\leq 0<t_{l}\}}\int_{0}^{t_{l}}(\chi_{M}\mu^{-\frac{1}{2}}\mathcal{K}_{w}h_{1}+K_{w}h_{2})(s)\,{\rm d}\Sigma_{l}(s){\rm d}s,\\
		 J_{12}+J_{13}&=\frac{e^{-\int _{t_{1}}^{t}\mathcal{A}(\mathscr{V}(\tau))\,{\rm d}\tau}}{\tilde{w}(\mathscr{V}(t_{1}))}\int_{\prod_{j=1}^{k-1}\mathcal{V}_{j}}\sum\limits_{l=1}^{k-1}\Big\{{\bf 1}_{\{t_{l+1}>0\}}\int_{t_{l+1}}^{t_{l}}+{\bf 1}_{\{t_{l+1}\leq 0<t_{l}\}}\int_{0}^{t_{l}}\Big\}(w^{\ell}\mathcal{S}_{2})(s)\,{\rm d}\Sigma_{l}(s){\rm d}s,\\
	 J_{14}&=\frac{e^{-\int _{t_{1}}^{t}\mathcal{A}(\mathscr{V}(\tau))\,{\rm d}\tau}}{\tilde{w}(\mathscr{V}(t_{1}))}\int_{\prod_{j=1}^{k-1}\mathcal{V}_{j}}\sum\limits_{l=1}^{k-1}{\bf 1}_{\{t_{l+1}\leq 0<t_{l}\}}h_{2}(\mathscr{X}(0),\mathscr{V}(0))\,{\rm d}\Sigma_{l}(0),\\
		  J_{15}&=\frac{e^{-\int _{t_{1}}^{t}\mathcal{A}(\mathscr{V}(\tau))\,{\rm d}\tau}}{\tilde{w}(\mathscr{V}(t_{1}))}\int_{\prod_{j=1}^{k-1}\mathcal{V}_{j}}{\bf 1}_{\{t_{k}>0\}}h_{2}(\mathscr{X}(t_{k}),\mathscr{V}(t_{k}))\,{\rm d}\Sigma_{k-1}(t_{k}),
	\end{align*}
	where
	\begin{align}
	&d\Sigma_{l}(s)=\prod_{j=l+1}^{k-1}\,{\rm d} \sigma_j \big(e^{-\int_s^{t_l} \mathcal{A}(\mathscr{V}(\tau))\,{\rm d}\tau} \tilde{w}(v_l)\,{\rm d} \sigma_l\big) \prod_{j=1}^{l-1} \Big(\frac{\tilde{w}(v_j)}{\tilde{w}(\mathscr{V}(t_{j+1})} e^{-\int_{t_{j+1}}^{t_j} \mathcal{A}(\mathscr{V}(\tau))d \tau} \,{\rm d} \sigma_j\Big).
    \label{3.18-2}
	\end{align}
	The validity of \eqref{3.18-3} can be checked by induction, and we omit it for simplicity.  
	
	Case 2. For the set  $\mathcal{C}_{2}: (\eta,v)\in (\bar{\mathcal{I}}\times \R^3\backslash (\gamma_{0}\cup \gamma_{-}))\cap \mathcal{A}_{2}$, it holds that
		\begin{align}\label{3.19-3}
		&h_{2}(\eta,v)={\bf 1}_{\{t_{1,*}\leq 0\}}\,e^{-\int_{0}^{t}\mathcal{A}(\mathscr{V}(\tau))\,{\rm d}\tau}\,h_{2}(\mathscr{X}(0),\mathscr{V}(0))\nonumber\\
			&\qquad \qquad \quad +\int_{\max\{t_{1,*},0\}}^{t}e^{-\int_{s}^{t}\mathcal{A}(\mathscr{V}(\tau))\,{\rm d}\tau}\big(\sigma\chi_{M}\mu^{-\frac{1}{2}}h_{1}+ \sigma K_{w}h_{2}+w^{\ell}\mathcal{S}_{2}\big)(s)\,{\rm d}s\nonumber\\
			&\qquad \qquad \quad
			+{\bf 1}_{\{t_{1,*}>0\}}e^{-\int_{t_{1,*}}^{t}\mathcal{A}(\mathscr{V}(\tau))\,{\rm d}\tau}h_{2}(\mathscr{X}(t_{1,*}),\mathscr{V}(t_{1,*})).
		\end{align}
		Recall $t_{1,*}$ is the moment when $\mathscr{V}_{\eta}(t_{1,*})=0$ and $\mathscr{X}(t_{1,*})=\eta_{+}$. Therefore we can further express $h_{2}(\mathscr{X}(t_{1,*}),\mathscr{V}(t_{1,*}))$ by using the representation formula in Case 1. 
	
	Case 3. For the set  $\mathcal{C}_{3}: (\eta,v)\in (\bar{\mathcal{I}}\times \R^3\backslash (\gamma_{0}\cup \gamma_{-}))\cap \mathcal{A}_{3}$, it holds that
		\begin{align}\label{3.19-4}
		&h_{2}(\eta,v)={\bf 1}_{\{t_{1}\leq 0\}}\,e^{-\int_{0}^{t}\mathcal{A}(\mathscr{V}(\tau))\,{\rm d}\tau}\,h_{2}(\mathscr{X}(0),\mathscr{V}(0))\nonumber\\
			&\qquad \qquad \quad +\int_{\max\{t_{1},0\}}^{t}e^{-\int_{s}^{t}\mathcal{A}(\mathscr{V}(\tau))\,{\rm d}\tau}\big(\sigma\chi_{M}\mu^{-\frac{1}{2}}h_{1}+ \sigma K_{w}h_{2}+w^{\ell}\mathcal{S}_{2}\big)(s)\,{\rm d}s\nonumber\\
			&\qquad \qquad \quad +{\bf 1}_{\{t_{1}>0\}}e^{-\int_{t_{1}}^{t}\mathcal{A}(\mathscr{V}(\tau))\,{\rm d}\tau}(w^{\ell}\mathcal{S}_{2,b}^{(2)})(\mathscr{V}(t_{1}))\nonumber\\
			&\qquad \qquad \quad +{\bf 1}_{\{t_{1}>0\}}e^{-\int_{t_{1}}^{t}\mathcal{A}(\mathscr{V}(\tau))\,{\rm d}\tau}\frac{\sigma}{\tilde{w}(\mathscr{V}(t_{1}))}\int_{\mathcal{V}_{1}}\tilde{w}(v_{1})\frac{h_{1}}{\sqrt{\mu}}(\eta_{1},v_{1})\,{\rm d}\sigma_{1}\nonumber\\
			&\qquad\qquad\quad  +{\bf 1}_{\{t_{1}>0\}}\,e^{-\int_{t_{1}}^{t}\mathcal{A}(\mathscr{V}(\tau))\,{\rm d}\tau}\frac{1}{\tilde{w}(\mathscr{V}(t_{1}))}\int_{\mathcal{V}_{1}}\tilde{w}(v_{1})h_{2}(\eta_{1},v_{1})\,{\rm d}\sigma_{1}.
		\end{align}
		Recall in this case $x_{1}=\eta_{1}$, then $v_{1}\in \mathcal{V}_{1}:=\{u\in \R^3:u_{\eta}>0\}$, which implies $v_{1,\eta}>0$. Therefore we can further express $h_{2}(\eta_{1},v_{1})$ by using the representation formula in Case 1.
	
	Step 4. {\it Estimates on $\|h_{2}\|_{L^{\infty}}$}. For $(\eta,v)\in \mathcal{C}_{1}$, we will control the right hand side of \eqref{3.18-3} term by term. It is clear that 
	\begin{align}\label{3.18-6}
	|J_{1}|\leq Ce^{-\nu_{0}t}\|h_{2}\|_{L^{\infty}}.
	\end{align}
	For $J_{2}$, $J_{4}$ and $J_{6}$, by using Lemma \ref{K}, one has
	\begin{align}\label{3.18-5}
	|J_{2}|+|J_{4}|+|J_{6}|\leq C\|h_{1}\|_{L^{\infty}}+C\|\nu^{-1}w^{\ell}\mathcal{S}_{2}\|_{L^{\infty}}+C|w^{\ell}\mathcal{S}_{2,b}^{(1)}|_{L^{\infty}(\gamma_{-})}.
	\end{align}
	For $J_{5}$, we obtain from Lemma \ref{lem2.1} that
	\begin{align}\label{3.18-7}
	|J_{5}|\leq C_{\ell}ke^{-\frac{|v_{1}|^2}{8}}\|h_{1}\|_{L^{\infty}},
	\end{align}
	where we have used the facts
	\begin{align*}
	 \frac{1}{\tilde{w}(v_{1})}=(\sqrt{2\pi}w^{\ell}\sqrt{\mu})(v_{1})\leq C_{\ell}e^{-\frac{|v_{1}|^2}{8}},\quad \int_{\mathcal{V}}\frac{h_{1}}{\sqrt{\mu}}\tilde{w}(u)\,{\rm d}\sigma=\frac{1}{\sqrt{2\pi}}\int_{\mathcal{V}}h_{1}w^{-\ell}|u_{\eta}|\,{\rm d}u\leq C_{\ell}\|h_{1}\|_{L^{\infty}},
	 \end{align*}
	 provided $\ell>4$. For $J_{7}$, by using Lemma \ref{lem2.1}, we have
	 \begin{align}\label{3.18-8}
	|J_{6}|+|J_{7}|\leq C_{\ell}ke^{-\frac{|v_{1}|^2}{8}}|w^{\ell}(\mathcal{S}_{2,b}^{(1)},\mathcal{S}_{2,b}^{(2)})|_{L^{\infty}(\gamma_{-})}.
	 \end{align}
	 Similarly, for $J_{8}$, $J_{10}$, $J_{12}$ and $J_{13}$, we have
	 \begin{align}\label{3.18-9}
	 |J_{8}|+|J_{10}|\leq C_{\ell}ke^{-\frac{|v_{1}|^2}{8}}\|h_{1}\|_{L^{\infty}},\quad |J_{12}|+|J_{13}|\leq C_{\ell}k\|\nu^{-1}w^{\ell}\mathcal{S}_{2}\|_{L^{\infty}}.
	 \end{align}
	  For $J_{9}$, a direct calculation shows that
   \begin{align}\label{C2}
   J_{9}&=\sigma\frac{e^{-\int_{t_{1}}^{t}\mathcal{A}(\mathscr{V}(\tau))\,{\rm d}\tau}}{\tilde{w}(\mathscr{V}(t_{1}))}\int_{\prod_{j=l+1}^{k-1}\mathcal{V}_{j}}\prod_{j=l+1}^{k-1}\,{\rm d}\sigma_{j}\sum\limits_{l=1}^{k-1}\int_{\mathcal{V}_{l}}{\bf{1}}_{\{t_{l+1}>0\}}\tilde{w}(v_{l})\nonumber\\
   & \times \int_{t_{l+1}}^{t_{l}}e^{-\int_{s}^{t_{l}}\mathcal{A}(\mathscr{V}(\tau))\,{\rm d}\tau}\int_{\R^3}k_{w}(\mathscr{V}(s),v')h_{2}(\mathscr{X}(s),v')\,{\rm d}v'\,{\rm d}s\,{\rm d}\sigma_{l}\int_{\prod_{j=1}^{l-1}}\Big(\frac{\tilde{w}(v_{j})}{\tilde{w}(\mathscr{V}_{t_{j+1}})}e^{-\int_{t_{j+1}}^{t_{j}}\mathcal{A}(\tau)\,{\rm d}\tau}\,{\rm d}\sigma_{j}\Big)\nonumber\\
   &\leq C_{\ell}k e^{-\frac{1}{8}|v_{1}|^2}\int_{\mathcal{V}_{l}}\int_{t_{l+1}}^{t_{l}}e^{-\int_{s}^{t_{l}}\mathcal{A}(\mathscr{V}(\tau))\,{\rm d}\tau}\int_{\R^3}k_{w}(\mathscr{V}(s),v')h_{2}(\mathscr{X}(s),v')\,{\rm d}v'\,{\rm d}s\tilde{w}(v_{l})\,{\rm d}\sigma_{l}.
   \end{align}
   We divide the estimate of $J_{9}$ into the following four cases.

{\it Case 1. $|v_{l}|\geq N\gg 1$}. Noting that
   $$
   \tilde{w}(v_{l})\,{\rm d}\sigma_{l}=w^{-\ell}(v)\mu^{\frac{1}{2}}(v)|v_{\eta}|\,{\rm d}v\leq e^{-\frac{1}{8}N^2}\mu^{\frac{1}{4}}|v_{\eta}|\,{\rm d}v_{\eta}{\rm d}v_{\phi},
   $$
   we have from \eqref{C2} and Lemma \ref{K} that the corresponding part in $J_{9}$ can be controlled by
   \begin{align}\label{C3}
    C_{\ell}k e^{-\frac{1}{8}|v_{1}|^2}e^{-\frac{1}{8}N^2}\|h_{2}\|_{L^{\infty}}.
   \end{align}
   
   {\it Case 2. $|v_{l}|\leq N$ but $|v'|\leq 3N$}. In this case, it holds that
   $$
   |\mathscr{V}(s)|\leq |v_{l}|+O(\alpha)\leq 2N,
   $$
   and $|v'-\mathscr{V}(s)|\geq |v'|-|\mathscr{V}(s)|\geq N$, which yields that
   $$
   \int_{\R^3}k_{w}(\mathscr{V}(s),v')\,{\rm d}v'\leq e^{-c_{0}N^2}\int_{\R^3}k_{w}(\mathscr{V}(s),v')e^{c_{0}|\mathscr{V}(s)-v'|^2}\,{\rm d}v'\leq Ce^{-c_{0}N^2},
   $$
   for some sufficiently small constant $c_{0}$, where we have used Lemma \ref{K}. Then the corresponding part in $J_{9}$ can be controlled by
   \begin{align}\label{C4}
    C_{\ell}ke^{-\frac{|v_{1}|^2}{8}}e^{-c_{0}N^2}\|h_{2}\|_{L^{\infty}}.
   \end{align}
   
   {\it Case 3. $|v_{l}|\leq N$, $|v'|\leq 3N$ and $|v_{l,\eta}|\leq m$ or $|v_{l,\phi}|\leq m$ where $m$ is a small constant}. In this case, the integral domain is small, and it follows from direct calculations that the corresponding part in $J_{9}$ can be bounded as
   \begin{align}\label{C5}
 C_{N,\ell}ke^{-\frac{|v_{1}|^2}{8}}m\|h_{2}\|_{L^{\infty}}.
   \end{align}
   
   {\it Case 4.  $|v_{l}|\leq N$, $|v'|\leq 3N$, $|v_{l,\eta}|\geq m$ and $|v_{l,\phi}|\geq m$}. We further divide this case into three parts: $(x_{l
   },v_{l})\in \mathcal{A}_{1}$, $(x_{l
   },v_{l})\in \mathcal{A}_{2}$ and $(x_{l
   },v_{l})\in \mathcal{A}_{3}$. 
   
   If $(x_{l
   },v_{l})\in \mathcal{A}_{1}$, then it follows that $v_{l,\eta}>0$ and $V_{\eta}(s)\geq v_{l,\eta}\geq m$ as $V_{\eta}(s)$ is increasing when $X(s)$ is decreasing. We make the change of variable: $s\mapsto X(s;t_{l},x_{l},v_{l})$ and use Lemma \ref{K} to get
   \begin{align}\label{C6}
   &\int_{\mathcal{V}_{l}\cap\{v_{l,\eta}\geq m,|v_{l,\phi}|\geq m,|v_{l}|\leq N\}}\int_{t_{l+1}}^{t_{l}}e^{-\int_{s}^{t_{l}}\mathcal{A}(\mathscr{V}(\tau))\,{\rm d}\tau}\int_{\{|v'|\leq 3N\}}k_{w}(\mathscr{V}(s),v')h_{2}\,{\rm d}v'\,{\rm d}s\tilde{w}(v_{l})\,{\rm d}\sigma_{l}\nonumber\\
   &=\int_{\mathcal{V}_{l}\cap\{v_{l,\eta}\geq m,|v_{l,\phi}|\geq m,|v_{l}|\leq N\}}\int_{0}^{\eta_{1}}\frac{e^{-\int_{X(s)}^{\eta_{1}}\frac{\mathcal{A}(V(X(\tau)))}{V_{\eta}(X(\tau))}\,{\rm d}X(\tau)}}{|V_{\eta}(X(s))|}\nonumber\\
   &\quad \times \int_{\{|v'|\leq 3N\}}k_{w}(\mathscr{V}(s),v')h_{2}(X(s),v')\,{\rm d}v'\,{\rm d}X(s)\tilde{w}(v_{l})\,{\rm d}\sigma_{l}\nonumber\\
   &\leq 
   C_{\ell,m,N}\|\frac{h_{2}}{w^{\ell}}\|_{L^{2}},
   \end{align}
   which yields the corresponding part in $J_{9}$ can be bounded by
    \begin{align}\label{C7}
   C_{\ell,m,N}e^{-\frac{1}{8}|v_{1}|^2}\|\frac{h_{2}}{w^{\ell}}\|_{L^{2}}.
   \end{align}
   
    If $(x_{l},v_{l})\in \mathcal{A}_{3}$, by changing the variable $s\mapsto X(s;t_{l},x_{l},v_{l})$, we have
   \begin{align}\label{C9}
   	&\int_{\mathcal{V}_{l}\cap\{|v_{l,\eta}|\geq m,|v_{l,\phi}|\geq m,|v_{l}|\leq N\}}\int_{t_{l+1}}^{t_{l}}e^{-\int_{s}^{t_{l}}\mathcal{A}(\mathscr{V}(\tau))\,{\rm d}\tau}\int_{\{|v'|\leq 3N\}}k_{w}(\mathscr{V}(s),v')h_{2}\,{\rm d}v'\,{\rm d}s\tilde{w}(v_{l})\,{\rm d}\sigma_{l}\nonumber\\
   	&=\int_{\mathcal{V}_{l}\cap\{|v_{l,\eta}|\geq m,|v_{l,\phi}|\geq m,|v_{l}|\leq N\}}\int_{0}^{\eta_{1}-\lambda}\frac{e^{-\int_{X(s)}^{\eta_{1}}\frac{\mathcal{A}(V(X(\tau)))}{V_{\eta}(X(\tau))}\,{\rm d}X(\tau)}}{|V_{\eta}(X(s))|}\nonumber\\
   	&\quad \times\int_{\{|v'|\leq 3N\}}k_{w}(\mathscr{V}(s),v')h_{2}(X(s),v')\,{\rm d}v'\,{\rm d}X(s)\tilde{w}(v_{l})\,{\rm d}\sigma_{l}\nonumber\\
   	&\quad +\int_{\mathcal{V}_{l}\cap\{|v_{l,\eta}|\geq m,|v_{l,\phi}|\geq m,|v_{l}|\leq N\}}\int_{\eta_{1}-\lambda}^{\eta_{1}}\frac{e^{-\int_{X(s)}^{\eta_{1}}\frac{\mathcal{A}(V(X(\tau)))}{V_{\eta}(X(\tau))}\,{\rm d}X(\tau)}}{|V_{\eta}(X(s))|}\nonumber\\
   	&\quad \times \int_{\{|v'|\leq 3N\}}k_{w}(\mathscr{V}(s),v')h_{2}(X(s),v')\,{\rm d}v'\,{\rm d}X(s)\tilde{w}(v_{l})\,{\rm d}\sigma_{l}\nonumber\\
   	&=:I_{1}+I_{2}.
   \end{align}
   For $I_{1}$, recall that
   $$
   V_{\eta}(X(s))^2=v_{l,\eta}^2-\int_{0}^{X(s)}\frac{2}{1-y}\Big[\frac{1}{1-y}v_{l,\phi}+\frac{\alpha}{2\eta_{1}}(1-y)-\frac{\alpha}{2\eta_{1}}\frac{1}{1-y}\Big]^2\,{\rm d}y.
   $$
   Noting that $(x_{l},v_{l})\in \mathcal{A}_{3}$, we have
   $$
   v_{l,\eta}^2\geq \int_{0}^{\eta_{1}}\frac{2}{1-y}\Big[\frac{1}{1-y}v_{l,\phi}+\frac{\alpha}{2\eta_{1}}(1-y)-\frac{\alpha}{2\eta_{1}}\frac{1}{1-y}\Big]^2\,{\rm d}y,
   $$
   which implies that
   \begin{align*}
   |V_{\eta}(X(s))|^2&\geq \int_{X(s)}^{\eta_{1}}\Big[\frac{1}{1-y}v_{l,\phi}+\frac{\alpha}{2\eta_{1}}(1-y)-\frac{\alpha}{2\eta_{1}}\frac{1}{1-y}\Big]^2\,{\rm d}y\nonumber\\
   &\geq \int_{\eta_{1}-\lambda}^{\eta_{1}}\Big[\frac{1}{1-y}v_{l,\phi}+\frac{\alpha}{2\eta_{1}}(1-y)-\frac{\alpha}{2\eta_{1}}\frac{1}{1-y}\Big]^2\,{\rm d}y\nonumber\\
   &\geq \lambda(\frac{1}{C}m^2-C\alpha^2)\geq \frac{\lambda}{2C}m^2,
   \end{align*}
   where we have used the fact that $\alpha\ll m$ and $|v_{l,\phi}|^2\geq m^2$. Then $|I_{1}|$ can be bounded by
   \begin{align}\label{C11}
   	|I_{1}|\leq  \frac{C_{\ell,N}}{m^2\lambda}\|\frac{h_{2}}{w^{\ell}}\|_{L^2}.
   \end{align}
   For $I_{2}$, by direct calculations, we have
      \begin{align*}
   	|V_{\eta}(X(s))|^2&\geq \int_{X(s)}^{\eta_{1}}\frac{2}{1-y}\Big[\frac{1}{1-y}v_{l,\phi}+\frac{\alpha}{2\eta_{1}}(1-y)-\frac{\alpha}{2\eta_{1}}\frac{1}{1-y}\Big]^2\,{\rm d}y\nonumber\\
   	&\geq |v_{l,\phi}|^2\int_{X(s)}^{\eta_{1}}\frac{1}{(1-y)^3}\,{\rm d}y-C\alpha^2(\eta_{1}-X(s))\nonumber\\
   	&=\frac{(\eta_{1}-X(s))}{2}\Big[\frac{2-X(s)-\eta_{1}}{(1-X(s))^2(1-\eta_{1})^2}|v_{l,\phi}|^2-C\alpha^2\Big]\nonumber\\
    &\geq Cm^2(\eta_{1}-X(s)),
   \end{align*}
   provided that $\alpha\ll m$, which implies that
  \begin{align*}
  	\int_{\eta_{1}-\lambda}^{\eta_{1}}\frac{1}{|V_{\eta}(X(s))|}\,{\rm d}X(s)\leq \frac{C}{m}\int_{\eta_{1}-\lambda}^{\eta_{1}}(\eta_{1}-X(s))^{-\frac{1}{2}}\,{\rm d}X(s)\leq \frac{C\sqrt{\lambda}}{m}.
  \end{align*}
  Then $I_{2}$ can be bounded by
  \begin{align}\label{C10}
  	|I_{2}|\leq \frac{C_{\ell,N}\sqrt{\lambda}}{m}\|h_{2}\|_{L^{\infty}}.
  \end{align}
  Combining \eqref{C9}, \eqref{C11} and \eqref{C10}, we can control the corresponding part in $J_{9}$ as
  \begin{align}\label{C10-1}
   C_{\ell,N}k e^{-\frac{1}{8}|v_{1}|^2}\frac{\sqrt{\lambda}}{m}\|h_{2}\|_{L^{\infty}}+C_{\ell,m,N,\lambda}ke^{-\frac{1}{8}|v_{1}|^2}\|\frac{h_{2}}{w^{\ell}}\|_{L^2}.
  \end{align}
 
   If $(x_{l},v_{l})\in \mathcal{A}_{2}$, then under the change of variable $s\mapsto X(s;t_{l},x_{l},v_{l})$, we have
    \begin{align}\label{C8}
   	&\int_{\mathcal{V}_{l}\cap\{|v_{l,\eta}|\geq m,|v_{l,\phi}|\geq m,|v_{l}|\leq N\}}\int_{t_{l+1}}^{t_{l}}e^{-\int_{s}^{t_{l}}\mathcal{A}(\mathscr{V}(\tau))\,{\rm d}\tau}\int_{\{|v'|\leq 3N\}}k_{w}(\mathscr{V}(s),v')h_{2}\,{\rm d}v'\,{\rm d}s\tilde{w}(v_{l})\,{\rm d}\sigma_{l}\nonumber\\
   	&=\int_{\mathcal{V}_{l}\cap\{|v_{l,\eta}|\geq m,|v_{l,\phi}|\geq m,|v_{l}|\leq N\}}\int_{t_{l+1}}^{t_{l,*}}e^{-\big(\int_{s}^{t_{l,*}}+\int_{t_{l,*}}^{t_{l}}\big)\mathcal{A}(\mathscr{V}(\tau))\,{\rm d}\tau}\int_{\R^3}k_{w}(\mathscr{V}(s),v')h_{2}\,{\rm d}v'\,{\rm d}s\tilde{w}(v_{l})\,{\rm d}\sigma_{l}\nonumber\\
   	&\quad +\int_{\mathcal{V}_{l}\cap\{|v_{l,\eta}|\geq m,|v_{l,\phi}|\geq m,|v_{l}|\leq N\}}\int_{t_{l,*}}^{t_{l}}e^{-\int_{s}^{t}\mathcal{A}(\mathscr{V}(\tau))\,{\rm d}\tau}\int_{\R^3}k_{w}(\mathscr{V}(s),v')h_{2}\,{\rm d}v'\,{\rm d}s\tilde{w}(v_{l})\,{\rm d}\sigma_{l}\nonumber\\
   	&=\int_{\mathcal{V}_{l}\cap\{|v_{l,\eta}|\geq m,|v_{l}|\leq N\}}\int_{0}^{\eta_{+}}\frac{1}{|V_{\eta}(X(s))|}e^{-\big(\int_{X(s)}^{\eta_{+}}\frac{\mathcal{A}(V(X(\tau))}{|V_{\eta}(X(\tau)|}\,{\rm d}X(\tau)+\int_{0}^{\eta_{+}}\frac{\mathcal{A}(V(X(\tau))}{|V_{\eta}(X(\tau))|}\,{\rm d}X(\tau)\big)}\nonumber\\
   	&\qquad \times \int_{\{|v'|\leq 3N\}}k_{w}(\mathscr{V}(s),v')h_{2}(X(s),v')\,{\rm d}v'\,{\rm d}X(s)\tilde{w}(v_{l})\,{\rm d}\sigma_{l}\nonumber\\
   	&\quad +\int_{\mathcal{V}_{l}\cap\{|v_{l,\eta}|\geq m,|v_{l,\phi}|\geq m,|v_{l}|\leq N\}}\int_{0}^{\eta_{+}}\frac{1}{|V_{\eta}(X(s))|}e^{-\int_{0}^{X(s)}\frac{\mathcal{A}(V(X(\tau))}{|V_{\eta}(X(\tau))|}
   	\,{\rm d}X(\tau)}\nonumber\\
   	&\qquad \times \int_{\{|v'|\leq 3N\}}k_{w}(\mathscr{V}(s),v')h_{2}(X(s),v')\,{\rm d}v'\,{\rm d}X(s)\tilde{w}(v_{l})\,{\rm d}\sigma_{l}\nonumber\\
   	&=:I_{1}+I_{2}.
   \end{align}
   For $I_{1}$, we divide it into two parts
   \begin{align*}
   I_{1}&=\int_{\mathcal{V}_{l}\cap\{v_{l,\eta}\geq m,|v_{l}|\leq N\}}\Big(\int_{0}^{\eta_{+}-\lambda}+\int_{\eta_{+}-\lambda}^{\eta_{+}}\Big)\frac{e^{-\big(\int_{X(s)}^{\eta_{+}}\frac{\mathcal{A}(V(X(\tau))}{|V_{\eta}(X(\tau))|}\,{\rm d}X(\tau)+\int_{0}^{\eta_{+}}\frac{\mathcal{A}(V(X(\tau))}{|V_{\eta}(X(\tau))|}\,{\rm d}X(\tau)\big)}}{V_{\eta}(X(s))}\nonumber\\
   	&\qquad \times \int_{\{|v'|\leq 3N\}}k_{w}(\mathscr{V}(s),v')h_{2}(X(s),v')\,{\rm d}v'{\rm d}X(s)\tilde{w}(v_{l})\,{\rm d}\sigma_{l}.
   \end{align*}
 Then, using similar calculations as in \eqref{C11}--\eqref{C10}, we have  
   \begin{align}\label{C13}
   |I_{1}|\leq C_{\ell,N}\frac{\sqrt{\lambda}}{m}\|h_{2}\|_{L^{\infty}}+C_{\ell,m,N,\lambda}\|\frac{h_{2}}{w^{\ell}}\|_{L^2}.
   \end{align}
   Similarly, $I_{2}$ can also be bounded by
   \begin{align}\label{C14}
   	|I_{2}|\leq C_{\ell,N}\frac{\sqrt{\lambda}}{m}\|h_{2}\|_{L^{\infty}}+C_{\ell,m,N,\lambda}\|\frac{h_{2}}{w^{\ell}}\|_{L^2}.
   \end{align}
   Therefore,  combining \eqref{C3}--\eqref{C5}, \eqref{C7} and \eqref{C10-1}--\eqref{C14} together, we get
   \begin{align}\label{C15}
    |J_{9}|\leq C_{\ell,N}k e^{-\frac{1}{8}|v_{1}|^2}\frac{\sqrt{\lambda}}{m}\|h_{2}\|_{L^{\infty}}+C_{\ell,m,N,\lambda}ke^{-\frac{1}{8}|v_{1}|^2}\|\frac{h_{2}}{w^{\ell}}\|_{L^2}.
   \end{align}
  Further, $|J_{11}|$ can be bounded in a similar way. In conclusion, $J_{9}$ and $J_{11}$ can be bounded as
  \begin{align}\label{C16}
  |J_{9}|+|J_{11}|&\leq C_{\ell} ke^{-\frac{1}{8}|v_{1}|^2}\Big(e^{-\varepsilon_{0}N^2}+C_{N}m+\frac{C\sqrt{\lambda}}{m}\Big)\|h_{2}\|_{L^{\infty}}+C_{\ell,m,N,\lambda}ke^{-\frac{1}{8}|v_{1}|^2}\|\frac{h_{2}}{w^{\ell}}\|_{L^2}.
  \end{align}
 For $J_{14}$, we obtain from Lemma \ref{lem2.1} that
  \begin{align}\label{C17}
    |J_{14}|\leq C_{\ell}ke^{-\frac{|v_{1}|^2}{8}}e^{-\nu_{0}t}\|h_{2}\|_{L^{\infty}}.  
  \end{align}
  For $J_{15}$, we first notice from the boundary conditions $\eqref{3.7-1}_{2}$--$\eqref{3.7-1}_{3}$ that
  \begin{align}\label{C17-1}
 |h_{2}|_{L^{\infty}(\gamma_{-})}\leq C_{\ell}|h_{2}|_{L^{\infty}(\gamma_{+})}+C_{\ell}|h_{1}|_{L^{\infty}(\gamma_{+})}+C|w^{\ell}(\mathcal{S}_{2,b}^{(1)},\mathcal{S}_{2,b}^{(2)})|_{L^{\infty}(\gamma_{-})},
 \end{align}
 which, together with \eqref{A2}, yields that
 \begin{align}\label{C18}
|J_{15}|&\leq C_{\ell}e^{-\frac{|v_{1}|^2}{8}}\Big(\frac{1}{2}\Big)^{C_{4}T_{0}^{\frac{5}{4}}}|h_{2}|_{L^{\infty}(\gamma_{-})}\nonumber\\
&\leq C_{\ell}e^{-\frac{|v_{1}|^2}{8}}\Big[\Big(\frac{1}{2}\Big)^{C_{4}T_{0}^{\frac{5}{4}}}|h_{2}|_{L^{\infty}(\gamma_{+})}+|h_{1}|_{L^{\infty}(\gamma_{+})}+|w^{\ell}(\mathcal{S}_{2,b}^{(1)},\mathcal{S}_{2,b}^{(2)})|_{L^{\infty}(\gamma_{-})}\Big].
 \end{align}
 Combining \eqref{3.18-6}--\eqref{3.18-9}, \eqref{C16}--\eqref{C17} and \eqref{C18} together, we have that for $(\eta,v)\in \mathcal{C}_{1}$,
 \begin{align}\label{C19}
 |h_{2}(\eta,v)|&\leq \sigma\int_{\max\{t_{1},0\}}^{t}e^{-\int_{s}^{t}\mathcal{A}(\mathscr{V}(\tau))\,{\rm d}\tau}K_{w}h_{2}(\mathscr{X}(s),\mathscr{V}(s))\,{\rm d}s+\mathcal{P}(t),
 \end{align}
 where
 \begin{align*}
     \mathcal{P}(t)&=C_{\ell}e^{-\nu_{0}t}\|h_{2}\|_{L^{\infty}}+C_{\ell}e^{-\frac{|v_{1}|^2}{8}}\big[ke^{-\nu_{0}t}+\big(\frac{1}{2}\big)^{C_{4}T_{0}^{\frac{5}{4}}}+ke^{-c_{0}N^2}+C_{\ell,N}mk+\frac{C_{\ell,N}k\sqrt{\lambda}}{m}\big]\|h_{2}\|_{L^{\infty}}\nonumber\\
     & +C_{\ell}k\|h_{1}\|_{L^{\infty}}+C_{\ell}k\|\nu^{-1}w^{\ell}\mathcal{S}_{2}\|_{L^{\infty}}+C_{\ell}e^{-\frac{|v_{1}|^2}{8}}k\|w^{\ell}[\mathcal{S}_{2,b}^{(1)},\mathcal{S}_{2,b}^{(2)}\|_{L^{\infty}(\gamma_{-})}+C_{\ell,m,N,\lambda}k\|\frac{h_{2}}{w^{\ell}}\|_{L^2}.
     \end{align*}

For $(\eta,v)\in \mathcal{C}_{2}$ or $(\eta,v)\in \mathcal{C}_{3}$, recalling the representation formulas in \eqref{3.19-3}--\eqref{3.19-4}, similar calculations shows that
\begin{align}\label{C20}
|h_{2}(\eta,v)|&\leq \sigma\int_{\max\{t_{1},0\}}^{t}e^{-\int_{s}^{t}\mathcal{A}(\mathscr{V}(\tau))\,{\rm d}\tau}K_{w}h_{2}(\mathscr{X}(s),\mathscr{V}(s))\,{\rm d}s+\mathcal{P}(t).
\end{align}
 Combining \eqref{C19}--\eqref{C20}, we obtain that for all $(\eta,v)\in (\bar{\mathcal{I}}\times \R^3)\backslash \big(\gamma_{0}\cup \gamma_{-}\big)$,
 \begin{align}\label{C21}
  \|h_{2}\|_{L^{\infty}}&\leq  \sigma\int_{\max\{t_{1},0\}}^{t}e^{-\int_{s}^{t}\mathcal{A}(\mathscr{V}(\tau))\,{\rm d}\tau}K_{w}h_{2}(\mathscr{X}(s),\mathscr{V}(s))\,{\rm d}s+\mathcal{P}(t)\nonumber\\
  &= \sigma \int_{\max\{t_{1},0\}}^{t}e^{-\int_{s}^{t}\mathcal{A}(\mathscr{V}(\tau))\,{\rm d}\tau}\int_{\R^3}k_{w}(\mathscr{V}(s),v')h_{2}(\mathscr{X}(s),\mathscr{V}(s))\,{\rm d}v'\,{\rm d}s+\mathcal{P}(t).
 \end{align}
 We denote
 \begin{align*}
     t_{i}'=t_{i}(s,\mathscr{X}(s),v'),\quad \mathscr{X}'(s')=\mathscr{X}(s';s,\mathscr{X}(s),v'),\quad \mathscr{V}'(s')=\mathscr{V}(s';s,\mathscr{X}(s),v').
 \end{align*}
 Then, by iterating \eqref{C21} again, we obtain
  \begin{align}\label{C21-1}
  \|h_{2}\|_{L^{\infty}}&\leq  \sigma \int_{\max\{t_{1},0\}}^{t}e^{-\int_{s}^{t}\mathcal{A}(\mathscr{V}(\tau))\,{\rm d}\tau}\int_{\R_{v'}^3}k_{w}(\mathscr{V}(s),v')\int_{\max\{t_{1}',0\}}^{s}e^{-\int_{s'}^{s}\mathcal{A}(\mathscr{V}(s'))\,{\rm d}y}\nonumber\\
  &\qquad \times \int_{\R_{v''}^3}k_{w}(\mathscr{V}'(s'),v'')h_{2}(\mathscr{X}(s'),v'')\,{\rm d}v''{\rm d}s'{\rm d}v'{\rm d}s
\nonumber\\
&\quad +\int_{\max\{t_{1},0\}}^{t}e^{-\int_{s}^{t}\mathcal{A}(\mathscr{V}(\tau))\,{\rm d}\tau}\int_{\R_{v'}^3}k_{w}(\mathscr{V}(s),v')\mathcal{P}(t)\,{\rm d}v'{\rm d}s\nonumber\\
&\leq \sigma \int_{\max\{t_{1},0\}}^{t}e^{-\int_{s}^{t}\mathcal{A}(\mathscr{V}(\tau))\,{\rm d}\tau}\int_{\R_{v'}^3}k_{w}(\mathscr{V}(s),v')\int_{\max\{t_{1}',0\}}^{s}e^{-\int_{s'}^{s}\mathcal{A}(\mathscr{V}(y))\,{\rm d}y}\nonumber\\
&\qquad \times \int_{\R_{v''}^3}k_{w}(\mathscr{V}'(s'),v'')h_{2}(\mathscr{X}'(s'),v'')\,{\rm d}v''{\rm d}s'{\rm d}v'{\rm d}s+C\mathcal{P}(t)\nonumber\\
&=:\mathcal{R}+C\mathcal{P}(t).
  \end{align}
Similar to the estimate of $J_{9}$, we shall divide the estimate of $\mathcal{R}$ on the right hand side of \eqref{C21-1} into several cases.
 
Case 1. {\it $|v|\geq N$.} We have
$$
|\mathscr{V}(s)|^2\geq |v|^2-C\alpha|v|\geq \frac{N^2}{4}\quad \text{provided $\alpha\ll N$},
$$
which implies that $|\mathscr{V}(s)|\geq \frac{N}{2}$. Then we obtain from Lemma \ref{K} that
$$
\int_{\R_{v'}^3}k_{w}(\mathscr{V}(s),v')\,{\rm d}v'\leq \frac{C}{1+|\mathscr{V}(s)|}\leq \frac{C}{N}.
$$
Then the corresponding part in $\mathcal{R}$ can be bounded as
\begin{align*}
\frac{C}{N}\|h_{2}\|_{L^{\infty}}.
\end{align*}

Case 2. {\it $|v|\leq N$, $|v'|\geq 3N$ or $|v'|\geq 3N$, $|v''|\geq 5N$.} Similar to \eqref{C4}, the corresponding part in $\mathcal{R}$ can be bounded as
\begin{align*}
 Ce^{-c_{0}N^2}\|h_{2}\|_{L^{\infty}}.
\end{align*}

Case 3. {\it $|v|\leq N$, $|v'|\leq 3N$, $|v''|\leq 5N$ but $|v_{\eta}'|\leq m$ or $|v_{\phi}'|\leq m$}. In this case, one has
\begin{align*}
&\Big\vert\int_{\{|v'|\leq 3N,|v_{\eta}'|\leq m\text{ or }|v_{\phi}'|\leq m\}}k_{w}(\mathscr{V}(s),v')\,{\rm d}v'\Big\vert\nonumber\\
&\leq C_{N}\Big\{\int_{\{|v'|\leq 3N,|v_{\eta}'|\leq m\text{ or }|v_{\phi}'|\leq m\}}\,{\rm d}v'\Big)^{\frac{1}{2}}\Big(\int_{\R_{v'}^3}|k_{w}(\mathscr{V}(s),v|^2\,{\rm d}v'\Big)^{\frac{1}{2}}\leq C_{\ell,N}m,
\end{align*}
which yields the corresponding part in $\mathcal{R}$ can be bounded as
\begin{align*}
C_{\ell,N}m\|h_{2}\|_{L^{\infty}}.
\end{align*}

Case 4. {\it $|v|\leq N$, $|v'|\leq 3N$, $|v''|\leq 3N$ and $|v_{\eta}'|\geq m$, $|v_{\phi}'|\geq m$.} If $(\mathscr{X}(s),v')\in \mathcal{A}_{1}$, then similar to \eqref{C6}, by making the change of variable $y\mapsto X(s';s,\mathscr{X}(s),v')$, we have that the corresponding part in $\mathcal{R}$ can be bounded as
\begin{align*}
C_{\ell,m,N}\|\frac{h_{2}}{w^{\ell}}\|_{L^2}.
\end{align*}
If $(\mathscr{X}(s),v')\in \mathcal{A}_{2}$ or $(\mathscr{X}(s),v')\in \mathcal{A}_{3}$, similar to \eqref{C10-1} and \eqref{C15}, the corresponding part in $\mathcal{R}$ can be bounded as
\begin{align*}
 C_{\ell,N}\frac{\sqrt{\lambda}}{m}\|h_{2}\|_{L^{\infty}}+C_{\ell,m,N,\lambda}\|\frac{h_{2}}{w^{\ell}}\|_{L^{2}}.
\end{align*}
Combining the above estimates 
together, we obtain
\begin{align*}
|\mathcal{R}|\leq C\Big(\frac{1}{N}+C_{\ell,N}m+C_{\ell,N,m}\sqrt{\lambda}\Big)\|h_{2}\|_{L^{\infty}}+C_{\ell,N,m,\lambda}\|\frac{h_{2}}{w^{\ell}}\|_{L^{2}},
\end{align*}
which, together with \eqref{C21}, yields that
\begin{align}\label{C27}
\|h_{2}\|_{L^{\infty}}&\leq C_{\ell}ke^{-\nu_{0}t}\|h_{2}\|_{L^{\infty}}+C_{\ell}\big[\big(\frac{1}{2}\big)^{C_{4}T_{0}^{\frac{5}{4}}}+\frac{k}{N}+C_{\ell,N}mk+\frac{C_{\ell,N}k\sqrt{\lambda}}{m}\big]\|h_{2}\|_{L^{\infty}}\nonumber\\
&\quad +C_{\ell}k\big[\|h_{1}\|_{L^{\infty}}+\|\nu^{-1}w^{\ell}\mathcal{S}_{2}\|_{L^{\infty}}+|w^{\ell}(\mathcal{S}_{2,b}^{(1)},\mathcal{S}_{2,b}^{(2)})|_{L^{\infty}(\gamma_{-})}+C_{\ell,m,N,\lambda}\|\frac{h_{2}}{w^{\ell}}\|_{L^2}\big].
\end{align}
Taking $k=C_{3}T_{0}^{\frac{5}{4}}$, $t=T_{0}$, and $T_{0}$ large enough, and then  $N$ large enough, next  $m$ small enough and finally $\lambda$ small enough, we obtain from \eqref{C27} that
\begin{align*}
\|h_{2}\|_{L^{\infty}}\leq C_{\ell}\|h_{1}\|_{L^{\infty}}+C_{\ell}\|\nu^{-1}w^{\ell}\mathcal{S}_{2}\|_{L^{\infty}}+C_{\ell}|w^{\ell}[\mathcal{S}_{2,b}^{(1)},\mathcal{S}_{2,b}^{(2)}|_{L^{\infty}(\gamma_{-})}+C_{\ell}\|\frac{h_{2}}{w^{\ell}}\|_{L^2},
\end{align*}
which concludes the proof of \eqref{3.17-2}. $\hfill\square$

\section{{\it A priori} $L^{2}$ estimates for approximate systems \texorpdfstring{\eqref{3.6-1}--\eqref{3.7-1}}{(2.9)--(2.10)}}\label{sec4}

In this section, to close the {\it a priori} $L^{\infty}-$estimate, we shall establish {\it a priori }$L^{2}-$estimate on $f_{2}$. Recall \eqref{3.6-1}--\eqref{3.7-1}. Then, $f_{1}=\frac{h_{1}}{w^{\ell}}$ and $f_{2}=\frac{h_{2}}{w^{\ell}}$ satisfy the following coupled systems:
	\begin{equation}\label{3.6-2}
		\left\{\begin{aligned}
			&\epsilon f_{1}+v_{\eta}\frac{\partial f_{1}}{\partial \eta}-\frac{1}{1-\eta}\Big(v_{\phi}^2\frac{\partial f_1}{\partial v_{\eta}}-v_{\eta}v_{\phi}\frac{\partial f_{1}}{\partial v_{\phi}}\Big)-\frac{\alpha}{\eta_{1}} v_{\eta}\frac{\partial f_{1}}{\partial v_{\phi}}+\frac{\alpha}{2\eta_{1}}v_{\eta}v_{\phi}\sqrt{\mu}f_{2}+\nu f_{1}\\
			&\quad =\sigma(1-\chi_{M})\mathcal{K}f_{1}+\mathcal{S}_{1},\\
			&f_{1}(0,v_{\eta},v_{\phi},v_{z})\vert_{v_{\eta}>0}=0,\quad f_{1}(\eta_{1},v_{\eta},v_{\phi},v_{z})\vert_{v_{\eta}<0}=0,
		\end{aligned}\right.
	\end{equation}
	and
	\begin{equation}\label{3.7-2}
		\left\{\begin{aligned}
			&\epsilon f_{2}+v_{\eta}\frac{\partial f_{2}}{\partial \eta}-\frac{1}{1-\eta}\Big(v_{\phi}^2\frac{\partial f_2}{\partial v_{\eta}}-v_{\eta}v_{\phi}\frac{\partial f_2}{\partial v_{\phi}}\Big)-\frac{\alpha}{\eta_{1}} v_{\eta}\frac{\partial f_2}{\partial v_{\phi}}+\nu f_2\\
			&\qquad =\sigma Kf_2+\sigma \chi_{M}\mu^{-\frac{1}{2}}\mathcal{K}f_1+\mathcal{S}_{2},\\
			&f_2(0,v_{\eta},v_{\phi},v_{z})\vert_{v_{\eta}>0}=\sqrt{2\pi\mu}\int_{u_{\eta}<0}(\sqrt{\mu}f_2+\sigma f_{1})(0,u)|u_{\eta}|\,\mathrm{d}u+\mathcal{S}_{2,b}^{(1)},\\
			&f_2(\eta_{1},v_{\eta},v_{\phi},v_{z})\vert_{v_{\eta}<0}=\sqrt{2\pi\mu}\int_{u_{\eta}>0}(\sqrt{\mu}f_2+\sigma f_{1})(\eta_{1},u)|u_{\eta}|\,\mathrm{d}u+\mathcal{S}_{2,b}^{(2)}.
		\end{aligned}\right.
	\end{equation}
	
\begin{lemma}
For any given $\epsilon\in (0,1]$ and $\sigma\in [0,1]$, let $(f_{1},f_{2})$ be the solution of the coupled systems \eqref{3.6-2}--\eqref{3.7-2}. There exists a positive constant $\alpha_{1}$ small enough, which is independent of $\epsilon$ and $\sigma$, such that if $\alpha\in (0,\alpha_{1})$, it holds that 
\begin{align}\label{B1}
\|f_{2}\|_{L^{2}}\leq C\|w^{\ell}f_{1}\|_{L^{\infty}}+C_{\epsilon}\|\nu^{-1}w^{\ell}(\mathcal{S}_{1},\mathcal{S}_{2})\|_{L^{\infty}}+C_{\epsilon}|w^{\ell}(\mathcal{S}_{2,b}^{(1)},\mathcal{S}_{2,b})|_{L^{\infty}(\gamma_{-})},
\end{align}
where $C$ is a positive constant independent of $\epsilon$ and $\sigma$, $C_{\epsilon}$ is a positive constant depending on $\epsilon$ but is independent of $\sigma$.
\end{lemma}

\noindent\textbf{Proof}. We divide the proof into three steps.

Step 1. Multiplying $\eqref{3.7-2}_{1}$ by $(1-\eta)f_{2}$ and integrating over over $\bar{\mathcal{I}}\times \R^3$ to get
\begin{align}\label{E1}
&\epsilon\int_{0}^{\eta_{1}}(1-\eta)\int_{\R^3}|f_{2}|^2\,{\rm d}v{\rm d}\eta+\frac{1}{2}(1-\eta_{1})\int_{\R^3}v_{\eta}|f_{2}(\eta_{1})|^2\,{\rm d}v-\frac{1}{2}\int_{\R^3}v_{\eta}|f_{2}(0)|^2\,{\rm d}v\nonumber\\
		&\quad +(1-\sigma)\int_{0}^{\eta_1}(1-\eta)\int_{\R^3}|f_{2}|^2\nu(v)\,{\rm d}v{\rm d}\eta+\sigma \int_{0}^{\eta_1}(1-\eta)\int_{\R^3}|(\mathbf{I-P})f_{2}|^2\nu(v)\,{\rm d}v{\rm d}\eta\nonumber\\
		&\leq \sigma \int_{0}^{\eta_1}(1-\eta)\int_{\R^3}f_{2}(\chi_{M}\mu^{-\frac{1}{2}}\mathcal{K}f_{1})\,{\rm d}v{\rm d}\eta+\int_{0}^{\eta_{1}}\int_{\R^3}f_{2}\mathcal{S}_{2}\,{\rm d}v{\rm d}\eta,
\end{align}
where we have used the fact that
\begin{align*}
    \sigma\int_{0}^{\eta_{1}}(1-\eta)\int_{\R^3}f_{2}(\nu-K)f_{2}\,{\rm d}v{\rm d}\eta&=\sigma\int_{0}^{\eta_{1}}(1-\eta)\int_{\R^3}f_{2}\mathbf{L}f_{2}\,{\rm d}v{\rm d}\eta\nonumber\\
    &\geq \sigma \int_{0}^{\eta_1}(1-\eta)\int_{\R^3}|(\mathbf{I-P})f_{2}|^2\nu(v)\,{\rm d}v{\rm d}\eta.
\end{align*}
Noting the boundary conditions in \eqref{3.7-2} and the definition of $P_{\gamma}f_{2}$ in \eqref{3.3}, a direct calculation shows that
	\begin{align}\label{E2}
		\int_{\R^3}v_{\eta}|f_{2}(\eta_{1})|^2\,{\rm d}v&=\int_{v_{\eta}>0}|v_{\eta}||f_{2}(\eta_{1})|^2\,{\rm d}v-\int_{v_{\eta}<0}|v_{\eta}||f_{2}(\eta_{1})|^2\,{\rm d}v\nonumber\\
		&=\int_{v_{\eta}>0}|v_{\eta}||f_{2}(\eta_{1})|^2\,{\rm d}v-\sqrt{2\pi}\Big(\int_{v_{\eta}>0}|v_{\eta}|(\sqrt{\mu}f_{2})(\eta_{1})\,{\rm d}v\Big)^2\nonumber\\
		&\quad -\sqrt{2\pi}\sigma^2\Big(\int_{v_{\eta}>0}f_{1}(\eta_{1})|v_{\eta}|\,{\rm d}v\Big)^2-\int_{v_{\eta}<0}|v_{\eta}|\,|\mathcal{S}_{2,b}^{(2)}|^2\,{\rm d}v\nonumber\\
        &\quad -2\sqrt{2\pi}\Big(\int_{v_{\eta}<0}|v_{\eta}|\mathcal{S}_{2,b}^{(2)}\,{\rm d}v\Big)\Big(\int_{v_{\eta}>0}|v_{\eta}|f_{1}(\eta_{1})\,{\rm d}v\Big)\nonumber\\
		&\quad -2\sqrt{2\pi}\Big(\int_{v_{\eta}<0}|v_{\eta}|\mathcal{S}_{2,b}^{(2)}\,{\rm d}v\Big)\Big(\int_{v_{\eta}>0}|v_{\eta}|(\sqrt{\mu}f_{2})(\eta_{1})\,{\rm d}v\Big)\nonumber\\
        &\quad -2\sqrt{2\pi}\sigma\Big(\int_{v_{\eta}>0}f_{1}(\eta_{1})|v_{\eta}|\,{\rm d}v\Big)\Big(\int_{v_{\eta}>0}(\sqrt{\mu}f_{2})(\eta_{1})\,{\rm d}v\Big)\nonumber\\
		&\geq |(I-P_{\gamma})f_{2}(\eta_{1})|_{L^2(\gamma_{+})}^2-\delta |P_{\gamma}f_{2}(\eta_{1})|_{L^2(\gamma_{+})}^2\nonumber\\
        &\quad -C_{\delta}(|w^{\ell}\mathcal{S}_{2,b}^{(2)}|_{L^{\infty}(\gamma_{-})}^2+|w^{\ell}f_{1}|_{L^{\infty}(\gamma_{+})}^2)
	\end{align}
	where $\delta>0$ is a small constant and will be chosen later. Similarly,
	\begin{align*}
		-\int_{\R^3}v_{\eta}|f_{2}(0)|^2\,{\rm d}v&\geq |(I-P_{\gamma})f_{2}(0)|_{L^2(\gamma_{+})}^2-\delta|P_{\gamma}f_{2}(0)|_{L^2(\gamma_{+})}^2-C_{\delta}|w^{\ell}\mathcal{S}_{2,b}^{(1)}|_{L^{\infty}(\gamma_{-})}^2-C_{\delta}|w^{\ell}f_{1}|_{L^{\infty}(\gamma_{+})}^2.
	\end{align*}
	For the terms on the right hand side of \eqref{E1}, it follows from the H\"{o}lder inequality and Lemma \ref{K} that
	\begin{align}\label{E4}
	 &\sigma \int_{0}^{\eta_1}(1-\eta)\int_{\R^3}f_{2}(\chi_{M}\mu^{-\frac{1}{2}}\mathcal{K}f_{1})\,{\rm d}v{\rm d}\eta+\int_{0}^{\eta_{1}}(1-\eta)\int_{\R^3}f_{2}\mathcal{S}_{2}\,{\rm d}v{\rm d}\eta\nonumber\\
	 &\leq \delta\|f_{2}\|_{L^2}^2+C_{\delta}\|w^{\ell}f_{1}\|_{L^{\infty}}^2+C_{\delta}\|\nu^{-1}w^{\ell}\mathcal{S}_{2}\|_{L^{\infty}}^2.
	\end{align}
Substituting \eqref{E2}--\eqref{E4} into \eqref{E1}, one has
\begin{align}\label{E5}
&\epsilon\int_{0}^{\eta_{1}}(1-\eta)\int_{\R^3}|f_{2}|^2\,{\rm d}v{\rm d}\eta+|(I-P_{\gamma})f_{2}|_{L^2(\gamma_{+})}^2\nonumber\\
&\quad +(1-\sigma)\int_{0}^{\eta_1}(1-\eta)\int_{\R^3}|f_{2}|^2\nu(v)\,{\rm d}v{\rm d}\eta+\sigma \int_{0}^{\eta_1}(1-\eta)\int_{\R^3}|(\mathbf{I-P})f_{2}|^2\nu(v)\,{\rm d}v{\rm d}\eta\nonumber\\
		&\leq \delta\big[\|f_{2}\|_{L^2}^2+|P_{\gamma}f_{2}|_{L^2(\gamma_{+})}^2\big]+C_{\delta}\|w^{\ell}f_{1}\|_{L^{\infty}}^2+C_{\delta}\|\nu^{-1}w^{\ell}\mathcal{S}_{2}\|_{L^{\infty}}^2+C_{\delta}|w^{\ell}[\mathcal{S}_{2,b}^{(1)},\mathcal{S}_{2,b}^{(2)}|_{L^{\infty}(\gamma_{-})}^2.
\end{align}
Note that
	\begin{align*}
		|P_{\gamma}f_{2}(\eta_{1})|_{L^2(\gamma_{+})}^2&=\sqrt{2\pi}\Big(\int_{v_{\eta}>0}f_{2}(\eta_{1},v)\sqrt{\mu}|v_{\eta}|\,{\rm d}v\Big)^2\\
		&=\sqrt{2\pi}\Big(\int_{V^{\tau}}f_{2}(\eta_{1},v)\sqrt{\mu}|v_{\eta}|\,{\rm d}v\Big)^2+\sqrt{2\pi}\Big(\int_{\{v_{\eta}>0\}\backslash V^{\tau}}f_{2}(\eta_{1},v)\sqrt{\mu}|v_{\eta}|\,{\rm d}v\Big)^2\\
		&=I_{1}+I_{2},
	\end{align*}
	where
	$$
	V^{\tau}=\{v_{\eta}>0: |v_{\eta}|<\tau\text{ or }|v_{\phi}|<\tau\text{ or }|v|>\frac{1}{\tau}\}.
	$$
	It is clear that
	\begin{align}\label{E7}
		I_{1}\leq C\Big(\int_{V^{\tau}}f_{2}(\eta_{1},v)\sqrt{\mu}|v_{\eta}|\,{\rm d}v\Big)&\leq  C\Big(\int_{V^{\tau}}\mu |v_{\eta}|\,{\rm d}v\Big)\Big(\int_{v_{\eta}>0}|f_{2}(\eta_{1},v)|^2|v_{\eta}|\,{\rm d}v\Big)\nonumber\\
		&\leq C \tau |P_{\gamma}f_{2}(\eta_{1})|_{L^2(\gamma_{+})}^2+C\tau|(I-P_{\gamma})f_{2}(\eta_{1})|_{L^2(\gamma_{+})}^2.
	\end{align}
	For $I_{2}$, we use Lemma \ref{tt} to get that
	\begin{align}\label{E8}
		I_{2}&\leq C \int_{\{v_{\eta}>0\}\backslash V^{\tau}}|f_{2}(\eta_{1},v)|^2|v_{\eta}|\,{\rm d}v\nonumber\\
		&\leq C_{\tau}\|f_{2}\|_{L^2}^2+C_{\tau}\|\{v_{\eta}\partial_{\eta}-\frac{1}{1-\eta}\Big(v_{\phi}^2\partial_{v_{\eta}}-v_{\eta}v_{\phi}\partial_{v_{\phi}}\Big)-\frac{\alpha}{\eta_{1}}v_{\eta}\partial_{v_{\phi}}\}f_{2}^2\|_{L^1}\nonumber\\
		&\leq C_{\tau}(1+\epsilon)\|f_{2}\|_{L^2}^2+C_{\tau}(1-\sigma)\|f_{2}\|_{L_{\nu}^2}^2 +C_{\tau}\sigma|\langle{\bf L}f_{2},f_{2}\rangle|+C_{\tau}\sigma|\langle f_{2},\chi_{M}\mu^{-\frac{1}{2}}\mathcal{K}f_{1}\rangle |+C_{\tau}|\langle f_{2},\mathcal{S}_{2}\rangle |\nonumber\\
		&\leq C_{\tau}\|f_{2}\|_{L_{\nu}^2}+C_{\tau}\|w^{\ell}f_{1}\|_{L^{\infty}}+C_{\tau}\|\nu^{-1}w^{\ell}\mathcal{S}_{2}\|_{L^{\infty}},
	\end{align}
where $\langle f,g\rangle=\int_{0}^{\eta_{1}}\int_{\R^3}fg\,{\rm d}v{\rm d}\eta$ represents the inner product.
Combining \eqref{E7} and \eqref{E8} together and taking $\tau$ small enough, we have
\begin{align}\label{E9}
		|P_{\gamma}f_{2}(\eta_{1})|_{L^2(\gamma_{+})}^2&\leq C\big[|(I-P_{\gamma})f_{2}(\eta_{1})|_{L^2(\gamma_{+})}^2+\|f_{2}\|_{L_{\nu}^2}^2+\|w^{\ell}f_{1}\|_{L^{\infty}}+\|\nu^{-1}w^{\ell}\mathcal{S}_{2}\|_{L^{\infty}}\big].
	\end{align}
Similar calculations show that
	\begin{align}\label{E10}
		|P_{\gamma}f_{2}(0)|_{L^2(\gamma_{+})}^2\leq C\big[|(I-P_{\gamma})f_{2}(0)|_{L^2(\gamma_{+})}^2+\|f_{2}\|_{L_{\nu}^2}^2+\|w^{\ell}f_{1}\|_{L^{\infty}}+\|\nu^{-1}w^{\ell}\mathcal{S}_{2}\|_{L^{\infty}}\big].
	\end{align}
Substituting \eqref{E9}--\eqref{E10} into \eqref{E5}, we have
\begin{align}\label{E5-1}
&\epsilon\|f_{2}\|_{L^2}^2+|(I-P_{\gamma})|_{L^2(\gamma_{+})}^2+(1-\sigma)\|f_{2}\|_{L_{\nu}^2}^2+\sigma\|(\mathbf{I-P})f_{2}\|_{L_{\nu}^2}^2\nonumber\\
&\leq C\delta\big[\|f_{2}\|_{L_{\nu}^2}^2+|(I-P_{\gamma})f_{2}|_{L^2(\gamma_{+})}^2\big] +C_{\delta}\big[\|w^{\ell}f_{1}\|_{L^{\infty}}^2+\|\nu^{-1}w^{\ell}\mathcal{S}_{2}\|_{L^{\infty}}^2+|w^{\ell}(\mathcal{S}_{2,b}^{(1)},\mathcal{S}_{2,b}^{(2)})|_{L^{\infty}(\gamma_{-})}^2\big].
\end{align}

We remark that, as $\sigma\to 1$, we can only obtain from \eqref{E5-1} that
\begin{align}\label{W}
\|f_{2}\|_{L^2}\leq C_{\epsilon}\big[\|w^{\ell}f_{1}\|_{L^{\infty}}^2+\|\nu^{-1}w^{\ell}\mathcal{S}_{2}\|_{L^{\infty}}^2+|w^{\ell}[\mathcal{S}_{2,b}^{(1)},\mathcal{S}_{2,b}^{(2)}|_{L^{\infty}(\gamma_{-})}^2\big],
\end{align}
by taking $\delta\ll \epsilon$. However, this estimate is not strong enough for us to close the {\it a priori} $L^{\infty}$ estimates \eqref{3.17-1}--\eqref{3.17-2}. Indeed, since \eqref{3.17-1}--\eqref{3.17-2} are coupled, by using the $\|f_{2}\|_{L^2}$ estimate in \eqref{W}, we are only able to close \eqref{3.17-1}--\eqref{3.17-2} by choosing $\alpha\ll \epsilon$, which will lead to $\alpha \to 0$ when we take the limit $\epsilon \to 0$. Therefore, we have to establish uniform estimate on $\|\mathbf{P}f_{2}\|_{L^2}$ with respect to $\epsilon$ and $\sigma$, especially when $\sigma$ approaches 1.

Step 2. Let $\sigma_{0}\in [\frac{1}{2},1)$ be a constant that is close to 1, which will be determined later. If $\sigma\in [0,\sigma_{0}]$, noting that
	$$
	(1-\sigma)\|f_{2}\|_{L_{\nu}^2}^2\geq (1-\sigma_{0})\|f_{2}\|_{L_{\nu}^2}^2\geq  (1-\sigma_{0})\|f_{2}\|_{L^2}^2,
	$$
we obtain from \eqref{E5-1} that
\begin{align*}
&\epsilon\|f_{2}\|_{L^2}^2+(1-\sigma_{0})\|f_{2}\|_{L_{\nu}^2}^2+|(I-P_{\gamma})f_{2}|_{L^2(\gamma_{+})}\nonumber\\
	&\leq C\delta\big[\|f_{2}\|_{L_{\nu}^2}^2+|(I-P_{\gamma})f_{2}|_{L^2(\gamma_{+})}^2\big] +C_{\delta}\big[\|w^{\ell}f_{1}\|_{L^{\infty}}^2+\|\nu^{-1}w^{\ell}\mathcal{S}_{2}\|_{L^{\infty}}^2+|w^{\ell}[\mathcal{S}_{2,b}^{(1)},\mathcal{S}_{2,b}^{(2)}|_{L^{\infty}(\gamma_{-})}^2\big].
\end{align*}
Then taking $\delta$ small enough such that $C\delta\leq \frac{1}{2}(1-\sigma_{0})$, we have
\begin{align}\label{E10-1}
&\|f_{2}\|_{L_{\nu}^2}^2+|(I-P_{\gamma})f_{2}|_{L^2(\gamma_{+})}^2
\leq C_{\sigma_{0}}\|w^{\ell}f_{1}\|_{L^{\infty}}^2+C_{\sigma_{0}}\|\nu^{-1}w^{\ell}\mathcal{S}_{2}\|_{L^{\infty}}^2+C_{\sigma_{0}}|w^{\ell}[\mathcal{S}_{2,b}^{(1)},\mathcal{S}_{2,b}^{(2)}|_{L^{\infty}(\gamma_{-})}^2.
\end{align}

Step 3. If $\sigma\in (\sigma_{0},1]$, then
\begin{align}\label{E10-2}
	\sigma\|(\mathbf{I-P})f_{2}\|_{L_{\nu}^2}^2\geq \sigma_{0}\|(\mathbf{I-P})f_{2}\|_{L_{\nu}^2}^2\geq \frac{1}{2}\|(\mathbf{I-P})f_{2}\|_{L_{\nu}^2}^2.
\end{align}
As mentioned above, we still need to establish uniform estimates on $\|\mathbf{P}f_{2}\|_{L^2}$. As in \cite{DL-2020,DL-2022,DLY-2021}, we denote
	$$
	\sqrt{\mu}f=f_{1}+\sqrt{\mu}f_{2}.
	$$
	Then $f$ satisfies
	\begin{align}\label{E11}
		\left\{
		\begin{aligned}
			&\epsilon f+v_{\eta}\partial_{\eta}f-\frac{1}{1-\eta}(v_{\phi}^2\partial_{v_{\eta}}f-v_{\eta}v_{\phi}\partial_{v_{\phi}}f)-\frac{\alpha}{\eta_{1}}v_{\eta}\partial_{v_{\phi}}f+\frac{\alpha}{2\eta_{1}}v_{\eta}v_{\phi}f+(1-\sigma)\nu f+\sigma {\bf L}f\\
			&\quad =\mu^{-\frac{1}{2}}\mathcal{S}_{1}+\mathcal{S}_{2},\\
		&f(0,v_{\eta},v_{\phi},v_{z})\vert_{v_{\eta}>0}=\sqrt{2\pi\mu}\int_{u_{\eta}<0}(\sqrt{\mu} f-(1-\sigma)f_{1})(0,u)|u_{\eta}|\,{\rm d}u+\mathcal{S}_{2,b}^{(1)},\\
			&f(\eta_{1},v_{\eta},v_{\phi},v_{z})|_{v_{\eta}<0}=\sqrt{2\pi\mu}\int_{u_{\eta}>0}(\sqrt{\mu} f-(1-\sigma)f_{1})(\eta_{1},u)|u_{\eta}|\,{\rm d}u+\mathcal{S}_{2,b}^{(1)}.
		\end{aligned}
		\right.
	\end{align}
	
	\smallskip
	
	We denote
	$$
	\mathbf{P}f=\{a(\eta)+b_{1}(\eta)v_{\eta}+b_{2}(\eta)v_{\phi}+b_{3}(\eta)v_{z}+c(\eta)(|v|^2-3)\}\sqrt{\mu}.
	$$
	Multiplying \eqref{E11} by $(1-\eta)\sqrt{\mu}$ and integrating the resultant equation over $[0,\eta_{1}]\times \R^3$, we have
	\begin{align}\label{E12}
		&[\epsilon+(1-\sigma)\beta_{0}]\int_{0}^{\eta_{1}}(1-\eta)a(\eta)\,{\rm d}\eta+(1-\eta_{1})b_{1}(\eta_{1})-b_{1}(0)\nonumber\\
		&\,\, +(1-\sigma)\beta_{\gamma}\int_{0}^{\eta_{1}}(1-\eta)c(\eta)\,{\rm d}\eta+(1-\sigma)\int_{0}^{\eta_{1}}(1-\eta)\int_{\R^3}\nu(v)(\mathbf{I-P})f\sqrt{\mu}\,{\rm d}v{\rm d}\eta\nonumber\\
		&=\int_{0}^{\eta_{1}}(1-\eta)\int_{\R^3}\big(\mathcal{S}_{1}+\sqrt{\mu}\mathcal{S}_{2}\big)\,{\rm d}v{\rm d}\eta,
	\end{align}
	where $\beta_{0}=\int_{\R^3}\nu(v)\mu\,{\rm d}v$ and $\beta_{\gamma}=\int_{\R^3}\nu(v)(|v|^2-3)\mu\,{\rm d}v$, and we have used the fact that
	$$
	\int_{\R^3}\nu(v)v_{\eta}\mu\,{\rm d}v=\int_{\R^3}\nu(v)v_{\phi}\mu\,{\rm d}v=\int_{\R^3}\nu(v)v_{z}\mu\,{\rm d}v=0,
	$$
	since $\nu(v)$ is an even function with respect to $v_{\eta}$, $v_{\phi}$ and $v_{z}$. Moreover, a direct calculation shows that
	\begin{align}\label{E13}
		b_{1}(\eta_{1})&=\int_{\R^3}v_{\eta}\sqrt{\mu}f(\eta_{1})\,{\rm d}v=\int_{v_{\eta}>0}v_{\eta}\sqrt{\mu}f(\eta_{1},v)\,{\rm d}v+\int_{v_{\eta}<0}v_{\eta}\sqrt{\mu}f(\eta_{1},v)\,{\rm d}v\nonumber\\
		&=(1-\sigma)\int_{v_{\eta}>0}f_{1}(\eta_{1},v)|v_{\eta}|\,{\rm d}v+\int_{v_{\eta}<0}v_{\eta}\sqrt{\mu}\mathcal{S}_{2,b}^{(2)}\,{\rm d}v.
	\end{align}
	Similarly, one has
	\begin{align}\label{E14}
		b_{1}(0)=-(1-\sigma)\int_{v_{\eta}<0}f_{1}(0,v)|v_{\eta}|\,{\rm d}v+\int_{v_{\eta}>0}v_{\eta}\sqrt{\mu}\mathcal{S}_{2,b}^{(1)}\,{\rm d}v.
	\end{align}
	Substituting \eqref{E13}--\eqref{E14} into \eqref{E12}, we have
	\begin{align}\label{E12-1}
		&[\epsilon+(1-\sigma)\beta_{0}]\int_{0}^{\eta_{1}}(1-\eta)a(\eta)\,{\rm d}\eta\nonumber\\
		&=-(1-\sigma)(1-\eta_{1})\int_{v_{\eta}>0}f_{1}(\eta_{1})|v_{\eta}|\,{\rm d}v-(1-\eta_{1})\int_{v_{\eta}<0}v_{\eta}\sqrt{\mu}\mathcal{S}_{2,b}^{(2)}\,{\rm d}v\nonumber\\
		&\quad -(1-\sigma)\int_{v_{\eta}<0}f_{1}(0)|v_{\eta}|\,{\rm d}v+\int_{v_{\eta}>0}v_{\eta}\sqrt{\mu}\mathcal{S}_{2,b}^{(1)}\,{\rm d}v\nonumber\\
		&\quad -(1-\sigma)\beta_{\gamma}\int_{0}^{\eta_{1}}(1-\eta)c(\eta)\,{\rm d}\eta-(1-\sigma)\int_{0}^{\eta_{1}}(1-\eta)\int_{\R^3}\nu(v)(\mathbf{I-P})f\sqrt{\mu}\,{\rm d}v{\rm d}\eta\nonumber\\
		&\quad +\int_{0}^{\eta_{1}}(1-\eta)\int_{\R^3}\big(\mathcal{S}_{1}+\sqrt{\mu}\mathcal{S}_{2}\big)\,{\rm d}v{\rm d}\eta,
	\end{align}
	which implies that
	\begin{align}\label{E12-2}
		\Big\vert\int_{0}^{\eta_{1}}(1-\eta)a(\eta)\,{\rm d}\eta\Big\vert&\lesssim \frac{1}{\beta_{0}}|w^{\ell}f_{1}|_{L^{\infty}(\gamma_{+})}+\frac{\beta_{\gamma}}{\beta_{0}}\|c(\eta)\|_{L^2}+\frac{1}{c_{0}}\|(\mathbf{I-P})f\|_{L_{\nu}^2}\nonumber\\
		&\quad +\frac{1}{\epsilon}|w^{\ell}(\mathcal{S}_{2,b}^{(1)},\mathcal{S}_{2,b}^{(2)})|_{L^{\infty}(\gamma_{-})}+\frac{1}{\epsilon}\|\nu^{-1}w^{\ell}(\mathcal{S}_{1},\mathcal{S}_{2})\|_{L^{\infty}}.
	\end{align}
	Here, it is worth pointing out that the coefficient in front of $|w^{\ell}f_{1}|_{L^{\infty}(\gamma_{+})}$ does not depend on $\epsilon$ and $\sigma$, and the coefficients in front of $|w^{\ell}(\mathcal{S}_{2,b}^{(1)},\mathcal{S}_{2,b}^{(2)})|_{L^{\infty}(\gamma_{-})}$ and $\|w^{\ell}(\mathcal{S}_{1},\mathcal{S}_{2})\|_{L^{\infty}}$ depends on $\epsilon$ but is independent of $\sigma$. In view of \eqref{E12-1}, the conservation of mass fails to hold for any given $\sigma\in [0,1)$.
	
		Let $\Psi(\eta,v)\in C^{\infty}([0,\eta_{1}]\times \R^3)$ be a test function, then we take the inner product of \eqref{E11} and $\Psi$ to get
	\begin{align}\label{E15}
		&\epsilon\langle f, \Psi\rangle -\langle v_{\eta}f,\partial_{\eta}\Psi\rangle +\int_{\R^3}v_{\eta}(\Psi f)(\eta_{1},v)\,{\rm d}v-\int_{\R^3}v_{\eta}(\Psi f)(0,v)\,{\rm d}v+\langle \frac{1}{1-\eta}v_{\phi}^2f,\partial_{v_{\eta}}\Psi\rangle\nonumber\\
		&\quad-\langle \frac{1}{1-\eta}v_{\eta}f,\Psi\rangle-\langle \frac{1}{1-\eta}v_{\eta}v_{\phi}f,\partial_{v_{\phi}}\Psi\rangle-\frac{\alpha}{\eta_{1}}\langle v_{\eta}f,\partial_{v_{\phi}}\Psi\rangle+\frac{\alpha}{2\eta_{1}}\langle v_{\eta}v_{\phi}f, \Psi\rangle\nonumber\\
		&\quad +\langle1-\sigma\rangle\langle \nu f,\Psi\rangle+\sigma\langle Lf,\Psi\rangle=\langle \mu^{-\frac{1}{2}}\mathcal{S}_{1}+\mathcal{S}_{2},\Psi\rangle.
	\end{align}
	Now, we shall establish $\|(a,b_{1},b_{2},b_{3},c)\|_{L^2}$ estimates by choosing appropriate test functions.
	
	\smallskip

	\noindent\underline{\textbf{Estimates on $\|a\|_{L^2}$}}: Let
	\begin{align}\label{a}
		\Psi_{a}=v_{\eta}(|v|^2-10)\sqrt{\mu}\phi_{a}'(\eta),
	\end{align}
    where $\phi_{a}(\eta)$ satisfies 
    $$
    \phi_{a}''(\eta)=(1-\eta)a(\eta),\quad \phi_{a}'(0)=0,\quad \phi_{a}(0)=0.
    $$
    Then, it is direct to get 
    $$
   \phi_{a}'(\eta)=\int_{0}^{\eta}(1-y)a(y)\,{\rm d}y,\quad \phi_{a}(\eta)=\int_{0}^{\eta}\int_{0}^{s}(1-y)a(y)\,{\rm d}y{\rm d}s,
    $$
	which, together with direct calculations and \eqref{E12-2}, yields that
	\begin{align}\label{E16}
    \begin{aligned}
		|\phi_{a}'(\eta_{1})|&\leq C|w^{\ell}f_{1}|_{L^{\infty}(\gamma_{+})}+C\|c(\eta)\|_{L^2}+C\|(\mathbf{I-P})f\|_{L_{\nu}^2}\\
		&\quad +C_{\epsilon}|w^{\ell}[\mathcal{S}_{2,b}^{(1)},\mathcal{S}_{2,b}^{(2)}|_{L^{\infty}(\gamma_{-})}+C_{\epsilon}\|\nu^{-1}w^{\ell}(\mathcal{S}_{1},\mathcal{S}_{2})\|_{L^{\infty}},\\
		\|\phi_{a}\|_{H^2}&\leq C\|a\|_{L^2}.
    \end{aligned}
	\end{align}

Taking $\Psi=\Psi_{a}$ in \eqref{E15}, we have
	\begin{align}\label{E17}
		-\la v_{\eta}f, \partial_{\eta}\Psi_{a}\ra=&-\epsilon\la f,\Psi_{a}\ra-(1-\sigma)\la \nu f,\Psi_{a}\ra-\int_{\R^3}v_{\eta}(\Psi_{a}f)(\eta_{1})\,{\rm d}v+\int_{\R^3}v_{\eta}(\Psi_{a}f)(0)\,{\rm d}v\nonumber\\
		&-\la \frac{1}{1-\eta}v_{\phi}^2f,\partial_{v_{\eta}}\Psi_{a}\ra+\la \frac{1}{1-\eta}v_{\eta}f,\Psi_{a}\ra+\la \frac{1}{1-\eta}v_{\eta}v_{\phi}f,\partial_{v_{\phi}}\Psi_{a}\ra\nonumber\\
		&-\frac{\alpha}{\eta_{1}}\la v_{\eta}f,\partial_{v_{\phi}}\Psi_{a}\ra-\frac{\alpha}{2\eta_{1}}\la v_{\eta}v_{\phi}f,\Psi_{a}\ra-\sigma\la {\bf L}f,\Psi_{a}\ra+\la \mathscr{S}, \Psi_{a}\ra=\sum\limits_{k=1}^{11}I_{k},
        \end{align}
	where $\mathscr{S}=:\mu^{-\frac{1}{2}}\mathcal{S}_{1}+\mathcal{S}_{2}$. It follows from \eqref{a} that
	\begin{align*}
		-\la v_{\eta}f,\partial_{\eta}\Psi_{a}\ra&=-\int_{0}^{\eta_{1}}\phi_{a}''(\eta)\int_{\R^3}v_{\eta}^2(|v|^2-10)\sqrt{\mu}f\,{\rm d}v{\rm d}\eta\nonumber\\
		&\geq 5\int_{0}^{\eta_{1}}(1-\eta)|a(\eta)|^2\,{\rm d}\eta-\kappa_{1}\|a\|_{L^2}^2-C_{\kappa_{1}}\|(\mathbf{I-P})f\|_{L_{\nu}^2}^2\nonumber\\
		&\geq 5(1-\eta_{1})\|a\|_{L^2}^2-\kappa_{1}\|a\|_{L^2}^2-C_{\kappa_{1}}\|(\mathbf{I-P})f\|_{L_{\nu}^2}^2,
	\end{align*}
	where $\kappa_{1}>0$ is a small parameter and will be chosen later.
	
	For $I_{1}$, $I_{2}$, $I_{10}$ and $I_{11}$, by using the H\"{o}lder inequality, we have
	\begin{align}
		I_{1}&=\epsilon\int_{0}^{\eta_{1}}\phi_{a}'(\eta)\int_{\R^3}v_{\eta}(|v|^2-10)\sqrt{\mu}f\,{\rm d}v{\rm d}\eta\leq  \kappa_{1}\|a\|_{L^2}^2+C_{\kappa_{1}}\|b_{1}\|_{L^2}^2+C_{\kappa_{1}}\|(\mathbf{I-P})f\|_{L_{\nu}^2}^2,\label{E19}\\
		I_{2}
        &
        \leq C(1-\sigma_{0})\|a\|_{L^2}\|b_{1}\|_{L^2}+\kappa_{1}\|a\|_{L^2}^2+C_{\kappa_{1}}\|(\mathbf{I-P})\|_{L_{\nu}^2}^2,\label{E20}\\
		I_{10}&\leq \kappa_{1}\|a\|_{L^2}^2+C_{\kappa_{1}}\|(\mathbf{I-P})f\|_{L_{\nu}^2}^2,\label{E21}\\
		I_{11}
		&\leq \kappa_{1}\|a\|_{L^2}^2+C_{\kappa_{1}}\|\nu^{-1}w^{\ell}(\mathcal{S}_{1},\mathcal{S}_{2})\|_{L^{\infty}}^2.\label{E22}
	\end{align}
	For the boundary terms $I_{3}$ and $I_{4}$, noting the boundary condition $\eqref{E11}_{3}$, $\phi_{a}'(0)=0$ and \eqref{E16}, we have
		\begin{align}\label{E23}
		I_{3}+I_{4}&=-\phi_{a}'(\eta_{1})\int_{\R^3}v_{\eta}^2(|v|^2-10)\sqrt{\mu}f(\eta_{1},v)\,{\rm d}v\nonumber\\
		&=-\phi_{a}'(\eta_{1})\Big[\int_{v_{\eta}>0}v_{\eta}^2(|v|^2-10)\sqrt{\mu}f(\eta_{1},v)\,{\rm d}v+\int_{v_{\eta}<0}v_{\eta}^2(|v|^2-10)\sqrt{\mu}f(\eta_{1},v)\,{\rm d}v\Big]\nonumber\\
		&=10\sqrt{2\pi}\phi_{a}'(\eta_{1})\int_{v_{\eta}>0}|v_{\eta}|\sqrt{\mu}f(\eta_{1},v)\,{\rm d}v
	 -\phi_{a}'(\eta_{1})\int_{v_{\eta}>0}v_{\eta}^2(|v|^2-10)\sqrt{\mu}(I-P_{\gamma})f(\eta_{1},v)\,{\rm d}v\nonumber\\
	 &\quad -5(1-\sigma)\phi_{a}'(\eta_{1})\int_{v_{\eta}>0}|v_{\eta}|f_{1}(\eta_{1},v)\,{\rm d}v-\phi_{a}'(\eta_{1})\int_{v_{\eta<0}}v_{\eta}^2(|v|^2-10)\sqrt{\mu}\mathcal{S}_{2,b}^{(2)}\,{\rm d}v\nonumber\\
	 &\leq \kappa_{1}|P_{\gamma}f|_{L^2(\gamma_{+})}^2+C_{\kappa_{1}}|\phi_{a}'(\eta_{1})|^2+C_{\kappa_{1}}\big[|(I-P_{\gamma})f|_{L^2(\gamma_{+})}^2+|w^{\ell}f_{1}|_{L^{\infty}(\gamma_{+})}+|w^{\ell}\mathcal{S}_{2,b}^{(2)}|_{L^{\infty}(\gamma_{-})}\big]\nonumber\\
	 &\leq \kappa_{1}|P_{\gamma}f|_{L^2(\gamma_{+})}^2+C_{\kappa_{1}}\|c\|_{L^2}^2+C_{\kappa_{1}}|(I-P_{\gamma})f|_{L^2(\gamma_{+})}^2+C_{\kappa_{1}}\|(\mathbf{I-P})f\|_{L_{\nu}^2}^2\nonumber\\
	 &\quad +C_{\kappa_{1}}|w^{\ell}f_{1}|_{L^{\infty}(\gamma_{+})}^2+C_{\epsilon,\kappa_{1}}|w^{\ell}(\mathcal{S}_{2,b}^{(1)},\mathcal{S}_{2,b}^{(2)})|_{L^{\infty}(\gamma_{-})}^2+C_{\epsilon,\kappa_{1}}\|\nu^{-1}w^{\ell}(\mathcal{S}_{1},\mathcal{S}_{2})\|_{L^{\infty}}^2.
	\end{align}
Furthermore, we obtain from \eqref{E9}--\eqref{E10} that
\begin{align*}
    |P_{\gamma}f|_{L^2(\gamma_{+})}^2&=|P_{\gamma}(\frac{f_{1}}{\sqrt{\mu}})+P_{\gamma}f_{2}|_{L^2(\gamma_{+})}^2\leq C|w^{\ell}f_{1}|_{L^{\infty}(\gamma_{+})}^2+C\|P_{\gamma}f_{2}\|_{L^2(\gamma_{+})}^2\nonumber\\
    &\leq C\big[|(I-P_{\gamma})f_{2}|_{L^2(\gamma_{+})}^2+\|f_{2}\|_{L_{\nu}^2}^2+\|w^{\ell}f_{1}\|_{L^{\infty}}^2+\|\nu^{-1}w^{\ell}(\mathcal{S}_{1},\mathcal{S}_{2})\|_{L^{\infty}}^2\big],
\end{align*}
which, together with \eqref{E23}, yields
\begin{align}\label{E24}
I_{3}+I_{4}&\lesssim \kappa_{1}\|f_{2}\|_{L_{\nu}^2}^2+C_{\kappa_{1}}\|c\|_{L^2}^2+C_{\kappa_{1}}|(I-P_{\gamma})f|_{L^2(\gamma_{+})}^2+C_{\kappa_{1}}\|(\mathbf{I-P})f\|_{L_{\nu}^2}^2\nonumber\\
	 &\quad +C_{\kappa_{1}}|w^{\ell}f_{1}|_{L^{\infty}(\gamma_{+})}^2+C_{\epsilon,\kappa_{1}}|w^{\ell}(\mathcal{S}_{2,b}^{(1)},\mathcal{S}_{2,b}^{(2)})|_{L^{\infty}(\gamma_{-})}^2+C_{\epsilon,\kappa_{1}}\|\nu^{-1}w^{\ell}(\mathcal{S}_{1},\mathcal{S}_{2})\|_{L^{\infty}}^2.
\end{align}
For the geometric terms $I_{5}+I_{6}+I_{7}$, direct calculations show that 
	\begin{align}\label{E25}
		I_{5}+I_{6}+I_{7}&=\int_{0}^{\eta_{1}}\frac{1}{1-\eta}\phi_{a}'(\eta)\int_{\R^3}(\mathbf{I-P})f(|v|^2-10)(v_{\phi}^2-v_{\eta}^2)\,{\rm d}v{\rm d}\eta\nonumber\\
		&\leq \kappa_{1}\|a\|_{L^2}^2+C_{\kappa_{1}}\|(\mathbf{I-P})f\|_{L_{\nu}^2}^2.
	\end{align}
For the terms $I_{8}$ and $I_{9}$ induced by shear force, we have
	\begin{align}\label{E26}
		I_{8}+I_{9}&=\int_{0}^{\eta_{1}}\frac{2\alpha}{\eta_{1}}\phi_{a}'(\eta)\int_{\R^3}v_{\eta}^2v_{\phi}\sqrt{\mu}f(\eta,v)\,{\rm d}v{\rm d}\eta\leq C\alpha\|a\|_{L^2}\|b_{2}\|_{L^2}+\kappa_{1}\|a\|_{L^2}^2+C_{\kappa_{1}}\|(\mathbf{I-P})f\|_{L_{\nu}^2}^2.
	\end{align}
Plugging all the above estimates into \eqref{E17}, we obtain
	\begin{align}\label{E27}
	\|a\|_{L^2}^2&\lesssim \kappa_{1}\|a\|_{L^2}^2+\alpha\|a\|_{L^2}\|b_{2}\|_{L^2}+(1-\sigma_{0})\|a\|_{L^2}\|b_{1}\|_{L^2}+\kappa_{1}\|f_{2}\|_{L_{\nu}^2}^2+C_{\kappa_{1}}\|(b_{1},c)\|_{L^2}^2\nonumber\\
	 &\quad +C_{\kappa_{1}}|(I-P_{\gamma})f|_{L^2(\gamma_{+})}^2+C_{\kappa_{1}}\|(\mathbf{I-P})f\|_{L_{\nu}^2}^2+C_{\kappa_{1}}\|w^{\ell}f_{1}\|_{L^{\infty}}^2\nonumber\\
	 &\quad +C_{\epsilon,\kappa_{1}}|w^{\ell}(\mathcal{S}_{2,b}^{(1)},\mathcal{S}_{2,b}^{(2)})|_{L^{\infty}(\gamma_{-})}^2+C_{\epsilon,\kappa_{1}}\|\nu^{-1}w^{\ell}(\mathcal{S}_{1},\mathcal{S}_{2})\|_{L^{\infty}}^2.
	\end{align}
	
\smallskip
	
\noindent\underline{\textbf{Estimates on $\|b_{1}\|_{L^2}$}}: Let
	\begin{align*}
		\Psi_{b_{1}}=-v_{\eta}^2(|v|^2-5)\sqrt{\mu}\phi_{b_{1}}'(\eta),
	\end{align*}
where $\phi_{b_{1}}(\eta)$ satisfies
$$
\phi_{b_{1}}''(\eta)=\frac{1}{1-\eta}b_{1}(y),\quad \phi_{b_{1}}'(0)=0,\,\,\phi_{b_{1}}(0)=0.
$$
Then it is direct to get
$$
\phi_{b_{1}}'(\eta)=\int_{0}^{\eta}\frac{1}{1-y}b_{1}(y)\,{\rm d}y,\quad \phi_{b_{1}}(\eta)=\int_{0}^{\eta}\int_{0}^{s}\frac{1}{1-y}b_{1}(y)\,{\rm d}y{\rm d}s.
$$
It follows from the H\"{o}lder inequality that
	\begin{align*}
		  \|\phi_{b_{1}}\|_{H^2}\lesssim \|b_{1}\|_{L^2}.
	\end{align*}
	Taking $\Psi=\Psi_{b_{1}}$ in \eqref{E15}, we have
	\begin{align}\label{E28}
		-\la v_{\eta}f,\partial_{\eta}\Psi_{b_{1}}\ra&=-\epsilon\la f,\Psi_{b_{1}}\ra-(1-\sigma)\la \nu f,\Psi_{b_{1}}\ra-\int_{\R^3}v_{\eta}(\Psi_{b_{1}}f)(\eta_{1},v)\,{\rm d}v+\int_{\R^3}v_{\eta}(\Psi_{b_{1}}f)(0,v)\,{\rm d}v\nonumber\\
		&\quad-\la \frac{1}{1-\eta}v_{\phi}^2f,\partial_{v_{\eta}}\Psi_{b_{1}}\ra+\la \frac{1}{1-\eta}v_{\eta}f,\Psi_{b_{1}}\ra+\la \frac{1}{1-\eta}v_{\eta}v_{\phi}f,\partial_{v_{\phi}}\Psi_{b_{1}}\ra\nonumber\\
		&\quad -\frac{\alpha}{\eta_{1}}\la v_{\eta}f,\partial_{v_{\phi}}\Psi_{b_{1}}\ra-\frac{\alpha}{2\eta_{1}}\la v_{\eta}v_{\phi}f,\Psi_{b_{1}}\ra-\sigma\la {\bf L}f,\Psi_{b_{1}}\ra+\la \mathscr{S},\Psi_{b_{1}}\ra=\sum\limits_{k=1}^{11}J_{k}.
	\end{align}
A direct calculation shows that
\begin{align*}
		-\la v_{\eta}f, \partial_{\eta}\Psi_{b_{1}}\ra
		&\geq (1-\kappa_{2})\|b_{1}\|_{L^2}^2-C_{\kappa_{2}}\|(\mathbf{I-P})f\|_{L_{\nu}^2}^2,
	\end{align*}
where $\kappa_{2}>0$ is a small parameter and will be chosen later. 


Similar to \eqref{E19}--\eqref{E22}, we have
\begin{align*}
|J_{1}|+|J_{2}|+|J_{10}|+|J_{11}|&\leq \kappa_{2}\|b_{1}\|_{L^2}^2+C_{\kappa_{2}}\|c\|_{L^2}^2+C(1-\sigma_{0})\|b_{1}\|_{L^2}\|(a,c)\|_{L^2}\nonumber\\
&\quad +C_{\kappa_{2}}\big[\|(\mathbf{I-P}f\|_{L_{\nu}^2}^2+\|\nu^{-1}w^{\ell}(\mathcal{S}_{1},\mathcal{S}_{2})\|_{L^{\infty}}\big].
\end{align*}
For the boundary terms $J_{3}$ and $J_{4}$, noting $\phi_{b_{1}}'(0)=0$, we have
\begin{align*}
J_{3}+J_{4}&=\phi_{b_{1}}'(\eta_{1})\int_{\R^3}v_{\eta}^3(|v|^2-5)\sqrt{\mu}f(\eta_{1},v)\,{\rm d}v\nonumber\\
&=\phi_{b_{1}}'(\eta_{1})\int_{v_{\eta}>0}v_{\eta}^3(|v|^2-5)\sqrt{\mu}f(\eta_{1},v)\,{\rm d}v+\phi_{b_{1}}'(\eta_{1})\int_{v_{\eta}<0}v_{\eta}^3(|v|^2-5)\sqrt{\mu}f(\eta_{1},v)\,{\rm d}v\nonumber\\
&=\phi_{b_{1}}'(\eta_{1})\int_{v_{\eta}>0}v_{\eta}^{3}(|v|^2-5)\sqrt{\mu}(I-P_{\gamma})f(\eta_{1},v)\,{\rm d}v+\phi_{b_{1}}'(\eta_{1})\int_{v_{\eta}<0}v_{\eta}^3(|v|^2-5)\sqrt{\mu}\mathcal{S}_{2,b}^{(2)}\,{\rm d}v\nonumber\\
&\quad -6(1-\sigma)\phi_{b_{1}}'(\eta_{1})\Big(\int_{v_{\eta}>0}f_{1}(\eta_{1},v)|v_{\eta}|\,{\rm d}v\Big)\nonumber\\
&\leq \kappa_{2}\|b_{1}\|_{L^2}^2+C_{\kappa_{2}}\big[|(I-P_{\gamma})f|_{L^2(\gamma_{+})}^2+|w^{\ell}f_{1}|_{L^{\infty}(\gamma_{+})}^2+|w^{\ell}\mathcal{S}_{2,b}^{(2)}|_{L^{\infty}(\gamma_{-})}^2\big].
\end{align*}
For the geometric terms $J_{5}$, $J_{6}$ and $J_{7}$, a direct calculation shows that
	\begin{align*}
		J_{5}+J_{6}+J_{7}&=-2\int_{0}^{\eta_{1}}\frac{1}{1-\eta}b_{1}(\eta)\phi_{b_{1}}'(\eta)\,{\rm d}\eta+\int_{0}^{\eta_{1}}\frac{2}{1-\eta}\phi_{b_{1}}'(\eta)\int_{\R^3}v_{\phi}^2v_{\eta}(|v|^2-5)\sqrt{\mu}(\mathbf{I-P})f\,{\rm d}v{\rm d}\eta\nonumber\\
		&\quad -\int_{0}^{\eta_{1}}\frac{1}{1-\eta}\phi_{b_{1}}'(\eta)\int_{\R^3}v_{\eta}^3(|v|^2-5)\sqrt{\mu}(\mathbf{I-P})f\,{\rm d}v{\rm d}\eta\nonumber\\
		&=-|\phi_{b_{1}}'(\eta_{1})|^2+\int_{0}^{\eta_{1}}\frac{2}{1-\eta}\phi_{b_{1}}'(\eta)\int_{\R^3}v_{\phi}^2v_{\eta}(|v|^2-5)\sqrt{\mu}(\mathbf{I-P})f\,{\rm d}v{\rm d}\eta\nonumber\\
		&\quad -\int_{0}^{\eta_{1}}\frac{1}{1-\eta}\phi_{b_{1}}'(\eta)\int_{\R^3}v_{\eta}^3(|v|^2-5)\sqrt{\mu}(\mathbf{I-P})f\,{\rm d}v{\rm d}\eta\nonumber\\
		&\leq -|\phi_{b_{1}}'(\eta_{1})|^2+\kappa_{2}\|b_{1}\|_{L^2}^2+C_{\kappa_{2}}\|(\mathbf{I-P})f\|_{L_{\nu}^2}^2.
	\end{align*}
For terms $J_{8}$ and $J_{9}$ induced by shear force, we have
	\begin{align}\label{E36}
	J_{8}+J_{9}&=-\frac{2\alpha}{\eta_{1}}\int_{0}^{\eta_{1}}\phi_{b_{1}}'(\eta)\int_{\R^3}v_{\eta}^3v_{\phi}\sqrt{\mu}(\mathbf{I-P})f\,{\rm d}v{\rm d}\eta\leq \kappa_{2}\|b_{1}\|_{L^2}^2+C_{\kappa_{2}}\|(\mathbf{I-P})f\|_{L_{\nu}^2}^2.
\end{align}
Plugging all the above estimates into \eqref{E28}, we obtain
\begin{align}\label{E37}
		\|b_{1}\|_{L^2}^2+|\phi_{b_{1}}'(\eta_{1})|^2&\lesssim \kappa_{2}\|b_{1}\|_{L^2}^2+(1-\sigma_{0})\|b_{1}\|_{L^2}\|(a,c)\|_{L^2}+C_{\kappa_{2}}\|c\|_{L^2}^2
		\nonumber\\
        &\quad +C_{\kappa_{2}}\big[\|(\mathbf{I-P})f\|_{L_{\nu}^2}^2
	  	 +|(I-P_{\gamma})f|_{L^2(\gamma_{+})}^2+\|w^{\ell}f_{1}\|_{L^{\infty}}^2\big]\nonumber\\
         &\quad +C_{\kappa_{2}}\big[|w^{\ell}\mathcal{S}_{2,b}^{(2)}|_{L^{\infty}(\gamma_{-})}^2+\|\nu^{-1}w^{\ell}(\mathcal{S}_{1},\mathcal{S}_{2})\|_{L^{\infty}}^2\big].
	\end{align}

\smallskip

\noindent\underline{\textbf{Estimates on $\|b_{2}\|_{L^2}$}}: Let 
	$$
\Psi_{b_{2}}=v_{\eta}v_{\phi}\sqrt{\mu}\phi_{b_{2}}'(\eta),
	$$
where $\phi_{b_{2}}(\eta)$ satisfies
$$
\phi_{b_{2}}''(\eta)=-\frac{1}{1-\eta}b_{2}(\eta),\quad \phi_{b_{2}}'(\eta_{1})=0,\,\,\phi_{b_{2}}(0)=0.
$$
Then it is direct to get
$$
\phi_{b_{2}}'(\eta)=\int_{\eta}^{\eta_{1}}\frac{1}{1-y}b_{2}(y)\,{\rm d}y,\quad \phi_{b_{2}}(\eta)=\int_{0}^{\eta}\int_{s}^{\eta_{1}}\frac{1}{1-y}b_{2}(y)\,{\rm d}y{\rm d}s.
$$
It follows from the H\"{o}lder inequality that
	\begin{align*}
		 \|\phi_{b_{2}}\|_{H^2}\lesssim \|b_{2}\|_{L^2}.
	\end{align*}
	Taking $\Psi=\Psi_{b_{2}}$ in \eqref{E15}, we have	
	\begin{align}\label{E38}
		-\la v_{\eta}f,\partial_{\eta}\Psi_{b_{2}}\ra&=-\epsilon\la f,\Psi_{b_{2}}\ra-(1-\sigma)\la \nu f,\Psi_{b_{2}}\ra-\int_{\R^3}v_{\eta}(\Psi_{b_{2}}f)(\eta_{1},v)\,{\rm d}v+\int_{\R^3}v_{\eta}(\Psi_{b_{2}}f)(0,v)\,{\rm d}v\nonumber\\
		&\,\,-\la \frac{1}{1-\eta}v_{\phi}^2f,\partial_{v_{\eta}}\Psi_{b_{2}}\ra+\la \frac{1}{1-\eta}v_{\eta}f,\Psi_{b_{2}}\ra+\la \frac{1}{1-\eta}v_{\eta}v_{\phi}f,\partial_{v_{\phi}}\Psi_{b_{2}}\ra\nonumber\\
		&\,\, -\frac{\alpha}{\eta_{1}}\la v_{\eta}f,\partial_{v_{\phi}}\Psi_{b_{2}}\ra-\frac{\alpha}{2\eta_{1}}\la v_{\eta}v_{\phi}f,\Psi_{b_{2}}\ra-\sigma\la {\bf L}f,\Psi_{b_{2}}\ra+\la \mathscr{S},\Psi_{b_{2}}\ra=\sum\limits_{k=1}^{11}H_{k}.
	\end{align}
	A direct calculation shows that
	\begin{align*}
		-\la v_{\eta}f,\partial_{\eta}\Psi_{b_{2}}\ra
		&\geq (1-\kappa_{3})\|b_{2}\|_{L^2}^2-C_{\kappa_{3}}\|(\mathbf{I-P})f\|_{L_{\nu}^2}^2,
	\end{align*}
where $\kappa_{3}>0$ is a small parameter, which will be determined later. 


Similar to \eqref{E19}--\eqref{E22}, we have
	\begin{align*}
    |H_{1}|+|H_{2}|+|H_{10}|+|H_{11}|\leq \kappa_{3}\|b_{2}\|_{L^2}+C_{\kappa_{3}}\big[\|(\mathbf{I-P})f\|_{L_{\nu}^2}+\|\nu^{-1}w^{\ell}(\mathcal{S}_{1},\mathcal{S}_{2})\|_{L^{\infty}}^2\big].
	\end{align*}
For the boundary terms $H_{3}$ and $H_{4}$,  noting $\phi_{b_{2}}'(\eta_{1})=0$, we have
\begin{align*}
	H_{3}+H_{4}
	&=\phi_{b_{2}}'(0)\int_{v_{\eta>0}}v_{\eta}^2v_{\phi}\sqrt{\mu}f(0,v)\,{\rm d}v+\phi_{b_{2}}'(0)\int_{v_{\eta}<0}v_{\eta}^2v_{\phi}\sqrt{\mu}f(0,v)\,{\rm d}v\nonumber\\
    &=\phi_{b_{2}}'(0)\int_{v_{\eta}>0}v_{\eta}^2v_{\phi}\sqrt{\mu}\mathcal{S}_{2,b}^{(1)}\,{\rm d}v+\phi_{b_{2}}'(0)\int_{v_{\eta}<0}v_{\eta}^2v_{\phi}\sqrt{\mu}(I-P_{\gamma})f(0,v)\,{\rm d}v\nonumber\\
	&\leq \kappa_{3}\|b_{2}\|_{L^2}^2+C_{\kappa_{3}}|(I-P_{\gamma})f|_{L^2(\gamma_{+})}^2+C_{\kappa_{3}}|w^{\ell}\mathcal{S}_{2,b}^{(1)}|_{L^{\infty}(\gamma_{-})}^2.
\end{align*}

For the geometric terms $H_{5}$, $H_{6}$ and $H_{7}$, a direct calculation shows that
\begin{align*}
H_{5}+H_{6}+H_{7}&=-\int_{0}^{\eta_{1}}\frac{1}{1-\eta}\phi_{b_{2}}'(\eta)\int_{\R^3}v_{\phi}^3\sqrt{\mu}f\,{\rm d}v{\rm d}\eta+2\int_{0}^{\eta_{1}}\frac{1}{1-\eta}\phi_{b_{2}}'(\eta)\int_{\R^3}v_{\phi}v_{\eta}^2\sqrt{\mu}f\,{\rm d}v{\rm d}\eta\nonumber\\
		&=-\int_{0}^{\eta_{1}}\frac{1}{1-\eta}\phi_{b_{2}}'(\eta)b_{2}(\eta)\,{\rm d}\eta-\int_{0}^{\eta_{1}}\frac{1}{1-\eta}\phi_{b_{2}}(\eta)\int_{\R^3}v_{\phi}^3\sqrt{\mu}(\mathbf{I-P})f\,{\rm d}v{\rm d}\eta\nonumber\\
		&\quad +2\int_{0}^{\eta_{1}}\frac{1}{1-\eta}\phi_{b_{2}}'(\eta)\int_{\R^3}v_{\phi}v_{\eta}^2\sqrt{\mu}(\mathbf{I-P})f\,{\rm d}v{\rm d}\eta\nonumber\\
		&\leq -\frac{1}{2}|\phi_{b_{2}}'(0)|^2+\kappa_{3}\|b_{2}\|_{L^2}^2+C_{\kappa_{3}}\|(\mathbf{I-P})f\|_{L_{\nu}^2}^2.
\end{align*}
For the terms $H_{8}$ and $H_{9}$ induced by shear force, we have
\begin{align}\label{E46}
	H_{8}+H_{9}&=\frac{\alpha}{\eta_{1}}\int_{0}^{\eta_{1}}\phi_{b_{2}}'(\eta)\int_{\R^3}v_{\eta}^2\sqrt{\mu}f\,{\rm d}v{\rm d}\eta\nonumber\\
		&\lesssim C\alpha\|a\|_{L^2}\|b_{2}\|_{L^2}+C\alpha\|c\|_{L^2}\|b_{2}\|_{L^2}+\kappa_{3}\|b_{2}\|_{L^2}^2+C_{\kappa_{3}}\|(\mathbf{I-P})f\|_{L_{\nu}^2}^2.
\end{align}
Plugging all the above estimates into \eqref{E38}, we obtain
\begin{align}\label{E47}
	\|b_{2}\|_{L^2}^2+|\phi_{b_{2}}'(0)|^2&\lesssim \kappa_{3}\|b_{2}\|_{L^2}^2+\alpha\|b_{2}\|_{L^2}\|(a,c)\|_{L^2}+C_{\kappa_{3}}\|(\mathbf{I-P})f\|_{L_{\nu}^2}^2
		\nonumber\\
	  	 &\quad +C_{\kappa_{3}}\big[|(I-P_{\gamma})f|_{L^2(\gamma_{+})}^2+|w^{\ell}\mathcal{S}_{2,b}^{(1)}|_{L^{\infty}(\gamma_{-})}^2+\|\nu^{-1}w^{\ell}(\mathcal{S}_{1},\mathcal{S}_{2})\|_{L^{\infty}}^2\big].
\end{align}

\smallskip

\noindent\underline{\textbf{Estimates on $\|b_{3}\|_{L^2}$}}: Let
$$
\Psi_{b_{3}}=v_{\eta}v_{z}\sqrt{\mu}\phi_{b_{3}}'(\eta),
$$
where $\phi_{b_{3}}(\eta)$ satisfies
$$
\phi_{b_{3}}''(\eta)=b_{3}(\eta),\quad \phi_{b_{3}}'(0)=0,\quad \phi_{b_{3}}(0)=0.
$$
Then it is direct to get
$$
\phi_{b_{3}}'(\eta)=-\int_{0}^{\eta}b_{3}(y)\,{\rm d}y,\quad \phi_{b_{3}}(\eta)=-\int_{0}^{\eta}\int_{0}^{s}b_{3}(y)\,{\rm d}y{\rm d}ds.
$$
It follows from the H\"{o}lder inequality that
$$
\|\phi_{b_{3}}\|_{H^2}\lesssim \|b_{3}\|_{L^2}.
$$
Taking $\Psi=\Psi_{b_{3}}$ in \eqref{E15}, we have	
	\begin{align}\label{z1}
		-\la v_{\eta}f,\partial_{\eta}\Psi_{b_{3}}\ra&=-\epsilon\la f,\Psi_{b_{3}}\ra-(1-\sigma)\la \nu f,\Psi_{b_{3}}\ra-\int_{\R^3}v_{\eta}(\Psi_{b_{3}}f)(\eta_{1},v)\,{\rm d}v+\int_{\R^3}v_{\eta}(\Psi_{b_{3}}f)(0,v)\,{\rm d}v\nonumber\\
		&\,\,-\la \frac{1}{1-\eta}v_{\phi}^2f,\partial_{v_{\eta}}\Psi_{b_{3}}\ra+\la \frac{1}{1-\eta}v_{\eta}f,\Psi_{b_{3}}\ra+\la \frac{1}{1-\eta}v_{\eta}v_{\phi}f,\partial_{v_{\phi}}\Psi_{b_{3}}\ra\nonumber\\
		&\,\, -\frac{\alpha}{\eta_{1}}\la v_{\eta}f,\partial_{v_{\phi}}\Psi_{b_{3}}\ra-\frac{\alpha}{2\eta_{1}}\la v_{\eta}v_{\phi}f,\Psi_{b_{3}}\ra-\sigma\la {\bf L}f,\Psi_{b_{3}}\ra+\la \mathscr{S},\Psi_{b_{3}}\ra=\sum\limits_{k=1}^{11}Q_{k}.
	\end{align}
A direct calculation shows that
\begin{align*}
-\langle v_{\eta}f,\partial_{\eta}\Psi_{b_{3}}\rangle\geq (1-\kappa_{4})\|b_{3}\|_{L^3}-C_{\kappa_{4}}\|(\mathbf{I-P}f\|_{L_{\nu}^2},
\end{align*}
where $\kappa_{4}>0$ is a small parameter, which will be determined later.

Similar to \eqref{E19}--\eqref{E22}, we have
	\begin{align*}
    |Q_{1}|+|Q_{2}|+|Q_{3}|+|Q_{4}|\leq \kappa_{4}\|b_{3}\|_{L^2}^2+C_{\kappa_{4}}\big[\|(\mathbf{I-P})f\|_{L_{\nu}^2}^2+\|\nu^{-1}w^{\ell}(\mathcal{S}_{1},\mathcal{S}_{2})\|_{L^{\infty}}\big].
	\end{align*}
For the boundary terms $Q_{3}$ and $Q_{4}$, noting $\phi_{b_{3}}'(0)=0$, we have
\begin{align*}
Q_{3}+Q_{4}
&=-\phi_{b_{3}}'(\eta_{1})\int_{v_{\eta}>0}v_{\eta}^2v_{z}\sqrt{\mu}f(\eta_{1},v)\,{\rm d}v-\phi_{b_{3}}'(\eta_{1})\int_{v_{\eta}<0}v_{\eta}^2v_{z}\sqrt{\mu}f(\eta_{1},v)\,{\rm d}v\nonumber\\
&=-\phi_{b_{3}}'(\eta_{1})\int_{v_{\eta}<0}v_{\eta}^2v_{z}\sqrt{\mu}\mathcal{S}_{2,b}^{(2)}\,{\rm d}v-\phi_{b_{3}}'(\eta_{1})\int_{v_{\eta}>0}v_{\eta}^2v_{z}\sqrt{\mu}(I-P_{\gamma})f(\eta_{1},v)\,{\rm d}v\nonumber\\
&\leq \kappa_{4}\|b_{3}\|_{L^2}^2+C_{\kappa_{4}}|(I-P_{\gamma})f|_{L^2(\gamma_{+})}^2+C_{\kappa_{4}}|w^{\ell}\mathcal{S}_{2,b}^{(2)}|_{L^{\infty}(\gamma_{-})}^2.
\end{align*}
For the geometric terms $Q_{5}$, $Q_{6}$ and $Q_{7}$, a direct calculation shows that
\begin{align}\label{z8}
    Q_{5}+Q_{6}+Q_{7}
    &=\int_{0}^{\eta_{1}}\frac{1}{1-\eta}\phi_{b_{3}}'(\eta)\int_{\R^3}(v_{\eta}^2-v_{\phi}^2)v_{z}\sqrt{\mu}(\mathbf{I-P})f(\eta,v)\,{\rm d}v{\rm d}\eta\nonumber\\
    &\leq \kappa_{4}\|b_{3}\|_{L^2}+C_{\kappa_{4}}\|(\mathbf{I-P})f\|_{L_{\nu}^2}^2.
\end{align}
For the terms $Q_{8}$ and $Q_{9}$ induced by shear force, a direct calculation shows that
\begin{align}\label{z9}
Q_{8}+Q_{9}=0.
\end{align}
Plugging all the above estimates 
into \eqref{z1}, we obtain
\begin{align}\label{z10}
\|b_{3}\|_{L^2}^2&\lesssim \kappa_{4}\|b_{3}\|_{L^2}^2+C_{\kappa_{4}}\big[\|(\mathbf{I-P})f\|_{L_{\nu}^2}^2+|(I-P_{\gamma})f|_{L^2(\gamma_{+})}^2\big]\nonumber\\
&\quad +C_{\kappa_{4}}\big[|w^{\ell}\mathcal{S}_{2,b}^{(2)}|_{L^{\infty}(\gamma_{-})}^2+\|\nu^{-1}w^{\ell}(\mathcal{S}_{1},\mathcal{S}_{2})\|_{L^{\infty}}^2\big].
\end{align}

\smallskip

\noindent\underline{\textbf{Estimates on $\|c\|_{L^2}$}}: Let
	\begin{align*}
		\Psi_{c}(\eta,v)=-v_{\eta}(|v|^2-5)\sqrt{\mu}\phi_{c}'(\eta),
	\end{align*}
where $\phi_{c}(\eta)$ satisfies
$$
\phi_{c}''(\eta)=(1-\eta)c(\eta),\quad \phi_{c}'(0)=0,\quad \phi_{c}(0)=0.
$$
Then it is direct to get
$$
\phi_{c}'(\eta)=\int_{0}^{\eta}(1-y)c(y)\,{\rm d}y,\quad \phi_{c}(\eta)=\int_{0}^{\eta}\int_{0}^{s}(1-y)c(y)\,{\rm d}y.
$$
It follows from the H\"{o}lder inequality that
	$$
	\|\phi_{c}\|_{H^2}\lesssim \|c\|_{L^2}.
	$$
Taking $\Psi=\Psi_{c}$ in \eqref{E15}, we have
	\begin{align}\label{E48}
		-\la v_{\eta}f,\partial_{\eta}\Psi_{c}\ra&=-\epsilon\la f,\Psi_{c}\ra-(1-\sigma)\la \nu f,\Psi_{c}\ra-\int_{\R^3}v_{\eta}(\Psi_{c}f)(\eta_{1},v)\,{\rm d}v+\int_{\R^3}v_{\eta}(\Psi_{c}f)(0,v)\,{\rm d}v\nonumber\\
		&\quad-\la \frac{1}{1-\eta}v_{\phi}^2f,\partial_{v_{\eta}}\Psi_{c}\ra+\la \frac{1}{1-\eta}v_{\eta}f,\Psi_{c}\ra+\la \frac{1}{1-\eta}v_{\eta}v_{\phi}f,\partial_{v_{\phi}}\Psi_{c}\ra\nonumber\\
		&\quad -\frac{\alpha}{\eta_{1}}\la v_{\eta}f,\partial_{v_{\phi}}\Psi_{c}\ra-\frac{\alpha}{2\eta_{1}}\la v_{\eta}v_{\phi}f,\Psi_{c}\ra-\sigma\la {\bf L}f,\Psi_{c}\ra+\la \mathscr{S},\Psi_{c}\ra=\sum\limits_{k=1}^{11}S_{k}.
	\end{align}
	A direct calculation shows that
	\begin{align*}
		-\la v_{\eta}f,\partial_{\eta}\Psi_{c}\ra
		&\geq (1-\kappa_{5})\|c\|_{L^2}^2-C_{\kappa_{5}}\|(\mathbf{I-P})f\|_{L_{\nu}^2}^2,
	\end{align*}
where $\kappa_{5}>0$ is a small parameter, which will be determined later.

\smallskip

Similar to \eqref{E19}--\eqref{E22}, we have
	\begin{align*}
    |S_{1}|+|S_{2}|+|S_{10}|+|S_{11}|&\leq \kappa_{5}\|c\|_{L^2}^2+C(1-\sigma_{0})\|c\|_{L^2}\|b_{1}\|_{L^2}\nonumber\\
    &\quad +C_{\kappa_{5}}\big[\|(\mathbf{I-P})f\|_{L_{\nu}^2}^2+\|\nu^{-1}w^{\ell}(\mathcal{S}_{1},\mathcal{S}_{2})\|_{L^{\infty}}^2\big].
	\end{align*}
For the boundary terms $S_{3}$ and $S_{4}$,  noting $\phi_{c}'(0)=0$, we have
\begin{align*}
	S_{3}+S_{4}
	&=-\phi_{c}'(\eta_{1})\int_{v_{\eta>0}}v_{\eta}^2(|v|^2-5)\sqrt{\mu}f(\eta_{1},v)\,{\rm d}v-\phi_{c}'(\eta_{1})\int_{v_{\eta}<0}v_{\eta}^2(|v|^2-5)\sqrt{\mu}f(\eta_{1},v)\,{\rm d}v\nonumber\\
    &=-\phi_{c}'(\eta_{1})\int_{v_{\eta}<0}v_{\eta}^2\sqrt{\mu}\mathcal{S}_{2,b}^{(2)}\,{\rm d}v-\phi_{c}'(\eta_{1})\int_{v_{\eta}<0}v_{\eta}^2v_{\phi}\sqrt{\mu}(I-P_{\gamma})f(\eta_{1},v)\,{\rm d}v\nonumber\\
	&\leq \kappa_{5}\|c\|_{L^2}^2+C_{\kappa_{5}}|(I-P_{\gamma})f|_{L^2(\gamma_{+})}^2+C_{\kappa_{5}}|w^{\ell}\mathcal{S}_{2,b}^{(2)}|_{L^{\infty}(\gamma_{-})}^2.
\end{align*}
For the geometric terms $S_{5}+S_{6}+S_{7}$, one has
	\begin{align}\label{E55}
		S_{5}+S_{6}+S_{7}&=\int_{0}^{\eta_{1}}\frac{1}{1-\eta}\phi_{c}'(\eta)\int_{\R^3}(v_{\eta}^2-v_{\phi}^2)(|v|^2-5)\sqrt{\mu}(\mathbf{I-P})f\,{\rm d}v{\rm d}\eta\nonumber\\
		&\lesssim \kappa_{5}\|c\|_{L^2}^2+C_{\kappa_{5}}\|(\mathbf{I-P})f\|_{L_{\nu}^2}^2.	
	\end{align}
For the terms $S_{8}$ and $S_{9}$ induced by shear force, we have
\begin{align}\label{E56}
	S_{8}+S_{9}&=\frac{2\alpha}{\eta_{1}}\int_{0}^{\eta_{1}}\phi_{c}'(\eta)\int_{\R^3}v_{\eta}^2v_{\phi}\sqrt{\mu}f\,{\rm d}v{\rm d}\eta\lesssim \alpha\|c\|_{L^2}\|b_{2}\|_{L^2}+\kappa_{5}\|c\|_{L^2}^2+C_{\kappa_{5}}\|(\mathbf{I-P})f\|_{L_{\nu}^2}^2.
\end{align}
Plugging all the above estimates into \eqref{E48}, 
we obtain
\begin{align}\label{E57}
\|c\|_{L^2}&\lesssim \kappa_{5}\|c\|_{L^2}^2+\alpha\|c\|_{L^2}\|b_{2}\|_{L^2}+(1-\sigma_{0})\|c\|_{L^2}\|b_{1}\|_{L^2}+C_{\kappa_{5}}\|(\mathbf{I-P})f\|_{L_{\nu}^2}^2\nonumber\\
&\quad +C_{\kappa_{5}}\big[|(I-P_{\gamma})f|_{L^2(\gamma_{+})}^2+|w^{\ell}\mathcal{S}_{2,b}^{(2)}|_{L^{\infty}(\gamma_{-})}^2+\|\nu^{-1}w^{\ell}(\mathcal{S}_{1},\mathcal{S}_{2})\|_{L^{\infty}}^2\big].
\end{align}

\medskip

\noindent\underline{\textbf{Conclusion on the estimates of $\|\mathbf{P}f_{2}\|_{L^2}$}}: Finally, combining \eqref{E27}, \eqref{E37}, \eqref{E47} and \eqref{E57}, we obtain
\begin{align}\label{E58}
\|(a,b_{1},b_{2},b_{3},c)\|_{L^2}^2
&\lesssim \kappa_{1}\|a\|_{L^2}^2+C_{\kappa_{1}}\kappa_{2}\|b_{1}\|_{L^2}^2+\kappa_{3}\|b_{2}\|_{L^2}^2+\kappa_{4}\|b_{3}\|_{L^2}^2+C_{\kappa_{1},\kappa_{2}}\kappa_{5}\|c\|_{L^2}^2+\kappa_{1}\|f_{2}\|_{L_{\nu}^2}^2\nonumber\\
&\quad +C_{\kappa}(\alpha+1-\sigma_{0})(\|a\|_{L^2}+\|c\|_{L^2})(\|b_{1}\|_{L^2}+\|b_{2}\|_{L^2})\nonumber\\
&\quad +C_{\kappa}\big[\|(\mathbf{I-P})f\|_{L_{\nu}^2}^2+|(I-P_{\gamma})f|_{L^2(\gamma_{+})}^2+|w^{\ell}f_{1}|_{L^{\infty}(\gamma_{+})}^2\big]\nonumber\\
&\quad +C_{\epsilon,\kappa}\big[|w^{\ell}(\mathcal{S}_{1,b}^{(1)},\mathcal{S}_{2,b}^{(2)})|_{L^{\infty}(\gamma_{-})}^2+\|\nu^{-1}w^{\ell}(\mathcal{S}_{1},\mathcal{S}_{2})\|_{L^{\infty}}^2\big],
\end{align}
where $C_{\kappa}$ means a positive constant depending on $\kappa_{1},\cdots,\kappa_{5}$. Furthermore, recalling $\sqrt{\mu}f=f_{1}+\sqrt{\mu}f$ and the definition of $P_{\gamma}f$ in \eqref{3.3}, we have
\begin{align}\label{E59}
\begin{aligned}
&\|\mathbf{P}f_{2}\|_{L^2}^2\lesssim \|\mathbf{P}f\|_{L^2}^2+\|w^{\ell}f_{1}\|_{L^{\infty}}^2\lesssim \|(a,b_{1},b_{2},b_{3},c)\|_{L^2}^2+\|w^{\ell}f_{1}\|_{L^{\infty}}^2,\\
&\|f_{2}\|_{L_{\nu}^2}^2\lesssim \|\mathbf{P}f_{2}\|_{L^2}^2+\|(\mathbf{I-P})f_{2}\|_{L_{\nu}^2}^2\\
&\qquad\,\,\,\lesssim \|(a,b_{1},b_{2},b_{3},c)\|_{L^2}^2+\|w^{\ell}f_{1}\|_{L^{\infty}}^2+\|(\mathbf{I-P})f_{2}\|_{L_{\nu}^2}^2,\\
&|(I-P_{\gamma})f|_{L^2(\gamma_{+})}^2\lesssim |(I-P_{\gamma})f_{2}|_{L^2(\gamma_{+})}^2+\|w^{\ell}f_{1}\|_{L^{\infty}}^2,\\
&\|(\mathbf{I-P})f\|_{L_{\nu}^2}^2\lesssim \|(\mathbf{I-P})f_{2}\|_{L_{\nu}^2}^2+\|w^{\ell}f_{1}\|_{L^{\infty}}^2,
\end{aligned}
\end{align}
provided that $\ell> 4$. Combining \eqref{E58} and \eqref{E59}, we have
\begin{align}\label{E60}
&\|(a,b_{1},b_{2},b_{3},c)\|_{L^2}^2\nonumber\\
&\lesssim \kappa_{1}\|a\|_{L^2}^2+(C_{\kappa_{1}}\kappa_{2}+\kappa_{1})\|b_{1}\|_{L^2}^2+(\kappa_{3}+\kappa_{1}+\kappa_{4})\|(b_{2},b_{3})\|_{L^2}^2+(C_{\kappa_{1},\kappa_{2}}\kappa_{5}+\kappa_{1})\|c\|_{L^2}^2\nonumber\\
&\quad +C_{\kappa}(\alpha+1-\sigma_{0})\|(a,b_{1},b_{2},c)\|_{L^2}^2+C_{\epsilon,\kappa}\big[|w^{\ell}(\mathcal{S}_{1,b}^{(1)},\mathcal{S}_{2,b}^{(2)})|_{L^{\infty}(\gamma_{-})}^2+\|\nu^{-1}w^{\ell}(\mathcal{S}_{1},\mathcal{S}_{2})\|_{L^{\infty}}^2\big]\nonumber\\
&\quad +C_{\kappa}\big[\|(\mathbf{I-P})f_{2}\|_{L_{\nu}^2}^2+|(I-P_{\gamma})f_{2}|_{L^2(\gamma_{+})}^2+\|w^{\ell}f_{1}\|_{L^{\infty}}^2\big].
\end{align}
Choosing $\kappa_{1},\kappa_{3},\kappa_{4}$ small enough such that $\kappa_{1}+\kappa_{3}+\kappa_{4}\ll \frac{1}{16}$ first, and then choosing $\kappa_{2}$ small enough such that $C_{\kappa_{1}}\kappa_{2}+\kappa_{1}\ll\frac{1}{16}$, and next choosing $\kappa_{5}$ small enough such that $C_{\kappa_{1},\kappa_{2}}\kappa_{5}+\kappa_{1}\ll \frac{1}{16}$, and finally taking $\alpha$ and $1-\sigma_{0}$ small enough such that $C_{\kappa}(\alpha+1-\sigma_{0})\ll \frac{1}{16}$, we obtain from
\eqref{E60} that
\begin{align}\label{E61}
\|(a,b_{1},b_{2},b_{3},c)\|_{L^2}^2
&\leq C\big[\|(\mathbf{I-P})f_{2}\|_{L_{\nu}^2}^2+|(I-P_{\gamma})f|_{L^2(\gamma_{+})}^2+\|w^{\ell}f_{1}\|_{L^{\infty}}^2\big]\nonumber\\
&\quad+C_{\epsilon}|w^{\ell}\big(\mathcal{S}_{1,b}^{(1)},\mathcal{S}_{2,b}^{(2)}\big)|_{L^{\infty}(\gamma_{-})}^2+\|\nu^{-1}w^{\ell}\big(\mathcal{S}_{1},\mathcal{S}_{2}\big)\|_{L^{\infty}}^2\big],
\end{align}
which, together with \eqref{E60}, yields that
\begin{align}\label{E61-1}
\|\mathbf{P}f_{2}\|_{L^2}^2&\leq C\big[\|(\mathbf{I-P})f_{2}\|_{L_{\nu}^2}^2+|(I-P_{\gamma})f_{2}|_{L^2(\gamma_{+})}^2+\|w^{\ell}f_{1}\|_{L^{\infty}}^2\big]\nonumber\\
&\quad+C_{\epsilon}\big[|w^{\ell}\big(\mathcal{S}_{1,b}^{(1)},\mathcal{S}_{2,b}^{(2)}\big)|_{L^{\infty}(\gamma_{-})}^2+\|\nu^{-1}w^{\ell}\big(\mathcal{S}_{1},\mathcal{S}_{2}\big)\|_{L^{\infty}}^2\big].
\end{align}
Combining \eqref{E5-1}, \eqref{E10-1} and \eqref{E61-1}, we obtain
\begin{align}\label{E62}
\|f_{2}\|_{L_{\nu}^2}^2+|(I-P_{\gamma})f_{2}|_{L^2(\gamma_{+})}^2&\lesssim \|\mathbf{P}f_{2}\|_{L^{2}}^2+\|(\mathbf{I-P})f_{2}\|_{L_{\nu}^2}^2+|(I-P_{\gamma})f_{2}|_{L^2(\gamma_{+})}^2\nonumber\\
&\lesssim \delta\big[\|f_{2}\|_{L_{\nu}^2}^2+|(I-P_{\gamma})f_{2}|_{L^2(\gamma_{+})}^2\big]+C_{\delta}\|w^{\ell}f_{1}\|_{L^{\infty}}^2\nonumber\\
&\quad +C_{\epsilon,\delta}\big[\|\nu^{-1}w^{\ell}\big(\mathcal{S}_{1},\mathcal{S}_{2}\big)\|_{L^{\infty}}^2+|w^{\ell}\big(\mathcal{S}_{2,b}^{(1)},\mathcal{S}_{2,b}^{(2)}\big)|_{L^{\infty}(\gamma_{-})}^2\big].
\end{align}
Taking $\delta$ small enough, we obtain from \eqref{E62} that
\begin{align}\label{E63}
\|f_{2}\|_{L_{\nu}^2}^2 \leq C\|w^{\ell}f_{1}\|_{L^{\infty}}^2+C_{\epsilon}\big[\|\nu^{-1}w^{\ell}\big(\mathcal{S}_{1},\mathcal{S}_{2}\big)\|_{L^{\infty}}^2+|w^{\ell}\big(\mathcal{S}_{2,b}^{(1)},\mathcal{S}_{2,b}^{(2)}\big)|_{L^{\infty}(\gamma_{-})}^2\big],
\end{align}
which, together with \eqref{E10-1}, concludes the proof of \eqref{B1}. $\hfill\square$

\begin{lemma}\label{lem2.5}
For any given $\epsilon\in (0,1]$ and $\sigma\in [0,1]$, let $(f_{1},f_{2})$ be the solution of the coupled systems \eqref{3.6-2}--\eqref{3.7-2}. There exists a positive constant $\alpha_{2}$ small enough, which is independent of $\epsilon$ and $\sigma$, such that if $\alpha\in (0,\alpha_{2})$, it holds that
    \begin{align}\label{E64}
  \|w^{\ell}\big(f_{1},f_{2}\big)\|_{L^{\infty}}\leq  C_{\epsilon,\ell}\big[\|\nu^{-1}w^{\ell}\big(\mathcal{S}_{1},\mathcal{S}_{2}\big)\|_{L^{\infty}}+|w^{\ell}\big(\mathcal{S}_{2,b}^{(1)},\mathcal{S}_{2,b}^{(2)}\big)|_{L^{\infty}(\gamma_{-})}\big],
   \end{align}
   where $C_{\epsilon,\ell}$ is a positive constant depending only $\epsilon$ and the weight $\ell$.
\end{lemma}

\noindent\textbf{Proof}. Substituting \eqref{B1} into \eqref{3.17-2}, we get
\begin{align}
&\|w^{\ell}f_{2}\|_{L^{\infty}}\leq C_{\ell}\|w^{\ell}f_{1}\|_{L^{\infty}}+C_{\epsilon,\ell}\|w^{\ell}\big(\mathcal{S}_{1},\mathcal{S}_{2}\big)\|_{L^{\infty}}+C_{\epsilon,\ell}|w^{\ell}\big(\mathcal{S}_{2,b}^{(1)},\mathcal{S}_{2,b}^{(2)}\big)|_{L^{\infty}(\gamma_{-})}.\label{E65}
\end{align}
Then, substituting \eqref{E65} into \eqref{3.17-1}, we obtain
\begin{align}\label{E66}
    \|w^{\ell}f_{1}\|_{L^{\infty}}\leq C_{\ell}\alpha\|w^{\ell}f_{1}\|_{L^{\infty}}+C_{\epsilon,\ell}\|w^{\ell}\big(\mathcal{S}_{1},\mathcal{S}_{2}\big)\|_{L^{\infty}}+C_{\epsilon,\ell}|w^{\ell}\big(\mathcal{S}_{2,b}^{(1)},\mathcal{S}_{2,b}^{(2)}\big)|_{L^{\infty}(\gamma_{-})}.
\end{align}
We conclude \eqref{E64} from \eqref{E65}--\eqref{E66} provided that $\alpha$ small enough. $\hfill\square$

\section{Existence of the linearized approximate system}\label{sec5}

In this section, we shall use the {\it a priori} $L^{\infty}$ uniform estimates established in Lemma \ref{lem2.5} to prove the existence of the linearized approximate systems \eqref{3.1}--\eqref{3.2} by taking the limits $\sigma\to 1$ and $\epsilon \to 0$ in \eqref{3.6-2}--\eqref{3.7-2}.

\subsection{The limit $\sigma\to 1$} 
With the help of {\it a priori} $L^{\infty}$ estimate, we are able to use a bootstrap argument to show the existence of solutions for the coupled systems \eqref{3.6-2}--\eqref{3.7-2} with $\sigma=1$. To this end, we introduce the following functional space
\begin{align*}
\mathbf{X}=\big\{(f_{1},f_{2})\,|\,&\|w^{\ell}(f_{1},f_{2})\|_{L^{\infty}}<\infty\big\},
\end{align*}
equipped with norm
\begin{align*}
\|(f_{1},f_{2})\|_{\mathbf{X}}=\|w^{\ell}(f_{1},f_{2})\|_{L^{\infty}}.
\end{align*}
We define the linear vector operator parameterized by $\sigma\in [0,1]$ as follows:
$$
\mathcal{L}_{\sigma}(f_{1},f_{2})=\big(\mathcal{L}_{\sigma}^{1}(f_{1},f_{2}),\mathcal{L}_{\sigma}^{2}(f_{1},f_{2})\big)^{t},
$$
with
\begin{align*}
&\mathcal{L}_{\sigma}^{1}(f_{1},f_{2})=\left\{
\begin{aligned}
&\epsilon f_{1}+v_{\eta}\frac{\partial f}{\partial\eta}-\frac{1}{1-\eta}\Big(v_{\phi}^2\frac{\partial f_{1}}{\partial v_{\eta}}-v_{\eta}v_{\phi}\frac{\partial f_{1}}{\partial v_{\phi}}\Big)-\frac{\alpha}{\eta_{1}}v_{\eta}\frac{\partial f_{1}}{\partial v_{\phi}}+\frac{\alpha}{2\eta_{1}}v_{\eta}v_{\phi}\sqrt{\mu}f_{2}\\
&\qquad +\nu f_{1}-\sigma(1-\chi_{M})\mathcal{K}f_{1},\qquad \eta\in (0,\eta_{1}),\\
&f_{1}(0,v_{\eta},v_{\phi},v_{z})\vert_{v_{\eta}<0},\qquad\,\,\eta=0,\\
& f_{1}(\eta_{1},v_{\eta},v_{\phi},v_{z})\vert_{v_{\eta}>0},\qquad\eta=\eta_{1},
\end{aligned}
\right.\\
&\mathcal{L}_{\sigma}^{2}(f_{1},f_{2})=\left\{
\begin{aligned}
&\epsilon f_{2}+v_{\eta}\frac{\partial f_{2}}{\partial\eta}-\frac{1}{1-\eta}\Big(v_{\phi}^2\frac{\partial f_{2}}{\partial v_{\eta}}-v_{\eta}v_{\phi}\frac{\partial f_{2}}{\partial v_{\phi}}\Big)-\frac{\alpha}{\eta_{1}}v_{\eta}\frac{\partial f_{2}}{\partial v_{\phi}}\\
&\qquad +\nu f_{2}-\sigma\chi_{M}\mathcal{K}f_{1}-\sigma Kf_{2},\qquad \eta\in (0,\eta_{1}),\\
&f_{2}(0,v_{\eta},v_{\phi},v_{z})\vert_{v_{\eta}>0}-\sqrt{2\pi\mu}\int_{u_{\eta}<0}(\sqrt{\mu}f_{2}+\sigma f_{1})(0,u)|u_{\eta}|\,{\rm d}u,\qquad\,\,\,\,\eta=0,\\
&f_{2}(\eta_{1},v_{\eta},v_{\phi},v_{z})\vert_{v_{\eta}<0}-\sqrt{2\pi\mu}\int_{u_{\eta}>0}(\sqrt{\mu}f_{2}+\sigma f_{1})(\eta_{1},u)|u_{\eta}|\,{\rm d}u,\qquad \eta=\eta_{1}.
\end{aligned}
\right.
\end{align*}
Hereafter, the superscript ``$t$" means the transpose of vectors. Then we can rewrite the coupled systems \eqref{3.6-2}--\eqref{3.7-2} as
\begin{align}\label{4.14-3}
\mathcal{L}_{\sigma}(f_{1},f_{2})=(\mathcal{S}_{1},0,0,\mathcal{S}_{2},\mathcal{S}_{2,b}^{(1)},\mathcal{S}_{2,b}^{(2)})^{t}.
\end{align}

\begin{lemma}\label{lem3.1}
Let $\epsilon\in (0,1]$ and $\alpha_{2}$ be the positive constant in Lemma \ref{lem2.5}. Then for $\alpha\in (0,\alpha_{2})$, there exists a unique solution $(f_{1},f_{2})\in \mathbf{X}$ to the system \eqref{4.14-3} with $\sigma=1$ satisfying
\begin{align}\label{4.4}
\|w^{\ell}(f_{1},f_{2})\|_{L^{\infty}}\leq C_{\epsilon,\ell}\big[\|\nu^{-1}w^{\ell}(\mathcal{S}_{1},\mathcal{S}_{2})\|_{L^{\infty}}+|w^{\ell}(\mathcal{S}_{2,b}^{(1)},\mathcal{S}_{2,b}^{(2)})|_{L^{\infty}(\gamma_{-})}\big].
\end{align}
\end{lemma}

\noindent\textbf{Proof}. We divide the proof into four steps.

\noindent Step 1. {\it $\sigma=0$}. We construct the approximate sequences $\{(f_{1}^{i},f_{2}^{i})\}_{i=1}^{\infty}$
by the following iterative systems
	\begin{equation*}
		\left\{\begin{aligned}
			&\epsilon f_{1}^{i}+v_{\eta}\frac{\partial f_{1}^{i}}{\partial \eta}-\frac{1}{1-\eta}\Big(v_{\phi}^2\frac{\partial f_1^{i}}{\partial v_{\eta}}-v_{\eta}v_{\phi}\frac{\partial f_{1}^{i}}{\partial v_{\phi}}\Big)-\frac{\alpha}{\eta_{1}} v_{\eta}\frac{\partial f_{1}^{i}}{\partial v_{\phi}}+\frac{\alpha}{2\eta_{1}}v_{\eta}v_{\phi}\sqrt{\mu}f_{2}^{i}+\nu f_{1}^{i}=\mathcal{S}_{1},\\
			&f_{1}^{i}(0,v_{\eta},v_{\phi},v_{z})\vert_{v_{\eta}>0}=0,\quad f_{1}^{i}(\eta_{1},v_{\eta},v_{\phi},v_{z})\vert_{v_{\eta}<0}=0,
		\end{aligned}\right.
	\end{equation*}
	and
	\begin{equation*}
		\left\{\begin{aligned}
			&\epsilon f_{2}^{i}+v_{\eta}\frac{\partial f_{2}^{i}}{\partial \eta}-\frac{1}{1-\eta}\Big(v_{\phi}^2\frac{\partial f_2^{i}}{\partial v_{\eta}}-v_{\eta}v_{\phi}\frac{\partial f_2^{i}}{\partial v_{\phi}}\Big)-\frac{\alpha}{\eta_{1}} v_{\eta}\frac{\partial f_2^{i}}{\partial v_{\phi}}+\nu f_2^{i}=\mathcal{S}_{2},\\
			&f_2^{i}(0,v_{\eta},v_{\phi},v_{z})\vert_{v_{\eta}>0}=\sqrt{2\pi\mu}\int_{u_{\eta}<0}(\sqrt{\mu}f_2^{i-1})(0,u)|u_{\eta}|\,\mathrm{d}u+\mathcal{S}_{2,b}^{(1)},\\
			&f_2^{i}(\eta_{1},v_{\eta},v_{\phi},v_{z})\vert_{v_{\eta}<0}=\sqrt{2\pi\mu}\int_{u_{\eta}>0}(\sqrt{\mu}f_2^{i-1})(\eta_{1},u)|u_{\eta}|\,\mathrm{d}u+\mathcal{S}_{2,b}^{(2)},
		\end{aligned}\right.
	\end{equation*}
with $[f_{1}^{0},f_{2}^{0}]=[0,0]$. Let $[h_{1}^{i},h_{2}^{i}]=w^{\ell}[f_{1}^{i},f_{2}^{i}]$, then $[h_{1}^{i},h_{2}^{i}]$ satisfies
\begin{equation*}
		\left\{\begin{aligned}
			&\epsilon h_{1}^{i}+v_{\eta}\frac{\partial h_{1}^{i}}{\partial \eta}-\frac{1}{1-\eta}\Big(v_{\phi}^2\frac{\partial h_1^{i}}{\partial v_{\eta}}-v_{\eta}v_{\phi}\frac{\partial h_{1}^{i}}{\partial v_{\phi}}\Big)-\frac{\alpha}{\eta_{1}} v_{\eta}\frac{\partial h_{1}^{i}}{\partial v_{\phi}}+\frac{\alpha}{2\eta_{1}}v_{\eta}v_{\phi}\sqrt{\mu}h_{2}^{i}+\nu h_{1}^{i}\\
			&\qquad =-\frac{\alpha}{\eta_{1}}\frac{\ell v_{\eta}v_{\phi}}{1+|v|^2}h_{1}^{i}+w^{\ell}\mathcal{S}_{1},\\
			&h_{1}^{i}(0,v_{\eta},v_{\phi})\vert_{v_{\eta}>0}=0,\quad h_{1}^{i}(\eta_{1},v_{\eta},v_{\phi})\vert_{v_{\eta}<0}=0,
		\end{aligned}\right.
	\end{equation*}
	and
	\begin{equation*}
		\left\{\begin{aligned}
			&\epsilon h_{2}^{i}+v_{\eta}\frac{\partial h_{2}^{i}}{\partial \eta}-\frac{1}{1-\eta}\Big(v_{\phi}^2\frac{\partial h_2^{i}}{\partial v_{\eta}}-v_{\eta}v_{\phi}\frac{\partial h_2^{i}}{\partial v_{\phi}}\Big)-\frac{\alpha}{\eta_{1}} v_{\eta}\frac{\partial h_2^{i}}{\partial v_{\phi}}+\nu h_2^{i}=-\frac{\alpha}{\eta_{1}}\frac{\ell v_{\eta}v_{\phi}}{1+|v|^2}h_{2}^{i}+w^{\ell}\mathcal{S}_{2},\\
			&h_2^{i}(0,v_{\eta},v_{\phi})\vert_{v_{\eta}>0}=\frac{1}{\tilde{w}(v)}\int_{u_{\eta}<0}h_2^{i-1}(0,u)\tilde{w}(u)\,\mathrm{d}\sigma+w^{\ell}\mathcal{S}_{2,b}^{(1)},\\
			&h_2^{i}(\eta_{1},v_{\eta},v_{\phi})\vert_{v_{\eta}<0}=\frac{1}{\tilde{w}(v)}\int_{u_{\eta}>0}h_2^{i-1}(\eta_{1},u)\tilde{w}(u)\,\mathrm{d}\sigma+w^{\ell}\mathcal{S}_{2,b}^{(2)},
		\end{aligned}\right.
	\end{equation*}
with $[h_{1}^{0},h_{2}^{0}]=[0,0]$. Once $[f_{1}^{i-1},f_{2}^{i-1}]\in \mathbf{X}$ is given, along the characteristic line defined in \eqref{3.16}, it is easy to check that 
\begin{align*}
   &\|h_{1}^{i}\|_{L^{\infty}}\leq C\alpha\|h_{2}^{i}\|_{L^{\infty}}+\|\nu^{-1}w^{\ell}\mathcal{S}_{1}\|_{L^{\infty}},\\ &\|h_{2}^{i}\|_{L^{\infty}}\leq C\|h_{2}^{i-1}\|_{L^{\infty}}+C\|\nu^{-1}w^{\ell}\mathcal{S}_{2}\|_{L^{\infty}}+C|w^{\ell}(\mathcal{S}_{2,b}^{(1)},\mathcal{S}_{2,b}^{(2)})|_{L^{\infty}(\gamma_{-})},
\end{align*}
which, together with induction method, yields that
\begin{align*}
    \|[h_{1}^{i},h_{2}^{i}]\|_{L^{\infty}}\leq C_{i}\|\nu^{-1}w^{\ell}(\mathcal{S}_{1},\mathcal{S}_{2})\|_{L^{\infty}}+C_{i}|w^{\ell}(\mathcal{S}_{2,b}^{(1)},\mathcal{S}_{2,b}^{(2)})|_{L^{\infty}(\gamma_{-})}<\infty,\quad \forall i\geq 1.
\end{align*}
In order to get uniform estimates with respect to $i$, we can represent $[h_{1}^{i},h_{2}^{i}]$ as those similar to \eqref{3.17-4}, \eqref{3.18-3}, and \eqref{3.19-3}--\eqref{3.19-4}. In fact, for example, for $(\eta,v)\in (\mathcal{I}\times \R^3\backslash(\gamma_{0}\cup \gamma_{-}))\cap \mathcal{A}_{1}$, we only need to set $\sigma=0$ and replace $h_{1}$ with $h_{1}^{i}$, $h_{2}$ with $h_{2}^{i}$ in \eqref{3.17-4}, $h_{2}$ on the left hand side of \eqref{3.18-3} with $h_{2}^{i}$, $h_{2}$ with $h_{2}^{i-1-l}$ in $J_{14}$, and $h_{2}$ with $h_{2}^{i-1-k}$ in $J_{15}$. Then, by similar calculations as in \eqref{3.18-4}--\eqref{C17-1}, we can obtain
\begin{align}
\|h_{1}^{i}\|_{L^{\infty}}&\leq C\alpha\|h_{2}^{i}\|_{L^{\infty}}+C\|w^{\ell}\mathcal{S}_{1}\|_{L^{\infty}},\label{4.6}\\
\|h_{2}^{i}\|_{L^{\infty}}&\leq C_{\ell}\big[ke^{-\nu_{0}t}+\Big(\frac{1}{2}\Big)^{C_{4}T_{0}^{\frac{5}{4}}}\big]\sup_{1\leq j\leq k}\|h_{2}^{i-1-j}\|_{L^{\infty}}\nonumber\\
&\quad +C_{\ell}k\big[\|\nu^{-1}w^{\ell}\mathcal{S}_{2}\|_{L^{\infty}}+|w^{\ell}(\mathcal{S}_{2,b}^{(1)},\mathcal{S}_{2,b}^{(2)})|_{L^{\infty}(\gamma_{-})}\big].
\label{4.7-1}
\end{align}
Then, taking $k=k_{0}:=C_{3}T_{0}^{\frac{5}{4}}$, $t=T_{0}$ and $T_{0}$ large enough such that $C_{\ell}\big[ke^{-\nu_{0}t}+(\frac{1}{2})^{C_{4}T_{0}^{\frac{5}{4}}}\big]\leq \frac{1}{8}$, we obtain from \eqref{4.7-1} and Lemma \ref{lemC1} that
\begin{align}\label{4.7}
\|h_{2}^{i}\|_{L^{\infty}}&\leq \Big(\frac{1}{8}\Big)^{\big[\frac{i}{k_{0}+1}\big]}\sup_{1\leq j\leq 2k_{0}}\|h_{2}^{j}\|_{L^{\infty}}+\frac{8+k_{0}}{7}C_{\ell}k_{0}\big[\|\nu^{-1}w^{\ell}\mathcal{S}_{2}\|_{L^{\infty}}+|w^{\ell}(\mathcal{S}_{2,b}^{(1)},\mathcal{S}_{2,b}^{(2)})|_{L^{\infty}(\gamma_{-})}\big]
\end{align}
for $i\geq k_{0}+1$. It then follows from \eqref{4.7} that
\begin{align}\label{4.8}
\|h_{2}^{i}\|_{L^{\infty}}\leq C_{\ell}\|\nu^{-1}w^{\ell}\mathcal{S}_{2}\|_{L^{\infty}}+C_{\ell}|w^{\ell}(\mathcal{S}_{2,b}^{(1)},\mathcal{S}_{2,b}^{(2)})|_{L^{\infty}(\gamma_{-})}\quad \text{for $ i\geq 1$}.
\end{align}
Substituting \eqref{4.8} into \eqref{4.6}, we have
\begin{align}\label{4.9}
\|h_{1}^{i}\|_{L^{\infty}}\leq C_{\ell}\|\nu^{-1}w^{\ell}[\mathcal{S}_{1},\mathcal{S}_{2}\|_{L^{\infty}}+C_{\ell}|w^{\ell}(\mathcal{S}_{2,b}^{(1)},\mathcal{S}_{2,b}^{(2)})|_{L^{\infty}(\gamma_{-})}\quad \text{for $ i\geq 1$}.
\end{align}
Considering the difference $[h_{1}^{i+1}-h_{1}^{i},h_{2}^{i+1}-h_{2}^{i}]$, and applying the similar calculations as in \eqref{4.6}--\eqref{4.8}, one has
\begin{align}
&\|h_{1}^{i+1}-h_{1}^{i}\|_{L^{\infty}}\leq 
C\alpha\|h_{2}^{i+1}-h_{2}^{i}\|_{L^{\infty}},\label{4.11}\\
&\|h_{2}^{i+1}-h_{2}^{i}\|_{L^{\infty}}\leq \frac{1}{8}\sup_{1\leq j\leq k_{0}}\|h_{2}^{i-j}-h_{2}^{i-1-j}\|_{L^{\infty}}.\label{4.10}
\end{align}
Applying Lemma \ref{lemC1} and \eqref{4.8} to \eqref{4.10} yields that
\begin{align*}
\|h_{2}^{i+1}-h_{2}^{i}\|_{L^{\infty}}&\leq \Big(\frac{1}{8}\Big)^{\big[\frac{i}{k_{0}+1}\big]}\sup_{1\leq j\leq 2k_{0}}\|h_{2}^{j}\|_{L^{\infty}}\leq \Big(\frac{1}{8}\Big)^{\big[\frac{i}{k_{0}+1}\big]}C_{\ell}\big[\|\nu^{-1}w^{\ell}\mathcal{S}_{2}\|_{L^{\infty}}+|w^{\ell}(\mathcal{S}_{2,b}^{(1)},\mathcal{S}_{2,b}^{(2)})|_{L^{\infty}(\gamma_{-})}\big],
\end{align*}
for any $i\geq k_{0}+1$, which implies immediately that $\{h_{2}^{i}\}$ is a Cauchy sequence in $L^{\infty}$. Furthermore, we obtain from \eqref{4.11} that $\{h_{1}^{i}\}$ is also a Cauchy sequence in $L^{\infty}$. Therefore, there exists a unique limit function pair $(h_{1},h_{2})$ such that
\begin{align}\label{4.12}
    \lim\limits_{i\to \infty}\|(h_{1}^{i}-h_{1},h_{2}^{i}-h_{2})\|_{L^{\infty}}=0.
\end{align}
Let $(f_{1},f_{2})=w^{-\ell}(h_{1},h_{2})$, then $f_{1}$ and $f_{2}$ will satisfy
\begin{align*}
\mathcal{L}_{0}(f_{1},f_{2})=(\mathcal{S}_{1},0,0,\mathcal{S}_{2},\mathcal{S}_{2,b}^{(1)},\mathcal{S}_{2,b}^{(2)})^{t}.
\end{align*}
Moreover, it follows from \eqref{4.8}--\eqref{4.9} and \eqref{4.12} that $(f_{1},f_{2})\in \mathbf{X}$ satisfying
\begin{align}\label{4.14}
\|(f_{1},f_{2})\|_{\mathbf{X}}&=\|\mathcal{L}_{0}^{-1}(\mathcal{S}_{1},0,0,\mathcal{S}_{2},\mathcal{S}_{2,b}^{(1)},\mathcal{S}_{2,b}^{(2)})^{t}\|_{\mathbf{X}}\leq C_{\ell}\big[\|\nu^{-1}w^{\ell}(\mathcal{S}_{1},\mathcal{S}_{2})\|_{L^{\infty}}+|w^{\ell}(\mathcal{S}_{2,b}^{(1)},\mathcal{S}_{2,b}^{(2)})|_{L^{\infty}(\gamma_{-})}\big],
\end{align}
where $C_{\ell}$ depends only $\ell$, which is independent of $\epsilon$.

Step 2. {\it $\sigma\in [0,\sigma_{*}]$} for some small $\sigma_{*}$. We rewrite
\begin{align}\label{4.14-1}
    \mathcal{L}_{\sigma}(f_{1},f_{2})=(\mathcal{S}_{1},0,0,\mathcal{S}_{2},\mathcal{S}_{2,b}^{(1)},\mathcal{S}_{2,b}^{(2)})^{t}
\end{align}
as
\begin{align*}
    \mathcal{L}_{0}(f_{1},f_{2})&=
    (\mathcal{S}_{1},0,0,\mathcal{S}_{2},\mathcal{S}_{2,b}^{(1)},\mathcal{S}_{2,b}^{(2)})^{t}\\
    &\quad +\sigma\big((1-\chi_{M})\mathcal{K}f_{1},0,0,\mu^{-\frac{1}{2}}\chi_{M}\mathcal{K}f_{1}+Kf_{2},P_{\gamma}\big(\frac{f_{1}}{\sqrt{\mu}}\big)(0),P_{\gamma}\big(\frac{f_{1}}{\sqrt{\mu}}\big)(\eta_{1})\big)^{t}.
\end{align*}
We construct the approximate sequence $\{(f_{1}^{i},f_{2}^{i})\}$ as follows:
\begin{align*}
    \mathcal{L}_{0}(f_{1}^{i+1},f_{2}^{i+1})&=
(\mathcal{S}_{1},0,0,\mathcal{S}_{2},\mathcal{S}_{2,b}^{(1)},\mathcal{S}_{2,b}^{(2)})^{t}\\
    &\quad +\sigma\big((1-\chi_{M})\mathcal{K}f_{1}^{i},0,0,\mu^{-\frac{1}{2}}\chi_{M}\mathcal{K}f_{1}^{i}+Kf_{2}^{i},P_{\gamma}\big(\frac{f_{1}^{i}}{\sqrt{\mu}}\big)(0),P_{\gamma}\big(\frac{f_{1}^{i}}{\sqrt{\mu}}\big)(\eta_{1})\big)^{t}.
\end{align*}
with $(f_{1}^{0},f_{2}^{0})=(0,0)$. It follows from \eqref{4.14} that
\begin{align}\label{4.15}
\|(f_{1}^{i+1},f_{2}^{i+1})\|_{\mathbf{X}}&\leq C_{\ell}\sigma\big[\|w^{\ell}(1-\chi_{M})\mathcal{K}f_{1}^{i}\|_{L^{\infty}}+\|w^{\ell}\mu^{-\frac{1}{2}}\chi_{M}f_{1}^{i}\|_{L^{\infty}}+\|Kf_{2}^{i}\|_{L^{\infty}}+\|w^{\ell}P_{\gamma}f_{1}^{i}\|_{L^{\infty}}\big]\nonumber\\
&\quad +C_{\ell}\|\nu^{-1}w^{\ell}(\mathcal{S}_{1},\mathcal{S}_{2})\|_{L^{\infty}}+C_{\ell}|w^{\ell}(\mathcal{S}_{2,b}^{(1)},\mathcal{S}_{2,b}^{(2)})|_{L^{\infty}(\gamma_{-})}\nonumber\\
&\leq C_{\ell}\sigma\|(f_{1}^{i},f_{2}^{i})\|_{\mathbf{X}}+C_{\ell}\|\nu^{-1}w^{\ell}(\mathcal{S}_{1},\mathcal{S}_{2})\|_{L^{\infty}}+C_{\ell}|w^{\ell}(\mathcal{S}_{2,b}^{(1)},\mathcal{S}_{2,b}^{(2)})|_{L^{\infty}(\gamma_{-})}.
\end{align}
Choosing $\sigma_{*}$ small enough such that $C_{\ell}\sigma_{*}\leq \frac{1}{2}$, we obtain from \eqref{4.15} and induction that
\begin{align*}
\|(f_{1}^{i+1},f_{2}^{i+1})\|_{\mathbf{X}}&\leq 2C_{\ell}\|\nu^{-1}w^{\ell}(\mathcal{S}_{1},\mathcal{S}_{2})\|_{L^{\infty}}+2C_{\ell}|w^{\ell}(\mathcal{S}_{2,b}^{(1)},\mathcal{S}_{2,b}^{(2)})|_{L^{\infty}(\gamma_{-})}\quad \text{for $i\geq 1$},
\end{align*}
when $\sigma\in [0,\sigma_{*}]$. Furthermore, taking the difference $[f_{1}^{i+1}-f_{1}^{i},f_{2}^{i+1}-f_{2}^{i}]$, we have
\begin{align}\label{4.17}
\|(f_{1}^{i+1}-f_{1}^{i},f_{2}^{i+1}-f_{2}^{i})\|_{\mathbf{X}}&\leq \frac{1}{2}\|(f_{1}^{i}-f_{1}^{i-1},f_{2}^{i}-f_{2}^{i-1})\|_{\mathbf{X}}\leq \cdots \leq \Big(\frac{1}{2}\Big)^{i}\|(f_{1}^{1},f_{2}^{1})\|_{\mathbf{X}}\nonumber\\
&\leq 2\Big(\frac{1}{2}\Big)^{i}C_{\ell}\big[\|\nu^{-1}w^{\ell}(\mathcal{S}_{1},\mathcal{S}_{2})\|_{L^{\infty}}+|w^{\ell}(\mathcal{S}_{2,b}^{(1)},\mathcal{S}_{2,b}^{(2)})|_{L^{\infty}(\gamma_{-})}\big],
\end{align}
which implies $\{(f_{1}^{i},f_{2}^{i})\}$ is a Cauchy sequence in $\mathbf{X}$. Therefore, there exists a limit function pair $(f_{1},f_{2})\in \mathbf{X}$ such that
$$
\lim\limits_{i\to \infty}\|(f_{1}^{i}-f_{1},f_{2}^{i}-f_{2})\|_{\mathbf{X}}=0.
$$
Moreover, $(f_{1},f_{2})$ solves \eqref{4.14-1}. Then applying the uniform $L^{\infty}-$estimates established in Lemma \ref{lem2.5}, we have
\begin{align}\label{4.18}
\|(f_{1},f_{2})\|_{\mathbf{X}}&=\|\mathcal{L}_{\sigma_{*}}^{-1}(\mathcal{S}_{1},0,0,\mathcal{S}_{2},\mathcal{S}_{2,b}^{(1)},\mathcal{S}_{2,b}^{(2)})^{t}\|_{\mathbf{X}}\leq C_{\epsilon,\ell}\big[\|\nu^{-1}w^{\ell}(\mathcal{S}_{1},\mathcal{S}_{2})\|_{L^{\infty}}+|w^{\ell}(\mathcal{S}_{2,b}^{(1)},\mathcal{S}_{2,b}^{(2)})|_{L^{\infty}(\gamma_{-})}\big],
\end{align}
where $C_{\epsilon,\ell}$ is independent of $\sigma$.

Step 3. {\it $\sigma\in [\sigma_{*},2\sigma_{*}]$}. We rewrite
\begin{align}\label{4.14-2}
\mathcal{L}_{\sigma}(f_{1},f_{2})=\big(\mathcal{S}_{1},0,0,\mathcal{S}_{2},\mathcal{S}_{2,b}^{(1)},\mathcal{S}_{2,b}^{(2)}\big)^{t}
\end{align}
as
\begin{align*}
\mathcal{L}_{\sigma_{*}}(f_{1},f_{2})&=\big(\mathcal{S}_{1},0,0,\mathcal{S}_{2},\mathcal{S}_{2,b}^{(1)},\mathcal{S}_{2,b}^{(2)}\big)^{t}\\
&\quad +\tilde{\sigma}\big((1-\chi_{M})\mathcal{K}f_{1},0,0,\mu^{-\frac{1}{2}}\mathcal{K}f_{1}+Kf_{2},P_{\gamma}\big(\frac{f_{1}}{\sqrt{\mu}}\big)(0),P_{\gamma}\big(\frac{f_{1}}{\sqrt{\mu}}\big)(\eta_{1})\big)^{t},
\end{align*}
where $\tilde{\sigma}=\sigma-\sigma_{*}$. With the help of \eqref{4.18}, by choosing $\sigma_{*}$ small enough such that $\sigma_{*}C_{\epsilon,\ell}=\frac{1}{2}$, we can apply similar arguments as in \eqref{4.15}--\eqref{4.17} to show that there exist a solution $(f_{1},f_{2})\in \mathbf{X}$ for \eqref{4.14-2} provided $\sigma\in [\sigma_{*},2\sigma_{*}]$. Moreover, applying Lemma \ref{lem2.5} to \eqref{4.14-2}, we can get
\begin{align*}
\|(f_{1},f_{2})\|_{\mathbf{X}}&=\|\mathcal{L}_{2\sigma_{*}}^{-1}\big(\mathcal{S}_{1},0,0,\mathcal{S}_{2},\mathcal{S}_{2,b}^{(1)},\mathcal{S}_{2,b}^{(2)}\big)^{t}\|_{\mathbf{X}}\nonumber\\
&\leq C_{\epsilon,\ell}\big[\|\nu^{-1}w^{\ell}(\mathcal{S}_{1},\mathcal{S}_{2})\|_{L^{\infty}}+|w^{\ell}\mathcal{S}_{2,b}^{(1)},\mathcal{S}_{2,b}^{(2)}|_{L^{\infty}(\gamma_{-})}\big].
\end{align*}

Step 4. {\it $\sigma\in [0,1]$}. Repeating the procedure in Step 3, we can establish the existence of solution for \eqref{4.14-3} for $\sigma\in [2\sigma_{*},3\sigma_{*}],[3\sigma_{*},4\sigma_{*}],\cdots$ until $\sigma=1$ step by step. Finally, \eqref{4.4} follows directly from Lemma \ref{lem2.5}. Therefore, the proof of Lemma \ref{lem3.1} is completed. $\hfill\square$

\medskip
In particular, by setting $\mathcal{S}_{1}=\mathcal{S}$, $\mathcal{S}_{2}=\mathcal{S}_{2,b}^{(1)}=\mathcal{S}_{2,b}^{(2)}=0$ in \eqref{4.14-3}, we have directly following Corollary from Lemma \ref{lem3.1}.
\begin{corollary}\label{Cor1}
Let $\epsilon\in (0,1]$ and $\alpha_{2}$ be the positive constant determined in Lemma \ref{lem2.5}. Then for $\alpha\in (0,\alpha_{2})$, there exist a unique solution $(f_{1},f_{2})\in \mathbf{X}$ of the following coupled systems
	\begin{equation}\label{3.1-1}
		\left\{\begin{aligned}
			&\epsilon f_{1}+v_{\eta}\frac{\partial f_{1}}{\partial \eta}-\frac{1}{1-\eta}\Big(v_{\phi}^2\frac{\partial f_1}{\partial v_{\eta}}-v_{\eta}v_{\phi}\frac{\partial f_{1}}{\partial v_{\phi}}\Big)-\frac{\alpha}{\eta_{1}} v_{\eta}\frac{\partial f_{1}}{\partial v_{\phi}}+\frac{\alpha}{2\eta_{1}}v_{\eta}v_{\phi}\sqrt{\mu}f_{2}+\nu f_{1}\\
			&\quad =(1-\chi_{M})\mathcal{K}f_{1}+\mathcal{S},\\
			&f_{1}(0,v_{\eta},v_{\phi},v_{z})\vert_{v_{\eta}>0}=0,\quad f_{1}(\eta_{1},v_{\eta},v_{\phi},v_{z})\vert_{v_{\eta}<0}=0,
		\end{aligned}\right.
	\end{equation}
	and
	\begin{equation}\label{3.2-1}
		\left\{\begin{aligned}
			&\epsilon f_{2}+v_{\eta}\frac{\partial f_{2}}{\partial \eta}-\frac{1}{1-\eta}\Big(v_{\phi}^2\frac{\partial f_2}{\partial v_{\eta}}-v_{\eta}v_{\phi}\frac{\partial f_2}{\partial v_{\phi}}\Big)-\frac{\alpha}{\eta_{1}} v_{\eta}\frac{\partial f_2}{\partial v_{\phi}}+\nu f_2\\
			&\quad =Kf_2+\chi_{M}\mu^{-\frac{1}{2}}\mathcal{K}f_1,\\
			&f_2(0,v_{\eta},v_{\phi},v_{z})\vert_{v_{\eta}>0}=\sqrt{2\pi\mu}\int_{u_{\eta}<0}|u_{\eta}|(f_{1}+\sqrt{\mu}f_2)(0,u)\,\mathrm{d}u,\\
			&f_2(\eta_{1},v_{\eta},v_{\phi},v_{z})\vert_{v_{\eta}<0}=\sqrt{2\pi\mu}\int_{u_{\eta}>0}|u_{\eta}|(f_{1}+\sqrt{\mu}f_2)(\eta_{1},u)\,\mathrm{d}u,
		\end{aligned}\right.
	\end{equation}
    for any given $\mathcal{S}$ with $\|\nu^{-1}w^{\ell}\mathcal{S}\|_{L^{\infty}}<\infty$. Moreover, it holds that
    $$
    \|w^{\ell}(f_{1},f_{2})\|_{L^{\infty}}\leq C_{\epsilon,\ell}\|\nu^{-1}w^{\ell}\mathcal{S}\|_{L^{\infty}}.
    $$
\end{corollary}

\subsection{The limit $\epsilon\to 0$}
In order to take the limit $\epsilon\to 0$ in \eqref{3.1-1}--\eqref{3.2-1}, we first need to establish the uniform estimate with respect to $\epsilon$. Note that the uniform estimate in Corollary \ref{Cor1} depends
on $\epsilon$. However, the dependence on $\epsilon$ arises from the invalidity of conservation of mass for any given $\sigma\neq 1$ and possibility of 
$$
\int_{0}^{\eta_{1}}(1-\eta)\int_{\R^3}\big(\mathcal{S}_{1}+\sqrt{\mu}\mathcal{S}_{2}\big)\,{\rm d}v{\rm d}\eta\neq 0,
$$
see \eqref{E12-1}--\eqref{E12-2} and \eqref{E23}--\eqref{E24} for details. So, in order to establish the uniform estimate with respect to $\epsilon$, we will further require $\mathcal{S}$ in \eqref{3.1-1} to satisfy
\begin{align}\label{3.1-2}
\int_{0}^{\eta_{1}}(1-\eta)\int_{\R^3}\mathcal{S}\,{\rm d}v{\rm d}\eta=0.
\end{align}

\begin{lemma}\label{lem3.2}
For any given $\epsilon\in [0,1]$ and $\mathcal{S}$ with $\|\nu^{-1}w^{\ell}\mathcal{S}\|_{L^{\infty}}<\infty$ satisfying \eqref{3.1-2}, let $(f_{1},f_{2})$ be the solution for coupled system \eqref{3.1-1}--\eqref{3.2-1}. There exists a positive constant $\alpha_{0}$ small enough, which is independent of $\epsilon$, such that if $\alpha\in (0,\alpha_{0})$, 
\begin{align}\label{E64-1}
  \|w^{\ell}(f_{1},f_{2})\|_{L^{\infty}}\leq  C_{\ell}\|\nu^{-1}w^{\ell}\mathcal{S}\|_{L^{\infty}},
 \end{align}
   where $C_{\ell}$ is a positive constant depending only on the weight $\ell$.
\end{lemma}

\noindent\textbf{Proof}. We denote $f=\mu^{-\frac{1}{2}}f_{1}+f_{2}$ and
$$
\mathbf{P}f=\{a(\eta)+b_{1}(\eta)v_{\eta}+b_{2}(\eta)v_{\phi}+b_{3}(\eta)v_{z}+(|v|^2-3)c(\eta)\}\sqrt{\mu},
$$
then $f$ satisfies
\begin{align*}
		\left\{
		\begin{aligned}
			&\epsilon f+v_{\eta}\partial_{\eta}f-\frac{1}{1-\eta}(v_{\phi}^2\partial_{v_{\eta}}f-v_{\eta}v_{\phi}\partial_{v_{\phi}}f)-\frac{\alpha}{\eta_{1}}v_{\eta}\partial_{v_{\phi}}f+\frac{\alpha}{2\eta_{1}}v_{\eta}v_{\phi}f+{\bf L}f=\mu^{-\frac{1}{2}}\mathcal{S},\\
		&f(0,v_{\eta},v_{\phi},v_{z})\vert_{v_{\eta}>0}=\sqrt{2\pi\mu}\int_{u_{\eta}<0}(\sqrt{\mu f})(0,u_{\eta},u_{\phi})|u_{\eta}|\,{\rm d}u,\\
			&f(\eta_{1},v_{\eta},v_{\phi},v_{z})|_{v_{\eta}<0}=\sqrt{2\pi\mu}\int_{u_{\eta}>0}(\sqrt{\mu} f)(\eta_{1},u_{\eta},u_{\phi})|u_{\eta}|\,{\rm d}u.
		\end{aligned}
		\right.
	\end{align*}
Noting \eqref{3.1-2}, by a similar calculation as in \eqref{E12}--\eqref{E12-1}, we have
\begin{align}\label{4.20}
    \epsilon\int_{0}^{\eta_{1}}(1-\eta)a(\eta)\,{\rm d}\eta=0.
\end{align}
With \eqref{4.20}, by choosing the same test functions as in \eqref{a}, the boundary terms $I_{3}+I_{4}$ in \eqref{E17} directly vanishes due to $\phi_{a}'(0)=\phi_{a}'(\eta_{1})=0$. Then by using similar calculations as in \eqref{E19}--\eqref{E26}, we can control $\|a\|_{L^2}^2$ as follows
\begin{align*}
	\|a\|_{L^2}^2&\lesssim \kappa_{1}\|a\|_{L^2}^2+\alpha\|a\|_{L^2}\|b_{2}\|_{L^2}+C_{\kappa_{1}}\|b_{1}\|_{L^2}^2+C_{\kappa_{1}}\|(\mathbf{I-P})f\|_{L_{\nu}^2}^2.
\end{align*}
Finally, applying the same calculations in \eqref{E37} to control $\|b_{1}\|_{L^2}^2$, \eqref{E47} to control $\|b_{2}\|_{L^2}^2$, \eqref{z10} to control $\|b_{3}\|_{L^2}^2$, \eqref{E57} to control $\|c\|_{L^2}^2$, we can obtain 
\begin{align*}
&\|(a,b_{1},b_{2},b_{3},c)\|_{L^2}^2\nonumber\\
&\lesssim \kappa_{1}\|a\|_{L^2}^2+C_{\kappa_{1}}\kappa_{2}\|b_{1}\|_{L^2}^2+\kappa_{3}\|b_{2}\|_{L^2}^2+\kappa_{4}\|b_{3}\|_{L^2}^2+C_{\kappa_{1},\kappa_{2}}\kappa_{5}\|c\|_{L^2}^2+C_{\kappa}\|(\mathbf{I-P})f_{2}\|_{L_{\nu}^2}^2\nonumber\\
&\quad +C_{\kappa}\alpha\|(a,b_{1},b_{2},c)\|_{L^2}^2+C_{\kappa}\|\nu^{-1}w^{\ell}\mathcal{S}\|_{L^{\infty}}^2+C_{\kappa}\|w^{\ell}f_{1}\|_{L^{\infty}}^2.
\end{align*}
Then applying the same calculations as in \eqref{E5-1}--\eqref{E10-2} to control $\|(\mathbf{I-P})f_{2}\|_{L_{\nu}^2}^2$ and similar arguments as in \eqref{E61}--\eqref{E63}, we can get
\begin{align}\label{E63-1}
\|f_{2}\|_{L_{\nu}^2}\leq C\|w^{\ell}f_{1}\|_{L^{\infty}}+C\|\nu^{-1}w^{\ell}\mathcal{S}\|_{L^{\infty}},
\end{align}
which, together with Lemma \ref{lem2.2}, yields \eqref{E64-1}. Hence, the proof of Lemma \ref{lem3.2} is completed. $\hfill\square$

\medskip
With the uniform estimates established in Lemma \ref{lem3.1}, we are able to take the limit $\epsilon \to 0$ in \eqref{3.1-1}--\eqref{3.2-1} to obtain the existence of solutions to the linearized systems \eqref{3.1}--\eqref{3.2}. 

\begin{lemma}\label{lem3.3}
There exists a positive constant $\alpha_{0}$ small enough such that if $\alpha\in (0,\alpha_{0})$, there exists a unique solution $(f_{1},f_{2})$ to the linearized coupled systems \eqref{3.1}--\eqref{3.2} satisfying
\begin{align}\label{E64-2}
  \|w^{\ell}(f_{1},f_{2})\|_{L^{\infty}}\leq  C_{\ell}\|\nu^{-1}w^{\ell}\mathcal{S}\|_{L^{\infty}},
 \end{align}
   where $C_{\ell}$ is a positive constant depending only on the weight $\ell$.
\end{lemma}
\noindent\textbf{Proof}. Let $[f_{1}^{\epsilon_{i}},f_{2}^{\epsilon_{i}}]\,\,(i=1,2)$ be two solutions of the coupled systems \eqref{3.1-1}--\eqref{3.2-1}. Considering the difference $[f_{1}^{\epsilon_{2}}-f_{1}^{\epsilon_{1}},f_{2}^{\epsilon_{2}}-f_{2}^{\epsilon_{1}}]$, we obtain from similar arguments as in \eqref{3.17-1}--\eqref{3.17-2} and \eqref{E63-1} that
\begin{align}
&\|w^{\ell}(f_{1}^{\epsilon_{2}}-f_{1}^{\epsilon_{1}})\|_{L^{\infty}}\leq C\alpha\|w^{\ell}(f_{2}^{\epsilon_{2}}-f_{2}^{\epsilon_{1}})\|_{L^{\infty}}+C\|w^{\ell}(\epsilon_{2}-\epsilon_{1})f_{1}^{\epsilon_{2}}\|_{L^{\infty}},\label{4.21}\\
&\|w^{\ell}(f_{2}^{\epsilon_{2}}-f_{2}^{\epsilon_{2}})\|_{L^{\infty}}\leq C_{\ell}\|w^{\ell}(f_{1}^{\epsilon_{2}}-f_{1}^{\epsilon_{1}})\|_{L^{\infty}}+C_{\ell}\|w^{\ell}(\epsilon_{2}-\epsilon_{1})f_{2}^{\epsilon_{2}}\|_{L^{\infty}}+C_{\ell}\|f_{2}^{\epsilon_{2}}-f_{2}^{\epsilon_{1}}\|_{L^2},\label{4.22}
\end{align}
and
\begin{align}\label{4.23}
\|f_{2}^{\epsilon_{2}}-f_{2}^{\epsilon_{1}}\|_{L^2}\leq C\|w^{\ell}(f_{1}^{\epsilon_{2}}-f_{1}^{\epsilon_{1}})\|_{L^{\infty}}+C(\epsilon_{2}+\epsilon_{1})\|w^{\ell}(f_{1}^{\epsilon_{2}},f_{2}^{\epsilon_{2}})\|_{L^{\infty}}.
\end{align}
Combining \eqref{4.21}, \eqref{4.22} and \eqref{4.23}, we get
\begin{align}\label{4.24}
&\|w^{\ell}(f_{1}^{\epsilon_{2}}-f_{1}^{\epsilon_{1}})\|_{L^{\infty}}\leq C_{\ell}\alpha\|w^{\ell}(f_{1}^{\epsilon_{2}}-f_{1}^{\epsilon_{1}})\|_{L^{\infty}}+C_{\ell}(\epsilon_{2}+\epsilon_{1})\|w^{\ell}(f_{1}^{\epsilon_{2}},f_{2}^{\epsilon_{2}})\|_{L^{\infty}}.
\end{align}
Taking $\alpha_{0}$ small enough such that $C_{\ell}\alpha\leq \frac{1}{2}$, we obtain from \eqref{4.24} by induction that
\begin{align*}
\|w^{\ell}(f_{1}^{\epsilon_{2}}-f_{1}^{\epsilon_{1}})\|_{L^{\infty}}&\leq 2C_{\ell}(\epsilon_{2}+\epsilon_{1})\big[\|w^{\ell}(f_{1}^{\epsilon_{2}},f_{2}^{\epsilon_{2}})\|_{L^{\infty}}\leq C_{\ell}(\epsilon_{2}+\epsilon_{1})\|\nu^{-1}w^{\ell}\mathcal{S}\|_{L^{\infty}}\rightarrow 0,
\end{align*}
as $\epsilon_{1},\epsilon_{2}\to 0$. Here we have used \eqref{E64-1} in the last inequality. Hence, $\{f_{1}^{\epsilon},f_{2}^{\epsilon}\}$ are both Cauchy sequences in $\mathbf{X}$, that is, there exist a unique solution $(f_{1},f_{2})\in \mathbf{X}$ for \eqref{3.1}--\eqref{3.2} satisfying
$$
\lim\limits_{\epsilon\to 0}\|(f_{1}^{\epsilon}-f_{1},f_{2}^{\epsilon}-f_{2})\|_{\mathbf{X}}=0,
$$
which, together with \eqref{E64-1}, yields \eqref{E64-2}. Therefore, the proof of Lemma \ref{lem3.3} is completed. $\hfill\square$

\section{Existence for the nonlinear system: Proof of Theorem \ref{thm1}}\label{sec6}
In this section, we are going to establish the existence of nonlinear coupled systems \eqref{2.3}--\eqref{2.4} and complete the proof of Theorem \ref{thm1}. We construct the approximate sequences for the coupled systems \eqref{2.3}--\eqref{2.4} by the following iterative systems
	\begin{equation*}
		\left\{\begin{aligned}
			&v_{\eta}\frac{\partial g_{1}^{k+1}}{\partial \eta}-\frac{1}{1-\eta}\Big(v_{\phi}^2\frac{\partial g_{1}^{k+1}}{\partial v_{\eta}}-v_{\eta}v_{\phi}\frac{\partial g_{1}^{k+1}}{\partial v_{\phi}}\Big)-\frac{\alpha}{\eta_{1}} v_{\eta}\frac{\partial g_{1}^{k+1}}{\partial v_{\phi}}+\frac{\alpha}{2\eta_{1}}v_{\eta}v_{\phi}\sqrt{\mu}g_{2}^{k+1}+\nu g_{1}^{k+1}\\
			&\quad =(1-\chi_{M})\mathcal{K}g_{1}^{k+1}-\frac{1}{\eta_{1}}v_{\eta}v_{\phi}\mu+\alpha Q(\sqrt{\mu}g^{k},\sqrt{\mu}g^{k}),\\
			&g_{1}^{k+1}(0,v_{\eta},v_{\phi},v_{z})\vert_{v_{\eta}>0}=0,\quad g_{1}^{k+1}(\eta_{1},v_{\eta},v_{\phi},v_{z})\vert_{v_{\eta}<0}=0,
		\end{aligned}\right.
	\end{equation*}
	and
	\begin{equation*}
		\left\{\begin{aligned}
			&v_{\eta}\frac{\partial g_{2}^{k+1}}{\partial \eta}-\frac{1}{1-\eta}\Big(v_{\phi}^2\frac{\partial g_{2}^{k+1}}{\partial v_{\eta}}-v_{\eta}v_{\phi}\frac{\partial g_{2}^{k+1}}{\partial v_{\phi}}\Big)-\frac{\alpha}{\eta_{1}} v_{\eta}\frac{\partial g_{2}^{k+1}}{\partial v_{\phi}}+\nu g_{2}^{k+1}\\
			&\qquad =Kg_{2}^{k+1}+\chi_{M}\mu^{-\frac{1}{2}}\mathcal{K}g_{1}^{k+1},\\
			&g_{2}^{k+1}(0,v_{\eta},v_{\phi})\vert_{v_{\eta}>0}=\sqrt{2\pi\mu}\int_{u_{\eta}<0}|u_{\eta}|(g_{1}^{k+1}+\sqrt{\mu}g_{2}^{k+1})(0,u)\,\mathrm{d}u,\\
			&g_{2}^{k+1}(\eta_{1},v_{\eta},v_{\phi},v_{z})\vert_{v_{\eta}<0}=\sqrt{2\pi\mu}\int_{u_{\eta}>0}|u_{\eta}|(g_{1}^{k+1}+\sqrt{\mu}g_{2}^{k+1})(\eta_{1},u)\,\mathrm{d}u,
		\end{aligned}\right.
	\end{equation*}
with $(g_{1}^{0},g_{2}^{0})=(0,0)$. Noting \eqref{3.1-2} holds, i.e.,
\begin{align*}
    -\int_{0}^{\eta_{1}}(1-\eta)\int_{\R^3}v_{\eta}v_{\phi}\mu\,{\rm d}v{\rm d}\eta=\int_{0}^{\eta_{1}}(1-\eta)\int_{\R^3}Q(\sqrt{\mu}g^{k},\sqrt{\mu}g^{k})\,{\rm d}v{\rm d}\eta=0,
\end{align*}
we obtain from Lemmas \ref{lem3.3}, and \ref{K} that $(g_{1}^{k+1},g_{2}^{k+1})\in \mathbf{X}$ once $(g_{1}^{k},g_{2}^{k})\in \mathbf{X}$ and
\begin{align*}
\|w^{\ell}(g_{1}^{k+1},g_{2}^{k+1})\|_{L^{\infty}}&\leq C_{\ell}\|\nu^{-1}w^{\ell}(v_{\eta}v_{\phi}\mu)\|_{L^{\infty}}+C_{\ell}\alpha\|\nu^{-1}w^{\ell}Q(\sqrt{\mu}g^{k},\sqrt{\mu}g^{k})\|_{L^{\infty}}\nonumber\\
&\leq C_{\ell}\alpha\|w^{\ell}(g_{1}^{k},g_{2}^{k})\|_{L^{\infty}}^2+C_{\ell}.
\end{align*}
Choosing $\alpha_{*}$ small enough such that $C_{\ell}^2\alpha_{*}\leq \frac{1}{8}$, we obtain from induction method that
\begin{align*}
\|w^{\ell}(g_{1}^{k},g_{2}^{k})\|_{L^{\infty}}\leq 2C_{\ell}\quad \text{for all $k\geq 1$},
\end{align*}
provided $\alpha\in (0,\alpha_{*})$. 
Now, taking the difference $(g_{1}^{k+1}-g_{1}^{k},g_{2}^{k+1}-g_{2}^{k})$ and applying the Lemma \ref{lem3.3}, one has
\begin{align*}
&\|w^{\ell}(g_{1}^{k+1}-g_{1}^{k},g_{2}^{k+1}-g_{2}^{k})\|_{L^{\infty}}\nonumber\\
&\leq C_{\ell}\alpha\|\nu^{-1}w^{\ell}\big[Q(\sqrt{\mu}g^{k}-\sqrt{\mu}g^{k-1},\sqrt{\mu}g^{k})+Q(\sqrt{\mu}g^{k-1},\sqrt{\mu}g^{k}-\sqrt{\mu}g^{k-1})\big]\|_{L^{\infty}}\nonumber\\
&\leq C_{\ell}\alpha\big[\|w^{\ell}(g_{1}^{k},g_{2}^{k})\|_{L^{\infty}}+\|w^{\ell}(g_{1}^{k-1},g_{2}^{k-1})\|_{L^{\infty}}\big]\cdot \|w^{\ell}(g_{1}^{k}-g_{1}^{k-1},g_{2}^{k}-g_{2}^{k-1})\|_{L^{\infty}}\nonumber\\
&\leq \frac{1}{2}\|w^{\ell}(g_{1}^{k}-g_{1}^{k-1},g_{2}^{k}-g_{2}^{k-1})\|_{L^{\infty}}\leq \cdots \leq \Big(\frac{1}{2}\Big)^{k}\|w^{\ell}(g_{1}^{1},g_{2}^{1})\|_{L^{\infty}}\leq C_{\ell}\Big(\frac{1}{2}\Big)^{k-1},
\end{align*}
which implies that $\{(g_{1}^{k},g_{2}^{k})\}$ is a Cauchy sequence in $\mathbf{X}$, that is, there exists a solution $(g_{1},g_{2})$ for the coupled systems \eqref{2.3}--\eqref{2.4} satisfying
\begin{align*}
\|w^{\ell}(g_{1},g_{2})\|_{L^{\infty}}= \lim\limits_{k\to \infty}\|w^{\ell}(g_{1}^{k},g_{2}^{k})\|_{L^{\infty}}\leq 2C_{\ell}.
\end{align*}
Therefore, the proof of Theorem \ref{thm1} is completed. $\hfill\square$

\section{Unsteady problem}\label{sec7}
	To show the non-negativity of steady solution
	$$
	G_{s}=\mu+\alpha\sqrt{\mu}g_{1}=\mu+\alpha(g_{1}+\sqrt{\mu}g_{2})
	$$
obtained in Theorem \ref{thm1}, we study the asymptotic behavior of the following unsteady problem:
	\begin{equation}\label{6.1}
		\left\{\begin{aligned}
			&\partial_{t}G+v_{\eta}\frac{\partial G}{\partial \eta}-\frac{1}{1-\eta}\Big(v_{\phi}^2\frac{\partial G}{\partial v_{\eta}}-v_{\eta}v_{\phi}\frac{\partial G}{\partial v_{\phi}}\Big)-\frac{\alpha}{\eta_{1}}v_{\eta}\frac{\partial G}{\partial v_{\phi}}=Q(G,G),\\
			&G(t,0,v_{\eta},v_{\phi},v_{z})\vert_{v_{\eta}>0}=\sqrt{2\pi}\mu\int_{u_{\eta}<0}G(t,0,u)|u_{\eta}|\,\mathrm{d}u,\\
			&G(t,\eta_{1},v_{\eta},v_{\phi},v_{z})\vert_{v_{\eta}<0}=\sqrt{2\pi}\mu\int_{u_{\eta}>0}G(t,\eta_{1},u)|u_{\eta}|\,\mathrm{d}u,\\
			&G(0,\eta,v_{\eta},v_{\phi},v_{z})\vert_{t=0}=G_{0}(\eta,v_{\eta},v_{\phi},v_{z}).
		\end{aligned}
		\right.
	\end{equation}
Note that the initial data in \eqref{6.1} is independent of the angle $\phi$.

Consider the perturbation $F=G-G_{s}$ and denote $F=\sqrt{\mu}f$, then $f=f(t,\eta,v)$ satisfies 
\begin{equation}\label{6.1-1}
	\left\{\begin{aligned}
		&\partial_{t}f+v_{\eta}\frac{\partial f}{\partial \eta}-\frac{1}{1-\eta}\Big(v_{\phi}^2\frac{\partial f}{\partial v_{\eta}}-v_{\eta}v_{\phi}\frac{\partial f}{\partial v_{\phi}}\Big)-\frac{\alpha}{\eta_{1}}v_{\eta}\frac{\partial f}{\partial v_{\phi}}+\frac{\alpha}{2\eta_{1}}v_{\eta}v_{\phi}f+\mathbf{L}f\\
		&\quad =\Gamma(f,f)+\alpha\{\Gamma(g_{1}+\sqrt{\mu}g_{2},f)+\Gamma(f,g_{1}+\alpha\sqrt{\mu}g_{2})\},\quad t>0, (\eta,v)\in (0,\eta_{1})\times \R^3,\\
		&f(t,0,v_{\eta},v_{\phi},v_{z})\vert_{v_{\eta}>0}=\sqrt{2\pi\mu}\int_{u_{\eta}<0}f(t,0,u)\sqrt{\mu}|u_{\eta}|\,\mathrm{d}u,\\
		&f(t,\eta_{1},v_{\eta},v_{\phi},v_{z})\vert_{v_{\eta}<0}=\sqrt{2\pi\mu}\int_{u_{\eta}>0}f(t,\eta_{1},u)\sqrt{\mu}|u_{\eta}|\,\mathrm{d}u,\\
		&\sqrt{\mu}f(0,\eta,v_{\eta},v_{\phi},v_{z})=G_{0}(\eta,v_{\eta},v_{\phi},v_{z})-G_{s}(\eta,v_{\eta},v_{\phi},v_{z}).
	\end{aligned}
	\right.
\end{equation}

\subsection{Local-in-time existence of \texorpdfstring{\eqref{6.1-1}}{（7.2）}}\label{sec7.1}
The goal of this subsection is to construct the local-in-time solution to the initial boundary value problem \eqref{6.1-1}. Similar to the steady problem, we decompose $\sqrt{\mu}f$ as
$$
\sqrt{\mu}f=f_{1}+\sqrt{\mu}f_{2},
$$
where $f_{1}$ and $f_{2}$ satisfy the following initial boundary value problems, respectively:
\begin{align}\label{6.1-2}
\left\{
\begin{aligned}
	&\partial_{t}f_{1}+v_{\eta}\frac{\partial f_{1}}{\partial \eta}-\frac{1}{1-\eta}\Big(v_{\phi}^2\frac{\partial f_{1}}{\partial v_{\eta}}-v_{\eta}v_{\phi}\frac{\partial f_{1}}{\partial v_{\phi}}\Big)-\frac{\alpha}{\eta_{1}}v_{\eta}\frac{\partial f_{1}}{\partial v_{\phi}}+\frac{\alpha}{2\eta_{1}}v_{\eta}v_{\phi}\sqrt{\mu}f_{2}+\nu f_{1}\\
	&\quad =(1-\chi_{M})\mathcal{K}f_{1}+H(f,f),\\
	&f_{1}(t,0,v_{\eta},v_{\phi},v_{z})\vert_{v_{\eta}>0}=f_{1}(t,\eta_{1},v_{\eta},v_{\phi},v_{z})\vert_{v_{\eta}<0}=0,\\
	&f_{1}(0,\eta,v_{\eta},v_{\phi},v_{z})=G_{0}(\eta,v_{\eta},v_{\phi},v_{z})-G_{s}(\eta,v_{\eta},v_{\phi},v_{z}),
\end{aligned}
\right.
\end{align}
and
\begin{align}\label{6.1-3}
	\left\{
	\begin{aligned}
		&\partial_{t}f_{2}+v_{\eta}\frac{\partial f_{2}}{\partial \eta}-\frac{1}{1-\eta}\Big(v_{\phi}^2\frac{\partial f_{2}}{\partial v_{\eta}}-v_{\eta}v_{\phi}\frac{\partial f_{2}}{\partial v_{\phi}}\Big)-\frac{\alpha}{\eta_{1}}v_{\eta}\frac{\partial f_{2}}{\partial v_{\phi}}+\nu f_{2} =\chi_{M}\mathcal{K}f_{1}+Kf_{2},\\
		&f_{2}(t,0,v_{\eta},v_{\phi},v_{z})\vert_{v_{\eta}>0}=\sqrt{2\pi\mu}\int_{u_{\eta}<0}(f_{1}+\sqrt{\mu}f_{2})(t,0,u)|u_{\eta}|\,\mathrm{d}u,\\
		&f_{2}(t,\eta_{1},v_{\eta},v_{\phi},v_{z})\vert_{v_{\eta}<0}=\sqrt{2\pi\mu}\int_{u_{\eta}>0}(f_{1}+\sqrt{\mu}f_{2})(t,\eta_{1},u)|u_{\eta}|\,\mathrm{d}u,\\
		&f_{2}(0,\eta,v_{\eta},v_{\phi},v_{z})=0,
	\end{aligned}
	\right.
\end{align}
where
$$
H(f,f)=Q(\sqrt{\mu}f,\sqrt{\mu}f)+\alpha[Q(g_{1}+\sqrt{\mu}g_{2},\sqrt{\mu}f)+Q(\sqrt{\mu}f,g_{1}+\alpha\sqrt{\mu}g_{2})].
$$
Note that the boundary conditions of $f_{1}$ and the initial data of $f_{2}$ are set to be zero.

We will look for solutions to \eqref{6.1-2}--\eqref{6.1-3} in the following functional space
\begin{align*}
	\mathbf{Y}_{T}=\big\{(\mathcal{G}_{1},\mathcal{G}_{2})\,\Big\vert\,\sup_{0\leq t\leq T}\|w^{\ell}(\mathcal{G}_{1},\mathcal{G}_{2})\|_{L^{\infty}}<+\infty\big\},
\end{align*}
equipped with the norm
\begin{align*}
	\|(\mathcal{G}_{1},\mathcal{G}_{2})\|_{\mathbf{Y}_{T}}=\sup_{0\leq t\leq T}\|w^{\ell}(\mathcal{G}_{1},\mathcal{G}_{2})\|_{L^{\infty}}.
\end{align*}

\begin{theorem}[Local-in-time existence]\label{thm4.1}
Under the assumptions of Theorem \ref{thm2}, there exists $T_{*}>0$ depending on $\alpha$ such that the coupled systems \eqref{6.1-2}--\eqref{6.1-3} admits a unique local-in-time solution $(f_{1},f_{2})(t,\eta,v)$ satisfying
\begin{align*}
\|(f_{1},f_{2})\|_{\mathbf{Y}_{T_{*}}}\leq C_{0}\varepsilon_{0},
\end{align*}
for some $C_{0}>0$. Moreover,
$$
G(t,\eta,v)=G_{s}(\eta,v)+\sqrt{\mu}f(t,\eta,v)=G_{s}(\eta,v)+(f_{1}+\sqrt{\mu}f_{2})(t,\eta,v)\geq 0
$$
provided that $G_{0}=G(0,\eta,v)\geq 0$.
\end{theorem}

\noindent\textbf{Proof}. Let us start with the following linearized iterative system for \eqref{6.1}:
\begin{align}\label{7.16}
\left\{
\begin{aligned}
&\partial_{t}G^{n+1}+v_{\eta}\frac{\partial G^{n+1}}{\partial \eta}-\frac{1}{1-\eta}\Big(v_{\phi}^2\frac{\partial G^{n+1}}{\partial v_{\eta}}-v_{\eta}v_{\phi}\frac{\partial G^{n+1}}{\partial v_{\phi}}\Big)-\frac{\alpha}{\eta_{1}}v_{\eta}\frac{\partial G^{n+1}}{\partial v_{\phi}}+G^{n+1}\mathcal{R}(G^{n})\\
&\quad =Q_{+}(G^{n},G^{n}),\quad t>0,\,\,\eta\in (0,\eta_{1}),\,\,v\in \R^3,\\
&G^{n+1}(t,0,v_{\eta},v_{\phi},v_{z})\vert_{v_{\eta}>0}=\sqrt{2\pi}\mu\int_{u_{\eta}<0}G^{n}(t,0,u)|u_{\eta}|\,{\rm d}u,\\
&G^{n+1}(t,\eta_{1},v_{\eta},v_{\phi},v_{z})\vert_{v_{\eta}<0}=\sqrt{2\pi}\mu\int_{u_{\eta}>0}G^{n}(t,\eta_{1},u)|u_{\eta}|\,{\rm d}u,\\
&G^{n+1}(0,\eta,v_{\eta},v_{\phi},v_{z})\vert_{t=0}=G_{0}(\eta,v_{\eta},v_{\phi},v_{z}),
\end{aligned}
\right.
\end{align}
where $G^{0}=G_{0}(\eta,v_{\eta},v_{\phi},v_{z})$, $Q_{+}(\cdot,\cdot)$ is the gain term of the collision operator $Q$ defined in \eqref{Q-op}, and
$$
\mathcal{R}(G^{n})=\int_{\R^3}\int_{\mathbb{S}^2}|v-v_{*}|^{\gamma}B_{0}(\cos\theta)G^{n}(v_{*})\,{\rm d}\omega{\rm d}v_{*}.
$$
It is direct to see $G^{n+1}\mathcal{R}(G^n)=Q_{-}(G^{n},G^{n+1})$, where $Q_{-}(\cdot,\cdot)$ is the loss term of the collision operator $Q$ defined in \eqref{Q-op}. Along the characteristic trajectory, by induction method, it follows from $G^{0}\geq 0$ that any solution $G^{n+1}$ of \eqref{7.16} is non-negative. 

Denote $G^{n+1}=G_{s}+\sqrt{\mu}f^{n+1}$ and $G^{n}=G_{s}+\sqrt{\mu}f^{n}$, and decompose $f^{n+1}$ and $f^{n}$ as 
$$
\sqrt{\mu}f^{n+1}=f_{1}^{n+1}+\sqrt{\mu}f_{2}^{n+1},\quad \sqrt{\mu}f^{n}=f_{1}^{n}+\sqrt{\mu}f_{2}^{n},
$$
where $(f_{1}^{n+1},f_{2}^{n+1})$ are required to satisfy
\begin{align*}
\left\{\begin{aligned}
&\partial_{t}f_{1}^{n+1}+v_{\eta}\frac{\partial f_{1}^{n+1}}{\partial \eta}-\frac{1}{1-\eta}\Big(v_{\phi}^2\frac{\partial f_{1}^{n+1}}{\partial v_{\eta}}-v_{\eta}v_{\phi}\frac{\partial f_{1}^{n+1}}{\partial v_{\phi}}\Big)-\frac{\alpha}{\eta_{1}}v_{\eta}\frac{\partial f_{1}^{n+1}}{\partial v_{\phi}}+\nu f_{1}^{n+1}\\
&\quad =-f_{1}^{n+1}\mathcal{R}(\alpha\sqrt{\mu}g+\sqrt{\mu}f^{n})-\frac{\alpha}{2\eta_{1}}v_{\eta}v_{\phi}\sqrt{\mu}f_{2}^{n+1}+(1-\chi_{M})\mathcal{K}f_{1}^{n}+\tilde{H}(f^{n},f^{n}),\\
&f_{1}^{n+1}(t,0,v_{\eta},v_{\phi},v_{z})\vert_{v_{\eta}>0}=f_{1}^{n+1}(t,\eta_{1},v_{\eta},v_{\phi},v_{z})\vert_{v_{\eta}<0}=0,\\
&f_{1}^{n+1}(0,\eta,v_{\eta},v_{\phi},v_{z})=G_{0}(\eta,v)-G_{s}(\eta,v),
\end{aligned}
\right.
\end{align*}
and
\begin{align*}
\left\{\begin{aligned}
&\partial_{t}f_{2}^{n+1}+v_{\eta}\frac{\partial f_{2}^{n+1}}{\partial \eta}-\frac{1}{1-\eta}\Big(v_{\phi}^2\frac{\partial f_{2}^{n+1}}{\partial v_{\eta}}-v_{\eta}v_{\phi}\frac{\partial f_{2}^{n+1}}{\partial v_{\phi}}\Big)-\frac{\alpha}{\eta_{1}}v_{\eta}\frac{\partial f_{2}^{n+1}}{\partial v_{\phi}}+\nu f_{2}^{n+1}\\
&\quad =-f_{2}^{n+1}\mathcal{R}(\alpha\sqrt{\mu}g+\sqrt{\mu}f^{n})+Kf_{2}^{n}+\chi_{M}\mathcal{K}f_{1}^{n},\\
&f_{2}^{n+1}(t,0,v_{\eta},v_{\phi},v_{z})\vert_{v_{\eta}>0}=\sqrt{2\pi\mu}\int_{u_{\eta}<0}(f_{1}^n+\sqrt{\mu}f_{2}^n)(t,0,u)|u_{\eta}|\,{\rm d}u,\\
&f_{2}^{n+1}(t,\eta_{1},v_{\eta},v_{\phi},v_{z})\vert_{v_{\eta}<0}=\sqrt{2\pi\mu}\int_{u_{\eta}>0}(f_{1}^{n}+\sqrt{\mu}f_{2}^{n})(t,\eta_{1},u)|u_{\eta}|\,{\rm d}u,\\
&f_{2}^{n+1}(0,\eta,v_{\eta},v_{\phi},v_{z})=0,
\end{aligned}
\right.
\end{align*}
respectively. Here we have denoted
$$
\tilde{H}(f^{n},f^{n})=-\alpha\sqrt{\mu}g\mathcal{R}(f^n)+Q_{+}(\alpha\sqrt{\mu}g,\sqrt{\mu}f^{n})+Q_{+}(\sqrt{\mu}f^{n},\alpha\sqrt{\mu}g)+Q_{+}(\sqrt{\mu}f^{n},\sqrt{\mu}f^{n}).
$$
Next, we shall show inductively that there exists a finite $T_{*}>0$ such that
\begin{align}\label{7.1-2}
\sup_{0\leq t\leq T_{*}}\|w^{\ell}\big(f_{1}^{m},f_{2}^{m}\big)\|_{L^{\infty}}\leq C_{0}\varepsilon_{0}<\frac{\nu_{0}}{16},\qquad \forall\,\,m\geq 1,
\end{align}
for some positive constant $C_{0}$ independent of $m$, provided that
\begin{align*}
\|w^{\ell}(f_{1}^{0},f_{2}^{0})\|_{L^{\infty}}=\|w^{\ell}(G_{0}-G_{s})\|_{L^{\infty}}\leq \varepsilon_{0}
\end{align*}
for $\varepsilon_{0}$ small enough.

Denote
\begin{align*}
(h_{1}^{n},h_{2}^{n})=w^{\ell}(f_{1}^{n},f_{2}^{n}),\quad \sqrt{\mu}h^{n}=w^{\ell}f^n=h_{1}^{n}+\sqrt{\mu}h_{2}^n,
\end{align*}
then $h_{1}^{n+1}$ and $h_{2}^{n+1}$ satisfy
\begin{align*}
\left\{
\begin{aligned}
&\partial_{t}h_{1}^{n+1}+v_{\eta}\frac{\partial h_{1}^{n+1}}{\partial \eta}-\frac{1}{1-\eta}\Big(v_{\phi}^2\frac{\partial h_{1}^{n+1}}{\partial v_{\eta}}-v_{\eta}v_{\phi}\frac{\partial h_{1}^{n+1}}{\partial v_{\phi}}\Big)-\frac{\alpha}{\eta_{1}}v_{\eta}\frac{\partial h_{1}^{n+1}}{\partial v_{\phi}}+\nu h_{1}^{n+1}+\frac{\alpha}{2\eta_{1}}\sqrt{\mu}v_{\eta}v_{\phi}h_{2}^{n+1}\\
    &\quad =-\frac{\alpha}{\eta_{1}}\frac{\ell v_{\eta}v_{\phi}}{1+|v|^2}h_{1}^{n+1}-h_{1}^{n+1}\mathcal{R}(\alpha \sqrt{\mu}g+\sqrt{\mu}f^{n})+(1-\chi_{M})\mathcal{K}_{w}h_{1}^{n}+w^{\ell}\tilde{H}(f^{n},f^{n}),\\
    &h_{1}^{n+1}(t,0,v_{\eta},v_{\phi},v_{z})\vert_{v_{\eta}>0}=h_{1}^{n+1}(t,\eta_{1},v_{\eta},v_{\phi},v_{z})\vert_{v_{\eta}<0}=0,\\
    &h_{1}^{n+1}(0,\eta,v_{\eta},v_{\phi},v_{z})=w^{\ell}\big(G_{0}(\eta,v)-G_{s}(\eta,v)\big),
\end{aligned}
\right.
\end{align*}
and
\begin{align*}
\left\{\begin{aligned}
&\partial_{t}h_{2}^{n+1}+v_{\eta}\frac{\partial h_{2}^{n+1}}{\partial \eta}-\frac{1}{1-\eta}\Big(v_{\phi}^2\frac{\partial h_{2}^{n+1}}{\partial v_{\eta}}-v_{\eta}v_{\phi}\frac{\partial h_{2}^{n+1}}{\partial v_{\phi}}\Big)-\frac{\alpha}{\eta_{1}}v_{\eta}\frac{\partial h_{2}^{n+1}}{\partial v_{\phi}}+\nu h_{2}^{n+1}\\
&\quad =-\frac{\alpha}{\eta_{1}}\frac{\ell v_{\eta}v_{\phi}}{1+|v|^2}h_{2}^{n+1}-h_{2}^{n+1}\mathcal{R}(\alpha\sqrt{\mu}g+\sqrt{\mu}f^{n})+\chi_{M}\mathcal{K}_{w}h_{1}^{n}+K_{w}h_{2}^{n},\\
&h_{2}^{n+1}(t,0,v_{\eta},v_{\phi},v_{z})\vert_{v_{\eta}>0}=\frac{1}{\tilde{w}(v)}\int_{u_{\eta}<0}\big(\frac{h_{1}^{n}}{\sqrt{\mu}}+\sqrt{\mu}h_{2}^{n}\big)(t,0,\tilde{u})\tilde{w}(u)\,{\rm d}\sigma,\\
&h_{2}^{n+1}(t,\eta_{1},v_{\eta},v_{\phi},v_{z})\vert_{v_{\eta}>0}=\frac{1}{\tilde{w}(v)}\int_{u_{\eta}>0}\big(\frac{h_{1}^{n}}{\sqrt{\mu}}+\sqrt{\mu}h_{2}^{n}\big)(t,\eta_{1},\tilde{u})\tilde{w}(u)\,{\rm d}\sigma,\\
&h_{2}^{n+1}(0,\eta,v_{\eta},v_{\phi},v_{z})=0,
\end{aligned}
\right.
\end{align*}
where $\tilde{w}(v)=(\sqrt{2\pi\mu}w^{\ell})^{-1}$ and ${\rm d}\sigma=\sqrt{2\pi}\mu|u_{\eta}|\,{\rm d}u$. Note that along the same characteristic line \eqref{3.16}, $s$ is no longer a parameter and is non-negative. We can represent $h_{1}$ as in \eqref{3.17-4}:
\begin{align}\label{7.21}
h_{1}^{n+1}(t,\eta,v)&={\bf 1}_{\{t_{1}\leq 0\}}h_{1}^{n+1}(0,\mathscr{X}(0),\mathscr{V}(0))\nonumber\\
&\quad+\int_{\max\{t_{1},0\}}^{t}e^{-\int_{s}^{t}\mathcal{B}(\mathscr{V}(\tau))\,{\rm d}\tau}\big[-\frac{\alpha}{2\eta_{1}}\mathscr{V}_{\eta}\mathscr{V}_{\phi}\sqrt{\mu}h_{2}^{n+1}+(1-\chi_{M})\mathcal{K}_{w}h_{1}^{n}\big](s,\mathscr{X}(s),\mathscr{V}(s))\,{\rm d}s\nonumber\\
&\quad +\int_{\max\{t_{1},0\}}^{t}e^{-\int_{s}^{t}\mathcal{B}(\mathscr{V}(\tau))\,{\rm d}\tau}w^{\ell}\tilde{H}(f^{n},f^{n})(s,\mathscr{X}(s),\mathscr{V}(s))\,{\rm d}s,
\end{align}
where
$$
\mathcal{B}(v)=\nu(v)+\frac{\alpha}{\eta_{1}}\frac{\ell v_{\eta}v_{\phi}}{|1+|v|^2}+\mathcal{R}(\alpha\sqrt{\mu}g+\sqrt{\mu}f^{n}).
$$
For $h_{2}^{n+1}(t,\eta,v)$, if $(t,\eta,v)\in [0,T]\times \bar{\mathcal{I}}\times \R^3\backslash \big((\gamma_{0}\cup \gamma_{-})\big)$ with $(\eta,v)\in \mathcal{A}_{1}$, then, similar to \eqref{3.18-3}, one has
\begin{align}\label{7.22}
h_{2}^{n+1}(t,\eta,v)=\sum\limits_{j=1}^{2}J_{j}+\sum\limits_{j=3}^{8}{\bf 1}_{\{t_{1}>0\}}J_{j},
\end{align}
with
\begin{align*}
			J_{1}+J_{2}&=\int_{\max\{t_{1},0\}}^{t}e^{-\int_{s}^{t}\mathcal{B}(\mathscr{V}(\tau))\,{\rm d}\tau}(\chi_{M}\mu^{-\frac{1}{2}} \mathcal{K}_{w}h_{1}^{n}+ K_{w}h_{2}^{n})(s,\mathscr{X}(s),\mathscr{V}(s))\,{\rm d}s,\\
		J_{3}&=\frac{e^{-\int_{t_{1}}^{t}\mathcal{B}(\mathscr{V}(\tau))\,{\rm d}\tau}}{\tilde{w}(\mathscr{V}(t_{1}))}\int_{\prod_{j=1}^{k-1}\mathcal{V}_{j}}\sum\limits_{l=1}^{k-1}{\bf 1}_{\{t_{l}>0\}}(\frac{h_{1}^{n+1-l}}{\sqrt{\mu}})(t_{l},\mathscr{X}(t_{l}),\mathscr{V}(t_{l}))\,{\rm d}\Sigma_{l}(t_{l}),\\
		J_{4}+J_{5}&=\frac{e^{-\int _{t_{1}}^{t}\mathcal{B}(\mathscr{V}(\tau))\,{\rm d}\tau}}{\tilde{w}(\mathscr{V}(t_{1}))}\int_{\prod_{j=1}^{k-1}\mathcal{V}_{j}}\sum\limits_{l=1}^{k-1}{\bf 1}_{\{t_{l+1}>0\}}\int_{t_{l+1}}^{t_{l}}(\chi_{M}\mu^{-\frac{1}{2}}\mathcal{K}_{w}h_{1}^{n-l}+ K_{w}h_{2}^{n-l})(s)\,{\rm d}\Sigma_{l}(s){\rm d}s,\\
		 J_{6}+J_{7}&=\frac{e^{-\int _{t_{1}}^{t}\mathcal{B}(\mathscr{V}(\tau))\,{\rm d}\tau}}{\tilde{w}(\mathscr{V}(t_{1}))}\int_{\prod_{j=1}^{k-1}\mathcal{V}_{j}}\sum\limits_{l=1}^{k-1}{\bf 1}_{\{t_{l+1}\leq 0<t_{l}\}}\int_{0}^{t_{l}}(\chi_{M}\mu^{-\frac{1}{2}}\mathcal{K}_{w}h_{1}^{n-l}+K_{w}h_{2}^{n-l})(s)\,{\rm d}\Sigma_{l}(s){\rm d}s,\\
		  J_{8}&=\frac{e^{-\int _{t_{1}}^{t}\mathcal{B}(\mathscr{V}(\tau))\,{\rm d}\tau}}{\tilde{w}(\mathscr{V}(t_{1}))}\int_{\prod_{j=1}^{k-1}\mathcal{V}_{j}}{\bf 1}_{\{t_{k}>0\}}h_{2}^{n+1-k}(t_{k},\mathscr{X}(t_{k}),\mathscr{V}(t_{k}))\,{\rm d}\Sigma_{k-1}(t_{k}).
	\end{align*}
Here, $d\Sigma_{l}(s)$ is similarly defined in \eqref{3.18-2} by replacing $\mathcal{A}(\mathscr{V}(\tau))$ with $\mathcal{B}(\mathscr{V}(\tau))$.

If $(t,\eta,v)\in [0,T]\times \bar{\mathcal{I}}\times \R^3\backslash \big((\gamma_{0}
\cup \gamma_{-})\big)$ with $(\eta,v)\in \mathcal{A}_{2}$, then, similar to \eqref{3.19-3}, one has
\begin{align*}
		h_{2}^{n+1}(\eta,v)&=
	 \int_{\max\{t_{1,*},0\}}^{t}e^{-\int_{s}^{t}\mathcal{B}(\mathscr{V}(\tau))\,{\rm d}\tau}\big(\chi_{M}\mu^{-\frac{1}{2}}h_{1}^{n}+ K_{w}h_{2}^{n}\big)(s)\,{\rm d}s\nonumber\\
		&\quad +{\bf 1}_{\{t_{1,*}>0\}}e^{-\int_{t_{1,*}}^{t}\mathcal{B}(\mathscr{V}(\tau))\,{\rm d}\tau}h_{2}^{n+1}(t_{1,*},\mathscr{X}(t_{1,*}),\mathscr{V}(t_{1,*})).
\end{align*}
Recall that $t_{1,*}$ is the time when $\mathscr{V}_{\eta}(t_{1,*})=0$. We can further express $h_{2}^{n+1}(t_{1,*},\mathscr{X}(t_{1,*}),\mathscr{V}(t_{1,*}))$ as in \eqref{7.22}.
        
If $(t,\eta,v)\in [0,T]\times \bar{\mathcal{I}}\times \R^3\backslash \big((\gamma_{0}
\cup \gamma_{-})\big)$ with $(\eta,v)\in \mathcal{A}_{3}$, similar to \eqref{3.19-4}, it holds that
	\begin{align}\label{7.24}
		h_{2}^{n+1}(\eta,v)&=
	\int_{\max\{t_{1},0\}}^{t}e^{-\int_{s}^{t}\mathcal{B}(\mathscr{V}(\tau))\,{\rm d}\tau}\big(\chi_{M}\mu^{-\frac{1}{2}}h_{1}^{n}+ K_{w}h_{2}^{n}\big)(s)\,{\rm d}s\nonumber\\
			&\quad +{\bf 1}_{\{t_{1}>0\}}e^{-\int_{t_{1}}^{t}\mathcal{B}(\mathscr{V}(\tau))\,{\rm d}\tau}\frac{1}{\tilde{w}(\mathscr{V}(t_{1}))}\int_{\mathcal{V}_{1}}\tilde{w}(v_{1})\frac{h_{1}^{n}}{\sqrt{\mu}}(t_{1},\eta_{1},v_{1,\eta},v_{1,\phi})\,{\rm d}\sigma_{1}\nonumber\\
			&\quad+{\bf 1}_{\{t_{1}>0\}}\,e^{-\int_{t_{1}}^{t}\mathcal{B}(\mathscr{V}(\tau))\,{\rm d}\tau}\frac{1}{\tilde{w}(\mathscr{V}(t_{1}))}\int_{\mathcal{V}_{1}}\tilde{w}(v_{1})h_{2}^{n}(t_{1},\eta_{1},v_{1,\eta},v_{1,\phi})\,{\rm d}\sigma_{1}.
		\end{align}
		Since $x_{1}=\eta_{1}$ in this case, it holds that $v_{1}\in \mathcal{V}_{1}:=\{u\in \R^3:u_{\eta}> 0\}$, which implies $v_{1,\eta}>0$. Therefore we can further express $h_{2}^{n}(t_{1},\eta_{1},v_{1})$ by using the representation formula in \eqref{7.22}.

For $n=0$, it is direct to see that
$$
\mathcal{B}\geq \nu-C\alpha-C\|w^{\ell}f^{0}\|_{L^{\infty}}\geq \nu_{0}-C(\alpha+\varepsilon_{0})\geq \frac{\nu_{0}}{2}>0,
$$
provided $\alpha$ and $\varepsilon_{0}$ suitably small. Then it follows form \eqref{7.21}--\eqref{7.24} and direct calculations that for any $t\geq 0$,
$$
\sup_{0\leq s\leq t}\|(h_{1}^{1},h_{2}^{1})(s)\|_{L^{\infty}}\leq C\|w^{\ell}f^{0}\|_{L^{\infty}}\leq C\varepsilon_{0}.
$$
Moreover, by induction, for any fixed finite $k_{0}>0$ and $t\geq 0$, we have
\begin{align}\label{7.8-0}
\sup_{0\leq s\leq t}\sup_{0\leq \ell \leq k_{0}}\|(h_{1}^{l},h_{2}^{l})(s)\|_{L^{\infty}}\leq C_{k_{0}}\|w^{\ell}f_{0}\|_{L^{\infty}}<\frac{1}{2}C_{0}\varepsilon_{0},
\end{align}
by choosing $C_{0}$ suitably large. Noting that $C_{0}$ depending on $k_{0}$, we can further choose $\varepsilon_{0}$ suitably small such that \eqref{7.1-2} holds for $0\leq m\leq k_{0}$.

In the following, we are going to prove \eqref{7.1-2} for $m=n+1$ under the assumption that it holds for $m\leq n$. It follows from the induction assumption that
$$
\mathcal{B}\geq \nu_{0}-C\alpha-C\|w^{\ell}f^{n}\|_{L^{\infty}}\geq \nu_{0}-C\alpha-CC_{0}\varepsilon_{0}\geq \frac{1}{2}\nu_{0}>0,
$$
provided that $\alpha$ and $\varepsilon_{0}$ small. Applying Lemma \ref{K}, we obtain from \eqref{7.21} that
\begin{align}\label{7.25}
\sup_{0\leq t\leq T_{*}}\|h_{1}^{n+1}\|_{L^{\infty}}&\leq C\|w^{\ell}f_{0}\|_{L^{\infty}}+C\big[{\bf 1}_{\{\gamma=0\}}\frac{1}{\ell}+{\bf 1}_{\{0\leq \gamma<1\}}[(1+M(\ell))^{-\gamma}+\zeta(\ell)]\big]\sup_{0\leq t\leq T_{*}}\|h_{1}^{n}\|_{L^{\infty}}\nonumber\\
&\quad +C\alpha\sup_{0\leq t\leq T_{*}}\|h_{2}^{n+1}\|_{L^{\infty}}+C_{\ell}\alpha\|(h_{1}^{n},h_{2}^{n})\|_{L^{\infty}}+C\|(h_{1}^{n},h_{2}^{n})\|_{L^{\infty}}^2\nonumber\\
&\leq C\varepsilon_{0}+\frac{1}{2}\sup_{0\leq t\leq T_{*}}\|h_{1}^{n}\|_{L^{\infty}}+C\alpha\|h_{2}^{n+1}\|_{L^{\infty}}+C(\alpha+C_{0}\varepsilon_{0})C_{0}\varepsilon_{0},
\end{align}
provided that $\ell$, $M(\ell)$ are suitably large and $\zeta(\ell)$ is suitably small. By iterating \eqref{7.25} $k_{0}$ times, one has
\begin{align}\label{7.25-1}
\sup_{0\leq t\leq T_{*}}\|h_{1}^{n+1}\|_{L^{\infty}}&\leq \Big(\frac{1}{2}\Big)^{k_{0}}\sup_{0\leq t\leq T_{*}}\|h_{1}^{n+1-k_{0}}\|_{L^{\infty}}+2C\alpha\sup_{0\leq l\leq k_{0}}\sup_{0\leq t\leq T_{*}}\|h_{2}^{n+1-l}\|_{L^{\infty}}\nonumber\\
&\quad +2C\varepsilon_{0}+2C(\alpha+C_{0}\varepsilon_{0})C_{0}\varepsilon_{0}.
\end{align}
Furthermore, by letting $t\leq T_{*}$ with $T_{*}>0$ being suitably small, we obtain from \eqref{7.22}--\eqref{7.24} that
\begin{align}\label{7.26}
\sup_{0\leq t\leq T_{*}}\|h_{2}^{n+1}\|_{L^{\infty}}&\leq C_{\ell}k_{0}\sup_{1\leq l\leq k}\sup_{0\leq t\leq T_{*}}\|h_{1}^{n+1-l}\|_{L^{\infty}}+C_{\ell}(k_{0}T_{*}+\tilde{\delta})\sup_{1\leq l\leq k_{0}}\sup_{0\leq t\leq T_{*}}\|\|h_{2}^{n+1-l}\|_{L^{\infty}}\nonumber\\
&\leq \frac{1}{8}\sup_{1\leq l\leq k_{0}}\sup_{0\leq t\leq T_{*}}\|h_{2}^{n+1-l}\|_{L^{\infty}}+C_{\ell}k_{0}\alpha\sup_{1\leq l\leq 2k_{0}}\|h_{2}^{n+1-l}\|_{L^{\infty}}\nonumber\\
&\quad +C_{\ell}k_{0}\Big(\frac{1}{2}\Big)^{k_{0}}\sup_{1\leq l\leq 2k_{0}}\sup_{0\leq t\leq T_{*}}\|h_{1}^{n+1-l}\|_{L^{\infty}}+C_{\ell}k_{0}\varepsilon_{0}+C_{\ell}k_{0}(\alpha+C_{0}\varepsilon_{0})C_{0}\varepsilon_{0}\nonumber\\
&\leq \frac{1}{8}\sup_{1\leq l\leq 2k_{0}}\sup_{0\leq t\leq T_{*}}\|h_{2}^{n+1-l}\|_{L^{\infty}} +C_{\ell}k_{0}\Big[\varepsilon_{0}+\Big(\Big(\frac{1}{2}\Big)^{k_{0}}+\alpha+C_{0}\varepsilon_{0}\Big)C_{0}\varepsilon_{0}\Big],
\end{align}
provided that $C_{\ell}(k_{0}T_{*}+k_{0}\alpha+\tilde{\delta})\leq \frac{1}{8}$. Here, we have used \eqref{7.25} in the second inequality, and the induction assumption in the last inequality, and $\tilde{\delta}$ is arbitrarily small constant, $k_{0}$ is the correspondingly large constant determined in Lemma \ref{lem2.1}. Applying Lemma \ref{lemC1} to \eqref{7.26}, we have
\begin{align}\label{7.27}
\sup_{0\leq t\leq T_{*}}\|h_{2}^{n+1}\|_{L^{\infty}}&\leq \Big(\frac{1}{8}\Big)^{\big[\frac{n}{k_{0}}\big]}\sup_{1\leq l\leq 4k_{0}}\sup_{0\leq t\leq T_{*}}\|h_{2}^{l}\|_{L^{\infty}}+\frac{8+k_{0}}{7}C_{\ell}k_{0}\Big[\varepsilon_{0}+\Big(\frac{1}{2}\Big)^{k_{0}}+\alpha+C_{0}\varepsilon_{0}\Big]C_{0}\varepsilon_{0}.
\end{align}
Then, noting \eqref{7.8-0} and taking $k_{0}$, $C_{0}$ suitably large, and $\alpha$ and $\varepsilon_{0}$ suitably small, we conclude \eqref{7.1-2} from \eqref{7.27} and \eqref{7.25-1}.

Now, we shall prove that $\{(f_{1}^{n},f_{2}^{n})\}$ is a Cauchy sequence in $\mathbf{Y}_{T_{*}}$. We denote
$$
\tilde{h}_{1}^{n+1}=h_{1}^{n+1}-h_{1}^{n},\quad \tilde{h}_{2}^{n+1}=h_{2}^{n+1}-h_{2}^{n}.
$$
Then $\tilde{h}_{1}^{n+1}$ and 
$\tilde{h}_{2}^{n+1}$ satisfy
\begin{align*}
\left\{
\begin{aligned}
&\partial_{t}\tilde{h}_{1}^{n+1}+v_{\eta}\frac{\partial \tilde{h}_{1}^{n+1}}{\partial \eta}-\frac{1}{1-\eta}\Big(v_{\phi}^2\frac{\partial \tilde{h}_{1}^{n+1}}{\partial v_{\eta}}-v_{\eta}v_{\phi}\frac{\partial \tilde{h}_{1}^{n+1}}{\partial v_{\phi}}\Big)-\frac{\alpha}{\eta_{1}}v_{\eta}\frac{\partial \tilde{h}_{1}^{n+1}}{\partial v_{\phi}}+\nu\tilde{h}_{1}^{n+1}+\frac{\alpha}{2\eta_{1}}\sqrt{\mu}v_{\eta}v_{\phi}\tilde{h}_{2}^{n+1}\\
&\quad =-\tilde{h}_{1}^{n+1}\mathcal{R}(\alpha\sqrt{\mu}g+\sqrt{\mu}f_{1}^{n})-\frac{\alpha}{\eta_{1}}\frac{\ell v_{\eta}v_{\phi}}{1+|v|^2}\tilde{h}_{1}^{n+1}+(1-\chi_{M})\mathcal{K}_{w}\tilde{h}_{1}^{n}+\tilde{\mathscr{H}},\\
&\tilde{h}_{1}^{n+1}(t,0,v_{\eta},v_{\phi},v_{z})\vert_{v_{\eta}>0}=\tilde{h}_{1}^{n+1}(t,\eta_{1},v_{\eta},v_{\phi},v_{z})\vert_{v_{\eta}<0}=0,\\
&\tilde{h}_{1}^{n+1}(0,\eta_{1},v_{\eta},v_{\phi},v_{z})=0,
\end{aligned}
\right.
\end{align*}
and
\begin{align*}
\left\{
\begin{aligned}
&\partial_{t}\tilde{h}_{2}^{n+1}+v_{\eta}\frac{\partial \tilde{h}_{2}^{n+1}}{\partial \eta}-\frac{1}{1-\eta}\Big(v_{\phi}^2\frac{\partial \tilde{h}_{2}^{n+1}}{\partial v_{\eta}}-v_{\eta}v_{\phi}\frac{\partial \tilde{h}_{2}^{n+1}}{\partial v_{\phi}}\Big)-\frac{\alpha}{\eta_{1}}v_{\eta}\frac{\partial \tilde{h}_{2}^{n+1}}{\partial v_{\phi}}+\nu\tilde{h}_{2}^{n+1}\\
    &\quad =-\frac{\alpha}{\eta_{1}}\frac{\ell v_{\eta}v_{\phi}}{1+|v|^2}\tilde{h}_{2}^{n+1}-\tilde{h}_{2}^{n+1}\mathcal{R}(\alpha\sqrt{\mu}g+\sqrt{\mu}f^{n})+\chi_{M}\mathcal{K}_{w}\tilde{h}_{1}^{n}+K_{w}\tilde{h}_{2}^{n}-h_{2}^{n}\mathcal{R}(\sqrt{\mu}\tilde{f}^{n}),\\
&\tilde{h}_{2}^{n+1}(t,0,v_{\eta},v_{\phi},v_{z})\vert_{v_{\eta}>0}=\frac{1}{\tilde{w}}\int_{u_{\eta}<0}\big(\frac{\tilde{h}_{1}^{n}}{\sqrt{\mu}}+\sqrt{\mu}\tilde{h}_{2}^{n}\big)(t,0,u)\tilde{w}(u)\,{\rm d}\sigma,\\
& \tilde{h}_{2}^{n+1}(t,\eta_{1},v_{\eta},v_{\phi},v_{z})\vert_{v_{\eta}<0}=\frac{1}{\tilde{w}}\int_{u_{\eta}>0}\big(\frac{\tilde{h}_{1}^{n}}{\sqrt{\mu}}+\sqrt{\mu}\tilde{h}_{2}^{n}\big)(t,\eta_{1},u)\tilde{w}(u)\,{\rm d}\sigma,
\end{aligned}
\right.
\end{align*}
where we have denoted
$$
\tilde{\mathscr{H}}=-h_{1}^{n}\mathcal{R}(\sqrt{\mu}(f^{n}-f^{n-1}))+w^{\ell}\big[\tilde{H}(f^{n},f^{n})-\tilde{H}(f^{n-1},f^{n-1})\big].
$$
Similar calculations as in \eqref{7.25} show that
\begin{align*}
\sup_{0\leq t\leq T^{*}}\|\tilde{h}_{1}^{n+1}\|_{L^{\infty}}&\leq C\big[\alpha+\varepsilon_{0}+{\bf 1}_{\{\gamma=0\}}\frac{1}{\ell}+{\bf 1}_{\{0<\gamma\leq 1\}}[(1+M(\ell))^{-\gamma}+\zeta(\ell)]\big]\sup_{0\leq t\leq T^{*}}\|\tilde{h}_{1}^{n}\|_{L^{\infty}}\nonumber\\
&\quad +C(\alpha+\varepsilon_{0})\sup_{0\leq l\leq 1}\sup_{0\leq t\leq T^{*}}\|\tilde{h}_{2}^{n+1-1}\|_{L^{\infty}}\nonumber\\
&\leq \frac{1}{2}\sup_{0\leq t\leq T^{*}}\|\tilde{h}_{1}^{n}\|_{L^{\infty}}+C(\alpha+\varepsilon_{0})\sup_{0\leq l\leq 1}\sup_{0\leq t\leq T^{*}}\|\tilde{h}_{2}^{n+1-1}\|_{L^{\infty}},
\end{align*}
which yields that
\begin{align}\label{7.31}
\sup_{0\leq t\leq T^{*}}\|\tilde{h}_{1}^{n+1}\|_{L^{\infty}}\leq \Big(\frac{1}{2}\Big)^{k_{0}}\sup_{0\leq t\leq T^{*}}\|\tilde{h}_{1}^{n+1-k_{0}}\|_{L^{\infty}}+2C(\alpha+\varepsilon_{0})\sup_{0\leq l\leq k_{0}}\sup_{0\leq t\leq T^{*}}\|\tilde{h}_{2}^{n+1-k_{0}}\|_{L^{\infty}}.
\end{align}
Similar to \eqref{7.26}, we have
\begin{align}\label{7.32}
\sup_{0\leq t\leq T^{*}}\|\tilde{h}_{2}^{n+1}\|_{L^{\infty}}&\leq C_{\ell}k_{0}\sup_{1\leq l\leq k_{0}}\sup_{0\leq t\leq T_{*}}\|\tilde{h}_{1}^{n+1-l}\|_{L^{\infty}}+C_{\ell}(k_{0}T_{*}+\tilde{\delta})\sup_{1\leq l\leq k_{0}}\sup_{0\leq t\leq T_{*}}\|\tilde{h}_{2}^{n+1-l}\|_{L^{\infty}}\nonumber\\
&\leq C_{\ell}k_{0}\Big(\frac{1}{2}\Big)^{k_{0}}\sup_{1\leq l\leq 2k_{0}}\sup_{0\leq t\leq T_{*}}\|\tilde{h}_{1}^{n+1-l}\|_{L^{\infty}}\nonumber\\
&\quad+2C_{\ell}(k_{0}\alpha+k_{0}\varepsilon_{0}+k_{0}T_{*}+\tilde{\delta})\sup_{1\leq l\leq 2k_{0}}\sup_{0\leq t\leq T_{*}}\|\tilde{h}_{2}^{n+1-l}\|_{L^{\infty}}.
\end{align}
Then, we obtain from \eqref{7.31}--\eqref{7.32} that
\begin{align}\label{7.33}
&\sup_{0\leq t\leq T_{*}}\|(\tilde{h}_{1}^{n+1},\tilde{h}_{2}^{n+1})\|_{L^{\infty}}\leq C_{\ell}\Big[k_{0}\Big(\frac{1}{2}\Big)^{k_{0}}+k_{0}(\alpha+\varepsilon_{0})+\tilde{\delta}\Big]\sup_{1\leq l\leq k_{0}}\sup_{0\leq t\leq T_{*}}\|(\tilde{h}_{1}^{n+1-l},\tilde{h}_{2}^{n+1-l})\|_{L^{\infty}}.
\end{align}
Taking $\tilde{\delta}$ suitably small and $k_{0}$ suitably large, and then $\alpha$, $\varepsilon_{0}$ and $T_{*}$ suitably small, we apply Lemma \ref{lemC1} to \eqref{7.33} to get
 \begin{align*}
\sup_{0\leq t\leq T_{*}}\|(\tilde{h}_{1}^{n+1},\tilde{h}_{2}^{n+1})\|_{L^{\infty}}&\leq \Big(\frac{1}{8}\Big)^{\big[\frac{n}{2k_{0}+1}\big]}\sup_{1\leq l\leq 2(2k_{0}+1)}\sup_{0\leq t\leq T_{*}}\|(\tilde{h}_{1}^{l},\tilde{h}_{2}^{l})\|_{L^{\infty}}\leq \Big(\frac{1}{8}\Big)^{\big[\frac{n}{2k_{0}+1}\big]}C_{0}\varepsilon_{0},
\end{align*}
which implies that $\{(f_{1}^n,f_{2}^n)\}$ is a Cauchy sequence in $\mathbf{Y}_{T_{*}}$. Therefore, there is a unique $(f_{1},f_{2})\in\mathbf{Y}_{T_{*}}$ such that $(f_{1}^{n},f_{2}^n)$ converges strongly to $(f_{1},f_{2})$ as $n\to \infty$, and $(f_{1},f_{2})$ is a desired a local-in-time solution of the coupled systems \eqref{6.1-2}--\eqref{6.1-3}. Furthermore, noting that $G^{n}=G_{s}+\sqrt{\mu}f^{n}\geq 0$ for any $n$, we have
$$
G=G_{s}+f_{1}+\sqrt{\mu}f_{2}=\lim\limits_{n\to \infty}(G_{s}+f_{1}^{n}+\sqrt{\mu}f_{2}^{n})\geq 0.
$$
Hence, the proof of Theorem \ref{thm4.1} is completed. $\hfill\square$

\subsection{{\it A priori} $L^{\infty}$ estimates}\label{sec7.2}
We are in position to show the global existence and large time behavior of solution to the initial boundary value problem \eqref{6.1-2}--\eqref{6.1-3}. Firstly, we focus on the uniform $L^{\infty}\cap L^2$ estimates under the {\it a priori} assumption
\begin{align}\label{8.1}
\sup_{s\geq 0}e^{\lambda_{0}s}\|w^{\ell}(f_{1},f_{2})\|_{L^{\infty}}\leq \tilde{\varepsilon},
\end{align}
for a constant $\tilde{\varepsilon}>0$ suitably small, where $\lambda_{0}>0$ independent of $\alpha$ is to be determined later.

\begin{lemma}\label{lem8.1}
Let $0<\lambda_{0}<\frac{\nu_{0}}{4}$, then under the assumption \eqref{8.1}, it holds that
\begin{align}
&\sup_{0\leq s\leq t}e^{\lambda_{0}s}\|w^{\ell}f_{1}(s)\|_{L^{\infty}}\leq C_{\ell}\|w^{\ell}f_{0}\|_{L^{\infty}}+C(\alpha+\tilde{\varepsilon})\sup_{0\leq s\leq t}e^{\lambda_{0}s}\|w^{\ell}f_{2}(s)\|_{L^{\infty}},\label{8.2}\\
&\sup_{0\leq s\leq t}e^{\lambda_{0}s}\|w^{\ell}f_{2}(s)\|_{L^{\infty}}\leq C_{\ell}\|w^{\ell}f_{0}\|_{L^{\infty}}+C_{\ell}\sup_{0\leq s\leq t}e^{\lambda_{0}s}\|f_{2}(s)\|_{L^2},
\label{8.3}
\end{align}
for any $t\geq 0$.
\end{lemma}

\noindent\textbf{Proof}. For simplicity of notation, we set
\begin{align}\label{def.8.4}
(\mathfrak{h}_{1},\mathfrak{h}_{2})=e^{\lambda_{0}t}w^{\ell}(f_{1},f_{2})(t,y,v)
\end{align}
with $\lambda_{0}$ to be chosen. Then, it follows from \eqref{6.1-2}--\eqref{6.1-3} that $(\mathfrak{h}_{1},\mathfrak{h}_{2})$ satisfies
\begin{align*}
\left\{
\begin{aligned}
	&\partial_{t}\mathfrak{h}_{1}+v_{\eta}\frac{\partial \mathfrak{h}_{1}}{\partial \eta}-\frac{1}{1-\eta}\Big(v_{\phi}^2\frac{\partial \mathfrak{h}_{1}}{\partial v_{\eta}}-v_{\eta}v_{\phi}\frac{\partial \mathfrak{h}_{1}}{\partial v_{\phi}}\Big)-\frac{\alpha}{\eta_{1}}v_{\eta}\frac{\partial \mathfrak{h}_{1}}{\partial v_{\phi}}+\frac{\alpha}{2\eta_{1}}v_{\eta}v_{\phi}\sqrt{\mu}\mathfrak{h}_{2}+\nu \mathfrak{h}_{1}-\lambda_{0}\mathfrak{h}_{1}\\
	&\qquad =-\frac{\alpha}{\eta_{1}}\frac{\ell v_{\eta}v_{\phi}}{1+|v|^2}\mathfrak{h}_{1}+(1-\chi_{M})\mathcal{K}_{w}\mathfrak{h}_{1}+w^{\ell}e^{\lambda_{0}t}H(f,f),\\
	&\mathfrak{h}_{1}(t,0,v_{\eta},v_{\phi},v_{z})\vert_{v_{\eta}>0}=\mathfrak{h}_{1}(t,\eta_{1},v_{\eta},v_{\phi},v_{z})\vert_{v_{\eta}<0}=0,\\
	&\mathfrak{h}_{1}(0,\eta,v_{\eta},v_{\phi},v_{z})=w^{\ell}\big(G_{0}(\eta,v_{\eta},v_{\phi},v_{z})-G_{s}(\eta,v_{\eta},v_{\phi},v_{z})\big),
\end{aligned}
\right.
\end{align*}
and
\begin{align*}
	\left\{
	\begin{aligned}
		&\partial_{t}\mathfrak{h}_{2}+v_{\eta}\frac{\partial \mathfrak{h}_{2}}{\partial \eta}-\frac{1}{1-\eta}\Big(v_{\phi}^2\frac{\partial \mathfrak{h}_{2}}{\partial v_{\eta}}-v_{\eta}v_{\phi}\frac{\partial \mathfrak{h}_{2}}{\partial v_{\phi}}\Big)-\frac{\alpha}{\eta_{1}}v_{\eta}\frac{\partial \mathfrak{h}_{2}}{\partial v_{\phi}}+\nu \mathfrak{h}_{2}-\lambda_{0}\mathfrak{h}_{2}\\
		&\qquad =-\frac{\alpha}{\eta_{1}}\frac{\ell v_{\eta}v_{\phi}}{1+|v|^2}\mathfrak{h}_{2}+\chi_{M}\mathcal{K}_{w}\mathfrak{h}_{1}+K_{w}\mathfrak{h}_{2},\\
		&\mathfrak{h}_{2}(t,0,v_{\eta},v_{\phi},v_{z})\vert_{v_{\eta}>0}=\frac{1}{\tilde{w}}\int_{u_{\eta}<0}(\frac{\mathfrak{h}_{1}}{\sqrt{\mu}}+\mathfrak{h}_{2})(t,0,u)\tilde{w}(u)\,{\rm d}\sigma,\\
		&\mathfrak{h}_{2}(t,\eta_{1},v_{\eta},v_{\phi},v_{z})\vert_{v_{\eta}<0}=\frac{1}{\tilde{w}}\int_{u_{\eta}>0}(\frac{\mathfrak{h}_{1}}{\sqrt{\mu}}+\mathfrak{h}_{2})(t,\eta_{1},u)\tilde{w}(u)\,{\rm d}\sigma,\\
		&\mathfrak{h}_{2}(0,\eta,v_{\eta},v_{\phi},v_{z})=0.
	\end{aligned}
	\right.
\end{align*}
Along the backward characteristic line, we can represent $\mathfrak{h}_{1}$ as follows
\begin{align*}
\mathfrak{h}_{1}(t,\eta,v)&={\bf 1}_{\{t_{1}\leq 0\}}\mathfrak{h}_{1}(0,\mathscr{X}(0),\mathscr{V}(0))\nonumber\\
&\quad+\int_{\max\{t_{1},0\}}^{t}e^{-\int_{s}^{t}\tilde{\mathcal{A}}(\mathscr{V}(\tau))\,{\rm d}\tau}\big[-\frac{\alpha}{2\eta_{1}}\mathscr{V}_{\eta}\mathscr{V}_{\phi}\sqrt{\mu}\mathfrak{h}_{2}+(1-\chi_{M})\mathcal{K}_{w}\mathfrak{h}_{1}\big](s,\mathscr{X}(s),\mathscr{V}(s))\,{\rm d}s\nonumber\\
&\quad +\int_{\max\{t_{1},0\}}^{t}e^{-\int_{s}^{t}\tilde{\mathcal{A}}(\mathscr{V}(\tau))\,{\rm d}\tau}w^{\ell}H(f,f)(s,\mathscr{X}(s),\mathscr{V}(s))\,{\rm d}s,
\end{align*}
where
\begin{align*}
    \tilde{\mathcal{A}}=\nu+\frac{\alpha}{\eta_{1}}\frac{\ell v_{\eta}v_{\phi}}{1+|v|^2}-\lambda_{0}\geq \frac{1}{4}\nu_{0}\quad \text{provided $\alpha\ell\ll 1$ and $\lambda_{0}\leq \frac{\nu_{0}}{4}$}.
\end{align*}
Then, by using Lemma \ref{K} and direct calculations, $\mathfrak{h}_{1}$ can be bounded as
\begin{align}\label{8.9}
\sup_{0\leq s\leq t}\|\mathfrak{h}_{1}(s)\|_{L^{\infty}}&\leq C\|w^{\ell}f_{0}\|_{L^{\infty}}+C(\alpha+\tilde{\varepsilon})\big[\sup_{0\leq s\leq t}\|\mathfrak{h}_{1}(s)\|_{L^{\infty}}+\sup_{0\leq s\leq t}\|\mathfrak{h}_{2}(s)\|_{L^{\infty}}\big]\nonumber\\
&\quad +C\big[{\bf 1}_{\{\gamma=0\}}\frac{1}{\ell}+{\bf 1}_{\{0\leq \gamma<1\}}\big[(1+M(\ell))^{-\gamma}+\zeta(\ell)\big]\big]\sup\limits_{0\leq s\leq t}\|\mathfrak{h}_{1}(s)\|_{L^{\infty}}.
\end{align}
Taking $\ell$ large enough and $\alpha$, $\tilde{\varepsilon}$ small enough, we obtain from \eqref{8.9} that
\begin{align*}
\sup_{0\leq s\leq t}\|\mathfrak{h}_{1}(s)\|_{L^{\infty}}&\leq C\|w^{\ell}f_{0}\|_{L^{\infty}}+C(\alpha+\tilde{\varepsilon})\sup_{0\leq s\leq t}\|\mathfrak{h}_{2}(s)\|_{L^{\infty}}.
\end{align*}
Recall \eqref{def.8.4}. The above estimate then concludes the proof of \eqref{8.2}.

For $\mathfrak{h}_{2}$, if $(t,\eta,v)\in [0,T]\times \bar{\mathcal{I}}\times \R^3\backslash \big((\gamma_{0}
\cup \gamma_{-})\big)$ with $(\eta,v)\in \mathcal{A}_{1}$, then, similar to \eqref{3.18-3}, one has
\begin{align}\label{8.11}
\mathfrak{h}_{2}(t,\eta,v)=\sum\limits_{j=1}^{2}J_{j}+\sum\limits_{j=3}^{8}{\bf 1}_{\{t_{1}>0\}}J_{j},
\end{align}
with
\begin{align*}
			J_{1}+J_{2}&=\int_{\max\{t_{1},0\}}^{t}e^{-\int_{s}^{t}\tilde{\mathcal{A}}(\mathscr{V}(\tau))\,{\rm d}\tau}(\chi_{M}\mu^{-\frac{1}{2}} \mathcal{K}_{w}\mathfrak{h}_{1}+ K_{w}\mathfrak{h}_{2})(s,\mathscr{X}(s),\mathscr{V}(s))\,{\rm d}s,\\
		J_{3}&=\frac{e^{-\int_{t_{1}}^{t}\tilde{\mathcal{A}}(\mathscr{V}(\tau))\,{\rm d}\tau}}{\tilde{w}(\mathscr{V}(t_{1}))}\int_{\prod_{j=1}^{k-1}\mathcal{V}_{j}}\sum\limits_{l=1}^{k-1}{\bf 1}_{\{t_{l}>0\}}(\frac{\mathfrak{h}_{1}}{\sqrt{\mu}})(t_{l},\mathscr{X}(t_{l}),\mathscr{V}(t_{l}))\,{\rm d}\Sigma_{l}(t_{l}),\\
		J_{4}+J_{5}&=\frac{e^{-\int _{t_{1}}^{t}\tilde{\mathcal{A}}(\mathscr{V}(\tau))\,{\rm d}\tau}}{\tilde{w}(\mathscr{V}(t_{1}))}\int_{\prod_{j=1}^{k-1}\mathcal{V}_{j}}\sum\limits_{l=1}^{k-1}{\bf 1}_{\{t_{l+1}>0\}}\int_{t_{l+1}}^{t_{l}}(\chi_{M}\mu^{-\frac{1}{2}}\mathcal{K}_{w}\mathfrak{h}_{1}+ K_{w}\mathfrak{h}_{2})(s)\,{\rm d}\Sigma_{l}(s){\rm d}s,\\
		 J_{6}+J_{7}&=\frac{e^{-\int _{t_{1}}^{t}\tilde{\mathcal{A}}(\mathscr{V}(\tau))\,{\rm d}\tau}}{\tilde{w}(\mathscr{V}(t_{1}))}\int_{\prod_{j=1}^{k-1}\mathcal{V}_{j}}\sum\limits_{l=1}^{k-1}{\bf 1}_{\{t_{l+1}\leq 0<t_{l}\}}\int_{0}^{t_{l}}(\chi_{M}\mu^{-\frac{1}{2}}\mathcal{K}_{w}\mathfrak{h}_{1}+K_{w}\mathfrak{h}_{1})(s)\,{\rm d}\Sigma_{l}(s){\rm d}s,\\
		  J_{8}&=\frac{e^{-\int _{t_{1}}^{t}\tilde{\mathcal{A}}(\mathscr{V}(\tau))\,{\rm d}\tau}}{\tilde{w}(\mathscr{V}(t_{1}))}\int_{\prod_{j=1}^{k-1}\mathcal{V}_{j}}{\bf 1}_{\{t_{k}>0\}}\mathfrak{h}_{2}(t_{k},\mathscr{X}(t_{k}),\mathscr{V}(t_{k}))\,{\rm d}\Sigma_{k-1}(t_{k}).
	\end{align*}
If $(t,\eta,v)\in [0,T]\times \bar{\mathcal{I}}\times \R^3\backslash \big((\gamma_{0}
\cup \gamma_{-})\big)$ with $(\eta,v)\in \mathcal{A}_{2}$, then, similar to \eqref{3.19-3}, one has
\begin{align*}
		\mathfrak{h}_{2}(\eta,v)&=
	 \int_{\max\{t_{1,*},0\}}^{t}e^{-\int_{s}^{t}\tilde{\mathcal{A}}(\mathscr{V}(\tau))\,{\rm d}\tau}\big(\chi_{M}\mu^{-\frac{1}{2}}\mathfrak{h}_{1}+ K_{w}\mathfrak{h}_{2}\big)(s)\,{\rm d}s\nonumber\\
		&\quad +{\bf 1}_{\{t_{1,*}>0\}}e^{-\int_{t_{1,*}}^{t}\tilde{\mathcal{A}}(\mathscr{V}(\tau))\,{\rm d}\tau}\mathfrak{h}_{2}(t_{1,*},\mathscr{X}(t_{1,*}),\mathscr{V}(t_{1,*})).
\end{align*}
		Recall that $t_{1,*}$ is the time when $\mathscr{V}_{\eta}(t_{1,*})=0$. We can further express $\mathfrak{h}_{2}(t_{1,*},\mathscr{X}(t_{1,*}),\mathscr{V}(t_{1,*}))$ as in \eqref{8.11}.
        
If $(t,\eta,v)\in [0,T]\times \bar{\mathcal{I}}\times \R^3\backslash \big((\gamma_{0}
\cup \gamma_{-})\big)$ with $(\eta,v)\in \mathcal{A}_{3}$, similar to \eqref{3.19-4}, it holds that
\begin{align*}
		\mathfrak{h}_{2}(\eta,v)&=
	\int_{\max\{t_{1},0\}}^{t}e^{-\int_{s}^{t}\tilde{\mathcal{A}}(\mathscr{V}(\tau))\,{\rm d}\tau}\big(\chi_{M}\mu^{-\frac{1}{2}}\mathfrak{h}_{1}+ K_{w}\mathfrak{h}_{2}\big)(s)\,{\rm d}s\nonumber\\
			&\quad +{\bf 1}_{\{t_{1}>0\}}e^{-\int_{t_{1}}^{t}\tilde{\mathcal{A}}(\mathscr{V}(\tau))\,{\rm d}\tau}\frac{1}{\tilde{w}(\mathscr{V}(t_{1}))}\int_{\mathcal{V}_{1}}\tilde{w}(v_{1})\frac{\mathfrak{h}_{1}}{\sqrt{\mu}}(t_{1},\eta_{1},v_{1})\,{\rm d}\sigma_{1}\nonumber\\
			&\quad+{\bf 1}_{\{t_{1}>0\}}\,e^{-\int_{t_{1}}^{t}\tilde{\mathcal{A}}(\mathscr{V}(\tau))\,{\rm d}\tau}\frac{1}{\tilde{w}(\mathscr{V}(t_{1}))}\int_{\mathcal{V}_{1}}\tilde{w}(v_{1})\mathfrak{h}_{2}(t_{1},\eta_{1},v_{1})\,{\rm d}\sigma_{1}.
\end{align*}
Since in this case $x_{1}=\eta_{1}$, then $v_{1}\in \mathcal{V}_{1}:=\{u\in \R^3:u_{\eta}> 0\}$, which implies $v_{1,\eta}>0$. Therefore we can further express $\mathfrak{h}_{2}(t_{1},\eta_{1},v_{1})$ by using the representation formula in \eqref{8.11}.
	
By using Lemmas \ref{K} and \ref{lem2.1} and similar arguments as in derivation of \eqref{C27}, we have
\begin{align}\label{8.14}
\sup_{0\leq s\leq t}\|\mathfrak{h}_{2}(s)\|_{L^{\infty}}&\leq C_{\ell}\big[\Big(\frac{1}{2}\Big)^{C_{4}T_{0}^{\frac{5}{4}}}+\frac{k}{N}+C_{\ell,N}mk+\frac{C_{\ell,N}k\sqrt{\lambda}}{m}\big]\sup_{0\leq s\leq t}\|\mathfrak{h}_{2}(s)\|_{L^{\infty}}\nonumber\\
&\quad +C_{\ell}k\sup_{0\leq s\leq t}\|\mathfrak{h}_{1}\|_{L^{\infty}}+C_{\ell,m,N}k\sup_{0\leq s\leq t}\|f_{2}\|_{L^2}.
\end{align}
Then taking $k= C_{3}T_{0}^{\frac{5}{4}}$ and $T_{0}$ suitably large, and then take $N$ large enough, and then $m$ small enough, and finally $\lambda$ small enough, we get from \eqref{8.14} that
\begin{align*}
\sup_{0\leq s\leq t}\|\mathfrak{h}_{2}(s)\|_{L^{\infty}}&\leq C_{\ell}\sup_{0\leq s\leq t}\|\mathfrak{h}_{1}\|_{L^{\infty}}+C_{\ell}\sup_{0\leq s\leq t}\|e^{\lambda_{0}s}f_{2}(s)\|_{L^2}.
\end{align*}
Recall \eqref{def.8.4}. The above estimate together with \eqref{8.9} and $C_{\ell}(\alpha+\tilde{\delta})$ small enough yields \eqref{8.3}. Therefore, the proof of Lemma \ref{lem8.1} is completed. $\hfill\square$

\subsection{{\it{A priori}} $L^2$ estimates}\label{sec7.3}
To close the estimates \eqref{8.2}--\eqref{8.3}, we still need to control the norm $\|e^{\lambda_{0}t}f_{2}(s)\|_{L^2}$. We denote
\begin{equation}\label{f12timew}
(\mathfrak{f}_{1},\mathfrak{f}_{2})=e^{\lambda_{0}t}(f_{1},f_{2}).
\end{equation}
Then, it follows from \eqref{6.1-2}--\eqref{6.1-3} that
\begin{align*}
\left\{
\begin{aligned}
	&\partial_{t}\mathfrak{f}_{1}+v_{\eta}\frac{\partial \mathfrak{f}_{1}}{\partial \eta}-\frac{1}{1-\eta}\Big(v_{\phi}^2\frac{\partial \mathfrak{f}_{1}}{\partial v_{\eta}}-v_{\eta}v_{\phi}\frac{\partial \mathfrak{f}_{1}}{\partial v_{\phi}}\Big)-\frac{\alpha}{\eta_{1}}v_{\eta}\frac{\partial \mathfrak{f}_{1}}{\partial v_{\phi}}+\frac{\alpha}{2\eta_{1}}v_{\eta}v_{\phi}\sqrt{\mu}\mathfrak{f}_{2}+\nu \mathfrak{f}_{1}-\lambda_{0}\mathfrak{f}_{1}\\
	&\qquad =(1-\chi_{M})\mathcal{K}\mathfrak{f}_{1}+e^{\lambda_{0}t}H(f,f),\\
	&\mathfrak{f}_{1}(t,0,v_{\eta},v_{\phi},v_{z})\vert_{v_{\eta}>0}=\mathfrak{f}_{1}(t,\eta_{1},v_{\eta},v_{\phi},v_{z})\vert_{v_{\eta}<0}=0,\\
	&\mathfrak{f}_{1}(0,\eta,v_{\eta},v_{\phi},v_{z})=G_{0}(\eta,v_{\eta},v_{\phi},v_{z})-G_{s}(\eta,v_{\eta},v_{\phi},v_{z}),
\end{aligned}
\right.
\end{align*}
and
\begin{align}\label{8.17}
	\left\{
	\begin{aligned}
		&\partial_{t}\mathfrak{f}_{2}+v_{\eta}\frac{\partial \mathfrak{f}_{2}}{\partial \eta}-\frac{1}{1-\eta}\Big(v_{\phi}^2\frac{\partial \mathfrak{f}_{2}}{\partial v_{\eta}}-v_{\eta}v_{\phi}\frac{\partial \mathfrak{f}_{2}}{\partial v_{\phi}}\Big)-\frac{\alpha}{\eta_{1}}v_{\eta}\frac{\partial \mathfrak{f}_{2}}{\partial v_{\phi}}+\nu \mathfrak{f}_{2}-\lambda_{0}\mathfrak{f}_{2} =\chi_{M}\mu^{-\frac{1}{2}}\mathcal{K}\mathfrak{f}_{1}+K\mathfrak{f}_{2},\\
		&\mathfrak{f}_{2}(t,0,v_{\eta},v_{\phi},v_{z})\vert_{v_{\eta}>0}=\sqrt{2\pi\mu}\int_{u_{\eta}<0}(\mathfrak{f}_{1}+\sqrt{\mu}\mathfrak{f}_{2})(t,0,u)|u_{\eta}|\,{\rm d}u\\
		&\mathfrak{f}_{2}(t,\eta_{1},v_{\eta},v_{\phi},v_{z})\vert_{v_{\eta}<0}=\sqrt{2\pi\mu}\int_{u_{\eta}>0}(\mathfrak{f}_{1}+\sqrt{\mu}\mathfrak{f}_{2})(t,\eta_{1},u)|u_{\eta}|\,{\rm d}u,\\
		&\mathfrak{f}_{2}(0,\eta,v_{\eta},v_{\phi},v_{z})=0.
	\end{aligned}
	\right.
\end{align}
Multiplying \eqref{8.17} by $(1-\eta)\mathfrak{f}_{2}$ and integrating the resultant equation over $\bar{\mathcal{I}}\times \R^3$, we have
\begin{align}\label{8.18}
&\frac{d}{dt}\int_{0}^{\eta_{1}}(1-\eta)\int_{\R^3}|\mathfrak{f}_{2}|^2\,{\rm d}v{\rm d}\eta+\frac{1}{2}(1-\eta_{1})\int_{\R^3}v_{\eta}|\mathfrak{f}_{2}(\eta_{1})|^2\,{\rm d}v-\frac{1}{2}\int_{\R^3}v_{\eta}|\mathfrak{f}_{2}(0)|^2\,{\rm d}v\nonumber\\
&\quad +\int_{0}^{\eta_{1}}(1-\eta)\int_{\R^3}|(\mathbf{I-P})\mathfrak{f}_{2}|^2\nu(v)\,{\rm d}v{\rm d}\eta\nonumber\\
&\leq \lambda_{0}\int_{0}^{\eta_{1}}(1-\eta)\int_{\R^3}|\mathfrak{f}_{2}|^2\,{\rm d}v{\rm d}\eta+\int_{0}^{\eta_{1}}(1-\eta)\int_{\R^3}\mathfrak{f}_{2}\chi_{M}\mu^{-\frac{1}{2}}\mathcal{K}\mathfrak{f}_{1}\,{\rm d}v{\rm d}\eta.
\end{align}
Using the H\"{o}lder inequality and Lemma \ref{K}, one has
\begin{align}\label{8.18-1}
\int_{0}^{\eta_{1}}(1-\eta)\int_{\R^3}\mathfrak{f}_{2}\chi_{M}\mu^{-\frac{1}{2}}\mathcal{K}\mathfrak{f}_{1}\,{\rm d}v{\rm d}\eta\leq \delta\|\mathfrak{f}_{2}\|_{L^2}^2+C_{\delta}\|w^{\ell}\mathfrak{f}_{1}\|_{L^{\infty}}^2.
\end{align}
Using the boundary conditions $\eqref{8.17}_{2}$--$\eqref{8.17}_{3}$, one has
\begin{align}\label{8.19}
\int_{\R^3}v_{\eta}|\mathfrak{f}_{2}(\eta_{1})|^2\,{\rm d}v&=\int_{v_{\eta}>0}v_{\eta}|\mathfrak{f}_{2}(\eta_{1})|^2\,{\rm d}v+\int_{v_{\eta}<0}v_{\eta}|\mathfrak{f}_{2}(\eta_{1})|^2\,{\rm d}v\nonumber\\
&\leq \int_{v_{\eta}>0}v_{\eta}|\mathfrak{f}_{2}(\eta_{1})|^2\,{\rm d}v-\sqrt{2\pi}\Big(\int_{v_{\eta}>0}|v_{\eta}|\sqrt{\mu}\mathfrak{f}_{2}(\eta_{1})\,{\rm d}v\Big)^2-\sqrt{2\pi}\Big(\int_{v_{\eta}>0}|v_{\eta}|\mathfrak{f}_{1}(\eta_{1})\,{\rm d}v\Big)^2\nonumber\\
&\quad -2\sqrt{2\pi}\Big(\int_{v_{\eta}>0}|v_{\eta}|\mathfrak{f}_{1}(\eta_{1})\,{\rm d}v\Big)\Big(\int_{v_{\eta}>0}|v_{\eta}|\sqrt{\mu}\mathfrak{f}_{2}(\eta_{1})\,{\rm d}v\Big)\nonumber\\
&\geq |(I-P_{\gamma})\mathfrak{f}_{2}(\eta_{1})|_{L^2(\gamma_{+})}^2-\delta |P_{\gamma}\mathfrak{f}_{2}(\eta_{1})|_{L^2(\gamma_{+})}^2-C_{\delta}\|w^{\ell}\mathfrak{f}_{1}\|_{L^{\infty}}^2.
\end{align}
Similarly,
\begin{align}\label{8.20}
-\int_{\R^3}v_{\eta}|\mathfrak{f}_{2}(0)|^2\,{\rm d}v\geq |(I-P_{\gamma})\mathfrak{f}_{2}(0)|_{L^2(\gamma_{+})}^2-\delta |P_{\gamma}\mathfrak{f}_{2}(0)|_{L^2(\gamma_{+})}^2-C_{\delta}\|w^{\ell}\mathfrak{f}_{1}\|_{L^{\infty}}^2.
\end{align}
Substituting \eqref{8.18-1}, \eqref{8.19} and \eqref{8.20} into \eqref{8.18}, we get
\begin{align}\label{8.21}
&\frac{d}{dt}\int_{0}^{\eta_{1}}(1-\eta)\int_{\R^3}|\mathfrak{f}_{2}|^2{\rm d}v{\rm d}\eta+|(I-P_{\gamma})\mathfrak{f}_{2}|_{L^2(\gamma_{+})}^2 +\|(\mathbf{I-P})\mathfrak{f}_{2}\|_{L_{\nu}^{2}}^2\nonumber\\
&\leq C(\lambda+\delta)\|\mathfrak{f}_{2}\|_{L^2}^2+C\delta|P_{\gamma}\mathfrak{f}_{2}|_{L^{2}(\gamma_{+})}^2+C_{\delta}\|w^{\ell}\mathfrak{f}_{1}\|_{L^{\infty}}^2.
\end{align}

We still need to control $\|\mathbf{P}\mathfrak{f}_{2}\|_{L^2}$. Since the conservation of mass fails to hold for the equation of $\mathfrak{f}_{2}$, we need to first consider the equation of $\mathfrak{f}:=\mu^{-\frac{1}{2}}\mathfrak{f}_{1}+\mathfrak{f}_{2}$, then obtain the estimates of $\|\mathbf{P}\mathfrak{f}\|_{L^2}$, and finally control $\|\mathbf{P}\mathfrak{f}_{2}\|_{L^2}$ by $\|\mathbf{P}\mathfrak{f}\|_{L^2}$ and $\|w^{\ell}\mathfrak{f}_{1}\|_{L^{\infty}}$ together. In fact, for $\mathfrak{f}=\mu^{-\frac{1}{2}}\mathfrak{f}_{1}+\mathfrak{f}_{2}$, it holds that
\begin{align}\label{8.22}
	\left\{
	\begin{aligned}
		&\partial_{t}\mathfrak{f}+v_{\eta}\frac{\partial \mathfrak{f}}{\partial \eta}-\frac{1}{1-\eta}\Big(v_{\phi}^2\frac{\partial \mathfrak{f}}{\partial v_{\eta}}-v_{\eta}v_{\phi}\frac{\partial \mathfrak{f}}{\partial v_{\phi}}\Big)-\frac{\alpha}{\eta_{1}}v_{\eta}\frac{\partial \mathfrak{f}}{\partial v_{\phi}}+\mathbf{L} \mathfrak{f}-\lambda_{0}\mathfrak{f}=\mu^{-\frac{1}{2}}e^{\lambda_{0}t}H(f,f),\\
		&\mathfrak{f}(t,0,v_{\eta},v_{\phi},v_{z})\vert_{v_{\eta}>0}=\sqrt{2\pi\mu}\int_{u_{\eta}<0}\sqrt{\mu}\mathfrak{f}(t,0,u)|u_{\eta}|\,{\rm d}u,\\
		&\mathfrak{f}(t,\eta_{1},v_{\eta},v_{\phi},v_{z})\vert_{v_{\eta}<0}=\sqrt{2\pi\mu}\int_{u_{\eta}>0}\sqrt{\mu}\mathfrak{f}(t,\eta_{1},u)|u_{\eta}|\,{\rm d}u,\\
		&\mathfrak{f}(0,\eta,v_{\eta},v_{\phi},v_{z})=\mu^{-\frac{1}{2}}\big(G_{0}(\eta,v_{\eta},v_{\phi},v_{z})-G_{s}(\eta,v_{\eta},v_{\phi},v_{z})\big).
	\end{aligned}
	\right.
\end{align}
We denote
\begin{align*}
    \mathbf{P}\mathfrak{f}=\{\mathfrak{a}(t,\eta)+\mathfrak{b}_{1}(t,\eta)v_{\eta}+\mathfrak{b}_{2}(t,\eta)v_{\phi}+\mathfrak{b}_{3}(t,\eta)v_{z}+\mathfrak{c}(t,\eta)(|v|^2-3)\}\sqrt{\mu}.
\end{align*}
Multiplying \eqref{8.22} by $(1-\eta)\sqrt{\mu}$ and integrating the resultant equation over $[0,\eta]\times \R^3$, one has
\begin{align}\label{8.23}
	\frac{d}{dt}\int_{0}^{\eta}(1-y)e^{-\lambda t}\mathfrak{a}(t,y)\,{\rm d}y+(1-\eta)e^{-\lambda t}\mathfrak{b}_{1}(t,\eta)-e^{-\lambda t}b_{1}(t,0)=0.	
\end{align}
Noting the boundary conditions $\eqref{8.22}_{2}$--$\eqref{8.22}_{3}$, one has
	\begin{align}\label{8.24}
		&\mathfrak{b}_{1}(t,\eta_{1})=\int_{\R^3}v_{\eta}\sqrt{\mu}\mathfrak{f}(t,\eta_{1},v)\,{\rm d}v=0,\quad \mathfrak{b}_{1}(t,0)=\int_{\R^3}v_{\eta}\sqrt{\mu}\mathfrak{f}(t,0,v)\,{\rm d}v=0.
	\end{align}
It follows from \eqref{8.23}--\eqref{8.24} that
	\begin{align}\label{8.25}
		e^{-\lambda t}\int_{0}^{\eta_{1}}(1-\eta)\mathfrak{a}(t,\eta)\,{\rm d}\eta=\int_{0}^{\eta_{1}}(1-\eta)\mathfrak{a}(0,\eta)\,{\rm d}\eta=0\quad \text{for any }t\geq 0,
	\end{align}
where from \eqref{M1}, we have used the fact that
	$$
\int_{0}^{\eta_{1}}\int_{\R^3}(1-\eta)[G_{0}(\eta,v)-G_{st}(\eta,v)]\,{\rm d}v{\rm d}\eta=0,
	$$
    or equivalently
$$
\int_{0}^{\eta_{1}}(1-\eta)\mathfrak{a}(0,\eta)\,{\rm d}\eta=0.
$$

As in Section 4, we shall prove the following result to capture the macroscopic dissipation of $\mathfrak{f}$.
\begin{lemma}\label{lem4.2}
Under the assumption \eqref{8.1}, there exists an instant functional $\mathcal{E}_{\text{\rm int}}(t)$ satisfying 
\begin{equation}\label{lem4.2.l1}
|\mathcal{E}_{\text{\rm int}}(t)|\lesssim \|\mathfrak{f}_{2}\|_{L^2}^2+\|w^{\ell}\mathfrak{f}_{1}\|_{L^{\infty}}^2,
\end{equation}
and a small positive constant $\kappa$ such that for any $t\geq 0$,
\begin{align}
\frac{d}{dt}\mathcal{E}_{\text{\rm int}}(t)+\kappa\|(\mathfrak{a},\mathfrak{b}_{1},\mathfrak{b}_{2},\mathfrak{b}_{3},\mathfrak{c})\|_{L^2}^2&\leq C\|(\mathbf{I-P})\mathfrak{f}_{2}\|_{L_{\nu}^2}^2+C|(I-P_{\gamma})\mathfrak{f}_{2}|_{L^2(\gamma_{+})}^2\nonumber\\
&\quad +C(\alpha+\tilde{\varepsilon})\|w^{\ell}\mathfrak{f}_{2}\|_{L^{\infty}}^2+C\|w^{\ell}\mathfrak{f}_{1}\|_{L^{\infty}}^2.
\label{lem4.2.l2}
\end{align}
\end{lemma}

\noindent\textbf{Proof}. Let $\Psi=\Psi(t,\eta,v)\in C^{\infty}([0,\infty)\times [0,\eta_{1}]\times \R^3)$ be a test function, then we take the inner product of \eqref{8.22} and $\Psi$ to see
	\begin{align}\label{8.26}
		&\frac{d}{dt}\la \mathfrak{f},\Psi\ra-\la \mathfrak{f},\partial_{t} \Psi\ra-\la v_{\eta}\mathfrak{f},\partial_{\eta}\Psi\ra+\int_{\R^3}v_{\eta}(\Psi \mathfrak{f})(\eta_{1},v)\,{\rm d}v-\int_{\R^3}v_{\eta}(\Psi \mathfrak{f})(0,v)\,{\rm d}v\nonumber\\
		&\quad +\la \frac{1}{1-\eta}v_{\phi}^2\mathfrak{f}, \partial_{v_{\eta}}\Psi\ra-\la \frac{1}{1-\eta}v_{\eta}\mathfrak{f},\Psi\ra-\la \frac{1}{1-\eta}v_{\eta}v_{\phi}\mathfrak{f},\partial_{v_{\phi}}\Psi\ra\nonumber\\
		&\quad +\frac{\alpha}{\eta_{1}}\la v_{\eta}\mathfrak{f},\partial_{v_{\phi}}\Psi\ra+\frac{\alpha}{2\eta_{1}}\la v_{\eta}v_{\phi}\mathfrak{f},\Psi\ra +\la (\mathbf{L}-\lambda)\mathfrak{f},\Psi\ra =\la \mathscr{H},\Psi\ra,
	\end{align}
where $\mathscr{H}=\mu^{-\frac{1}{2}}e^{\lambda_{0}t}H(f,f)$.

\vspace{0.2cm}
Now, we shall construct appropriate test functions for \eqref{8.26} in order to establish the macroscopic dissipation of $\mathfrak{f}$.

\smallskip

\noindent\underline{\textbf{Estimates of $\|a\|_{L^2}$}}: Let 
	$$
	\Psi_{a}=v_{\eta}(|v|^2-10)\sqrt{\mu}\int_{0}^{\eta}(1-y)\mathfrak{a}(t,y)\,{\rm d}y=:v_{\eta}(|v|^2-10)\sqrt{\mu}\phi_{a}(t,\eta).
	$$
	Then it follows from \eqref{8.25} that
	\begin{align}\label{8.27}
		\phi_{a}(t,\eta_{1})=\phi_{a}(t,0)=0,\quad \|\phi_{a}\|_{H^1}\leq C\|a\|_{L^2}.
	\end{align}
	Taking $\Psi=\Psi_{a}$ in \eqref{8.26} and using \eqref{8.27}, one has
	\begin{align}\label{8.28}
		\frac{d}{dt}\la \mathfrak{f},\Psi_{a}\ra -\la v_{\eta}\mathfrak{f},\partial_{\eta}\Psi_{a}\ra&=\la \mathfrak{f},\partial_{t}\Psi_{a}\ra-\la \frac{1}{1-\eta}v_{\phi}^2\mathfrak{f}, \partial_{v_{\eta}}\Psi_{a}\ra+\la \frac{1}{1-\eta}v_{\eta}\mathfrak{f},\Psi_{a}\ra+\la \frac{1}{1-\eta}v_{\eta}v_{\phi}\mathfrak{f},\partial_{v_{\phi}}\Psi_{a}\ra\nonumber\\
		&\quad -\frac{\alpha}{\eta_{1}}\la v_{\eta}\mathfrak{f},\partial_{v_{\phi}}\Psi_{a}\ra-\frac{\alpha}{2\eta_{1}}\la v_{\eta}v_{\phi}\mathfrak{f},\Psi_{a}\ra -\la (\mathbf{L}-\lambda)\mathfrak{f},\Psi_{a}\ra+\la \mathscr{H},\Psi_{a}\ra\nonumber\\
		& =:\sum\limits_{k=1}^{8}I_{j}.
	\end{align}
	A direct calculation shows that
	\begin{align*}
		-\la v_{\eta}\mathfrak{f},\partial_{\eta}\Psi_{a}\ra\geq [4(1-\eta_{1})-\kappa_{1}]\|\mathfrak{a}\|_{L^2}^2-C_{\kappa_{1}}\|(\mathbf{I-P})\mathfrak{f}\|_{L_{\nu}^2}^2.
	\end{align*}	
For $I_{1}$, we have 
\begin{align}\label{8.29}
		I_{1}=\la \mathfrak{f},\partial_{t}\Psi_{a}\ra&=-\int_{0}^{\eta_{1}}\partial_{t}\phi_{a}(t,\eta)\int_{\R^3}v_{\eta}(|v|^2-10)\mathfrak{f}\sqrt{\mu}\,{\rm d}v{\rm d}\eta\nonumber\\
		&=-5\int_{0}^{\eta_{1}}\mathfrak{b}_{1}(t,\eta)\Big(\int_{0}^{\eta}(1-y)\partial_{t}\mathfrak{a}\,{\rm d}y\Big){\rm d}\eta\nonumber\\
		&\quad +\int_{0}^{\eta_{1}}\Big(\int_{0}^{\eta}(1-y)\partial_{t}\mathfrak{a}\,{\rm d}y\Big)\int_{\R^3}v_{\eta}(|v|^2-10)\sqrt{\mu}(\mathbf{I-P})\mathfrak{f}\,{\rm d}v{\rm d}\eta.
	\end{align}
It follows from \eqref{8.23} that
	$$
	\int_{0}^{\eta}(1-y)\partial_{t}\mathfrak{a}\,{\rm d}y=\lambda\int_{0}^{\eta}(1-y)\mathfrak{a}\,{\rm d}y-(1-\eta)\mathfrak{b}_{1}(t,\eta),
	$$
	which, together with \eqref{8.29}, implies that
	\begin{align*}
		\vert\la \mathfrak{f},\partial_{t}\Psi\ra\vert\leq C\|\mathfrak{b}_{1}\|_{L^2}^2+C\lambda\|\mathfrak{a}\|_{L^2}^2+C\|(\mathbf{I-P})\mathfrak{f}\|_{L_{\nu}^2}^2.
	\end{align*}
    For $I_{2}$ to $I_{6}$, similar to \eqref{E25}--\eqref{E26}, one has
\begin{align*}
I_{2}+I_{3}+I_{4}
&\leq \kappa_{1}\|\mathfrak{a}\|_{L^2}^2+C_{\kappa_{1}}\|(\mathbf{I-P})\mathfrak{f}\|_{L_{\nu}^2}^2.
\\
I_{5}+I_{6}
		&\leq C\alpha\|\mathfrak{a}\|_{L^2}\|\mathfrak{b}_{2}\|_{L^2}+\kappa_{1}\|\mathfrak{a}\|_{L^2}+C_{\kappa_{1}}\|(\mathbf{I-P})\mathfrak{f}\|_{L_{\nu}^2}^2.
	\end{align*}
For $I_{7}$ and $I_{8}$, using \eqref{8.1}, we have
\begin{align*}
I_{7}+I_{8}	&\leq C(\alpha+\kappa_{1}+\tilde{\v}+\lambda)\|\mathfrak{a}\|_{L^2}^2+C_{\kappa_{1}}\|(\mathbf{I-P})\mathfrak{f}\|_{L_{\nu}^2}^2+C(\alpha+\tilde{\v})\|w^{\ell}(\mathfrak{f}_{1},\mathfrak{f}_{2})\|_{L^{\infty}}^2.
\end{align*}
Plugging all the above estimates into \eqref{8.28}, 
we obtain
	\begin{align}\label{8.34}
		\frac{d}{dt}\la \mathfrak{f},\Psi_{a}\ra+\|\mathfrak{a}\|_{L^2}^2&\lesssim (\alpha+\tilde{\v}+\lambda+\kappa_{1})\|\mathfrak{a}\|_{L^2}^2+\alpha\|\mathfrak{b}_{2}\|_{L^2}\|\mathfrak{a}\|_{L^2}+\|\mathfrak{b}_{1}\|_{L^2}^2\nonumber\\
		&\quad +C_{\kappa_{1}}\|(\mathbf{I-P})\mathfrak{f}\|_{L_{\nu}^2}^2+(\alpha+\tilde{\v})\|w^{\ell}(\mathfrak{f}_{1},\mathfrak{f}_{2})\|_{L^{\infty}}^2.
	\end{align}
	
\smallskip

\noindent\underline{\textbf{Estimates of $\|b_{1}\|_{L^2}$}}: Let
	$$
	\Psi_{b_{1}}=-v_{\eta}^2(|v|^2-5)\sqrt{\mu}\,\partial_{\eta}\phi_{b_{1}}(t,\eta),
	$$
	where $\phi_{b_{1}}(t,\eta)$ satisfies
	\begin{align}\label{8.35}
		\partial_{\eta\eta}^2\phi_{b_{1}}(t,\eta)+\frac{1}{3(1-\eta)}\partial_{\eta}\phi_{b_{1}}(t,\eta)=\mathfrak{b}_{1}(t,\eta),\quad \phi_{b_{1}}(t,0)=\phi_{b_{1}}(t,\eta_{1})=0.
	\end{align}
	A straightforward calculation shows that
	\begin{align*}
		\phi_{b_{1}}(t,\eta)=\int_{0}^{\eta}(1-y)^{\frac{1}{3}}\int_{0}^{y}(1-z)^{-\frac{1}{3}}\mathfrak{b}_{1}(t,z)\,{\rm d}z{\rm d}y+C_{1}(t)\int_{0}^{\eta}(1-y)^{\frac{1}{3}}\,{\rm d}y,
	\end{align*}
	with 
	$$
	C_{1}(t)=-\frac{\int_{0}^{\eta_{1}}(1-y)^{\frac{1}{3}}\int_{0}^{y}(1-z)^{-\frac{1}{3}}\mathfrak{b}_{1}(t,z)\,{\rm d}z{\rm d}y}{\int_{0}^{\eta_{1}}(1-y)^{\frac{1}{3}}\,{\rm d}y}.
	$$
Furthermore, one has
	$$
	\|\phi_{b_{1}}(t,\cdot)\|_{H^2}\lesssim \|\mathfrak{b}_{1}(t,\cdot)\|_{L^2}.
	$$
	Taking $\Psi=\Psi_{b_{1}}$ in \eqref{8.26}, one has
	\begin{align}\label{8.37}
		\frac{d}{dt}\la \mathfrak{f},\Psi_{b_{1}}\ra-\la v_{\eta}\mathfrak{f},\partial_{\eta}\Psi_{b_{1}}\ra&=\la \mathfrak{f},\partial_{t}\Psi_{b_{1}}\ra-\int_{\R^3}v_{\eta}(\Psi_{b_{1}} \mathfrak{f})(\eta_{1},v)\,{\rm d}v+\int_{\R^3}v_{\eta}(\Psi_{b_{1}} \mathfrak{f})(0,v)\,{\rm d}v\nonumber\\
		&\quad -\la \frac{1}{1-\eta}v_{\phi}^2\mathfrak{f}, \partial_{v_{\eta}}\Psi_{b_{1}}\ra+\la \frac{1}{1-\eta}v_{\eta}\mathfrak{f},\Psi_{b_{1}}\ra+\la \frac{1}{1-\eta}v_{\eta}v_{\phi}\mathfrak{f},\partial_{v_{\phi}}\Psi_{b_{1}}\ra\nonumber\\
		&\quad -\frac{\alpha}{\eta_{1}}\la v_{\eta}\mathfrak{f},\partial_{v_{\phi}}\Psi_{b_{1}}\ra-\frac{\alpha}{2\eta_{1}}\la v_{\eta}v_{\phi}\mathfrak{f},\Psi_{b_{1}}\ra-\la (\mathbf{L}-\lambda)\mathfrak{f},\Psi_{b_{1}}\ra +\la \mathscr{H},\Psi_{b_{1}}\ra\nonumber\\
	&=:\sum\limits_{k=1}^{10}J_{k}.
	\end{align}
For $J_{1}$, a direct calculation shows that
	\begin{align}\label{8.38}
		J_{1}&=-\int_{0}^{\eta_{1}}\partial_{t\eta}\phi_{b_{1}}(t,\eta)\Big(\int_{\R^3}v_{\eta}^2(|v|^2-5)\mathfrak{f}\sqrt{\mu}\,{\rm d}v\Big){\rm d}\eta\nonumber\\
		&=10\int_{0}^{\eta_{1}}\mathfrak{c}(t,\eta)\partial_{t\eta}\phi_{b_{1}}(t,\eta){\rm d}\eta -\int_{0}^{\eta_{1}}\partial_{t\eta}\phi_{b_{1}}(t,\eta)\Big(\int_{\R^3}v_{\eta}^2(|v|^2-5)\sqrt{\mu}(\mathbf{I-P})\mathfrak{f}\,{\rm d}v\Big){\rm d}\eta.
	\end{align}
We need to control $\|\partial_{t\eta}\phi_{b_{1}}(t,\cdot)\|_{L^2}$. Motivated by \cite{EGKM-13}, we resort to the weak formulation of the equation of $\partial_{t}\mathfrak{b}_{1}$. Letting $\Psi=\psi(\eta)v_{\eta}\sqrt{\mu}$ in \eqref{8.26}, we have
	\begin{align*}
		&\int_{0}^{\eta_{1}}(\mathfrak{b}_{1}(t+\varepsilon)-\mathfrak{b}_{1}(t))\psi(\eta)\,{\rm d}\eta+\int_{t}^{t+\varepsilon}\int_{\R^3}\psi(\eta_{1})v_{\eta}^2\sqrt{\mu}\mathfrak{f}(s,\eta_{1},v)\,{\rm d}v{\rm d}s\nonumber\\
		&\quad-\int_{t}^{t+\varepsilon}\int_{\R^3}\psi(0)v_{\eta}^2\sqrt{\mu}\mathfrak{f}(s,0,v)\,{\rm d}v{\rm d}s-\int_{t}^{t+\varepsilon}\int_{0}^{\eta_{1}}\int_{\R^3}\psi'(\eta)v_{\eta}^2\sqrt{\mu}\mathfrak{f}(s,\eta,v)\,{\rm d}v{\rm d}\eta{\rm d}s\nonumber\\
		&\quad-\int_{t}^{t+\varepsilon}\int_{0}^{\eta_{1}}\int_{\R^3}\frac{\psi(\eta)}{1-\eta}(v_{\phi}^2-v_{\eta}^2)\sqrt{\mu}(\mathbf{I-P})\mathfrak{f}\,{\rm d}v{\rm d}\eta{\rm d}s\nonumber\\
		&=\lambda \int_{t}^{t+\varepsilon}\int_{0}^{\eta_{1}}\mathfrak{b}_{1}(s,\eta)\psi(\eta)\,{\rm d}\eta{\rm d}s.
	\end{align*}
	Taking the difference quotient, we obtain
	\begin{align}\label{J2}
		&\int_{0}^{\eta_{1}}\partial_{t}\mathfrak{b}_{1}\psi(\eta)\,{\rm d}\eta+\psi(\eta_{1})\int_{\R^3}v_{\eta}^2\sqrt{\mu}\mathfrak{f}(t,\eta_{1},v)\,{\rm d}v\nonumber\\
		&\quad-\psi(0)\int_{\R^3}v_{\eta}^2\sqrt{\mu}\mathfrak{f}(t,0,v)\,{\rm d}v-\int_{0}^{\eta_{1}}\psi'(\eta)(\mathfrak{a}(t,\eta)+2\mathfrak{c}(t,\eta))\,{\rm d}v{\rm d}\eta\nonumber\\
		&\quad-\int_{0}^{\eta_{1}}\int_{\R^3}\psi'(\eta)v_{\eta}^2\sqrt{\mu}(\mathbf{I-P})\mathfrak{f}(t,\eta,v)\,{\rm d}v{\rm d}\eta-\int_{0}^{\eta_{1}}\int_{\R^3}\frac{1}{1-\eta}\psi(\eta)(v_{\phi}^2-v_{\eta}^2)\sqrt{\mu}\mathfrak{f}\,{\rm d}v{\rm d}\eta\nonumber\\
		&=\lambda \int_{0}^{\eta_{1}}\mathfrak{b}_{1}(t,\eta)\psi(\eta)\,{\rm d}\eta.
	\end{align}
	For any fixed $t>0$, we choose $\psi(\eta)=-(1-\eta)^{-\frac{1}{3}}\partial_{t}\phi_{b_{1}}(t,\eta)$, where $\phi_{b_{1}}(t,\eta)$ satisfies \eqref{8.35}. Then $\psi(0)=\psi(\eta_{1})=0$, 
	$$
	\psi'(\eta)=-\frac{1}{3}(1-\eta)^{-\frac{4}{3}}\partial_{t}\phi_{b_{1}}(t,\eta)-(1-\eta)^{-\frac{1}{3}}\partial_{t\eta}\phi_{b_{1}}(t,\eta),
	$$
	and
	$$
	(1-\eta)^{\frac{1}{3}}\partial_{\eta}[(1-\eta)^{-\frac{1}{3}}\partial_{t\eta}\phi_{b_{1}}(t,\eta)]=\partial_{t}\mathfrak{b}_{1}(t,\eta).
	$$
 Using the Poincar\'{e} inequality, one has directly 
 $$
 \|\partial_{t}\phi_{b_{1}}(t,\cdot)\|_{L^2}\leq \|\partial_{t\eta}\phi_{b_{1}}(t,\eta)\|_{L^2},\qquad \|\psi\|_{H^1}\lesssim \|\psi'(\eta)\|_{L^2}\leq \|\partial_{t\eta}\phi_{b_{1}}(t,\eta)\|_{L^2}.
 $$
	Then we obtain from \eqref{J2} and the H\"{o}lder inequality that
	\begin{align*}
		\int_{0}^{\eta_{1}}(1-\eta)^{-\frac{1}{3}}|\partial_{t\eta}\phi_{b_{1}}(t,\eta)|^2\,{\rm d}\eta&\lesssim \delta\|\psi'(\eta)\|_{L^2}^2+C_{\delta}\|(\mathfrak{a},\mathfrak{c})\|_{L^2}^2+C_{\delta}\|(\mathbf{I-P})\mathfrak{f}\|_{L_{\nu}^2}^2+C_{\delta}\lambda\|\mathfrak{b}_{1}\|_{L^2}\nonumber\\
		&\lesssim \delta\|\partial_{t\eta}\phi_{b_{1}}(t,\eta)\|_{L^2}^2+\|+C_{\delta}\big(\|(\mathfrak{a},\mathfrak{c})\|_{L^2}^2+\lambda\|\mathfrak{b}_{1}\|_{L^2}^2+\|(\mathbf{I-P})\mathfrak{f}\|_{L_{\nu}^2}^2\big).
	\end{align*}	
Taking $\delta$ small enough, we obtain 
	\begin{align*}
		\|\partial_{t\eta}\phi_{b_{1}}(t,\eta)\|_{L^2}^2&\lesssim \|(\mathfrak{a},\mathfrak{c})\|_{L^2}^2+\|(\mathbf{I-P})\mathfrak{f}\|_{L_{\nu}^2}^2+\lambda\|\mathfrak{b}_{1}\|_{L^2}^2.
	\end{align*}	
	Then \eqref{8.38} can be bounded as
	\begin{align*}
		J_{1}\leq \kappa_{2}\|a\|_{L^2}^2+C_{\kappa_{2}}\|c\|_{L^2}^2+C_{\kappa_{2}}\|(\mathbf{I-P})\mathfrak{f}\|_{L_{\nu}^2}^2+C\lambda\|b_{1}\|_{L^2}^2.
	\end{align*}	
For the boundary terms $J_{2}$ and $J_{3}$, by using the boundary conditions $\eqref{8.22}_{2}$--$\eqref{8.22}_{3}$, we have
\begin{align*}
	J_{2}+J_{3}	&=\partial_{\eta}\phi_{b_{1}}(t,\eta_{1})\int_{v_{\eta}>0}v_{\eta}^3(|v|^2-5)\sqrt{\mu}(I-P_{\gamma})\mathfrak{f}(t,\eta_{1},v)\,{\rm d}v\nonumber\\
	&\quad -\partial_{\eta}\phi_{b_{1}}(t,0)\int_{v_{\eta}<0}v_{\eta}^3(|v|^2-5)\sqrt{\mu}(I-P_{\gamma})\mathfrak{f}(t,0,v)\,{\rm d}v\nonumber\\
	&\leq \kappa_{2}\|b_{1}\|_{L^2}^2+C_{\kappa_{2}}|(I-P_{\gamma})\mathfrak{f}|_{L^2(\gamma_{+})}^2.
	\end{align*}
For the geometric terms $J_{4}$, $J_{5}$ and $J_{6}$, we need to combine them with $-\la v_{\eta}\mathfrak{f},\partial_{\eta}\Psi_{b_{1}}\ra$ to get
\begin{align*}
		&-\la v_{\eta}\mathfrak{f},\partial_{\eta}\Psi_{b_{1}}\ra-J_{4}-J_{5}-J_{6}\nonumber\\
		&=\int_{0}^{\eta_{1}}
		6{\mathfrak{b}_{1}(t,\eta)[\partial_{\eta\eta}^2\phi_{b_{1}}+\frac{1}{3(1-\eta)}\partial_{\eta}\phi_{b_{1}}]}\,{\rm d}\eta
		+\int_{0}^{\eta_{1}}\partial_{\eta\eta}^2\phi_{b_{1}}(t,\eta)\int_{\R^3}v_{\eta}^3(|v|^2-5)\sqrt{\mu}(\mathbf{I-P})\mathfrak{f}\,{\rm d}v
		\nonumber\\
		&\quad -\int_{0}^{\eta_{1}}\frac{2}{1-\eta}\partial_{\eta}\phi_{b_{1}}\int_{\R^3}v_{\phi}^2v_{\eta}(|v|^2-5)\sqrt{\mu}(\mathbf{I-P})\mathfrak{f}\,{\rm d}v{\rm d}\eta\nonumber\\
		&\quad +\int_{0}^{\eta_{1}}\frac{1}{1-\eta}\partial_{\eta}\phi_{b_{1}}\int_{\R^3}v_{\eta}^3(|v|^2-5)\sqrt{\mu}(\mathbf{I-P})\mathfrak{f}\,{\rm d}v{\rm d}\eta\nonumber\\
		&\geq (6-\kappa_{2})\|\mathfrak{b}_{1}\|_{L^2}^2-C_{\kappa_{2}}\|(\mathbf{I-P})\mathfrak{f}\|_{L_{\nu}^2}^2.
	\end{align*}
For $J_{7}$ and $J_{8}$, similar to \eqref{E36}, we have
	\begin{align*}
J_{7}+J_{8}\leq \kappa_{2}\|\mathfrak{b}_{1}\|_{L^2}^2+C_{\kappa_{2}}\|(\mathbf{I-P})\mathfrak{f}\|_{L_{\nu}^2}^2.
	\end{align*}
For $J_{9}$ and $J_{10}$, it follows from the H\"{o}lder inequality and \eqref{8.1} that
\begin{align*}
J_{9}+J_{10}	&\leq C(\alpha+\kappa_{2}+\tilde{\v}+\lambda)\|\mathfrak{b}_{1}\|_{L^2}^2+C_{\kappa_{2}}\|(\mathbf{I-P})\mathfrak{f}\|_{L_{\nu}^2}^2+C(\alpha+\tilde{\v})\|w^{\ell}(\mathfrak{f}_{1},\mathfrak{f}_{2})\|_{L^{\infty}}^2.
\end{align*}
Plugging all the above estimates into \eqref{8.37}, 
we obtain
	\begin{align}\label{8.44}
		\frac{d}{dt}\la \mathfrak{f},\Psi_{b_{1}}\ra+\|\mathfrak{b}_{1}\|_{L^2}^2&\lesssim (\alpha+\tilde{\v}+\lambda+\kappa_{2})\|\mathfrak{b}_{1}\|_{L^2}^2+\kappa_{2}\|\mathfrak{a}\|_{L^2}^2+C_{\kappa_{2}}\|\mathfrak{c}\|_{L^2}^2+C_{\kappa_{2}}\|(\mathbf{I-P})\mathfrak{f}\|_{L^2}^2\nonumber\\
		&\quad
		+C_{\kappa_{2}}|(I-P_{\gamma})\mathfrak{f}|_{L^2(\gamma_{+})}^2+(\alpha+\tilde{\v})\|w^{\ell}(\mathfrak{f}_{1},\mathfrak{f}_{2})\|_{L^{\infty}}^2.
	\end{align}
	
\smallskip

\noindent\underline{\textbf{Estimates of $\|b_{2}\|_{L^2}$}}: Let
	$$
	\Psi_{b_{2}}=-v_{\eta}v_{\phi}\sqrt{\mu}\partial_{\eta}\phi_{b_{2}}(t,\eta),
	$$
	where $\phi_{b_{2}}(t,\eta)$ satisfies
	\begin{align}\label{8.45}
		\partial_{\eta\eta}^2\phi_{b_{2}}-\frac{1}{1-\eta}\partial_{\eta}\phi_{b_{2}}=\mathfrak{b}_{2}(t,\eta),\quad \phi_{b_{2}}(t,\eta_{1})=\phi_{b_{2}}(t,0)=0.
	\end{align}
	A direct calculation shows that
	\begin{align*}
		\phi_{b_{2}}(t,\eta)=\int_{0}^{\eta}(1-y)^{-1}\int_{0}^{y}(1-z)\mathfrak{b}_{2}(t,z)\,{\rm d}z{\rm d}y+C_{2}(t)\int_{0}^{\eta}(1-y)^{-1}\,{\rm d}y,
	\end{align*}
	where
	$$
	C_{2}(t)=-\frac{\int_{0}^{\eta_{1}}(1-y)^{-1}\int_{0}^{y}(1-z)\mathfrak{b}_{2}(t,z)\,{\rm d}z{\rm d}y}{\int_{0}^{\eta_{1}}(1-y)^{-1}\,{\rm d}y}.
	$$
	Furthermore, one has
	$$
	\|\phi_{b_{2}}(t,\cdot)\|_{H^2}\lesssim \|\mathfrak{b}_{2}(t,\cdot)\|_{L^2}.
	$$
Taking $\Psi=\Psi_{b_{2}}$ in \eqref{8.26}, one has
	\begin{align}\label{8.47}
		\frac{d}{dt}\la \mathfrak{f},\Psi_{b_{2}}\ra-\la v_{\eta}\mathfrak{f},\partial_{\eta}\Psi_{b_{2}}\ra&=\la \mathfrak{f},\partial_{t}\Psi_{b_{2}}\ra-\int_{\R^3}v_{\eta}(\Psi_{b_{2}} \mathfrak{f})(\eta_{1},v)\,{\rm d}v+\int_{\R^3}v_{\eta}(\Psi_{b_{2}} \mathfrak{f})(0,v)\,{\rm d}v\nonumber\\
		&\quad -\la \frac{1}{1-\eta}v_{\phi}^2\mathfrak{f}, \partial_{v_{\eta}}\Psi_{b_{2}}\ra+\la \frac{1}{1-\eta}v_{\eta}\mathfrak{f},\Psi_{b_{2}}\ra+\la \frac{1}{1-\eta}v_{\eta}v_{\phi}\mathfrak{f},\partial_{v_{\phi}}\Psi_{b_{2}}\ra\nonumber\\
		&\quad -\frac{\alpha}{\eta_{1}}\la v_{\eta}\mathfrak{f},\partial_{v_{\phi}}\Psi_{b_{2}}\ra-\frac{\alpha}{2\eta_{1}}\la v_{\eta}v_{\phi}\mathfrak{f},\Psi_{b_{2}}\ra -\la (\mathbf{L}-\lambda)\mathfrak{f},\Psi_{b_{2}}\ra +\la \mathscr{H},\Psi_{b_{2}}\ra\nonumber\\
		&=:\sum\limits_{k=1}^{10}H_{k}.
	\end{align}
For $H_{1}$, a direct calculation shows that
	\begin{align}\label{8.48}
		J_{1}
&=\int_{0}^{\eta_{1}}\partial_{t\eta}\phi_{b_{2}}(t,\eta)\,{\rm d}y\int_{\R^3}v_{\eta}v_{\phi}\sqrt{\mu}(\mathbf{I-P})\mathfrak{f}\,{\rm d}v{\rm d}\eta.
	\end{align}
We need to control $\|\partial_{t\eta}\phi_{b_{2}}\|_{L^2}$. Letting $\Psi=\psi(\eta)v_{\phi}\sqrt{\mu}$ in \eqref{8.26}, we have
	\begin{align*}
		&\int_{0}^{\eta_{1}}(\mathfrak{b}_{2}(t+\varepsilon)-\mathfrak{b}_{2}(t))\psi(\eta)\,{\rm d}\eta+\int_{t}^{t+\varepsilon}\int_{\R^3}\psi(\eta_{1})v_{\eta}v_{\phi}\sqrt{\mu}\mathfrak{f}(s,\eta_{1},v)\,{\rm d}v{\rm d}s\nonumber\\
		&\quad-\int_{t}^{t+\varepsilon}\int_{\R^3}\psi(0)v_{\eta}v_{\phi}\sqrt{\mu}\mathfrak{f}(s,0,v)\,{\rm d}v{\rm d}s-\int_{t}^{t+\varepsilon}\int_{0}^{\eta_{1}}\int_{\R^3}\psi'(\eta)v_{\eta}v_{\phi}\sqrt{\mu}(\mathbf{I-P})\mathfrak{f}(s,\eta,v)\,{\rm d}v{\rm d}\eta{\rm d}s\nonumber\\
		&\quad-\int_{t}^{t+\varepsilon}\int_{0}^{\eta_{1}}\int_{\R^3}\frac{2}{1-\eta}\psi(\eta)v_{\eta}v_{\phi}\sqrt{\mu}\mathfrak{f}\,{\rm d}v{\rm d}\eta{\rm d}s+\alpha\int_{t}^{t+\varepsilon}\int_{0}^{\eta_{1}}\mathfrak{b}_{1}(s,\eta)\psi(\eta)\,{\rm d}\eta{\rm d}s\nonumber\\
		&=\lambda \int_{t}^{t+\varepsilon}\int_{0}^{\eta_{1}}\mathfrak{b}_{2}(s,\eta)\psi(\eta)\,{\rm d}\eta{\rm d}s.
	\end{align*}
	Taking the difference quotient, we obtain
	\begin{align}\label{J5}
		&\int_{0}^{\eta_{1}}\partial_{t}\mathfrak{b}_{2}\psi(\eta)\,{\rm d}\eta+\int_{\R^3}\psi(\eta_{1})v_{\eta}v_{\phi}\sqrt{\mu}\mathfrak{f}(t,\eta_{1},v)\,{\rm d}v\nonumber\\
		&\quad-\int_{\R^3}\psi(0)v_{\eta}v_{\phi}\sqrt{\mu}\mathfrak{f}(t,0,v)\,{\rm d}v{\rm d}s-\int_{0}^{\eta_{1}}\int_{\R^3}\psi'(\eta)v_{\eta}v_{\phi}\sqrt{\mu}(\mathbf{I-P})\mathfrak{f}(t,\eta,v)\,{\rm d}v{\rm d}\eta\nonumber\\
		&\quad-\int_{0}^{\eta_{1}}\int_{\R^3}\frac{2}{1-\eta}\psi(\eta)v_{\eta}v_{\phi}\sqrt{\mu}\mathfrak{f}\,{\rm d}v{\rm d}\eta+\alpha\int_{0}^{\eta_{1}}\mathfrak{b}_{1}(t,\eta)\psi(\eta)\,{\rm d}\eta\nonumber\\
		&=\lambda \int_{t}^{t+\varepsilon}\int_{0}^{\eta_{1}}\mathfrak{b}_{2}(s,\eta)\psi(\eta)\,{\rm d}\eta.
	\end{align}
	For fixed $t>0$, we choose $\psi(\eta)=-(1-\eta)\partial_{t}\phi_{b_{2}}(t,\eta)$, where $\phi_{b_{2}}$ satisfies \eqref{8.45}. Then $\psi(0)=\psi(\eta_{1})=0$,
	$$
	\psi'(\eta)=\partial_{t}\phi_{b_{2}}-(1-\eta)\partial_{t\eta}\phi_{b_{2}}(t,\eta),
	$$
	and
	\begin{align*}
		\int_{0}^{\eta_{1}}\partial_{t}\mathfrak{b}_{2}(t,\eta)\psi(\eta)\,{\rm d}\eta&=-\int_{0}^{\eta_{1}}\partial_{\eta}[(1-\eta)\partial_{t\eta}\phi_{b_{2}}]\partial_{t}\phi_{b_{2}}\,{\rm d}\eta=\int_{0}^{\eta_{1}}(1-\eta)|\partial_{t\eta}\phi_{b_{2}}(t,\eta)|^2\,{\rm d}\eta.
	\end{align*}
	Using the Poincar\'{e} inequality, we have $\|\psi\|_{H^1}\lesssim \|\psi'(\eta)\|_{L^2}\lesssim \|\partial_{t\eta}\phi_{b_{2}}\|_{L^2}$. Then we obtain from \eqref{J5} that
	\begin{align*}
		\|\partial_{t\eta}\phi_{b_{2}}\|_{L^2}^2\leq
		 C\|(\mathbf{I-P})\mathfrak{f}\|_{L_{\nu}^2}^2+C\alpha\|\mathfrak{b}_{1}\|^2+C\lambda\|\mathfrak{b}_{2}\|_{L^2}^2.
	\end{align*}
	Hence, \eqref{8.48} can be bounded as
	\begin{align*}
	H_{1}\leq C\|(\mathbf{I-P})\mathfrak{f}\|_{L_{\nu}^2}^2+C\alpha\|\mathfrak{b}_{1}\|_{L^2}^2+C\lambda\|b_{2}\|_{L^2}^2.
\end{align*}
For the boundary terms $H_{2}$ and $H_{3}$, using the boundary conditions $\eqref{8.22}_{2}$--$\eqref{8.22}_{3}$, we have
\begin{align*}
	H_{2}+H_{3}&=-\partial_{\eta}\phi_{b_{2}}(t,0)\int_{v_{\eta}<0}v_{\eta}^2v_{\phi}\sqrt{\mu}(I-P_{\gamma})\mathfrak{f}(t,0,v)\,{\rm d}v\nonumber\\
	&\quad +\partial_{\eta}\phi_{b_{2}}(t,\eta_{1})\int_{v_{\eta}>0}v_{\eta}^2v_{\phi}\sqrt{\mu}(I-P_{\gamma})\mathfrak{f}(t,\eta_{1},v)\,{\rm d}v\nonumber\\
		&\leq \kappa_{3}\|\mathfrak{b}_{2}\|_{L^2}^2+C_{\kappa_{3}}|(I-P_{\gamma})\mathfrak{f}|_{L^2(\gamma_{+})}^2.
	\end{align*}
For the geometric terms $H_{4}$, $H_{5}$ and $H_{6}$, we need to combine them with $-\la v_{\eta}\mathfrak{f},\partial_{\eta}\Psi_{b_{2}}\ra$ to get
\begin{align*}
-\la v_{\eta}\mathfrak{f},\partial_{\eta}\Psi_{b_{2}}\ra-H_{4}-H_{5}-H_{6}&={\int_{0}^{\eta_{1}}\mathfrak{b}_{2}(t,\eta)[\partial_{\eta\eta}^2\phi_{b_{2}}-\frac{1}{1-\eta}\phi_{b_{2}}(t,\eta)]\,{\rm d}\eta}\nonumber\\
&\quad +\int_{0}^{\eta_{1}}\partial_{\eta\eta}^2\phi_{b_{2}}\int_{\R^3}v_{\eta}^2v_{\phi}\sqrt{\mu}(\mathbf{I-P})\mathfrak{f}\,{\rm d}v{\rm d}\eta\nonumber\\
&\quad +\int_{0}^{\eta_{1}}\frac{1}{1-\eta}\partial_{\eta}\phi_{b_{2}}(t,\eta)\int_{\R^3}(2v_{\eta}^2-v_{\phi}^2)v_{\phi}\sqrt{\mu}(\mathbf{I-P})\mathfrak{f}\,{\rm d}v{\rm d}\eta\nonumber\\
		&\geq (1-\kappa_{3})\|\mathfrak{b}_{2}\|_{L^2}^2-C_{\kappa_{3}}\|(\mathbf{I-P})\mathfrak{f}\|_{L_{\nu}^2}^2.
	\end{align*}
For $H_{7}$ and $H_{8}$, similar to \eqref{E46}, we have
\begin{align*}
		H_{7}+H_{8}
		&\leq C\alpha\|\mathfrak{b}_{2}\|_{L^2}\|(\mathfrak{a},\mathfrak{c})\|_{L^2}^2+\kappa_{3}\|\mathfrak{b}_{2}\|_{L^2}^2+C_{\kappa_{3}}\|(\mathbf{I-P})\mathfrak{f}\|_{L_{\nu}^2}^2.
	\end{align*}
For $H_{9}$ and $H_{10}$, it follows from the H\"{o}lder inequality and \eqref{8.1} that
\begin{align*}
H_{9}+H_{10}&\leq C(\alpha+\kappa_{3}+\tilde{\v}+\lambda)\|\mathfrak{b}_{2}\|_{L^2}^2+C_{\kappa_{3}}\|(\mathbf{I-P})\mathfrak{f}\|_{L_{\nu}^2}^2+C(\alpha+\tilde{\v})\|w^{\ell}(\mathfrak{f}_{1},\mathfrak{f}_{2})\|_{L^{\infty}}^2.
\end{align*}
Plugging all the above estimates to \eqref{8.47}, 
we obtain
	\begin{align}\label{8.53}
		&\frac{d}{dt}\la \mathfrak{f},\Psi_{b_{2}}\ra+\|b_{2}\|_{L^2}^2\lesssim (\alpha+\kappa_{3}+\tilde{\v}+\lambda)\|\mathfrak{b}_{2}\|_{L^2}^2+\alpha\|\mathfrak{b}_{2}\|_{L^2}\|(\mathfrak{a},\mathfrak{c})\|_{L^2}^2+\alpha\|\mathfrak{b}_{1}\|_{L^2}^2\nonumber\\
		&\qquad\qquad\qquad
		+C_{\kappa_{3}}\|(\mathbf{I-P})\mathfrak{f}\|_{L_{\nu}^2}^2 +C_{\kappa_{3}}|(I-P_{\gamma})\mathfrak{f}|_{L^2(\gamma_{+})}^2+(\alpha+\tilde{\v})\|w^{\ell}(\mathfrak{f}_{1},\mathfrak{f}_{2})\|_{L^{\infty}}^2.
	\end{align}

\smallskip

\noindent\underline{\textbf{Estimates of $\|b_{3}\|_{L^2}$}}: Let
$$
\Psi_{b_{3}}=-v_{\eta}v_{z}\sqrt{\mu}\partial_{\eta}\phi_{b_{3}}(t,\eta),
$$
where $\phi_{b_{3}}(t,\eta)$ satisfies
$$
\partial_{\eta\eta}^2\phi_{b_{3}}(t,\eta)=\mathfrak{b}_{3}(t,\eta),\quad \phi_{b_{3}}(t,0)=\phi_{b_{3}}(t,\eta_{1})=0.
$$
A direct calculation shows that
$$
\phi_{b_{3}}(t,\eta)=\int_{0}^{\eta}\int_{0}^{y}\phi_{b_{3}}(t,z)\,{\rm d}z{\rm d}y+C_{3}(t)\eta,
$$
where
$$
C_{3}(t)=-\frac{1}{\eta_{1}}\int_{0}^{\eta_{1}}\int_{0}^{y}\mathfrak{b}_{3}(t,z)\,{\rm d}z{\rm d}y.
$$
Furthermore, we have
$\|\phi_{b_{3}}(t,\cdot)\|_{H^2}\lesssim \|\mathfrak{b}_{3}(t,\cdot)\|_{L^2}$. Taking $\Psi=\Psi_{b_{3}}$ in \eqref{8.26}, one has
	\begin{align}\label{8.55-1}
		\frac{d}{dt}\la \mathfrak{f},\Psi_{b_{3}}\ra-\la v_{\eta}\mathfrak{f},\partial_{\eta}\Psi_{b_{3}}\ra&=\la \mathfrak{f},\partial_{t}\Psi_{b_{2}}\ra-\int_{\R^3}v_{\eta}(\Psi_{b_{3}} \mathfrak{f})(\eta_{1},v)\,{\rm d}v+\int_{\R^3}v_{\eta}(\Psi_{b_{3}} \mathfrak{f})(0,v)\,{\rm d}v\nonumber\\
		&\quad -\la \frac{1}{1-\eta}v_{\phi}^2\mathfrak{f}, \partial_{v_{\eta}}\Psi_{b_{3}}\ra+\la \frac{1}{1-\eta}v_{\eta}\mathfrak{f},\Psi_{b_{3}}\ra+\la \frac{1}{1-\eta}v_{\eta}v_{\phi}\mathfrak{f},\partial_{v_{\phi}}\Psi_{b_{3}}\ra\nonumber\\
		&\quad -\frac{\alpha}{\eta_{1}}\la v_{\eta}\mathfrak{f},\partial_{v_{\phi}}\Psi_{b_{3}}\ra-\frac{\alpha}{2\eta_{1}}\la v_{\eta}v_{\phi}\mathfrak{f},\Psi_{b_{3}}\ra -\la (\mathbf{L}-\lambda)\mathfrak{f},\Psi_{c}\ra +\la \mathscr{H},\Psi_{b_{3}}\ra\nonumber\\
	&\quad=:\sum\limits_{k=1}^{10}Q_{k}.
	\end{align}
A direct calculation shows that 
\begin{align*}
-\langle v_{\eta}\mathfrak{f},\partial_{\eta}\Psi_{b_{3}}\rangle\geq (1-\kappa_{4})\|\mathfrak{b}_{3}\|_{L^2}^2-C_{\kappa_{4}}\|(\mathbf{I-P})\mathfrak{f}\|_{L_{\nu}^2}^2.
\end{align*}
For $Q_{1}$, we have
\begin{align}\label{8.56-1}
Q_{1}=-\int_{0}^{\eta_{1}}\partial_{t\eta}\phi_{b_{3}}\int_{\R^3}v_{\eta}v_{z}\sqrt{\mu}(\mathbf{I-P})\mathfrak{f}\,{\rm d}v{\rm d}\eta.
\end{align}
We need to control $\|\partial_{t\eta}\phi_{b_{3}}\|_{L^2}$. Letting $\Psi=\psi(\eta)v_{z}\sqrt{\mu}$ in \eqref{8.26}, we have
\begin{align*}
&\int_{0}^{\eta_{1}}(\mathfrak{b}_{3}(t+\varepsilon)-\mathfrak{b}_{3}(t))\psi(\eta)\,{\rm d}\eta+\int_{t}^{t+\varepsilon}\int_{\R^3}\psi(\eta_{1})v_{\eta}v_{z}\sqrt{\mu}(\mathbf{I-P})\mathfrak{f}(s,\eta_{1},v)\,{\rm d}v{\rm d}s\nonumber\\
&\quad -\int_{t}^{t+\varepsilon}\int_{\R^3}\psi(0)v_{\eta}v_{z}\sqrt{\mu}(\mathbf{I-P})\mathfrak{f}(s,0,v)\,{\rm d}v{\rm d}s\nonumber\\
&\quad-\int_{t}^{t+\varepsilon}\int_{0}^{\eta_{1}}\int_{\R^3}\psi'(\eta)v_{\eta}v_{z}\sqrt{\mu}(\mathbf{I-P})\mathfrak{f}(s,\eta,v)\,{\rm d}v{\rm d}\eta{\rm d}s\nonumber\\
&\quad +\int_{t}^{t+\varepsilon}\int_{0}^{\eta_{1}}\int_{\R^3}\frac{\psi(\eta)}{1-\eta}v_{\eta}v_{z}\sqrt{\mu}(\mathbf{I-P})\mathfrak{f}(s,\eta,v)\,{\rm d}v{\rm d}\eta{\rm d}s\nonumber\\
&\quad =\lambda \int_{t}^{t+\varepsilon}\int_{0}^{\eta_{1}}\mathfrak{b}_{3}(s,\eta)\psi(\eta)\,{\rm d}\eta{\rm d}s.
\end{align*}
Taking the difference quotient, we obtain
\begin{align}\label{8.58-2}
&\int_{0}^{\eta_{1}}\partial_{t}\mathfrak{b}_{3}\psi(\eta)\,{\rm d}\eta+\psi(\eta_{1})\int_{\R^3}v_{\eta}v_{z}\sqrt{\mu}(\mathbf{I-P})\mathfrak{f}(s,\eta_{1},v)\,{\rm d}v\nonumber\\
&\quad-\psi(0)\int_{\R^3}\int_{\R^3}v_{\eta}v_{z}\sqrt{\mu}(\mathbf{I-P})\mathfrak{f}(s,0,v)\,{\rm d}v -\int_{0}^{\eta_{1}}\int_{\R^3}\psi'(\eta)v_{\eta}v_{z}\sqrt{\mu}(\mathbf{I-P})\mathfrak{f}(t,\eta,v)\,{\rm d}v{\rm d}\eta\nonumber\\
&\quad +\int_{0}^{\eta_{1}}\int_{\R^3}\frac{\psi(\eta)}{1-\eta}v_{\eta}v_{z}\sqrt{\mu}(\mathbf{I-P})\mathfrak{f}(t,\eta,v)\,{\rm d}v{\rm d}\eta=\lambda\int_{0}^{\eta_{1}}\mathfrak{b}_{3}(t,\eta)\psi(\eta)\,{\rm d}\eta{\rm d}s.
\end{align}
For fixed $t>0$, we denote $\psi(\eta)=-\partial_{t}\phi_{b_{3}}(t,\eta))$, where $\phi_{b_{3}}$ satisfies \eqref{8.54}. Then $\psi(0)=\psi(\eta_{1})=0$ and
$$
\int_{0}^{\eta_{1}}\partial_{t}\mathfrak{b}_{3}\psi(\eta)\,{\rm d}\eta=\int_{0}^{\eta_{1}}|\partial_{t\eta}\phi_{b_{3}}|^2\,{\rm d}\eta,\quad \|\psi\|_{H^1}\lesssim \|\partial_{t\eta}\phi_{b_{3}}\|_{L^2}.
$$
Then we obtain from \eqref{8.58-2} that
\begin{align*}
\|\partial_{t\eta}\phi_{b_{3}}\|_{L^2}\lesssim C\|(\mathbf{I-P})\mathfrak{f}\|_{L_{\nu}^2}+C\lambda\|\mathfrak{b}_{3}\|_{L^2}^2.
\end{align*}
Hence, \eqref{8.56-1} can be controlled by
\begin{align*}
    Q_{1}\leq \kappa_{4}\|b_{3}\|_{L^2}^2+C_{\kappa_{4}}\|(\mathbf{I-P})\mathfrak{f}\|_{L_{\nu}^2}.
\end{align*}
For the boundary terms $Q_{2}$ and $Q_{3}$, using the boundary conditions $\eqref{8.22}_{2}$--$\eqref{8.22}_{3}$, we have
\begin{align*}
Q_{2}+Q_{3}&=\partial_{\eta}\phi_{b_{3}}(t,\eta_{1})\int_{v_{\eta}>0}v_{\eta}^2v_{z}\sqrt{\mu}(I-P_{\gamma})\mathfrak{f}(t,\eta_{1},v)\,{\rm d}v\nonumber\\
&\quad -\partial_{\eta}\phi_{b_{3}}(t,0)\int_{v_{\eta}>0}v_{\eta}^2v_{z}\sqrt{\mu}(I-P_{\gamma})\mathfrak{f}(t,0,v)\,{\rm d}v\nonumber\\
&\leq \kappa_{4}\|\mathfrak{b}_{3}\|_{L^2}^2+C_{\kappa_{4}}\|(I-P_{\gamma})\mathfrak{f}\|_{L^2(\gamma_{+})}^2.
\end{align*}
For the geometric terms $Q_{4}$, $Q_{5}$ and $Q_{6}$, similar to \eqref{z8}, we have
\begin{align*}
Q_{4}+Q_{5}+Q_{6}
&\leq \kappa_{4}\|\mathfrak{b}_{3}\|_{L^2}^2+C_{\kappa_{4}}\|(\mathbf{I-P})\mathfrak{f}\|_{L_{\nu}^2}^2.
\end{align*}
For $Q_{7}$ and $Q_{8}$, similar to \eqref{z9}, one has
\begin{align*}
Q_{7}+Q_{8}=0.
\end{align*}
For $Q_{9}$ and $Q_{10}$, it follows from the H\"{o}lder inequality and \eqref{8.1} that
\begin{align*}
Q_{9}+Q_{10}	&\leq C(\alpha+\kappa_{4}+\tilde{\v}+\lambda)\|\mathfrak{b}_{3}\|_{L^2}^2+C_{\kappa_{4}}\|(\mathbf{I-P})\mathfrak{f}\|_{L_{\nu}^2}^2+C(\alpha+\tilde{\v})\|w^{\ell}(\mathfrak{f}_{1},\mathfrak{f}_{2})\|_{L^{\infty}}^2.
\end{align*}
Plugging all the above estimates to \eqref{8.55-1}, 
we obtain
\begin{align}\label{8.58-9}
\frac{d}{dt}\langle \mathfrak{f},\Psi_{b_{3}}\rangle+\|\mathfrak{b}_{3}\|_{L^2}^2&\lesssim (\alpha+\tilde{\v}+\kappa_{4}+\lambda)\|\mathfrak{b}_{3}\|_{L^2}^2+C_{\kappa_{4}}\|(\mathbf{I-P})\mathfrak{f}\|_{L_{\nu}^2}^2\nonumber\\
&\quad +C_{\kappa_{4}}|(I-P_{\gamma})\mathfrak{f}|_{L^2(\gamma_{+})}^2+(\alpha+\tilde{\v})\|w^{\ell}(\mathfrak{f}_{1},\mathfrak{f}_{2})\|_{L^{\infty}}^2.
\end{align}

\noindent\underline{\textbf{Estimates of $\|c\|_{L^2}$}}: Let
	$$
	\Psi_{c}=-v_{\eta}(|v|^2-5)\sqrt{\mu}\partial_{\eta}\phi_{c}(t,\eta),
	$$
	where $\phi_{c}(t,\eta)$ satisfies
	\begin{align}\label{8.54}
	\partial_{\eta\eta}^2\phi_{c}(t,\eta)=\mathfrak{c}(t,\eta),\quad \phi_{c}(t,0)=\phi_{c}(t,\eta_{1})=0.
	\end{align}
A direct calculation shows that
$$
		\phi_{c}(t,\eta)=\int_{0}^{\eta}\int_{0}^{y}\mathfrak{c}(t,z)\,{\rm d}z{\rm d}y+C_{3}(t)\eta,
$$
	where
	$$
	C_{3}(t)=-\frac{1}{\eta_{1}}\int_{0}^{\eta_{1}}\int_{0}^{y}\mathfrak{c}(t,z)\,{\rm d}z{\rm d}y.
	$$
	Taking $\Psi=\Psi_{c}$ in \eqref{8.26}, one has
	\begin{align}\label{8.55}
		\frac{d}{dt}\la \mathfrak{f},\Psi_{c}\ra-\la v_{\eta}\mathfrak{f},\partial_{\eta}\Psi_{c}\ra&=\la \mathfrak{f},\partial_{t}\Psi_{c}\ra-\int_{\R^3}v_{\eta}(\Psi_{c} \mathfrak{f})(\eta_{1},v)\,{\rm d}v+\int_{\R^3}v_{\eta}(\Psi_{c} \mathfrak{f})(0,v)\,{\rm d}v\nonumber\\
		&\quad -\la \frac{1}{1-\eta}v_{\phi}^2\mathfrak{f}, \partial_{v_{\eta}}\Psi_{c}\ra+\la \frac{1}{1-\eta}v_{\eta}\mathfrak{f},\Psi_{c}\ra+\la \frac{1}{1-\eta}v_{\eta}v_{\phi}\mathfrak{f},\partial_{v_{\phi}}\Psi_{c}\ra\nonumber\\
		&\quad -\frac{\alpha}{\eta_{1}}\la v_{\eta}\mathfrak{f},\partial_{v_{\phi}}\Psi_{c}\ra-\frac{\alpha}{2\eta_{1}}\la v_{\eta}v_{\phi}\mathfrak{f},\Psi_{c}\ra -\la (\mathbf{L}-\lambda)\mathfrak{f},\Psi_{c}\ra +\la \mathscr{H},\Psi_{c}\ra\nonumber\\
		&=:\sum\limits_{k=1}^{10}S_{k}.
	\end{align}
A direct calculation shows that
\begin{align*}
		-\la v_{\eta}\mathfrak{f},\partial_{\eta}\Psi_{c}\ra\geq (1-\kappa_{5})\|\mathfrak{c}\|_{L^2}^2-C_{\kappa_{5}}\|(\mathbf{I-P})\mathfrak{f}\|_{L_{\nu}^2}^2.
\end{align*}	
For $S_{1}$, a straightforward calculation shows that
	\begin{align}\label{8.56}
		S_{1}
		&=-\int_{0}^{\eta_{1}}\partial_{t\eta}\phi_{c}(t,\eta)\int_{\R^3}v_{\eta}(|v|^2-5)\sqrt{\mu}(\mathbf{I-P})\mathfrak{f}\,{\rm d}v{\rm d}\eta.
	\end{align}
We need to control $\|\partial_{t\eta}\phi_{c}(t,\eta)\|_{L^2}$. Letting $\Psi=\psi(\eta)(|v|^2-3)\sqrt{\mu}$ in \eqref{8.26}, we have
	\begin{align*}
		&2\int_{0}^{\eta_{1}}(\mathfrak{c}(t+\varepsilon)-\mathfrak{c}(t))\psi(\eta)\,{\rm d}\eta+2\int_{t}^{t+\varepsilon}\psi(\eta_{1})\mathfrak{b}_{1}(s,\eta_{1})\,{\rm d}s-2\int_{t}^{t+\varepsilon}\psi(0)\mathfrak{b}_{1}(s,0)\,{\rm d}s\nonumber\\
		&\quad+\int_{t}^{t+\varepsilon}\int_{\R^3}\psi(\eta_{1})v_{\eta}(|v|^2-3)\sqrt{\mu}(\mathbf{I-P})\mathfrak{f}(s,\eta_{1},v)\,{\rm d}v{\rm d}s\nonumber\\
		&\quad-\int_{t}^{t+\varepsilon}\int_{\R^3}\psi(0)v_{\eta}(|v|^2-3)\sqrt{\mu}(\mathbf{I-P})\mathfrak{f}(s,0,v)\,{\rm d}v{\rm d}s-2\int_{t}^{t+\varepsilon}\int_{0}^{\eta_{1}}\psi'(\eta)\mathfrak{b}_{1}(s,\eta)\,{\rm d}\eta{\rm d}s\nonumber\\
		&\quad-	\int_{t}^{t+\varepsilon}\int_{0}^{\eta_{1}}\int_{\R^3}\psi'(\eta)v_{\eta}(|v|^2-3)\sqrt{\mu}(\mathbf{I-P})\mathfrak{f}(s,\eta,v)\,{\rm d}v{\rm d}\eta{\rm d}s-\int_{t}^{t+\varepsilon}\int_{0}^{\eta_{1}}\frac{2}{1-\eta}\psi(\eta)\mathfrak{b}_{1}(t,\eta)\,{\rm d}\eta{\rm d}s\nonumber\\
		&\quad-\int_{t}^{t+\varepsilon}\int_{0}^{\eta_{1}}\int_{\R^3}\frac{2}{1-\eta}\psi(\eta)v_{\eta}(|v|^2-3)\sqrt{\mu}(\mathbf{I-P})\mathfrak{f}\,{\rm d}v{\rm d}\eta{\rm d}s\nonumber\\
		&\quad+\frac{2\alpha}{\eta_{1}}\int_{t}^{t+\varepsilon}\int_{0}^{\eta_{1}}\int_{\R^3}\psi(\eta)v_{\eta}v_{\phi}\sqrt{\mu}(\mathbf{I-P})\mathfrak{f}(s,\eta,v)\,{\rm d}v{\rm d}\eta{\rm d}s\nonumber\\
		&=\lambda \int_{t}^{t+\varepsilon}\int_{0}^{\eta_{1}}\mathfrak{c}(s,\eta)\psi(\eta)\,{\rm d}\eta{\rm d}s.
	\end{align*}
	Taking the difference quotient, we obtain
	\begin{align}\label{H5-1}
		&2\int_{0}^{\eta_{1}}\partial_{t}\mathfrak{c}\psi(\eta)\,{\rm d}\eta+2\psi(\eta_{1})\mathfrak{b}_{1}(t,\eta_{1})-2\psi(0)\mathfrak{b}_{1}(t,0)\nonumber\\
		&\quad+\int_{\R^3}\psi(\eta_{1})v_{\eta}(|v|^2-3)\sqrt{\mu}(\mathbf{I-P})\mathfrak{f}(t,\eta_{1},v)\,{\rm d}v\nonumber\\
		&\quad-\int_{\R^3}\psi(0)v_{\eta}(|v|^2-3)\sqrt{\mu}(\mathbf{I-P})\mathfrak{f}(t,0,v)\,{\rm d}v-2\int_{0}^{\eta_{1}}\psi'(\eta)\mathfrak{b}_{1}(t,\eta)\,{\rm d}\eta\nonumber\\
		&\quad-	\int_{0}^{\eta_{1}}\int_{\R^3}\psi'(\eta)v_{\eta}(|v|^2-3)\sqrt{\mu}(\mathbf{I-P})\mathfrak{f}(t,\eta,v)\,{\rm d}v{\rm d}\eta-\int_{0}^{\eta_{1}}\frac{2}{1-\eta}\psi(\eta)\mathfrak{b}_{1}(t,\eta)\,{\rm d}\eta\nonumber\\
		&\quad-\int_{0}^{\eta_{1}}\int_{\R^3}\frac{2}{1-\eta}\psi(\eta)v_{\eta}(|v|^2-3)\sqrt{\mu}(\mathbf{I-P})\mathfrak{f}\,{\rm d}v{\rm d}\eta\nonumber\\
		&\quad+\frac{2\alpha}{\eta_{1}}\int_{0}^{\eta_{1}}\int_{\R^3}\psi(\eta)v_{\eta}v_{\phi}\sqrt{\mu}(\mathbf{I-P})\mathfrak{f}(t,\eta,v)\,{\rm d}v{\rm d}\eta\nonumber\\
		&=\lambda \int_{0}^{\eta_{1}}\mathfrak{c}(t,\eta)\psi(\eta)\,{\rm d}\eta.
	\end{align}
	For fixed $t>0$, we choose $\psi(\eta)=-\partial_{t}\phi_{c}(t,\eta)$, where $\phi_{c}$ satisfies \eqref{8.54}. Then $\psi(0)=\psi(\eta_{1})=0$ and
	$$
	\int_{0}^{\eta_{1}}\partial_{t}c\psi(\eta)\,{\rm d}\eta=\|\partial_{t\eta}\phi_{c}(t,\eta)\|_{L^2}^2,\quad \|\psi\|_{H^1}\lesssim \|\partial_{t\eta}\phi_{c}\|_{L^2}.
	$$
	Then we obtain from \eqref{H5-1} that
	\begin{align*}
		\|\partial_{t\eta}\phi_{c}(t,\eta)\|_{L^2}^2\lesssim \|\mathfrak{b}_{1}\|_{L^2}^2+\|(\mathbf{I-P})\mathfrak{f}\|_{L_{\nu}^2}^2+\lambda\|\mathfrak{c}\|_{L^2}^2.
	\end{align*}
Hence, \eqref{8.56} can be controlled by
	\begin{align*}
S_{1}\leq \kappa_{5}\|(\mathfrak{b}_{1},\mathfrak{c})\|_{L^2}^2+C_{\kappa_{5}}\|(\mathbf{I-P})\mathfrak{f}\|_{L_{\nu}^2}^2.
	\end{align*}
For the boundary terms $S_{2}$ and $S_{3}$, using the boundary conditions $\eqref{8.22}_{2}$--$\eqref{8.22}_{3}$, we have
	\begin{align*}
		S_{2}+S_{3}&=\partial_{\eta}\phi_{c}(t,\eta_{1})\int_{v_{\eta}>0}v_{\eta}^2(|v|^2-5)\sqrt{\mu}(I-P_{\gamma})\mathfrak{f}(t,\eta_{1},v)\,{\rm d}v\nonumber\\
		&\quad -\partial_{\eta}\phi_{c}(t,0)\int_{v_{\eta}<0}v_{\eta}^2(|v|^2-5)\sqrt{\mu}(I-P_{\gamma})\mathfrak{f}(t,0,v)\,{\rm d}v\nonumber\\
		&\leq \kappa_{5}\|\mathfrak{c}\|_{L^2}^2+C_{\kappa_{5}}|(I-P_{\gamma})\mathfrak{f}|_{L^2(\gamma_{+})}^2.
	\end{align*}
For the geometric terms $S_{4}$, $S_{5}$ and $S_{6}$, similar to \eqref{E55}, we have
	\begin{align*}
S_{4}+S_{5}+S_{6}
		&\leq \kappa_{5}\|\mathfrak{c}\|_{L^2}^2+C_{\kappa_{5}}\|(\mathbf{I-P})\mathfrak{f}\|_{L_{\nu}^2}^2.
	\end{align*}
For $S_{7}$ and $S_{8}$, similar to \eqref{E56}, we have
\begin{align*}
S_{7}+S_{8}
&\leq C\alpha\|\mathfrak{c}\|_{L^2}\|\mathfrak{b}_{2}\|_{L^2}+\kappa_{5}\|\mathfrak{c}\|_{L^2}^2+C_{\kappa_{5}}\|(\mathbf{I-P})\mathfrak{f}\|_{L_{\nu}^2}^2.
\end{align*}
For $S_{9}$ and $S_{10}$, it follows from the H\"{o}lder inequality and \eqref{8.1} that
\begin{align*}
S_{9}+S_{10}	&\leq C(\alpha+\kappa_{5}+\tilde{\v}+\lambda)\|\mathfrak{c}\|_{L^2}^2+C_{\kappa_{5}}\|(\mathbf{I-P})\mathfrak{f}\|_{L_{\nu}^2}^2+C(\alpha+\tilde{\v})\|w^{\ell}(\mathfrak{f}_{1},\mathfrak{f}_{2})\|_{L^{\infty}}^2.
\end{align*}
Applying the above estimates to \eqref{8.55}, we obtain
\begin{align}\label{8.62}
	\frac{d}{dt}\la \mathfrak{f},\Psi_{c}\ra+\|\mathfrak{c}\|_{L^2}^2&\lesssim (\alpha+\tilde{\v}+\kappa_{5}+\lambda)\|\mathfrak{c}\|_{L^2}^2+\alpha\|\mathfrak{c}\|_{L^2}\|\mathfrak{b}_{2}\|_{L^2}+\kappa_{5}\|\mathfrak{b}_{1}\|_{L^2}^2\nonumber\\
		&\quad
		+C_{\kappa_{5}}\|(\mathbf{I-P})\mathfrak{f}\|_{L_{\nu}^2}^2 +C_{\kappa_5}|(I-P_{\gamma})\mathfrak{f}|_{L^2(\gamma_{+})}^2+C(\alpha+\tilde{\v})\|w^{\ell}(\mathfrak{f}_{1},\mathfrak{f}_{2})\|_{L^{\infty}}^2.
\end{align}

\medskip

\noindent\underline{\textbf{Conclusions on $\|\mathbf{P}\mathfrak{f}\|_{L^2}$}}: Combining \eqref{8.34}, \eqref{8.44}, \eqref{8.53}, \eqref{8.58-9} and \eqref{8.62} together and taking $0< \kappa_{5}\ll \kappa_{2}\ll \kappa_{1}\ll 1$ and $\kappa_{3}$, $\kappa_{4}$, $\alpha$ and $\lambda$ small, we get
\begin{align*}
&\frac{d}{dt}\big[\kappa_{2}\kappa_{5}\langle \mathfrak{f},\Psi_{a}\rangle+\kappa_{5}\langle \mathfrak{f},\Psi_{b_{1}}\rangle+\langle \mathfrak{f},\Psi_{b_{2}}\rangle+\langle \mathfrak{f},\Psi_{b_{3}}\rangle+\langle \mathfrak{f},\Psi_{c}\rangle\big]+\kappa \|(\mathfrak{a},\mathfrak{b}_{1},\mathfrak{b}_{2},\mathfrak{b}_{3},\mathfrak{c})\|_{L^2}^2\nonumber\\
&\leq C\alpha\|\mathfrak{b}_{2}\|_{L^2}\|(\mathfrak{a},\mathfrak{c})\|_{L^2}+C\alpha\|b_{1}\|_{L^2}^2+C_{\kappa}\big(\|(\mathbf{I-P})\mathfrak{f}\|_{L_{\nu}^2}^2+|(I-P_{\gamma})\mathfrak{f}|_{L^2(\gamma_{+})}^2+(\alpha+\tilde{\varepsilon})\|w^{\ell}(\mathfrak{f}_{1},\mathfrak{f}_{2})\|_{L^{\infty}}^2\big),
\end{align*}
which, together with H\"{o}lder inequality and smallness of $\alpha$, implies that
\begin{align}\label{8.64}
&\frac{d}{dt}\mathcal{E}_{\text{int}}(t)+\kappa \|(\mathfrak{a},\mathfrak{b}_{1},\mathfrak{b}_{2},\mathfrak{b}_{3},\mathfrak{c})\|_{L^2}^2 \leq C_{\kappa}\big(\|(\mathbf{I-P})\mathfrak{f}\|_{L_{\nu}^2}^2+|(I-P_{\gamma})\mathfrak{f}|_{L^2(\gamma_{+})}^2+(\alpha+\tilde{\varepsilon})\|w^{\ell}(\mathfrak{f}_{1},\mathfrak{f}_{2})\|_{L^{\infty}}^2\big),
\end{align}
where $\mathcal{E}_{\text{int}}$ is the quadratic functional defined as
\begin{align*}
\mathcal{E}_{\text{int}}=\kappa_{2}\kappa_{5}\langle \mathfrak{f},\Psi_{a}\rangle+\kappa_{5}\langle \mathfrak{f},\Psi_{b_{1}}\rangle+\langle \mathfrak{f},\Psi_{b_{2}}\rangle+\langle \mathfrak{f},\Psi_{b_{3}}\rangle+\langle \mathfrak{f},\Psi_{c}\rangle.
\end{align*}
Note that it holds
\begin{align*}
    |\mathcal{E}_{\text{int}}|\lesssim \|\mathfrak{f}_{2}\|_{L^2}^2+\|w^{\ell}\mathfrak{f}_{1}\|_{L^{\infty}}^2.
\end{align*}
Thus, \eqref{lem4.2.l2} together with \eqref{lem4.2.l1} follows. Therefore, the proof of Lemma \ref{lem4.2} is completed. $\hfill\square$

\subsection{Proof of Theorem \ref{thm2}} Recall the basic energy estimate \eqref{8.21} and the macroscopic dissipation estimate \eqref{lem4.2.l2}. Let $\tilde{C}>0$ be a constant sufficiently large. Then, by taking the summation $\tilde{C}\times \eqref{8.21}$ and \eqref{lem4.2.l2}, we get
\begin{align}\label{8.65}
&\frac{d}{dt}\mathcal{E}(t)+\kappa\|\mathfrak{f}_{2}\|_{L^2}^2+|(I-P_{\gamma})\mathfrak{f}_{2}|_{L^2(\gamma_{+})}^2+\|(\mathbf{I-P})\mathfrak{f}_{2}\|_{L_{\nu}^2}^2\nonumber\\
&\leq C_{\delta}\|w^{\ell}\mathfrak{f}_{1}\|_{L^{\infty}}^{2}+C_{\kappa}(\alpha+\tilde{\varepsilon})\|w^{\ell}\mathfrak{f}_{2}\|_{L^{\infty}}^2+C\delta\big[|P_{\gamma}\mathfrak{f}_{2}|_{L^2(\gamma_{+})}^2\big],
\end{align}
where
$$
\mathcal{E}(t):=\tilde{C}\|(1-\eta)^{1/2}\mathfrak{f}_{2}\|_{L^2}^2+\mathcal{E}_{\text{int}}.
$$
Since $\tilde{C}$ can be large enough, it follows from \eqref{lem4.2.l1} that
\begin{align}\label{8.66}
    \frac{\tilde{C}}{2}\|\mathfrak{f}_{2}\|_{L^2}^2-C\|w^{\ell}\mathfrak{f}_{1}\|_{L^{\infty}}^2\leq \mathcal{E}(t)\leq 2\tilde{C}\|\mathfrak{f}_{2}\|_{L^2}^2+C\|w^{\ell}\mathfrak{f}_{1}\|_{L^{\infty}}^2.
\end{align}
Applying \eqref{8.66} to \eqref{8.65}, one can replace $\|\mathfrak{f}_{2}\|_{L^2}^2$ by $\mathcal{E}(t)$, to obtain  
\begin{align}\label{8.67}
&\frac{d}{dt}\mathcal{E}(t)+\frac{\kappa}{2\tilde{C}}\mathcal{E}(t)+|(I-P_{\gamma})\mathfrak{f}_{2}|_{L^2(\gamma_{+})}^2+\|(\mathbf{I-P})\mathfrak{f}_{2}\|_{L_{\nu}^2}^2\nonumber\\
&\leq C_{\delta}\|w^{\ell}\mathfrak{f}_{1}\|_{L^{\infty}}^2+C_{\kappa}(\alpha+\tilde{\varepsilon})\|w^{\ell}\mathfrak{f}_{2}\|_{L^{\infty}}^2+C\delta|P_{\gamma}\mathfrak{f}_{2}|_{L^{2}(\gamma_{+})}^2.
\end{align}
It follows from \eqref{8.67} that
\begin{align}\label{8.68}
&\mathcal{E}(t)+\int_{0}^{t}e^{-\frac{\kappa}{2\tilde{C}}(t-s)}\big[|(I-P_{\gamma})\mathfrak{f}_{2}(s)|_{L^2(\gamma_{+})}^2\,{\rm d}s+\int_{0}^{t}e^{-\frac{\kappa}{2\tilde{C}}(t-s)}\|(\mathbf{I-P})\mathfrak{f}_{2}(s)\|_{L_{\nu}^2}^2\,{\rm d}s\nonumber\\
&\leq \mathcal{E}(0)e^{-\frac{\kappa}{2\tilde{C}}t}+C_{\delta}\int_{0}^{t}e^{-\frac{\kappa}{2\tilde{C}}(t-s)}\|w^{\ell}\mathfrak{f}_{1}(s)\|_{L^{\infty}}^2\,{\rm d}s+C_{\kappa}(\alpha+\tilde{\varepsilon})\int_{0}^{t}e^{-\frac{\kappa}{2\tilde{C}}(t-s)}\|w^{\ell}\mathfrak{f}_{2}(s)\|_{L^{\infty}}^2\,{\rm d}s\nonumber\\
&\quad + C\delta\int_{0}^{t}e^{-\frac{\kappa}{2\tilde{C}}(t-s)}|P_{\gamma}\mathfrak{f}_{2}(s)|_{L^{2}(\gamma_{+})}^2\,{\rm d}s\nonumber\\
&\leq \mathcal{E}(0)+C_{\delta}    \sup_{0\leq s\leq t} \|w^{\ell}\mathfrak{f}_{1}(s)\|_{L^{\infty}}^2+C_{\kappa}(\alpha+\tilde{\varepsilon})   \sup_{0\leq s\leq t} \|w^{\ell}\mathfrak{f}_{2}\|_{L^{\infty}}^2+ C\delta\int_{0}^{t}e^{-\frac{\kappa}{2\tilde{C}}(t-s)}|P_{\gamma}\mathfrak{f}_{2}(s)|_{L^{2}(\gamma_{+})}^2\,{\rm d}s.
\end{align}
Next, we estimate the last term on the right hand side of \eqref{8.68}. Note that at the boundary point $\eta=\eta_1$,
\begin{align}\label{8.21-1}
|P_{\gamma}\mathfrak{f}_{2}(\eta_{1})|_{L^2(\gamma_{+})}^2&=\sqrt{2\pi}\Big(\int_{v_{\eta}>0}\mathfrak{f}_{2}(\eta_{1})\sqrt{\mu}|v_{\eta}|\,{\rm d}v\Big)^2\nonumber\\
&=\sqrt{2\pi}\Big(\int_{V^{\tau}}\mathfrak{f}_{2}(\eta_{1})\sqrt{\mu}|v_{\eta}|\,{\rm d}v\Big)^2+\sqrt{2\pi}\Big(\int_{\{v_{\tau}>0\}\backslash V^{\tau}}\mathfrak{f}_{2}(\eta_{1})\sqrt{\mu}|v_{\eta}|\,{\rm d}v\Big)^2\nonumber\\
&=:I_{1}+I_{2},
\end{align}
where for any small enough constant $\tau>0$, we have denoted
$$
V^{\tau}=\{v\in \R^3\,:0<v_{\eta}<\tau\,\,\text{ or }|v_{\phi}|<\tau\,\,\text{ or }|v|>\frac{1}{\tau}\}.
$$
We need to estimate $I_1$ and $I_2$ in \eqref{8.21-1}. For $I_1$, from the H\"{o}lder inequality, it follows that
\begin{align*}
I_{1}&\leq \sqrt{2\pi}\Big(\int_{V^{\tau}}\mu|v_{\eta}|\,{\rm d}v\Big)\Big(\int_{V^{\delta}}|\mathfrak{f}_{2}|^2(\eta_{1})|v_{\eta}|\,{\rm d}v\Big)\leq C\tau\big[|P_{\gamma}\mathfrak{f}_{2}(\eta_{1})|_{L^2(\gamma_{+})}^2+|(I-P_{\gamma})\mathfrak{f}_{2}(\eta_{1})|_{L^2(\gamma_{+})}^2\big],
\end{align*}
which implies that
\begin{align}\label{8.69}
   &\int_{0}^{t}e^{-\frac{\kappa}{2\tilde{C}}(t-s)}I_{1}\,{\rm d}s\leq C\tau  \int_{0}^{t}e^{-\frac{\kappa}{2\tilde{C}}(t-s)}\big[|P_{\gamma}\mathfrak{f}_{2}(s,\eta_{1})|_{L^2(\gamma_{+})}^2+|(I-P_{\gamma})\mathfrak{f}_{2}(s,\eta_{1})|_{L^2(\gamma_{+})}^2\big]\,{\rm d}s.
\end{align}
For $I_{2}$, we obtain from Lemma \ref{tt} that
\begin{align}\label{8.70}
&\int_{0}^{t}e^{-\frac{\kappa}{2\tilde{C}}(t-s)}I_{2}\,{\rm d}s\leq C\int_{0}^{t}\int_{\{v_{\eta}\}>0\backslash V^{\tau}}e^{-\frac{\kappa}{2\tilde{C}}(t-s)}|\mathfrak{f}_{2}(s,\eta_{1}|^2|v_{\eta}|\,{\rm d}v{\rm d}s\nonumber\\
&\leq C_{\tau}\int_{0}^{t}e^{-\frac{\kappa}{2\tilde{C}}(t-s)}\|\mathfrak{f}_{2}(s,\eta)\|_{L^2}^2\,{\rm d}s\nonumber\\
&\quad +C_{\tau}\int_{0}^{t}\|\{\partial_{s}+v_{\eta}\partial_{\eta}-\frac{1}{1-\eta}\big(v_{\phi}^2\partial_{v_{\phi}}-v_{\eta}v_{\phi}\partial_{v_{\eta}}\big)-\frac{\alpha}{\eta_{1}}v_{\eta}\partial_{v_{\phi}}\}(e^{-\frac{\kappa}{2\tilde{C}}(t-s)}\mathfrak{f}_{2}^2)\|_{L^1}\,{\rm d}s\nonumber\\
&\leq C_{\tau}    \sup_{0\leq s\leq t} \|\mathfrak{f}_{2}(s)\|_{L_{\eta,v}^{2}}^2+C_{\tau}    \sup_{0\leq s\leq t} \|w^{\ell}\mathfrak{f}_{1}(s)\|_{L^{\infty}}^2+ C_{\tau}\int_{0}^{t}e^{-\frac{\kappa}{2\tilde{C}}(t-s)}\|(\mathbf{I-P})\mathfrak{f}_{2}(s)\|_{L_{\nu}^2}^2\,{\rm d}s,
\end{align}
where we have used \eqref{8.17}. Combining \eqref{8.69}--\eqref{8.70} and taking $\tau$ small enough, we conclude
\begin{align}\label{8.72}
&\int_{0}^{t}e^{-\frac{\kappa}{2\tilde{C}}(t-s)}|P_{\gamma}\mathfrak{f}_{2}(s,\eta_{1})|_{L^2(\gamma_{+})}^2\,{\rm d}s\nonumber\\
&\leq C\int_{0}^{t}e^{-\frac{\kappa}{2\tilde{C}}(t-s)}|(I-P_{\gamma})\mathfrak{f}_{2}(s,\eta_{1})|_{L^2(\gamma_{+})}^2\,{\rm d}s+C_{\tau}\int_{0}^{t}e^{-\frac{\kappa}{2\tilde{C}}(t-s)}\|(\mathbf{I-P})\mathfrak{f}_{2}(s)\|_{L_{\nu}^2}^2\,{\rm d}s\nonumber\\
&\quad +C_{\tau}    \sup_{0\leq s\leq t} \|w^{\ell}\mathfrak{f}_{1}(s)\|_{L^{\infty}}^2+C_{\tau}   \sup_{0\leq s\leq t} \|\mathfrak{f}_{2}(s)\|_{L_{\eta,v}^{2}}^2.
\end{align}
Similarly, at the other boundary point $\eta=0$, we also have
\begin{align}\label{8.71}
&\int_{0}^{t}e^{-\frac{\kappa}{2\tilde{C}}(t-s)}|P_{\gamma}\mathfrak{f}_{2}(s,0)|_{L^2(\gamma_{+})}^2\,{\rm d}s\nonumber\\
&\leq C\int_{0}^{t}e^{-\frac{\kappa}{2\tilde{C}}(t-s)}|(I-P_{\gamma})\mathfrak{f}_{2}(s,0)|_{L^2(\gamma_{+})}^2\,{\rm d}s+C_{\tau}\int_{0}^{t}e^{-\frac{\kappa}{2\tilde{C}}(t-s)}\|(\mathbf{I-P})\mathfrak{f}_{2}(s)\|_{L_{\nu}^2}^2\,{\rm d}s\nonumber\\
&\quad +C_{\tau}   \sup_{0\leq s\leq t} \|w^{\ell}\mathfrak{f}_{1}(s)\|_{L^{\infty}}^2+C_{\tau}   \sup_{0\leq s\leq t} \|\mathfrak{f}_{2}(s)\|_{L^{2}}^2.
\end{align}
Substituting \eqref{8.72}--\eqref{8.71} into \eqref{8.68} and taking $\delta$ suitably small, we obtain
\begin{align}\label{8.73}
&\sup_{0\leq s\leq t}  \mathcal{E}(s)
+\int_{0}^{t}e^{-\frac{\kappa}{2\tilde{C}}(t-s)}\big[|(I-P_{\gamma})\mathfrak{f}_{2}(s)|_{L^2(\gamma_{+})}^2+\int_{0}^{t}e^{-\frac{\kappa}{2\tilde{C}}(t-s)}\|(\mathbf{I-P})\mathfrak{f}_{2}(s)\|_{L_{\nu}^2}^2\,{\rm d}s\nonumber\\
&\leq \mathcal{E}(0) 
+C_{\delta}\sup_{0\leq s\leq t} \|w^{\ell}\mathfrak{f}_{1}(s)\|_{L^{\infty}}^2+C_{\kappa} (\alpha+\tilde{\varepsilon})\sup_{0\leq s\leq t} \|w^{\ell}\mathfrak{f}_{2}(s)\|_{L^{\infty}}^2+C\delta\sup_{0\leq s\leq t} \|w^{\ell}\mathfrak{f}_{2}(s)\|_{L^{\infty}}^2.
\end{align}
Therefore, using \eqref{8.66} and recalling \eqref{f12timew}, the above estimate \eqref{8.73} implies that 
\begin{align}\label{8.74}
\sup_{0\leq s\leq t}e^{2\lambda_0 s}\|f_{2}(s)\|_{L^{2}}^2& \leq \mathcal{E}(0)+C_{\delta}\sup_{0\leq s\leq t}e^{2\lambda_0 s}\|w^{\ell}f_{1}(s)\|_{L^{\infty}}^2 \nonumber\\
&\quad +C_{\kappa,\delta}(\alpha+\tilde{\varepsilon})\sup_{0\leq  s\leq t}e^{2\lambda_0 s}\|w^{\ell}f_{2}(s)\|_{L^{\infty}}^2+C\delta\sup_{0\leq  s\leq t}e^{2\lambda_0 s}\|w^{\ell}f_{2}(s)\|_{L^{\infty}}^2.
\end{align}
Plugging the above estimate \eqref{8.74} back to \eqref{8.3} and using the smallness of $\alpha$, $\tilde{\varepsilon}$ and $\delta$, one has
\begin{align}\label{8.75}
\sup_{0\leq s\leq t}e^{2\lambda_0 s}\|f_{2}(s)\|_{L^{\infty}}^2
&\leq C_{\ell}\|w^{\ell}f_{0}\|_{L^{\infty}}^2+C_{\ell}\sup_{0\leq s\leq t}e^{2\lambda_0 s}\|w^{\ell}f_{1}(s)\|_{L^{\infty}}^2.
\end{align}
Moreover, plugging \eqref{8.75} into \eqref{8.2}, and using the smallness of $\alpha$ and $\tilde{\varepsilon}$ as well as \eqref{8.75} again, we have
\begin{align*}
\sup_{0\leq s\leq t}e^{2\lambda_0 s}\|w^{\ell}[f_{1},f_{2}\|_{L^{\infty}}^2\leq C\|w^{\ell}f_{0}\|_{L^{\infty}}^2.
\end{align*}
Since $\|w^{\ell}f_{0}\|_{L^{\infty}}$ is sufficiently small, the {\it a priori} assumption \eqref{8.1} is closed.

Finally, combining the existence of non-negativity of local-in-time solution established in Theorem \ref{thm4.1} with exponentially asymptotic behavior \eqref{7.31}, we obtain the existence of non-negative global solution to \eqref{2.6} satisfying \eqref{M2}. Therefore, the proof of Theorem \ref{thm2} is completed. $\hfill\square$

\bigskip

\appendix
\renewcommand{\appendixname}{Appendix~\Alph{section}}

\section{Some useful estimates}\label{AppendixA}
In this section, we collect some useful estimates, which have been used extensively in the previous sections.

\begin{lemma}[Properties of $K$ and $\mathcal{K}$ {\cite[Section 5]{DL-2022}}]\label{K}

		{\rm 1)} Let $K_{w}$ be defined as \eqref{3.9} with $\ell\geq 0$, then it holds that
		\begin{align*}
			&\int_{\R^3}k_{w}(v,u)e^{\frac{c_{0}|v-u|^2}{8}}\,{\rm d}u\leq \frac{C}{1+|v|},\quad \int_{\R^3}|k(v,u)|^2\,{\rm d}u<C,
		\end{align*}
		for $c_{0}=0$ or $c_{0}>0$ small enough. Moreover, for any $\ell>0$, one has
		\begin{align*}
			\|w^{\ell}Kf\|_{L^{\infty}}\leq C\|w^{\ell-1}f\|_{L^{\infty}}.
		\end{align*}
		
{\rm 2)} Let $\mathcal{K}_{w}$ be defined in \eqref{3.10}. 
		\begin{itemize}
			\item For $\gamma=0$, there exists $C>0$ such that for any arbitrarily $\ell>0$, one can choose $M\geq \ell^2$ such that it holds
			$$
			(1-\chi_{M})|\mathcal{K}_{w}h|\leq \frac{C}{\ell}\|h\|_{L^{\infty}}.
			$$
			\item For $0<\gamma\leq 1$, then there is a constant $C>0$ and for any arbitrarily large $\ell>4$, there are sufficiently large $M=M(\ell)>0$ and suitably small $\zeta=\zeta(\ell)>0$ such that it holds
			$$
			(1-\chi_{M})\nu^{-1}|\mathcal{K}_{w}h|\leq C\{(1+M)^{-\gamma}+\zeta\}\|h\|_{L^{\infty}}.
			$$
			\item Let $\ell>4$ and $0\leq \gamma\leq 1$, the it holds that
			\begin{align*}
				|w^{\ell}\nu^{-1}Q_{+}(F_{1},F_{2})|+|w^{\ell}\nu^{-1}Q_{-}(F_{1},F_{2})|\leq C\|w^{\ell}F_{1}\|_{L^{\infty}}\|w^{\ell}F_{2}\|_{L^{\infty}}.
			\end{align*}
			In particular, recall that
			$$
			\mathcal{K}f=Q(f,\mu)+Q_{+}(\mu,f)=\sqrt{\mu}K\big(\frac{f}{\sqrt{\mu}}\big),
			$$
		    it holds that
		    \begin{align*}
		    	\|w^{\ell}\nu^{-1}\mathcal{K}f\|_{L^{\infty}}\leq C\|w^{\ell}f\|_{L^{\infty}}.
		    \end{align*}
		\end{itemize}
	\end{lemma}

	\begin{lemma}[A trace lemma]\label{tt}
		Let $\tau\ll 1$ and $\eta\in [0,\eta_{1}]$, and denote the near-grazing set of $\gamma_{+}$ 
		$$
		\gamma_{+}^{\tau}=\Big\{(\eta,v_{\eta},v_{\phi},v_{z})\in \gamma_{+}:|v_{\eta}|\leq \tau\,\,\text{or }\,\,|v_{\phi}|\leq \tau\,\,\text{ or }\,\,|v|\geq \frac{1}{\tau}\Big\}.
		$$
		Then, there exists a constant $C_{\tau}$  depending on $\tau$ such that
		\begin{align}\label{t1}
		|f\mathbf{1}_{\gamma_{+}\backslash \gamma_{+}^{\tau}}|_{L^1}\leq C_{\tau}\big\{\|f\|_{L^1}+\|\big\{v_{\eta}\partial_{\eta}-\frac{1}{1-\eta}\big(v_{\phi}^2\partial_{v_{\phi}}-v_{\eta}v_{\phi}\partial_{v_{\eta}}\big)-\frac{\alpha}{\eta_{1}}v_{\eta}\partial_{v_{\phi}}\big\}f\|_{L^1}\big\}.
		\end{align}
		Moreover, it also holds that
		\begin{align}\label{t2}
		&\int_{0}^{T}|f{\bf 1}_{\gamma_{+}\backslash \gamma_{+}^{\tau}}|_{L_{v}^1}\,{\rm d}t\leq C_{\tau}\Big[
		\|f(0)\|_{L^1}+\int_{0}^{T}\|f(t)\|_{L^1}\,{\rm d}t\Big]\nonumber\\
		&\qquad \qquad +C_{\tau}\int_{0}^{T}\|\{\partial_{t}+v_{\eta}\partial_{\eta}-\frac{1}{1-\eta}\big(v_{\phi}^2\partial_{v_{\phi}}-v_{\eta}v_{\phi}\partial_{v_{\eta}}\big)-\frac{\alpha}{\eta_{1}}v_{\eta}\partial_{v_{\phi}}\}f(t)\|_{L^1}\big]\,{\rm d}t,
		\end{align}
		for any $T\geq 0$.
	\end{lemma}
	
\noindent\textbf{Proof.}
1. We first focus on the steady case, {\it i.e.}, \eqref{t1}. We define a characteristic line 
$$
[y,\mathscr{V}_{\eta}(y;\eta,v), \mathscr{V}_{\phi}(y;\eta,v),\mathscr{V}_{z}(y;\eta,v)]
$$
passing through $(\eta,v)$ such that
	\begin{align*}
		\left\{
		\begin{aligned}
			&\frac{d \mathscr{V}_{\eta}}{d y}=-\frac{1}{(1-y)\mathscr{V}_{\eta}}\mathscr{V}_{\phi}^2,\quad \frac{d \mathscr{V}_{\phi}}{d y}=\frac{1}{1-y}\mathscr{V}_{\phi}-\frac{\alpha}{\eta_{1}},\quad \frac{d\mathscr{V}_{z}}{dy}=0,\\
			&\mathscr{V}_{\eta}\vert_{y=\eta}=v_{\eta},\quad \mathscr{V}_{\phi}\vert_{y=\eta}=v_{\phi},\quad \mathcal{V}_{z}\vert_{y=\eta}=v_{z}.
		\end{aligned}
		\right.
	\end{align*}
	Then it follows from direct calculations that
	\begin{align}\notag
		\begin{aligned}
			&\mathscr{V}_{\phi}(y;\eta,v_{\eta},v_{\phi})=\frac{1-\eta}{1-y}v_{\phi}+\frac{\alpha}{2\eta_{1}}(1-y)-\frac{\alpha}{2\eta_{1}}\frac{(1-\eta)^2}{1-y},\quad \mathscr{V}_{z}(y;\eta,v)=v_{z},\\
			&|\mathscr{V}_{\eta}(y;\eta,v_{\eta},v_{\phi})|^2=v_{\eta}^2-\int_{\eta}^{y}\frac{2}{1-z}|\mathscr{V}_{\phi}(z;\eta,v_{\eta},v_{\phi})|^2\,{\rm d}z,
		\end{aligned}
	\end{align}
	for $(\eta,v)\in \gamma_{+}\backslash \gamma_{+}^{\tau}$. Furthermore,
	\begin{align}\notag
		|\mathscr{V}_{\eta}(y;\eta,v)|^2+|\mathscr{V}_{\phi}(y;\eta,v)|^2+|\mathscr{V}_{z}(y;\eta,v)|^2=v_{\eta}^2+v_{\phi}^2+v_{z}^2+O(\alpha)(1+|v_{\phi}|),
	\end{align}
	which yields
	\begin{align*}
		|\mathscr{V}_{\eta}(y;\eta,v)|^2+|\mathscr{V}_{\phi}(y;\eta,v)|^2+|\mathscr{V}_{z}(y;\eta,v)|^2\leq C(|v|^2+\alpha)\leq \frac{C}{\tau}\quad \text{for }(\eta,v)\in \gamma_{+}\backslash \gamma_{+}^{\tau},
	\end{align*}
	provided $\alpha<1$ and $\tau\ll 1$. We also define $\eta_{+}(\eta,v)$ by
	\begin{align*}
	  \mathscr{V}_{\eta}(\eta_{+};\eta,v)=0\Leftrightarrow v_{\eta}^2=\int_{\eta}^{\eta_{+}}\big[\frac{1-\eta}{1-y}v_{\phi}+\frac{\alpha}{2\eta_{1}}(1-y)-\frac{\alpha}{2\eta_{1}}\frac{(1-\eta)^2}{1-y}\big]^2\,{\rm d}y.
	\end{align*}
	Along the trajectory, one has the identity
	\begin{align*}
		f(\eta,v)=f(y,\mathscr{V}_{\eta}(y;\eta,v), \mathscr{V}_{\phi}(y;\eta,v),\mathscr{V}_{z}(y;\eta,v))-\int_{\eta}^{y}\frac{d}{d z}f(z,\mathscr{V}_{\eta}(z;\eta,v), \mathscr{V}_{\phi}(z;\eta,v),\mathscr{V}_{z}(z;\eta,v))\,{\rm d}z.
	\end{align*}
If $\eta=\eta_{1}$ or $\eta=0$ but $\eta_{+}(0,v_{\eta},v_{\phi})>\eta_{1}$, then
	\begin{align}\label{CL-1}
		&\int_{\gamma_{+}\backslash \gamma_{+}^{\tau}} f(\eta,v)|v_{\eta}|\,{\rm d}v\nonumber\\
		&\leq \frac{C}{\eta_{1}}\int_{\gamma_{+} \backslash \gamma_{+}^{\tau}}\int_{0}^{\eta_{1}}|f(y,\mathscr{V}_{\eta}(y;\eta,v), \mathscr{V}_{\phi}(y;\eta,v),\mathscr{V}_{z}(y;\eta,v))|\cdot |v_{\eta}|\,{\rm d}y{\rm d}v\nonumber\\
		&\,\,\,\,+\frac{C}{\eta_{1}}\int_{\gamma_{+} \backslash \gamma_{+}^{\tau}}\int_{0}^{\eta_{1}}\frac{|v_{\eta}|}{|\mathscr{V}_{\eta}(y)|}\big|\big\{\mathscr{V}_{\eta}\partial_{y}-\frac{1}{1-y}\big(\mathscr{V}_{\phi}^2\partial_{\mathscr{V}_{\eta}}-\mathscr{V}_{\eta}\mathscr{V}_{\phi}\partial_{\mathscr{V}_{\phi}}\big)-\frac{\alpha}{\eta_{1}}\mathscr{V}_{\eta}\partial_{\mathscr{V}_{\phi}}\big\}f(y,\mathscr{V})\big|\,{\rm d}y{\rm d}v.
	\end{align}
If $\eta=0$ but  $\eta_{+}(0,v)\leq \eta_{1}$, then one has
\begin{align*}
   \tau^2\leq v_{\eta}^2=\int_{0}^{\eta_{+}}\big[\frac{1}{1-y}v_{\phi}+\frac{\alpha}{2\eta_{1}}(1-y)-\frac{\alpha}{2\eta_{1}}\frac{1}{1-y}\big]^2\,{\rm d}y\leq \frac{C}{\tau^2}\eta_{+}\quad \text{for }(0,v)\in \gamma_{+}\backslash \gamma_{+}^{\tau},
\end{align*}
which implies $\eta_{+}\geq \frac{\tau^4}{C}$. Then
	\begin{align}\label{CL-2-0}
		&\int_{\gamma_{+}\backslash \gamma_{+}^{\tau}} f(\eta,v_{\eta},v_{\phi},v_{z})|v_{\eta}|\,{\rm d}v\nonumber\\
		&\leq \frac{C_{\tau}}{\eta_{1}}\int_{\gamma_{+} \backslash \gamma_{+}^{\tau}}\int_{0}^{\eta_{+}}|f(y,\mathscr{V}_{\eta}(y;\eta,v), \mathscr{V}_{\phi}(y;\eta,v),\mathscr{V}_{z}(y;\eta,v))|\cdot |v_{\eta}|\,{\rm d}y{\rm d}v\nonumber\\
		&\,\,\,\,+C_{\tau}\int_{\gamma_{+} \backslash \gamma_{+}^{\tau}}\int_{0}^{\eta_{+}}\frac{|v_{\eta}|}{|\mathscr{V}_{\eta}(y)|}\big|\big\{\mathscr{V}_{\eta}\partial_{y}-\frac{1}{1-y}\big(\mathscr{V}_{\phi}^2\partial_{\mathscr{V}_{\eta}}-\mathscr{V}_{\eta}\mathscr{V}_{\phi}\partial_{\mathscr{V}_{\phi}}\big)-\frac{\alpha}{\eta_{1}}\mathscr{V}_{\eta}\partial_{\mathscr{V}_{\phi}}\big\}f(y,\mathscr{V})\big|\,{\rm d}y{\rm d}v.
	\end{align}

Noting the Jacobian
	\begin{align*}
		\frac{\partial (y,\mathscr{V}_{\eta}(y;\eta,v), \mathscr{V}_{\phi}(y;\eta,v),\mathscr{V}_{z}(y;\eta,v))}{\partial (y,v_{\eta},v_{\phi},v_{z})}=\left\vert
		\begin{matrix}
			1&0&0&0\\
		\frac{\partial \mathscr{V}_{\eta}(y;\eta,v_{\eta},v_{\phi})}{\partial y}&\frac{v_{\eta}}{\mathscr{V}_{\eta}(y)}&	\frac{\partial \mathscr{V}_{\eta}(y;\eta,v_{\eta},v_{\phi})}{\partial v_{\phi}}&0\\
				\frac{\partial \mathscr{V}_{\phi}(y;\eta,v_{\eta},v_{\phi})}{\partial y}&0&\frac{1-\eta}{1-y}&0\\
                0&0&0&1
		\end{matrix}
		\right\vert=\frac{v_{\eta}}{\mathscr{V}_{\eta}(y)}\frac{1-\eta}{1-y},
	\end{align*}
	we make the following change of variable:
	$$
	[z,u_{\eta},u_{\phi},u_{z}]=[y,\mathscr{V}_{\eta}(y;\eta,v),\mathscr{V}_{\phi}(y;\eta,v),\mathscr{V}_{z}(y;\eta,v)],
	$$
	to obtain
	\begin{align}\label{CL-5}
		&\int_{\gamma_{+} \backslash \gamma_{+}^{\tau}}\int_{0}^{\eta_{1}}|f(y,\mathscr{V}_{\eta}(y;\eta,v), \mathscr{V}_{\phi}(y;\eta,v),\mathscr{V}_{z}(y;\eta,v))|\,|v_{\eta}|\,{\rm d}y{\rm d}v\nonumber\\
		&\leq C\int_{0}^{\eta_{1}}\int_{B(0,\sup_{\{z\in [0,\eta_{1}],(\eta,v)\in \gamma_{+}\backslash \gamma_{+}^{\tau}\}}|u|)}|u_{\eta}|\, |f(w,u_{\eta},u_{\phi},u_{z})|\,{\rm d}u{\rm d}z\nonumber\\
		&\leq C_{\tau}\|f\|_{L^1},
	\end{align}
	and
	\begin{align}\label{CL-6}
		&\int_{\gamma_{+} \backslash \gamma_{+}^{\tau}}\int_{0}^{\eta_{1}}\frac{|v_{\eta}|}{|\mathscr{V}_{\eta}(y)|}\big|\big\{\mathscr{V}_{\eta}\partial_{y}-\frac{1}{1-y}\big(\mathscr{V}_{\phi}^2\partial_{\mathscr{V}_{\eta}}-\mathscr{V}_{\eta}\mathscr{V}_{\phi}\partial_{\mathscr{V}_{\phi}}\big)-\frac{\alpha}{\eta_{1}}\mathscr{V}_{\eta}\partial_{\mathscr{V}_{\phi}}\big\}f(y,\mathscr{V}_{\eta},\mathscr{V}_{\phi},\mathscr{V}_{z})\big|\,{\rm d}y{\rm d}v\nonumber\\
		&\leq C\int_{0}^{\eta_{1}}\int_{\R^3} \big|\big\{u_{\eta}\partial_{w}-\frac{1}{1-z}\big(u_{\phi}^2\partial_{u_{\eta}}-u_{\eta}u_{\phi}\partial_{u_{\phi}}\big)-\frac{\alpha}{\eta_{1}}u_{\eta}\partial_{u_{\phi}}\big\}f(z,u)\big|\,{\rm d}u{\rm d}z\nonumber\\
		&\leq C\|\big\{v_{\eta}\partial_{\eta}-\frac{1}{1-\eta}\big(v_{\phi}^2\partial_{v_{\phi}}-v_{\eta}v_{\phi}\partial_{v_{\eta}}\big)-\frac{\alpha}{\eta_{1}}v_{\eta}\partial_{v_{\phi}}\big\}f\|_{L^1}.
	\end{align}
	Substituting \eqref{CL-5}--\eqref{CL-6} into \eqref{CL-1}--\eqref{CL-2-0}, we conclude the proof of \eqref{t1}. 
	
	2. For the non-steady case \eqref{t2}, we define the characteristic line
	$$
	[\mathscr{T}(y;t,\eta,v),y,\mathscr{V}_{\eta}(y;t,\eta,v),\mathscr{V}_{\phi}(y;t,\eta,v),\mathscr{V}_{z}(y;t,\eta,v)],
	$$
	passing through $(t,\eta,v)$, in terms of
	\begin{align*}
	    \left\{
	    \begin{aligned}
	    &\frac{d\mathscr{T}}{dy}=\frac{1}{\mathscr{V}_{\eta}},\,\,\,\frac{d\mathscr{V}_{\eta}}{ dy}=-\frac{1}{1-y}\mathscr{V}_{\phi}^2,\,\,\,\frac{d\mathscr{V}_{\phi}}{dy}=\frac{1}{1-y}\mathscr{V}_{\phi}-\frac{\alpha}{\eta_{1}},\,\,\,\frac{d\mathscr{V}_{z}}{dy}=0\\
	    &\mathscr{T}\vert_{y=\eta}=t,\,\,\,\mathscr{V}_{\eta}\vert_{y=\eta}=v_{\eta},\,\,\,\mathscr{V}_{\phi}\vert_{y=\eta}=v_{\phi},\,\,\,\mathscr{V}_{z}\vert_{y=\eta}=v_{z}.
	   \end{aligned}
	    \right.
	\end{align*}
It follows from direct calculations that
\begin{align*}
&\mathscr{V}_{\phi}(y;t,\eta,v)=\frac{1-\eta}{1-y}v_{\phi}+\frac{\alpha}{2\eta_{1}}(1-y)-\frac{\alpha}{2\eta_{1}}\frac{(1-\eta)^2}{1-y},\quad \mathscr{V}_{z}(y;t,\eta,v)=v_{z},\\
&|\mathscr{V}_{\eta}(y;t,\eta,v)|^2=v_{\eta}^2+\int_{y}^{\eta}\frac{2}{1-z}|\mathscr{V}_{\phi}(z;t,\eta,v)|^2\,{\rm d}z,\\
    &\mathscr{T}(y;t,\eta,v)=t+\int_{\eta}^{y}\frac{1}{\mathscr{V}_{\eta}(z;t,\eta,v)}\,{\rm d}z=\int_{\eta}^{y}\frac{1}{\mathscr{V}_{\eta}(z;\eta,v)}\,{\rm d}z,
\end{align*}
	which yields that
	\begin{align}\label{t3}
		|\mathscr{V}_{\eta}(y;t,\eta,v)|^2+|\mathscr{V}_{\phi}(y;t,\eta,v)|^2+|\mathscr{V}_{z}(y;t,\eta,v)|^2\leq C(|v|^2+\alpha)\leq \frac{C}{\tau}\quad (\eta,v)\in \gamma_{+}\backslash \gamma_{+}^{\tau},
	\end{align}
	provided $\alpha<1$ and $\tau\ll 1$. Moreover, along the trajectory, we have
	\begin{align}\label{t4-1}
	 f(t,\eta_{1},v)&=f(\mathscr{T}(y;t,\eta,v),y,\mathscr{V}_{\eta}(y;t,\eta,v),\mathscr{V}_{\phi}(y;t,\eta,v),\mathscr{V}_{z}(y;t,\eta,v))\nonumber\\
	 &\quad +\int_{y}^{\eta_{1}}\frac{d}{dz}f(\mathscr{T}(z;t,\eta,v),z,\mathscr{V}_{\eta}(z;t,\eta,v),\mathscr{V}_{\phi}(z;t,\eta,v),\mathscr{V}_{z}(z;t,\eta,v))\,{\rm d}z.
	\end{align}
	For each $(t,\eta,v)\in (0,T]\times \big(\gamma_{+}\backslash \gamma_{+}^{\tau}\big)$ and $\tilde{t}\in (0,t]$, we define $\mathscr{Y}(\tilde{t};\eta,v)$ as the position where
	$$
	\mathscr{T}(\mathscr{Y}(\tilde{t});t,\eta,v)=t+\int_{\eta}^{\mathscr{Y}(\tilde{t})}\frac{1}{\mathscr{V}_{\eta}(z;t,\eta,v)}\,{\rm d}z=t-\tilde{t}\Rightarrow \int_{\mathscr{Y}(\tilde{t})}^{\eta}\frac{1}{\mathscr{V}_{\eta}(z;t,\eta,v)}\,{\rm d}z=\tilde{t}.
	$$
	{\it Case 1}. If $(\eta,v)\in \mathcal{A}_{1}\cap (\gamma_{+}\backslash \gamma_{+}^{\tau})$, i.e., $\eta=\eta_{1}$ and $v_{\eta}>\tau>0$, then noting $\mathscr{V}_{\eta}(z;\eta_{1},v)\geq v_{\eta}\geq \tau$ in this case, one has
	\begin{align}\label{t5-1}
	\tilde{t}=\int_{\mathscr{Y}(\tilde{t})}^{\eta_{1}}\frac{1}{\mathscr{V}_{\eta}(z;t,\eta,v)}\,{\rm d}z\leq \frac{C}{\tau}(\eta_{1}-\mathscr{Y}(\tilde{t}))\Rightarrow |\eta_{1}-\mathscr{Y}(\tilde{t})|\geq \frac{\tilde{t}\tau}{C}.
	\end{align}
    We then get from \eqref{t4-1} that
    \begin{align}\label{t5-2}
     &\min\{\eta_{1},\eta_{1}-\mathscr{Y}(\tilde{t})\}\int_{\tilde{t}}^{T}\int_{A_{\tau}}|f(t,\eta_{1},v)|\,|v_{\eta}|\,{\rm d}v{\rm d}t\nonumber\\
     &\leq \int_{\tilde{t}}^{T}\int_{A_{\tau}}\int_{\max\{0,\mathscr{Y}(\tilde{t})\}}^{\eta_{1}}|f(\mathscr{T}(y;t,\eta_{1},v),y,\mathscr{V}_{\eta}(y;t,\eta_{1},v),\mathscr{V}_{\phi}(y;t,\eta_{1},v),\mathscr{V}_{z}(y;t,\eta_{1},v))|\,|v_{\eta}|\,{\rm d}y{\rm d}v{\rm d}t\nonumber\\
     &\quad +\int_{\tilde{t}}^{T}\int_{A_{\tau}}\int_{\max\{0,\mathscr{Y}(\tilde{t})\}}^{\eta_{1}}\int_{y}^{\eta_{1}}|\mathscr{D}_{z}f(\mathscr{T}(z),z,\mathscr{V}_{\eta}(z),\mathscr{V}_{\phi}(z),\mathscr{V}_{z}(z))|\times \frac{|v_{\eta}|}{|\mathscr{V}_{\eta}(z)|}\,{\rm d}z{\rm d}y{\rm d}v{\rm d}t\nonumber\\
     &\leq \int_{\tilde{t}}^{T}\int_{A_{\tau}}\int_{\max\{0,\mathscr{Y}(\tilde{t})\}}^{\eta_{1}}|f(\mathscr{T}(y;t,\eta_{1},v),y,\mathscr{V}_{\eta}(y;t,\eta_{1},v),\mathscr{V}_{\phi}(y;t,\eta_{1},v),\mathscr{V}_{z}(y;t,\eta_{1},v))|\,|v_{\eta}|\,{\rm d}y{\rm d}v{\rm d}t\nonumber\\
     &\quad +\eta_{1}\int_{\tilde{t}}^{T}\int_{A_{\tau}}\int_{\max\{0,\mathscr{Y}(\tilde{t})\}}^{\eta_{1}}|\mathscr{D}_{y}f(\mathscr{T}(y),y,\mathscr{V}_{\eta}(y),\mathscr{V}_{\phi}(y),\mathscr{V}_{z}(y))|\times \frac{|v_{\eta}|}{|\mathscr{V}_{\eta}(y)|}\,{\rm d}y{\rm d}v{\rm d}t,
    \end{align}
	where we have denoted
	$$
	A_{\tau}=\{(v_{\eta},v_{\phi},v_{z})\in \R^3\,|\,\tau\leq |v_{\phi}|\leq |v|\leq \frac{1}{\tau},\,\,|v_{\eta}|\geq \tau\}
	$$
	and
	$$
	\mathscr{D}_{z}=:\partial_{t}-\mathscr{V}_{\eta}\partial_{z}-\frac{1}{1-z}(\mathscr{V}_{\phi}^2\partial_{\mathscr{V}_{\eta}}-\mathscr{V}_{\eta}\mathscr{V}_{\phi}\partial_{\mathscr{V}_{\phi}})-\frac{\alpha}{\eta_{1}}\mathscr{V}_{\eta}\partial_{\mathscr{V}_{\phi}}.
	$$
Noting the Jacobian
\begin{align*}
&\frac{\partial(\mathscr{T}(y;t,\eta,v),y,\mathscr{V}_{\eta}(y;t,\eta,v),\mathscr{V}_{\phi}(y;t,\eta,v),\mathscr{V}_{z}(y;t,\eta,v))}{\partial(t,y,v_{\eta},v_{\phi},v_{z})}\\&=\left\vert
\begin{matrix}
1&\frac{\partial\mathscr{T}}{\partial y}&\frac{\partial\mathscr{T}}{\partial v_{\eta}}&\frac{\partial\mathscr{T}}{\partial v_{\phi}}&0\\
0&1&0&0&0\\
0&\frac{\partial \mathscr{V}_{\eta}}{\partial y}&\frac{v_{\eta}}{\mathscr{V}_{\eta}}&\frac{\partial \mathscr{V}_{\eta}}{\partial v_{\phi}}&0\\
0&\frac{\partial \mathscr{V}_{\phi}}{\partial y}&0&\frac{1-y}{1-\eta}&0\\
0&0&0&0&1
\end{matrix}
\right\vert
=\frac{v_{\eta}}{\mathscr{V}_{\eta}(y)}\frac{1-\eta}{1-y},
\end{align*}
we make the following change of variable:
\begin{align*}
    [s,z,u_{\eta},u_{\phi},u_{z}]=[\mathscr{T}(y;t,\eta,v),y,\mathscr{V}_{\eta}(y;t,\eta,v),\mathscr{V}_{\phi}(y;t,\eta,v),\mathscr{V}_{z}(y;t,\eta,v)],
\end{align*}
to get
\begin{align}\label{t5-3}
&\int_{\tilde{t}}^{T}\int_{A_{\tau}}\int_{\max\{0,\mathscr{Y}(\tilde{t}\}}^{\eta_{1}}|f(\mathscr{T}(y;t,\eta_{1},v),y,\mathscr{V}_{\eta}(y;t,\eta_{1},v),\mathscr{V}_{\phi}(y;t,\eta_{1},v),\mathscr{V}_{z}(y;t,\eta_{1},v))|\,|v_{\eta}|\,{\rm d}y{\rm d}v{\rm d}t\nonumber\\
&\leq C_{\tau}\int_{0}^{T}\int_{\R^3}\int_{0}^{\eta_{1}}|f(s,z,u)|\,{\rm d}z{\rm d}u{\rm d}s,
\end{align}
and
\begin{align}\label{t5-4}
&\eta_{1}\int_{\tilde{t}}^{T}\int_{A_{\tau}}\int_{\max\{0,\mathscr{Y}(\tilde{t}\}}^{\eta_{1}}|\mathscr{D}_{y}f(\mathscr{T}(y),y,\mathscr{V}_{\eta}(y),\mathscr{V}_{\phi}(y),\mathscr{V}_{z}(y))|\times \frac{|v_{\eta}|}{|\mathscr{V}_{\eta}(y)|}\,{\rm d}y{\rm d}v{\rm d}t\nonumber\\
&\leq C\int_{0}^{T}\int_{\R^3}\int_{0}^{\eta_{1}}|\{\partial_{s}-u_{\eta}\partial_{w}-\frac{1}{1-z}(u_{\phi}^2\partial_{u_{\eta}}-u_{\eta}u_{\phi}\partial_{u_{\phi}})-\frac{\alpha}{\eta_{1}}u_{\eta}\partial_{u_{\phi}}\}f(s,z,u)|\,{\rm d}z{\rm d}u{\rm d}s,
\end{align}
where we have used \eqref{t3} together with the fact that
$$
s=\mathscr{T}(y;t,\eta,v_{\eta},v_{\phi})\in [0,T]\quad \text{for }t\in [\tilde{t},T]\text{ and }y\in [\max\{0,\mathscr{Y}(\tilde{t})\},\eta_{1}].
$$
Combining \eqref{t5-1}--\eqref{t5-4} together, we get
\begin{align}\label{t5-5}
&\int_{\tilde{t}}^{T}\int_{A_{\tau}}|f(t,\eta_{1},v)|\,|v_{\eta}|\,{\rm d}v{\rm d}t\nonumber\\
&\leq C_{\tau,\tilde{t}}\Big[\int_{0}^{T}\|f(t)\|_{L^1}\,{\rm d}t
+\int_{0}^{T}\|\{\partial_{t}+v_{\eta}\partial_{v_{\eta}}-\frac{1}{1-\eta}(v_{\phi}^2\partial_{v_{\eta}}-v_{\eta}v_{\phi}\partial_{v_{\phi}})-\frac{\alpha}{\eta_{1}}v_{\eta}\partial_{v_{\phi}}\}f(t)\|_{L^1}\,{\rm d}t\Big].
\end{align}
It remains to control 
$$
\int_{0}^{\tilde{t}}\int_{A_{\tau}}|f(t,\eta_{1},v_{\eta},v_{\phi})|\,|v_{\eta}|\,{\rm d}v{\rm d}t.
$$
In fact, if we choose $\tilde{t}$ to be small enough such that $\tilde{t}\lesssim \eta_{1}\tau$, at this stage, the backward trajectory hits the initial plane first before hitting the lower boundary $\{\eta=0\}$ due to the fact:
$$
\int_{0}^{\eta_{1}}\frac{1}{\mathscr{V}_{\eta}(z;t,\eta_{1},v)}\,{\rm d}z\gtrsim \eta_{1}\tau.
$$
So, there must exist a unique $\mathscr{Y}(t,\eta_{1},v)>0$ satisfying
\begin{align*}
t=\int_{\mathscr{Y}(t)}^{\eta_{1}}\frac{1}{\mathscr{V}_{\eta}(z;t,\eta_{1},v)}\,{\rm d}z\quad \text{for each $t\in (0,\tilde{t})$}.
\end{align*}
Direct calculations show that
\begin{align*}
\frac{\partial \mathscr{Y}(t,\eta_{1},v)}{\partial t}=-\mathscr{V}_{\eta}(\mathscr{Y}(t);t,\eta_{1},v).
\end{align*}
Therefore,
\begin{align}\label{t6-1}
&\int_{0}^{\tilde{t}}\int_{A_{\tau}}|f(t,\eta_{1},v_{\eta},v_{\phi})|\,|v_{\eta}|\,{\rm d}v{\rm d}t\nonumber\\
&\leq \int_{0}^{\tilde{t}}\int_{A_{\tau}}|f(0,\mathscr{Y}(t),\mathscr{V}_{\eta}(\mathscr{Y}(t);t,\eta_{1},v),\mathscr{V}_{\phi}(\mathscr{Y}(t);t,\eta_{1},v),\mathscr{V}_{z}(\mathscr{Y}(t);t,\eta_{1},v))\,{\rm d}v{\rm d}t\nonumber\\
&\quad +\int_{0}^{\tilde{t}}\int_{A_{\tau}}\int_{\mathscr{Y}(t)}^{\eta_{1}}|\mathscr{D}_{y}f(\mathscr{T}(y),y,\mathscr{V}_{\eta}(y),\mathscr{V}_{\phi}(y),\mathscr{V}_{z}(y))|\times \frac{|v_{\eta}|}{|\mathscr{V}_{\eta}(y)|}\,{\rm d}y{\rm d}v{\rm d}t\nonumber\\
&\leq \int_{0}^{\tilde{t}}\int_{A_{\tau}}|f(0,\mathscr{Y}(t),\mathscr{V}_{\eta}(\mathscr{Y}(t);t,\eta_{1},v),\mathscr{V}_{\phi}(\mathscr{Y}(t);t,\eta_{1},v),\mathscr{V}_{z}(\mathscr{Y}(t);t,\eta_{1},v))\,{\rm d}v{\rm d}t\nonumber\\
&\quad +C\int_{0}^{T}\|\{\partial_{t}+v_{\eta}\partial_{v_{\eta}}-\frac{1}{1-\eta}(v_{\phi}^2\partial_{v_{\eta}}-v_{\eta}v_{\phi}\partial_{v_{\phi}})-\frac{\alpha}{\eta_{1}}v_{\eta}\partial_{v_{\phi}}\}f(t)\|_{L^1}\,{\rm d}t.
\end{align}
Noting the Jacobian
\begin{align*}
&\frac{\partial(\mathscr{Y}(t,\eta_{1},v),\mathscr{V}_{\eta}(\mathscr{Y}(t);t,\eta_{1},v),\mathscr{V}_{\phi}(\mathscr{Y}(t);t,\eta_{1},v),\mathscr{V}_{z}(\mathscr{Y}(t);t,\eta_{1},v)}{\partial(t,v_{\eta},v_{\phi},v_{z})}\nonumber\\
&=
\left\vert
\begin{matrix}
\frac{\partial\mathscr{Y}}{\partial t}&\frac{\partial\mathscr{Y}}{\partial v_{\eta}}&\frac{\partial\mathscr{Y}}{\partial v_{\phi}}&0\\
\frac{\partial\mathscr{V}_{\eta}}{\partial \mathscr{Y}}\times \frac{\partial\mathscr{Y}}{\partial t}&\frac{v_{\eta}}{\mathscr{V}_{\eta}}+\frac{\partial\mathscr{V}_{\eta}}{\partial \mathscr{Y}}\times \frac{\partial\mathscr{Y}}{\partial v_{\eta}}&\frac{\mathscr{V}_{\eta}}{v_{\phi}}+\frac{\partial\mathscr{V}_{\eta}}{\partial \mathscr{Y}}\times \frac{\partial\mathscr{Y}}{\partial v_{\phi}}&0\\
\frac{\partial\mathscr{V}_{\phi}}{\partial \mathscr{Y}}\times \frac{\partial\mathscr{Y}}{\partial t}&\frac{\partial\mathscr{V}_{\phi}}{\partial \mathscr{Y}}\times \frac{\partial\mathscr{Y}}{\partial v_{\eta}}&\frac{1-\eta_{1}}{1-\mathscr{Y}(t)}+\frac{\partial\mathscr{V}_{\phi}}{\partial \mathscr{Y}}\times \frac{\partial\mathscr{Y}}{\partial v_{\phi}}&0\\
0&0&0&1
\end{matrix}
\right\vert=\frac{1-\eta_{1}}{1-\mathscr{Y}(t)}|v_{\eta}|,
\end{align*}
we make the following change of variable
$$
(y,u_{\eta},u_{\phi},u_{z})\mapsto (\mathscr{Y}(t,\eta_{1},v),\mathscr{V}_{\eta}(\mathscr{Y}(t);t,\eta_{1},v),\mathscr{V}_{\phi}(\mathscr{Y}(t);t,\eta_{1},v),\mathscr{V}_{z}(\mathscr{Y}(t);t,\eta_{1},v))
$$
to get
\begin{align}\label{t6-2}
&\int_{0}^{\tilde{t}}\int_{A_{\tau}}|f(0,\mathscr{Y}(t),\mathscr{V}_{\eta}(\mathscr{Y}(t);t,\eta_{1},v),\mathscr{V}_{\phi}(\mathscr{Y}(t);t,\eta_{1},v),\mathscr{V}_{z}(\mathscr{Y}(t);t,\eta_{1},v))\,{\rm d}v{\rm d}t\nonumber\\
&\leq C\int_{0}^{\eta_{1}}\int_{\R^3}|f(0,y,u)|\,{\rm d}u{\rm d}y.
\end{align}
Substituting \eqref{t6-2} into \eqref{t6-1}, we have
\begin{align}\label{t6-3}
&\int_{0}^{\tilde{t}}\int_{A_{\tau}}|f(t,\eta_{1},v)|\,|v_{\eta}|\,{\rm d}v{\rm d}t\nonumber\\
&\leq C\Big[\|f(0)\|_{L^1}+\int_{0}^{T}\|\{\partial_{t}+v_{\eta}\partial_{v_{\eta}}-\frac{1}{1-\eta}(v_{\phi}^2\partial_{v_{\eta}}-v_{\eta}v_{\phi}\partial_{v_{\phi}})-\frac{\alpha}{\eta_{1}}v_{\eta}\partial_{v_{\phi}}\}f(t)\|_{L^1}\,{\rm d}t\Big].
\end{align}
Combining \eqref{t6-3} with \eqref{t5-4}, we conclude \eqref{t2} for the case $(\eta,v)\in \mathcal{A}_{1}\cap (\gamma_{+}\cup \gamma_{+}^{\tau})$.

{\it Case 2}. If $(\eta,v)\in \mathcal{A}_{2}\cap (\gamma_{+}\backslash \gamma_{+}^{\tau})$, i.e., $\eta=0$, $\tau<|v_{\phi}|\leq \frac{1}{\tau}$, $v_{\eta}<0$ and $\eta_{+}(0,v_{\eta},v_{\phi})\leq \eta_{1}$. Noting that
	$$
	\tau^2\leq v_{\eta}^2=\int_{0}^{\eta_{+}}\frac{2}{1-y}\big[\frac{1}{1-y}v_{\phi}+\frac{\alpha}{2\eta_{1}}-\frac{\alpha}{2\eta_{1}}\frac{1}{1-y}\big]^2\,{\rm d}y\leq \frac{C}{\tau^2}\eta_{+}\Rightarrow \eta_{+}\geq \frac{\tau^4}{C},
	$$
	we have
	\begin{align*}
	|\mathscr{V}_{\eta}(z;t,\eta,v)|^2\geq \int_{z}^{\eta_{+}}\frac{2}{1-y}\big[\frac{1}{1-y}v_{\phi}+\frac{\alpha}{2\eta_{1}}-\frac{\alpha}{2\eta_{1}}\frac{1}{1-y}\big]^2\,{\rm d}y\geq C\tau^2(\eta_{+}-z),
	\end{align*}
	and
	\begin{align*}
	\tilde{t}&=\int_{0}^{\mathscr{Y}(\tilde{t})}\frac{1}{|\mathscr{V}_{\eta}(z;t,\eta,v)|}\,{\rm d}z\leq \frac{1}{C\tau}\int_{0}^{\mathscr{Y}(\tilde{t})}(\eta_{+}-z)^{-\frac{1}{2}}\,{\rm d}z
    =\frac{2}{C\tau}\frac{\mathscr{Y}(\tilde{t})}{\eta_{+}^{\frac{1}{2}}+(\eta_{+}-\mathscr{Y}(\tilde{t}))^{\frac{1}{2}}}\leq \frac{2}{C\tau^3}\mathscr{Y}(\tilde{t}),
	\end{align*}
	which yields that $\mathscr{Y}(\tilde{t})\geq \frac{C\tau^3\tilde{t}}{2}$. It follows from \eqref{t4-1} that 
	 \begin{align}\label{t7}
     &\min\{\eta_{+},\mathscr{Y}(t)\}\int_{\tilde{t}}^{T}\int_{A_{\tau}}|f(t,0,v)|\,|v_{\eta}|\,{\rm d}v{\rm d}t\nonumber\\
     &\leq \int_{\tilde{t}}^{T}\int_{A_{\tau}}\int_{0}^{\min\{\eta_{+},\mathscr{Y}(\tilde{t})\}}|f(\mathscr{T}(y;t,0,v),y,\mathscr{V}_{\eta}(y;t,0,v),\mathscr{V}_{\phi}(y;t,0,v),\mathscr{V}_{z}(y;t,0,v))|\,|v_{\eta}|\,{\rm d}y{\rm d}v{\rm d}t\nonumber\\
     &\quad +\int_{\tilde{t}}^{T}\int_{A_{\tau}}\int_{0}^{\min\{0,\mathscr{Y}(\tilde{t})\}}\int_{0}^{y}|\mathscr{D}_{z}f(\mathscr{T}(z),z,\mathscr{V}_{\eta}(z),\mathscr{V}_{\phi}(z),\mathscr{V}_{z})|\times \frac{|v_{\eta}|}{|\mathscr{V}_{\eta}(z)|}\,{\rm d}z{\rm d}y{\rm d}v{\rm d}t\nonumber\\
     &\leq \int_{\tilde{t}}^{T}\int_{A_{\tau}}\int_{0}^{\min\{\eta_{+},\mathscr{Y}(\tilde{t})\}}|f(\mathscr{T}(y;t,0,v_{\eta},v_{\phi}),y,\mathscr{V}_{\eta}(y;0,v_{\eta},v_{\phi}),\mathscr{V}_{\phi}(y;0,v_{\eta},v_{\phi})|\,|v_{\eta}|\,{\rm d}y{\rm d}v{\rm d}t\nonumber\\
     &\quad +\eta_{1}\int_{\tilde{t}}^{T}\int_{A_{\tau}}\int_{0}^{\min\{\eta_{+},\mathscr{Y}(\tilde{t})\}}|\mathscr{D}_{y}f(\mathscr{T}(y),y,\mathscr{V}_{\eta}(y),\mathscr{V}_{\phi}(y))|\times \frac{|v_{\eta}|}{|\mathscr{V}_{\eta}(y)|}\,{\rm d}y{\rm d}v{\rm d}t.
    \end{align}
Then applying similar calculations as in \eqref{t5-3}--\eqref{t5-5} to \eqref{t7}, we have
\begin{align}\label{t7-1}
&\int_{\tilde{t}}^{T}\int_{A_{\tau}}|f(t,0,v)|\,|v_{\eta}|\,{\rm d}v{\rm d}t\nonumber\\
&\leq C_{\tau,\tilde{t}}\Big[\int_{0}^{T}\|f(t)\|_{L^1}\,{\rm d}t
+\int_{0}^{T}\|\{\partial_{t}+v_{\eta}\partial_{v_{\eta}}-\frac{1}{1-\eta}(v_{\phi}^2\partial_{v_{\eta}}-v_{\eta}v_{\phi}\partial_{v_{\phi}})-\frac{\alpha}{\eta_{1}}v_{\eta}\partial_{v_{\phi}}\}f(t)\|_{L^1}\,{\rm d}t\Big].
\end{align}
It remains to control
$$
\int_{0}^{\tilde{t}}\int_{A_{\tau}}|f(t,0,v_{\eta},v_{\phi})|\,|v_{\eta}|\,{\rm d}v{\rm d}t.
$$
In fact, if we choose $\tilde{t}$ to be small enough such that $\tilde{t}\lesssim \tau^{5}$, at this stage, the backward trajectory hits the initial plane first before hitting the plane $\{\eta=\eta_{1}\}$ due to the fact
$$
\int_{0}^{\eta_{+}}\frac{1}{|\mathscr{V}(z;t,0,v)|}\,{\rm d}z\gtrsim \eta_{+}\tau\geq \frac{\tau^5}{C}.
$$
So, there must exist a unique $\mathscr{Y}(t,0,v_{\eta},v_{\phi})\in (0,\eta_{+})$ satisfying 
$$
t=\int_{0}^{\mathscr{Y}(t)}\frac{1}{|\mathscr{V}_{\eta}(z;t,0,v)|}\,{\rm d}z\quad \text{for each }t\in (0,\tilde{t}).
$$
A direct calculation shows that
$$
\frac{\partial \mathscr{Y}(t,0,v)}{\partial t}=|\mathscr{V}_{\eta}(\mathscr{Y}(t);t,0,v)|.
$$
Therefore, similar as in \eqref{t6-1}--\eqref{t6-3}, we get
\begin{align}\label{t7-2}
&\int_{0}^{\tilde{t}}\int_{A_{\tau}}|f(t,0,v)|\,|v_{\eta}|\,{\rm d}v{\rm d}t\nonumber\\
&\leq C\Big[\|f(0)\|_{L^1}+\int_{0}^{T}\|\{\partial_{t}+v_{\eta}\partial_{v_{\eta}}-\frac{1}{1-\eta}(v_{\phi}^2\partial_{v_{\eta}}-v_{\eta}v_{\phi}\partial_{v_{\phi}})-\frac{\alpha}{\eta_{1}}v_{\eta}\partial_{v_{\phi}}\}f(t)\|_{L^1}\,{\rm d}t\Big].
\end{align}
Combining \eqref{t7-1}--\eqref{t7-2} together, we conclude \eqref{t2} in the case of $(\eta,v)\in \mathcal{A}_{2}\cap (\gamma_{+}\backslash \gamma_{+}^{\tau})$.

{\it Case 3}. If $(\eta,v)\in \mathcal{A}_{3}\cap (\gamma_{+}\backslash \gamma_{+}^{\tau})$, i.e., $\eta=0$, $|v_{\phi}|\geq \tau$, $v_{\eta}<-\tau<0$ and
	$$
	v_{\eta}^2\geq \int_{0}^{\eta_{1}}\frac{2}{1-y}\big[\frac{1}{1-y}v_{\phi}+\frac{\alpha}{2\eta_{1}}-\frac{\alpha}{2\eta_{1}}\frac{1}{1-y}\big]^2\,{\rm d}y,
	$$
	then noting
	\begin{align*}
	|\mathscr{V}_{\eta}(z;t,\eta,v)|^2\geq \int_{z}^{\eta_{1}}\frac{2}{1-y}\big[\frac{1}{1-y}v_{\phi}+\frac{\alpha}{2\eta_{1}}-\frac{\alpha}{2\eta_{1}}\frac{1}{1-y}\big]^2\,{\rm d}y\geq C\tau^2(\eta_{1}-z),
	\end{align*}
	provided that $\alpha\ll \tau$, one has
	\begin{align*}
	\tilde{t}&=\int_{0}^{\mathscr{Y}(\tilde{t})}\frac{1}{|\mathscr{V}_{\eta}(z;t,\eta,v)|}\,{\rm d}z\leq \frac{1}{C\tau}\int_{0}^{\mathscr{Y}(\tilde{t})}(\eta_{1}-z)^{-\frac{1}{2}}\,{\rm d}z\nonumber\\
	&=\frac{2}{C\tau}\big[\eta_{1}^{\frac{1}{2}}-(\eta_{1}-\mathscr{Y}(\tilde{t}))^{\frac{1}{2}}\big]=\frac{2}{C\tau}\frac{\mathscr{Y}(\tilde{t})}{\eta_{1}^{\frac{1}{2}}+(\eta_{1}-\mathscr{Y}(\tilde{t}))^{\frac{1}{2}}}\leq \frac{2}{C\tau}\mathscr{Y}(\tilde{t}),
	\end{align*}
	which implies 
	$\mathscr{Y}(\tilde{t})\geq \frac{C\tau\tilde{t}}{2}$. Then applying almost the same arguments as in Case 2 (just replacing $\eta_{+}$ in Case 2 with $\eta_{1}$), we conclude \eqref{t2} in the case of $(\eta,v)\in \mathcal{A}_{3}\cap (\gamma_{+}\backslash \gamma_{+}^{\tau})$. Therefore, the proof of Lemma \ref{tt} is completed. $\hfill\square$

\medskip
Motivated by \cite{DLY-2021,Guo-2010}, we have the following lemma on controlling the iterative integral.
	\begin{lemma}\label{lem2.1}
	For any $\ell>0$ in the weight function $w_{\ell}(v)$, there exist positive constants $C_{1}$ and $C_{2}$ independent of $k$ and $\ell$ such that if $\alpha\ll \frac{1}{\ell k}$,
		\begin{align}\label{A3}
		\int_{\prod_{j=1}^{k-1}\mathcal{V}_{j}}\sum\limits_{l=1}^{k-1}\mathbf{1}_{\{t_{l+1}\leq 0<t_{l}\}}\int_{0}^{t_{l}}\,{\rm d}\Sigma_{l}^{i}(s)\,{\rm d}s\leq C_{1}k,
		\end{align}
		and
		\begin{align}\label{A4}
			\int_{\prod_{j=1}^{k-1}\mathcal{V}_{j}}\sum\limits_{l=1}^{k-1}\mathbf{1}_{\{t_{l+1}>0\}}\int_{0}^{t_{l}}\,{\rm d}\Sigma_{l}^{i}(s)\,{\rm d}s\leq C_{2}k.
		\end{align}
		
	Furthermore, for any $\tilde{\delta}>0$, $\ell\gg 4$ and any $T_{0}>0$, there exists an integer $k_{0}=k_{0}(\tilde{\delta},T_{0})$ such that for any integer $k\geq k_{0}$ and $\alpha\ll \frac{1}{\ell}$, and for all $(t,\eta,v)\in [0,T_{0}]\times [0,\eta_{1}]\times \R^3$, it holds that
		\begin{align}\label{A1}
	\int_{\prod_{l=1}^{k-1}\mathcal{V}_{l}}\mathbf{1}_{\{t_{k}>0\}}\,{\rm d}\Sigma_{k-1}^{i}(t_{k})<\tilde{\delta}.
		\end{align}
		In particular, let $T_{0}>0$ large enough, there exist positive constants $C_{3}$ and $C_{4}$ independent of $T_{0}$ such that for $k=C_{3}T_{0}^{\frac{5}{4}}$ with a suitable choice of $C_{3}$ such that $k$ is an integer and for all $(t,\eta,v)\in [0,T_{0})\times [0,\eta_{1}]\times \R^3$, it holds
		\begin{align}\label{A2}
			\int_{\prod_{l=1}^{k-1}\mathcal{V}_{l}}\mathbf{1}_{\{t_{k}>0\}}\,{\rm d}\Sigma_{k-1}^{i}(t_{k})\leq \Big(\frac{1}{2}\Big)^{C_{4}T_{0}^{\frac{5}{4}}}.
		\end{align}
		
\end{lemma}

\noindent\textbf{Proof}.
Recall the definition \eqref{3.18-2}. We have
\begin{align*}
	&\int_{\prod_{j=1}^{k-1}\mathcal{V}_{j}}\sum\limits_{l=1}^{k-1}\mathbf{1}_{\{t_{l+1}\leq 0<t_{l}\}}\int_{0}^{t_{l}}\,{\rm d}\Sigma_{l}^{i}(s)\,{\rm d}s\nonumber\\
	&\leq \int_{\prod_{j=1}^{k-1}\mathcal{V}_{j}}\sum\limits_{l=1}^{k-1}{\bf 1}_{\{t_{l+1}\leq 0<t_{l}\}}\int_{0}^{t_{l}}\prod_{j=l+1}^{k-1}\,{\rm d}\sigma_{j}e^{-\int_{s}^{t_{l}}\mathcal{A}(\tau)\,{\rm d}\tau}\tilde{w}(v_{l})\,{\rm d}\sigma_{l}\prod_{j=1}^{l-1}\Big(\frac{\tilde{w}(v_{j})}{\tilde{w}(\mathscr{V}(t_{j+1}))}\Big)e^{-\int_{t_{l+1}}^{t_{l}}\mathcal{A}(\tau)\,{\rm d}\tau}\,{\rm d}\sigma_{j}\,{\rm d}s.
\end{align*}
Using the fact that $\int_{\mathcal{V}_{j}}\,{\rm d}\sigma_{j}=1$, it can be bounded by
\begin{align}\label{A3-2}
	\sum\limits_{l=1}^{k-1}\int_{\prod_{j=1}^{l}\mathcal{V}_{j}}\int_{0}^{t_{l}}e^{-\frac{\nu_{0}}{2}(t_{l}-s)}\tilde{w}(v_{l})\,{\rm d}\sigma_{l}\prod_{j=1}^{l-1}\Big(\frac{\tilde{w}(v_{j})}{\tilde{w}(\mathscr{V}(t_{j+1})}\Big)e^{-\int_{t_{l+1}}^{t_{l}}\mathcal{A}(\tau)\,{\rm d}\tau}\,{\rm d}\sigma_{j}\,{\rm d}s.
\end{align}
Note that
$$
\int_{\mathcal{V}_{l}}\tilde{w}(v_{l})\,{\rm d}\sigma_{l}\leq C,
$$
and
\begin{align}\label{A3-5}
	\frac{\tilde{w}(v_{j})}{\tilde{w}(\mathscr{V}(t_{j+1}))}&=\frac{(1+|\mathscr{V}(t_{j+1})|^2)^{\frac{\ell}{2}}}{(1+|v_{j}|^2)^{\frac{\ell}{2}}}\frac{\sqrt{\mu}(\mathscr{V}(t_{j+1}))}{\sqrt{\mu}(v_{j})}=\frac{(1+|\mathscr{V}(t_{j+1})|^2)^{\frac{\ell}{2}}}{(1+|v_{j}|^2)^{\frac{\ell}{2}}}e^{\frac{|v_{j}|^2-|\mathscr{V}(t_{j+1})|^2}{4}}\nonumber\\
	&\leq (1+C\alpha)^{\frac{\ell}{2}}e^{\frac{|v_{j}|^2-|\mathscr{V}(t_{j+1})|^2}{4}}\leq (1+C\alpha)^{\frac{\ell}{2}}e^{C\alpha+\frac{\alpha}{2}|v_{j,\phi}|^2},
\end{align}
where we have used the fact that
$$
|v_{j}|^2-|\mathscr{V}(t_{j+1})|^2=-\frac{\alpha^2}{2\eta_{1}}(1-\eta_{1})^2+\frac{\alpha^2}{2\eta_{1}}+\frac{\alpha^2}{\eta_{1}^2}\ln (1-\eta_{1})-\frac{2\alpha}{\eta_{1}}v_{j,\phi}\ln (1-\eta_{1})\leq C\alpha+\alpha|v_{k,\phi}|^2.
$$
Then \eqref{A3-2} can be further bounded by
\begin{align}\label{A3-4}
	C\sum\limits_{l=1}^{k-1}(1+C\alpha)^{(l-1)\frac{\ell}{2}}e^{C\alpha (l-1)}(1-\alpha)^{-\frac{l-1}{2}}\leq Ck(1+C\alpha)^{(k-2)\frac{\ell}{2}}e^{C\alpha(k-2)}(1-\alpha)^{-\frac{k-2}{2}},
\end{align}
where we have used the fact that
\begin{align*}
\int_{\mathcal{V}_{j}}e^{\frac{\alpha}{2}|v_{j,\phi}|^2}\,{\rm d}\sigma_{j}=\int_{\mathcal{V}_{j}}\frac{1}{\sqrt{2\pi}}e^{\frac{(\alpha-1)}{2}|v_{j,\phi}|^2}e^{-\frac{1}{2}\big(|v_{j,\eta}|^2+|v_{j,z}|^2\big)}|v_{j,\eta}|\,{\rm d}v_{j}=(1-\alpha)^{-\frac{1}{2}}.
\end{align*}
By choosing $C\alpha\ll \frac{1}{\ell k}$, one has from \eqref{A3-4} that \eqref{A3} can be bounded by $Ce^{3}k$. Similarly, \eqref{A4} can be proved.


Now, we are going to prove \eqref{A1}--\eqref{A2}. Choosing $0<\delta<1$ sufficiently small, we further define non-grazing sets
	\begin{align*}
	\mathcal{V}_{j}^{\delta}=\{v_{j}\in \mathcal{V}_{j}\,|\, v_{j}\cdot \vec{n}(x_{j})\geq \delta\}\cap \Big\{v_{j}\in \mathcal{V}_{j}\,|\,|v_{j}|\leq \frac{1}{\delta}\Big\}.
	\end{align*}
By direct calculations, we have 
	\begin{align*}
		\int_{\mathcal{V}_{j}\backslash \mathcal{V}_{j}^{\delta}}e^{\frac{\alpha}{2}|v_{j,\phi}|^2}\,{\rm d}\sigma_{j}&\leq \int_{v_{j}\cdot \vec{n}(x_{j})\leq \delta}e^{\frac{\alpha}{2}|v_{j,\phi}|^2}\,{\rm d}\sigma_{j}+\int_{|v_{j,\phi}|\geq \frac{1}{2\delta}}e^{\frac{\alpha}{2}|v_{j,\phi}|^2}\,{\rm d}\sigma_{j}+\int_{|v_{j,\eta}|\geq \frac{1}{2\delta}}e^{\frac{\alpha}{2}|v_{j,\phi}|^2}\,{\rm d}\sigma_{j}\nonumber\\
		&\leq C\delta\Big(1+\frac{1}{1-\alpha}\Big)\leq C\delta(1-\alpha)^{-1},
	\end{align*}
where $C$ is independent of $j$. On the other hand, if $v_{j}\in \mathcal{V}_{j}^{\delta}$, then
	$$
	\delta\leq |v_{j,\eta}|\leq \frac{1}{\delta},\quad |v_{j,\phi}|\leq \frac{1}{\delta}.
	$$
	Then, it follows from \eqref{3.12}--\eqref{3.13} that
    \begin{align*}
	&|V_{\phi}(s;t_{j},x_{j},v_{j})|\leq C(\frac{1}{\delta}+\alpha)\leq \frac{C}{\delta},\quad |V_{\eta}(s;t_{j},x_{j},v_{j})|\leq \frac{C}{\delta}\quad \text{for }s\in [t_{j+1},t_{j}).
	\end{align*}
	Furthermore, if $(x_{j},v_{j})\in  \mathcal{A}_{1}$ or $(x_{j},v_{j})\in \mathcal{A}_{3}$, one has
	\begin{align}\label{B2}
	\eta_{1}=|x_{j+1}-x_{j}|\leq \int_{t_{j+1}}^{t_{j}}|V_{\eta}(\tau)|\,{\rm d}\tau\leq C(t_{j}-t_{j+1})\frac{1}{\delta}\Rightarrow t_{j}-t_{j+1}\geq C\delta.
    \end{align}
	If $(x_{j},v_{j})\in \mathcal{A}_{2}$, one has
	$$
	\delta^2\leq v_{j,\eta}^2=\int_{0}^{\eta_{+}(x_{j},v_{j})}\frac{2}{1-y}\Big[\frac{1}{1-y}v_{j,\phi}+\frac{\alpha}{2\eta_{1}}(1-y)-\frac{\alpha}{2\eta_{1}}\frac{1}{1-y}\Big]^2\,{\rm d}y\leq C\frac{\eta_{+}}{\delta^2},
	$$
	which implies $\eta_{+}\geq C\delta^4$. Note that 
    \begin{align*}
    \eta_{+}-0=\int_{t_{j,*}}^{t_{j}}|V_{\eta}(\tau)|\,{\rm d}\tau+\int_{t_{j+1}}^{t_{j,*}}|V_{\eta}(\tau)|\,{\rm d}\tau\Rightarrow t_{j}-t_{j+1}\geq C\delta^{5},
    \end{align*}
    which, together with \eqref{B2}, yields that
    $$
    t_{j}-t_{j+1}\geq C\delta^5.
    $$
	Therefore, if $t_{k}(t,\eta,v,v_{1},v_{2},\cdots,v_{k-1})>0$, then there can be at most $\Big[\frac{CT_{0}}{\delta^5}\Big]+1$ number of $v_{j}\in \mathcal{V}_{j}^{\delta}$ for $1\leq j\leq k-1$. We therefore have
	\begin{align*}
	&\int_{\mathcal{V}_{1}}\cdots\Big\{\int_{\mathcal{V}_{k-1}}\mathbf{1}_{\{t_{k}>0\}}e^{\frac{\alpha}{2}|v_{k-1,\phi}|^2}\,{\rm d}\sigma_{k-1}\Big\}e^{\frac{\alpha}{2}|v_{k-2,\phi}|^2}\,{\rm d}\sigma_{k-2}\cdots e^{\frac{\alpha}{2}|v_{1,\phi}|^2}\,{\rm d}\sigma_{1}\nonumber\\
	&\leq \sum\limits_{j=0}^{\Big[\frac{CT_{0}}{\delta^5}\Big]+1}\int_{\{\text{There are exactly $j$ of $v_{l_{j}}\in \mathcal{V}_{{l_{j}}}^{\delta}$, and $k-1-j$ of $v_{l_{j}}\notin \mathcal{V}_{l_{j}}^{\delta}$}\}}\prod_{l=1}^{k-1}e^{\frac{\alpha}{2}|v_{l,\phi}|^2}{\rm d}\sigma_{l}\nonumber\\
	&\leq \sum\limits_{j=0}^{\Big[\frac{CT_{0}}{\delta^5}\Big]+1}\binom{k-1}{j}\Big|\sup_{l}\int_{V_{l}^{\delta}} e^{\frac{\alpha}{2}|v_{l,\phi}|^2}\,{\rm d}\sigma_{l}\Big|^{j}\Big|\sup_{l}\int_{\mathcal{V}_{l}\backslash \mathcal{V}_{l}^{\delta}} e^{\frac{\alpha}{2}|v_{l,\phi}|^2}\,{\rm d}\sigma_{l}\Big\vert^{k-j-1}.
	\end{align*}
	Note that
	$$
	\int_{V_{l}^{\delta}}e^{\frac{\alpha}{2}|v_{l,\phi}|^2}\,{\rm d}\sigma_{l}\leq (1-\alpha)^{-\frac{1}{2}}\leq (1-\alpha)^{-1},
	$$
	and
	$$
\Big|\int_{\mathcal{V}_{l}\backslash \mathcal{V}_{l}^{\delta}} e^{\frac{\alpha}{2}|v_{l,\phi}|^2}\,{\rm d}\sigma_{l}\Big\vert^{k-j-1}\leq \big[C\delta(1-\alpha)^{-1}\big]^{k-j-1}\leq (1-\alpha)^{-(k-j-1)}(C\delta)^{k-2-\Big[\frac{CT_{0}}{\delta^5}\Big]}.
$$
But it holds that
$$
\sum\limits_{j=0}^{\Big[\frac{CT_{0}}{\delta^{5}}\Big]+1}\binom{k-1}{j}\leq \sum\limits_{j=0}^{\Big[\frac{CT_{0}}{\delta^{5}}\Big]+1}(k-1)^{j}\leq (k-1)^{\Big[\frac{CT_{0}}{\delta^{5}}\Big]+2}.
$$
We then deduce that
\begin{align}\label{A1-1}
	\int_{\prod_{l=1}^{k-1}\mathcal{V}_{l}}\mathbf{1}_{\{t_{k}>0\}}\prod_{l=1}^{k-1}e^{\frac{\alpha}{2}|v_{l,\phi}|^2}\,{\rm d}\sigma_{l}&\leq (k-1)^{\Big[\frac{CT_{0}}{\delta^5}\Big]+2}(1-\alpha)^{-(k-1)}(C\delta)^{k-2-\Big[\frac{CT_{0}}{\delta^5}\Big]}\nonumber\\
	&\leq (k-1)^{\Big[\frac{CT_{0}}{\delta^5}\Big]+2}(\frac{C\delta}{1-\alpha})^{k-1}(C\delta)^{-1-\Big[\frac{CT_{0}}{\delta^5}\Big]}.
\end{align}
Recalling the definitions of ${\rm d}\Sigma_{k-1}^i(t_{k})$ in \eqref{3.18-2} and using \eqref{A3-5}, we have
\begin{align*}
	\int_{\prod_{l=1}^{k-1}\mathcal{V}_{l}}\mathbf{1}_{\{t_{k}>0\}}\,{\rm d}\Sigma_{k-1}^{i}(t_{k})&\leq(1+C\alpha)^{\frac{\ell}{2}(k-1)}e^{C\alpha(k-1)}
	\int_{\prod_{l=1}^{k-1}\mathcal{V}_{l}}\mathbf{1}_{\{t_{k}>0\}}\prod_{l=1}^{k-1}e^{\frac{\alpha}{2}|v_{l,\phi}|^2}\,{\rm d}\sigma_{l}\nonumber\\
	&\leq (k-1)^{\Big[\frac{CT_{0}}{\delta^5}\Big]+2}(\frac{C\delta(1+C\alpha)^{\frac{\ell}{2}}(e^{C\alpha})}{1-\alpha})^{k-1}(C\delta)^{-1-\Big[\frac{CT_{0}}{\delta^5}\Big]}.
\end{align*}
Then for any $\tilde{\delta}>0$ and $C\alpha\ll \frac{1}{\ell}<\frac{1}{2}$, \eqref{A1} follows by choosing $\delta$ such that $\frac{C\delta(1+C\alpha)^{\frac{\ell}{2}}(e^{C\alpha})}{1-\alpha}<1$, and then taking $k\gg \Big[\frac{CT_{0}}{\delta^5}\Big]+1$.

In particular, we let $k-1=N\big\{\Big[\frac{CT_{0}}{\delta^5}\Big]+2\big\}\geq \frac{1}{C\delta}$ so that if $\frac{CT_{0}}{\delta^5}\geq 1$, then \eqref{A1-1} can be further bounded by
\begin{align*}
	\Big\{N\big(\frac{CT_{0}}{\delta^5}+2\big)\big(\frac{C\delta(1+C\alpha)^{\frac{\ell}{2}}(e^{C\alpha})}{1-\alpha}\big)^{N}\frac{1}{C\delta}\Big\}^{\Big[\frac{CT_{0}}{\delta^5}\Big]+2}\leq \Big(C_{N}T_{0}\delta^{N-6}\Big)^{\Big[\frac{CT_{0}}{\delta^5}\Big]+2}\quad \text{if $C\alpha\ll \frac{1}{\ell}<\frac{1}{2}$}.
\end{align*}
We choose $C_{N}T_{0}\delta^{N-6}=\frac{1}{2}$, so that $\delta=\big(\frac{1}{2C_{N}T_{0}}\big)^{\frac{1}{N-6}}$ is small for $T_{0}$ large and $N>6$. Moreover,
$$
\Big[\frac{CT_{0}}{\delta^5}\Big]+2\sim C_{N}T_{0}^{1+\frac{5}{N-6}},
$$
and $\frac{CT_{0}}{\delta^5}\geq 2$ for $T_{0}$ large so that we can close our estimate. Finally, we choose $N=26$. For $T_{0}$  sufficiently large, $k=N\big\{\Big[\frac{CT_{0}}{\delta^5}\Big]+2\big\}+1\sim CT_{0}^{\frac{5}{4}}$ and \eqref{A2} follows.\qed

\medskip
Finally, we recall a useful lemma in \cite{DHWZ-19}, which has been frequently used in obtaining $L^{\infty}$ estimates for approximate solutions.

\begin{lemma}[\cite{DHWZ-19}]\label{lemC1}
	Consider a sequence $\{a_{i}\}_{i=0}^{\infty}$ with each $a_{i}\geq 0$. For any fixed $k\in \mathbb{N}_{+}$. We denote
	$$
	A_{i}^{k}=\max\{a_{i},a_{i+1},\cdots,a_{i+k}\}.
	$$
	Assume $D\geq 0$. If $a_{i+1+k}\leq \frac{1}{8}A_{i}^{k}+D$ for $i=0,1,\cdots$, then it holds that
	\begin{align*}
	A_{i}^{k}\leq \Big(\frac{1}{8}\Big)^{\Big[\frac{i}{k+1}\Big]}\cdot \max\{A_{0}^{k},A_{1}^{k},\cdots,A_{k}^{k}\}+\frac{8+k}{7}D,\quad \text{for }i\geq k+1.
	\end{align*}
\end{lemma}

\medskip
\noindent\textbf{Acknowledgments}: The research was partially supported by National Key R\&D Program of China grant No. 2021YFA1000800. Renjun Duan's research was partially supported by the General Research Fund (Project No.~14303523) from RGC of Hong Kong and also by the grant from the National Natural Science Foundation of China (Project No.~12425109). Yong Wang's research was partially supported by the National Natural Science Foundation of China grants No. 12421001 \& 12288201, and CAS Project for Young Scientists in Basic Research, grant No. YSBR-031.

\medskip
\noindent\textbf{Data Availability Statement:}
Data sharing is not applicable to this article as no datasets were generated or analysed during the current study.

\noindent\textbf{Conflict of Interest:}
The authors declare that they have no conflict of interest.


\end{document}